\documentclass[11pt]{amsart}
\usepackage{amssymb}
\usepackage{graphicx}
\usepackage{csquotes}
\usepackage{mathrsfs}
\usepackage[all]{xy}
\usepackage{imakeidx}
\usepackage[colorinlistoftodos]{todonotes}
\usepackage{bm}
\usepackage[colorlinks=true]{hyperref}
\makeindex[name=terms, title=Subject index, columns=1]
\makeindex[name=not, title=Index of notation, columns=2]
\newcommand{\Ft}{F^{*} \otimes_{\Z_{2}}\mathrm{o}\left(TY\right)}

\newcommand{\Ou}{\scriptsize\bullet}

\newcommand{\Aut}{\mathrm{Aut}}

\newcommand{\TsX}{T^{*}X}

\newcommand{\ad}{\mathrm{ad}}
\newcommand{\pa}{\partial}

\newcommand{\we}{\wedge}
\newcommand{\ho}{\widehat{\otimes}}
\newcommand{\n}{\nabla}
\newcommand{\N}{\mathbf{N}}
\newcommand{ \Z}{\mathbf{Z}}
\newcommand{ \R}{\mathbf{R}}
\newcommand{ \C}{\mathbf{C}}

\def \Hom{\mathrm {Hom}}
\def \Trs{\mathrm {Tr_{s}}}
\def \Tr{\mathrm{Tr}}

\def \End{\mathrm {End}}

\def \Re{\mathrm{Re}}
\def\deg{\mathrm{deg}}
\def\Im{\mathrm{Im}}

\newcommand{\cD}{\mathcal{D}}

\newcommand{\cH}{\mathcal{H}}

\newcommand{\bR}{\mathbf{R}}
\newcommand{\bC}{\mathbf{C}}

\renewcommand{\Re}{\mathrm{Re}\,}

\DeclareMathOperator{\Dom}{\mathrm Dom}

\newcommand{\<}{\langle}
\renewcommand{\>}{\rangle}
\newcommand{\ol}{\overline}
\newcommand{\ul}{\underline}

\renewcommand{\(}{\left(}
\renewcommand{\)}{\right)}

\renewcommand{\]}{\right]}
\renewcommand{\l}{\leqslant}

\newcommand{\g}{\geqslant}

\newcommand{\e}{\epsilon}

\theoremstyle{plain}
\newtheorem{theorem}{Theorem}

\newtheorem{proposition}[theorem]{Proposition}

\theoremstyle{definition}
\newtheorem{definition}[theorem]{Definition}

\newtheorem*{conjecture}{A modified Fried conjecture}

\theoremstyle{remark}
\newtheorem{remark}[theorem]{Remark}

\numberwithin{equation}{section}
\numberwithin{theorem}{section}
\begin{document} 
\bibliographystyle{amsalpha} 
  
\title{Anosov vector fields and Fried sections}


\author{Jean-Michel \textsc{Bismut}}
\address{Institut de Mathématique d'Orsay \\ Université Paris-Saclay
\\ Bâtiment 307 \\ 91405 Orsay \\ France} 
\curraddr{}
\email{jean-michel.bismut@universite-paris-saclay.fr}
\author{Shu \textsc{Shen}}
\address{Institut de Mathématiques de Jussieu-Paris Rive Gauche\\
Sorbonne Université\\ Case Courrier 247\\4 place Jussieu\\75252 Paris 
Cedex 05\\ France}
\email{shu.shen@imj-prg.fr}

\thanks{}
\date{November $19^{\rm{th}}$ 2024}


\dedicatory{}

\begin{abstract}
	The purpose of this paper is to prove that if $Y$ is a compact 
	manifold, if $Z$ is an Anosov vector field on $Y$,  and if $F$ is a flat vector bundle,  there is a 
	corresponding canonical nonzero section 
	$\tau_{\nu}\left(i_{Z}\right)$  of the 
	determinant line $\nu=\det H\left(Y,F\right)$. In families, this 
	section is $C^{1}$ with respect to the canonical smooth structure on $\nu$. When $F$ is flat on the total space of the 
	corresponding fibration, our section is flat with respect to the 
	Gauss-Manin connection on $\nu$. 
 \end{abstract}
\maketitle
\tableofcontents
\section{Introduction}
The purpose of this paper is to develop the  theory of 
determinant line bundles associated with a family of Anosov vector 
fields, and with a given flat vector bundle.

Our paper builds on three different and related threads:
\begin{enumerate}
	\item  The construction by Fried \cite[Section 5]{Fried87}   of the zeta function of a 
	hyperbolic vector field $Z$ on a compact manifold $Y$ equipped with 
	a flat vector bundle $F$.  The Fried zeta function 
	$R_{Z,F}\left(\sigma\right)$ is obtained by a principle similar to the Ruelle zeta 
	function \cite{Ruelle76a,Ruelle76b}, and led Fried to formulate the Fried conjecture, that 
	states that this zeta function has a meromorphic continuation, 
	and that the corresponding constant term at $0$ may be related 
	to the Reidemeister torsion \cite{Reidemeister35}  of the manifold.
\item  The proofs by Giulietti-Liverani-Pollicott \cite[Corollary 2.2]{GiuliettiLiveraniPollicott13}  and by Dyatlov-Zworski \cite[p. 543]{DyatlovZworski16} of the fact 
	that when the vector field $Z$  is Anosov, the Fried zeta function can 
	be extended by analytic continuation. In particular, in a series 
	of fundamental papers, Dyatlov-Zworski \cite{DyatlovZworski16}, Faure-Roy-Sjöstrand \cite{FaureRoySjostrand08}, and Faure-Sjöstrand \cite{FaureSjostrand11} 
	developed the proper functional analytic theory for  an Anosov vector 
	field $Z$, introduced associated anisotropic Sobolev spaces,  proved that the 
	spectrum of the Lie derivative operator $L_{Z}$ is discrete, and showed that the zeros and 
	poles of the Fried zeta function are contained in this spectrum.
	\item  The proof  by Dang, Guillarmou, Rivière and Shen \cite{DangShen20} that when $0$ is not included in  the 
	spectrum of $L_{Z}$, the value at $0$ of the Fried zeta function 
	(which is a complex number)
	does not change under small deformations of $Z$.
\end{enumerate}

The purpose of the present paper is to get rid of the invertibility 
assumption of $L_Z$ that was made in \cite{DangShen20}. Let $H\left(Y,F\right)$ be the cohomology of $Y$ with coefficients in the flat vector bundle $F$.\footnote{If $\left(\Omega\left(Y,F\right),d^{Y}\right)$ is the de Rham complex on $Y$ with coefficients in $F$, $H\left(Y,F\right)$ is just the cohomology of this complex.} The invertibility assumption on $L_{Z}$ 
implies in particular that $H\left(Y,F\right)=0$.

Set
\index[not]{nu@$\nu$}%
\begin{equation}\label{eq:defi1}
\nu=\det H\left(Y,F\right).
\end{equation}
Then $\nu$ is a complex line. 

Let $\pi:M\to S$ be a proper 
submersion of smooth manifolds with compact fiber $Y$, and let $\left(F,\n^{F}\right)$ be a smooth complex vector bundle with connection on $M$. Assume that $\n^{F}$ is flat along the fiber $Y$. By a construction inspired from Quillen \cite{Quillen85b}, using the family of  de Rham complexes $\Omega\left(Y,F\vert_{Y}\right)$, $\nu$ can be turned into a smooth complex line bundle on $S$. If $\n^{F}$ is flat on $M$, $\nu$ inherits the flat Gauss-Manin connection.\footnote{In Section 
\ref{subsec:flam}, we will recall the construction of the Gauss-Manin connection.}

We are thus 
forced  to  construct explicitly the determinant 
line bundle $\nu$  on $S$ in the context of Anosov vector fields, and to 
compare it  with the corresponding construction via the 
family of de Rham complexes.  

When no invertibility assumption is made on $L_{Z}$, the 
question arises of how to state the proper version of the Fried 
conjecture \cite{Fried87}.  
\subsection{Our main result}%
\label{subsec:main}
In the present paper, we prove the following result.
\begin{theorem}\label{thm:funda}
	There is a  canonical nonzero section 
	\index[not]{tniZ@$\tau_{\nu}\left(i_{Z}\right)$}%
	$\tau_{\nu}\left(i_{Z}\right)$ of $\nu$,  whose 
	construction is done using the Fried zeta function and the 
	spectral theory of $L_{Z}$. In families, 
	$\tau_{\nu}\left(i_{Z}\right)$
	 is a $C^{1}$ section of $\nu$ with respect to its canonical smooth structure. When $\n^{F}$ is flat on $M$, 
	 $\tau_{\nu}\left(i_{Z}\right)$ is 
	a flat section of $\nu$ with respect to the Gauss-Manin connection, and so it is  smooth.
\end{theorem}

Let us give the proper perspective to the above result. When $L_{Z}$ 
is invertible, then $H\left(Y,F\right)=0$, and  $\nu$ is 
canonically trivial. In this case, 
$\tau_{\nu} \left(i_{Z}\right)=R_{Z,F}\left(0\right)$. 
When $F$ is flat on $M$, our result reduces to the main result of \cite{DangShen20} 
that 
says that $R_{Z,F}\left(0\right)$ is invariant by a small deformation of $Z$. 
When $L_{Z}$ is non invertible, special care has to be given to the 
fact that, in general,  the dimension of the space of resonances\footnote{The space of resonances associated with an eigenvalue of $L_Z$ is just the characteristic subspace associated with this eigenvalue.} associated with the $0$ 
eigenvalue do not necessarily form a vector bundle on $S$, even when $F$ is flat 
on $M$. As we shall see later, a related result was established by \cite{ChaubetDang24} when $Y$ is contact.
\subsection{The methods used in the paper}%
\label{subsec:methods}
Let us briefly describe our techniques. As explained before, the results by Dyatlov-Zworski \cite{DyatlovZworski16}, by Faure-Roy-Sjöstrand \cite{FaureRoySjostrand08} and Faure-Sjöstrand \cite{FaureSjostrand11} play a fundamental role in the present paper. More specifically, we use in  an extensive way their construction of the spectral theory of the operator $L_{Z}$, and of the corresponding microlocal Sobolev spaces. For $a>0$,  the finite dimensional   vector spaces of currents $\mathcal{D}_{Z,<a}\left(Y,F\right)$ associated with eigenvalues $\lambda$ of $L_{Z}$ such that $\left\vert \lambda\right\vert<a$ play the same role as  in the study by \cite{Quillen85b,BismutFreed86a,BismutGilletSoule88c} of the determinant of the cohomology. A  conceptually important point is that we ultimately view the resolvent of $L_{Z}$ or the spectral projectors as classical currents, that are naturally dual to the spaces of smooth forms. Their intersection traces or supertraces are  obtained via their intersection with the diagonal, once the proper wave front set conditions are satisfied. An other important point is that the classical Poincaré duality of currents is used systematically in our constructions. In particular, $\mathcal{D}_{Z,<a}\left(Y,F\right)$ and $\mathcal{D}_{-Z,<a}\left(Y,\Ft\right)$ are shown to be naturally Poincaré dual.

Let us explain  the role of Poincaré duality. In section \ref{sec:opcur}, we develop in detail the theory of currents on $Y\times Y$ and of the corresponding operators that are obtained by integration along the fiber. If $n=\dim Y$, the objects we deal with are essentially currents of degree $n$ on $Y\times Y$. As an example, the resolvent of $L_{Z}$ is  viewed as a closed current on $Y\times Y$ instead of an operator that commutes with $d^{Y}$. This allows us  to ignore questions of orientability of $Y$, and to deal with complicate questions of signs. Both $L_{Z}$ and its proper adjoint are dealt with simultaneously. This is especially important in light of the fact that as is well-known from earlier work \cite{BismutLott95}, analytic torsion is a secondary invariant with respect to  the Euler class, that itself ignores orientability, and corresponds  to the self-intersection of the diagonal in $Y\times Y$.

To express the intersection of a current of degree $n$ with the diagonal in $Y\times Y$, we decided to call it its intersection supertrace. The reason is that when considering smooth currents, this is just the supertrace of the corresponding operator. In the literature, this is sometimes called a flat supertrace, a terminology we preferred to avoid here.

To establish our results in an  infinite dimensional context,  formally 
similar results have to be established first in a finite dimensional context, that are later applied to the 
spaces of truncated resonances.

To study the dependence of the section $\tau_{\nu}\left(i_{Z}\right)$ on a  family of Anosov vector fields parametrized by $s\in S$, like in \cite{DangShen20}, we take into account the fact that our microlocal Sobolev spaces vary  with $s$. As was shown by Chaubet-Dang \cite{ChaubetDang24}, some  arguments can be  simplified using a method given by Guedes
Bonthonneau \cite{Bonthonneau20}, who was able to use instead locally constant microlocal Sobolev spaces. In this situation in families, where the cohomology of $F$ may jump from fiber to fiber, there is a corresponding smooth structure on the line bundle $\nu$, that is obtained via the corresponding family of de Rham complexes. In this context, we prove that $\tau_{\nu}\left(i_{Z}\right)$ is a $C^{1}$ section of $\nu$. As mentioned before, when $F$ is a flat vector bundle on the total space of the deformation, then $\nu$ is a smooth flat line bundle on $S$, and $ \tau_{\nu}\left(i_{Z}\right) $ is a smooth flat section of $\nu$. 

\subsection{Earlier work on the subject}%
\label{subsec:ear}
Our paper should be compared with earlier work by Chaubet and Dang \cite{ChaubetDang24}, where similar results were obtained in the case where the Anosov vector field $Z$ is  the Reeb vector field on $Y$ associated with a contact form, which forces $Y$ to be odd-dimensional and oriented. In their paper,  a canonical nonzero section $\tau_{\nu }^{\mathrm{CD}}$ of $\nu$ is obtained.  Our results are valid in full generality, and do not require  $Y$ to be contact. In Subsection \ref{subsec:chada}, we will show that in this special case, our section $\tau_{\nu}\left(i_{Z}\right)$ coincides with $\tau_{\nu}^{\mathrm{CD}}$.  
\subsection{The organization of the paper}%
\label{subsec:org}
This paper is organized as follows. In Section \ref{sec:twdi}, we 
develop the proper finite dimensional algebraic machinery, that will 
be later used in our infinite dimensional constructions. In 
particular, when  $E=\bigoplus_{i=p}^{q}  E^{i}$ is a finite dimensional complex $\Z$-graded vector 
space equipped with differentials $d,\delta$  of degree $1$ and $-1$, if 
$d,\delta$ are exact, we construct the canonical nonzero sections 
$\tau\left(d\right),\tau\left(\delta\right)$ of 
$\det E$, and we compare these two sections. We also extend the 
comparison formula to the case where $d$ is not exact. 

In Section \ref{sec:frizet}, if $Z$ is an Anosov vector field, we 
recall the construction of the Fried zeta function 
$R_{Z,F}\left(\sigma\right), \Re\,\sigma\gg 1$. Also we prove an 
elementary identity that describes the behavior of the Fried zeta 
function under Poincaré duality\footnote{More precisely, we establish an identity that will be later related to Poincaré duality.}. 

In Section \ref{sec:opcur}, we develop the theory of currents on 
$Y\times Y$ and of 
the associated operators, and we introduce the intersection traces and  
supertraces. This allows us to reformulate statements on a class of operators on $Y$  as being associated with their behavior as currents on $Y\times Y$. The $\Z$-grading of this space of currents plays a fundamental role in the formulation  of the results. 

In Section \ref{sec:avf}, following earlier work by Faure-Roy-Sjöstrand \cite{FaureRoySjostrand08}, Faure-Sjöstrand \cite{FaureSjostrand11}, and Dyatlov-Zworski \cite{DyatlovZworski16}, we introduce the microlocal anisotropic Sobolev spaces 
associated with an Anosov vector field $Z$, and we prove that in a 
suitable sense, the resolvent $\left(L_{Z}-z\right)^{-1}$ is a 
meromorphic family of Fredholm operators. In particular the spectrum 
of $L_{Z}$ is discrete.

In Section \ref{sec:frize},  we describe in detail the microlocal properties of the 
spectral projectors $P_{Z,<a}$ corresponding to eigenvalues $z$ of 
$L_{Z}$ with $\left\vert  z\right\vert<a$, and we construct a 
truncated version $R_{Z,F,>a}\left(\sigma\right)$ of the Fried zeta 
function.

In Section \ref{sec:detli}, we construct a complex line $\lambda$ 
that is canonically isomorphic to $\nu=\det H\left(Y,F\right)$, and a 
corresponding canonical nonzero section 
$\tau_{\lambda}\left(i_{Z}\right)$ of $\lambda$.  It defines a corresponding 
section $\tau_{\nu}\left(i_{Z}\right)$ of $\nu$.

In Section \ref{sec:para}, in the context of families of 
Anosov vector fields, we prove that 
$\tau_{\nu}\left(i_{Z}\right)$ is a $C^{1}$ 
section of $\nu$, and that when $F$ is flat on $M$, that it is a flat section of $\nu$.

Finally, in Section \ref{sec:rsmfr}, we give a proper formulation of the Fried conjecture. Recall that by results of Reidemeister \cite{Reidemeister35}, Cheeger \cite{Cheeger79}, Müller \cite{Muller78,Muller93b}, Bismut-Zhang \cite{BismutZhang92}, when $F$ is unimodular,  a canonical metric $\left\Vert \,\right\Vert_{\nu}$ on the complex line $\nu$ can be defined.  This metric can either be obtained combinatorially or via the Ray-Singer analytic torsion \cite{RaySinger71}. Our  conjecture says that  $\left\Vert \tau_{\nu}\left(i_{Z}\right)\right\Vert_{\nu}$   is equal to $1$.
\section{Complexes equipped with two differentials}%
\label{sec:twdi}
The purpose of this section is to develop in some detail the theory of finite dimensional determinant lines
 associated with a finite dimensional complex  $\Z$-graded vector space $E$ equipped 
with two differentials $d,\delta$ of respective degree $1$ and $-1$. 
When the differentials are exact, we give an explicit formula 
comparing corresponding nonzero canonical sections $\tau\left(d\right), \tau\left(\delta\right)$ of $\det E$ in 
terms of an algebraic version of analytic torsion. 

This Section is organized as follows. In Subsection 
\ref{subsec:twdif}, we consider the spaces of such differentials and 
their tangent bundle.

In Subsection \ref{subsec:siid}, we establish a simple identity on 
$2$-forms over such spaces.

In Subsection \ref{subsec:clof}, we construct a canonical closed $1$-form on 
the above spaces.

In Subsection \ref{subsec:sipid}, we establish an elementary identity 
of linear algebra.

In Subsection \ref{subsec:det}, we briefly recall basic results of \cite{KnudsenMumford} on the determinant functor.

In Subsection \ref{subsec:exad}, given exact differentials $d,\delta$ 
of degree $1,-1$, we 
recall the definition of the corresponding canonical sections $\tau\left(d\right),\tau\left(\delta\right)$  of $\det E$, and we  
establish a comparison formula in terms of a version of analytic 
torsion for $\left[d,\delta\right]$ when $\left[d,\delta\right]$ is 
 invertible.

 In Subsection \ref{subsec:acau}, we consider the action of  the automorphism group $\mathrm{Aut}^{0}\left(E\right)$ on $\det E$ and on the section $\tau\left(d\right)$.

 In Subsection \ref{subsec:exap}, when $\left[d,\delta\right]$ is invertible,  we show that $\tau\left(d\right) \otimes \tau\left(\delta\right)^{-1}$ is just the analytic torsion of $\left[d,\delta\right]$.

In Subsection \ref{subsec:canco}, we show that the canonical $1$-form of 
Subsection \ref{subsec:clof} is just the connection $1$-form for the 
obvious connection on $\det E$.

In Subsection \ref{subsec:hefo}, when $E$ is equipped with a Hermitian 
metric,  we specialize the results of Subsection 
\ref{subsec:exad} to the case where  the two differentials are adjoint to each 
other.

In Subsection \ref{subsec:deex}, we extend the results of Subsection 
\ref{subsec:exad} when only $\delta$ is exact.

In Subsection \ref{subsec:dual}, when replacing $E$ by $E^{*}$, we compare the canonical 
sections of the determinants of $E$ and $E^{*}$.

In Subsection \ref{subsec:eveca}, following \cite{KnudsenMumford}, we consider the effect of shifts in the determinant line.

In Subsection \ref{subsec:fix}, we define the determinant of a linear map as a section of a complex line.

Finally, in Subsection \ref{subsec:spe}, we consider a special case where our exact complex is equipped with an odd automorphism $\Gamma$, and we expressed the canonical section  as the determinant of this automorphism. These results  will be used in Subsection \ref{subsec:chada}, where we compare the  results of \cite{ChaubetDang24} with ours.
\subsection{Differentials and their tangent bundle}%
\label{subsec:twdif}
Let $p,q\in \Z,p\le q$. Let $E=\bigoplus_{i=p}^{n}E^{q}$ be a finite dimensional complex 
$\Z$-graded vector space.  Then $\End\left(E\right)$ is a $\Z$-graded algebra, the elements 
of degree $i\in \Z$, $\End^{i}\left(E\right)$,  increasing the degree by $i$ 
in $E$. Let $\End^{\ge 0}\left(E\right)$ denote the subalgebra of 
elements of nonnegative degree.   In the sequel, $\left[\ \right]$ denotes the
\index[terms]{supercommutator}%
supercommutator\footnote{If $A$ is a $\Z_2$-graded algebra, if $a,b\in \mathcal{A}$, the supercommutator $\left[a,b\right]$ is given by $\left[a,b\right]=ab-\left(-1\right)^{\mathrm{deg}a\mathrm{deg}b}ba$.} 
in $\End\left(E\right)$. Let 
\index[not]{N@$N$}%
$N\in \End^{0}\left(E\right)$ define the 
$\Z$-grading, i.e., $N$ acts by multiplication by $i$ on $E^{i}$. This operator will also be called the 
\index[terms]{number operator}%
number operator of $E$.

Let $\Aut^{0}\left(E\right)$ be the Lie group of degree preserving 
automorphisms of $E$. Then $\End^{0}\left(E\right)$ is the 
corresponding Lie algebra. 

Let
\index[not]{D@$\mathcal{D}$}%
 $\mathcal{D}$ denote the space of the  differentials  $d$ on $E$, i.e., $d$ has degree $1$, and $d^{2}=0$.   If 
$d\in \mathcal{D}$, we denote by $H\left(E,d\right)$ the cohomology 
of $\left(E,d\right)$.
Let $\mathcal{D}$ denote the space of such differentials. If 
$d\in\mathcal{D}$,  and if $\ad\left(d\right)$ denotes the 
supercommutator with $d$,  $\ad\left(d\right)$ 
is a differential on  $\End\left(E\right)$.

Then 
$\Aut^{0}\left(E\right)$ acts on $\mathcal{D}$ by conjugation.  If 
$g\in \Aut^{0}\left(E\right),d\in \mathcal{D}$, we use the notation
\begin{equation}\label{eq:par0}
g.d=gdg^{-1}.
\end{equation}

If $d\in 
\mathcal{D}$, let $\Aut^{0}_{d}\left(E\right)$ denote the Lie group of 
automorphisms of the complex $\left(E,d\right)$.  Equivalently 
$\Aut_{d}^{0}\left(E\right)$ is the stabilizer of $d$ in 
$\Aut^{0}\left(E\right)$. Then $\left(\End^{\ge 
0}E,\ad\left(d\right) \right) $ is also a complex. The  Lie algebra of 
$\Aut^{0}_{d}\left(E\right)$ coincides with 
$H^{0}\left(\End^{\ge 
0}E,\ad\left(d\right) \right)$. The orbit $O_{d} \subset \mathcal{D}$ of 
$d$ under $\Aut^{0}\left(E\right)$ can be identified with 
$\Aut^{0}\left(E\right)/\Aut^{0}_{d}\left(E\right)$. The normal bundle to 
$O_{d}$ in $\mathcal{D}$ is given by  
$H^{1}\left(\End^{\ge 0}\left(E\right),\ad\left(d\right)\right)$. Also 
$\Aut^{0}_{d}\left(E\right)$ acts  on 
$H^{1}\left(\End^{\ge 0}\left(E\right),\ad\left(d\right)\right)$. 
Near $d$, the quotient $\Aut^{0}\left(E\right)\setminus \mathcal{D}$ is 
modeled on the quotient $\Aut_{d}^{0}\left(E\right)\setminus 
H^{1}\left(\End^{\ge 0}\left(E\right),\ad\left(d\right)\right)$.

Let
\index[not]{De@$\mathcal{D}_{\mathrm{e}}$}%
$\mathcal{D}_{\mathrm{e}}$ denote the subspace of exact 
differentials, i.e., the space of $d\in\mathcal{D}$ such that 
$H\left(E,d\right)=0$.  In the sequel, we assume that 
$\mathcal{D}_{\mathrm{e}}\neq \emptyset$. Then 
$\Aut^{0}\left(E\right)$ acts transitively on $\mathcal{D}_{\mathrm{e}}$, so 
that if $d\in \mathcal{D}_{\mathrm{e}}$, 
\begin{equation}\label{eq:dif1}
\mathcal{D}_{\mathrm{e}} = \Aut^{0}\left(E\right)/\Aut^{0}_{d}\left(E\right).
\end{equation}

If $d\in \mathcal{D}_{\mathrm{e}}$, let $T$ be a
\index[terms]{homotopy}%
homotopy for $d$, 
i.e., $T\in \End^{-1}\left(E\right)$ is such that
\begin{equation}\label{eq:dif2}
\left[d,T\right]=1.
\end{equation}
Let $L_T$ denote   left multiplication by $T$. Then $L_{T}$ acts on 
$\End\left(E\right)$ with degree $-1$,  and $L_{T}$ is a homotopy for the 
differential $\ad\left(d\right)$.
\subsection{An identity of $2$-forms}%
\label{subsec:siid}
Let 
\index[not]{Dm@$\mathcal{D}^{-}$}%
$\mathcal{D}^{-}$ be the space of differentials of degree $-1$, 
and let  
\index[not]{Dme@$\mathcal{D}^{-}_{\mathrm{e}}$}
$\mathcal{D}^{-}_{\mathrm{e}}$ denote the space of exact 
differentials of degree $-1$. Since $\mathcal{D}_{\mathrm{e}}\neq \emptyset$, 
 we have $\mathcal{D}^{-}_{\mathrm{e}}\neq \emptyset$.  If $\delta\in \mathcal{D}_{\mathrm{e}}^{-}$, let 
$\alpha\in \mathcal{D}_{\mathrm{e}}$ be a differential such that
\begin{equation}\label{eq:clev1}
\left[\delta,\alpha\right]=1.
\end{equation}
In particular, $\alpha$ is a homotopy with respect to $\delta$.  If $g^{E}$ is 
a Hermitian metric on $E$, and if $\delta^{*}\in \mathcal{D}_{\mathrm{e}}$ is 
the adjoint of $\delta$, we can choose 
\begin{equation}\label{eq:clev2}
\alpha=\left[\delta,\delta^{*}\right]^{-1}\delta^{*}.
\end{equation}
The above shows that we
can take $\alpha$ depending smoothly on $d$. In the sequel, we will 
not assume that (\ref{eq:clev2}) holds.

If $\alpha'\in \mathcal{D}_{\mathrm{e}}$ is also such that 
$\left[\delta,\alpha'\right]=1$, observe that
\begin{equation}\label{eq:clev3}
\left[\delta,\alpha\alpha'\right]=\alpha'-\alpha.
\end{equation}

If $\mathcal{A},\mathcal{A}'$ are $\Z_2$-graded algebras, we denote by  $\mathcal{A}\ho\mathcal{A}'$ the algebra which is the $\Z_2$-graded tensor product of $\mathcal{A},\mathcal{A}'$. In particular $\Lambda\left(T^{*}\mathcal{D}_{\mathrm{e}}^{-}\right)\ho\End\left(E\right)$ is a $\Z_2$-graded bundle of algebras.

Let 
\index[not]{d@$\mathbf{d}$}
$\mathbf{d}$ denote the de Rham operator on 
$\mathcal{D}^{-}_{\mathrm{e}}$. Then  $\mathbf{d}\delta$ is a smooth 
section of $T^{*}\mathcal{D}^{-}_{\mathrm{e}}\ho \End^{-1}\left(E\right)$. 
Since $\delta^{2}=0$, we get
\begin{equation}\label{eq:clev4}
\left[\delta,\mathbf{d}\delta\right]=0.
\end{equation}
In the sequel, we will allow ourselves to pick $\alpha$ depending 
smoothly on $\delta$, so that $\mathbf{d}\alpha$ is well-defined. 
Ultimately, the expressions we will consider will not depend on this 
choice.
\begin{proposition}\label{prop:pidf}
	The following identity holds:
	\begin{equation}\label{eq:idab}
\mathbf{d}\left(\alpha\mathbf{d}\delta\right)=\left[\delta,\left(\mathbf{d}\alpha\right)\alpha\mathbf{d}\delta\right]+\left(\alpha\mathbf{d}\delta\right)^{2}.
\end{equation}
\end{proposition}
\begin{proof}
	Using (\ref{eq:clev4}), we get
\begin{equation}\label{eq:clev8}
\mathbf{d}\delta=\left[\delta,\alpha\mathbf{d}\delta\right].
\end{equation}
By (\ref{eq:clev8}), we obtain
\begin{equation}\label{eq:clev9}
\mathbf{d} \left( \alpha\mathbf{d}\delta \right) =\left[\delta,\left(\mathbf{d}\alpha\right)\alpha\mathbf{d}\delta\right]-\left[\delta,\mathbf{d}\alpha\right]\alpha\mathbf{d}\delta.
\end{equation}
The identity $\left[\delta,\alpha\right]=1$ gives
\begin{equation}\label{eq:clev11}
\left[\mathbf{d\delta},\alpha\right]-\left[\delta,\mathbf{d}\alpha\right]=0.
\end{equation}
Therefore (\ref{eq:clev9}) can be rewritten in the form,
\begin{equation}\label{eq:clev9a1}
\mathbf{d}\left( \alpha\mathbf{d}\delta \right) =\left[\delta,\left(\mathbf{d}\alpha\right)\alpha\mathbf{d}\delta\right]-\left[\mathbf{d}\delta,\alpha\right]\alpha\mathbf{d}\delta.
\end{equation}
Since $\alpha^{2}=0$, (\ref{eq:idab}) follows from (\ref{eq:clev9a1}).
 The proof of our proposition is complete. 
\end{proof}
\subsection{A closed $1$-form on the space of exact differentials}%
\label{subsec:clof}
Let $\alpha$ be a homotopy for $\delta$, that depends smoothly on 
$\delta$. Namely $\alpha$ is a smooth section of 
$\End^{1}\left(E\right)$, and 
\begin{equation}\label{eq:comi1}
\left[\alpha,\delta\right]=1.
\end{equation}
Observe  that we do not assume  $\alpha$ to be a differential.

Let $\tau$ denote the involution defining the $\Z_2$-grading of $E$. If $A\in\End\left(E\right)$, let 
\index[not]{Trs@$\Trs$}%
$\Trs\left[A\right]$ denote the supertrace
\index[terms]{supertrace}%
of $A$, i.e., 
\index[not]{Trs@$\Trs$}%
\begin{equation}\label{eq:comi1z1}
\Trs\left[A\right]=\Tr\left[\tau A\right].
\end{equation}
Then $\Trs$ extends to a map from $\Lambda\left(T^{*}\mathcal{D}_{\mathrm{e}}^{-}\right)\ho\End\left(E\right)$  into $\Lambda\left(T^{*}\mathcal{D}^{-}_{\mathrm{e}}\right)$ so that if $\beta\in \Lambda\left(T^{*}\mathcal{D}_{\mathrm{e}}^{-}\right),A\in\End\left(E\right)$,
\begin{equation}\label{eq:comilz2}
\Trs\left[\beta A\right]=\beta\Trs\left[A\right].
\end{equation}
Then $\Trs$ is a supertrace, in the sense that it vanishes on supercommutators.

\begin{definition}\label{def:form}
	Let 
	\index[not]{k@$\kappa$}%
	$\kappa$ be the $1$-form on $\mathcal{D}^{-}_{\mathrm{e}}$, 
	\begin{equation}\label{eq:clev6}
\kappa=\Trs\left[\alpha\mathbf{d}\delta\right].
\end{equation}
\end{definition}
\begin{proposition}\label{prop:clo}
	The $1$-form $\kappa$ does not depend on the choice of 
	$\alpha$, and it is closed. If 
	$f\in T^{*}\mathcal{D}_{\mathrm{e}}^{-}\ho\End^{0}\left(E\right)$, then
	\begin{equation}\label{eq:gap1}
\Trs\left[\alpha\left[f,\delta\right]\right]=\Trs\left[f\right].
\end{equation}
\end{proposition}
\begin{proof}
	Let $\alpha'$ have the same properties as $\alpha$, so that if 
	$\beta=\alpha'-\alpha$,
	\begin{equation}\label{eq:comi2}
\left[\delta,\beta\right]=0.
\end{equation}
Then 
\begin{equation}\label{eq:comi3}
\beta=\left[\delta,\alpha\beta\right].
\end{equation}
By (\ref{eq:comi3}), we obtain
\begin{equation}\label{eq:comi4}
\Trs\left[\beta\mathbf{d}\delta\right]=\Trs\left[\left[\delta,\alpha\beta\right]\mathbf{d}\delta\right].
\end{equation}
Using (\ref{eq:clev4}), (\ref{eq:comi4}), since supertraces vanish on 
supercommutators, we get
\begin{equation}\label{eq:comi5}
\Trs\left[\beta\mathbf{d}\delta\right]=\Trs\left[\left[\delta,\alpha\beta\mathbf{d}\delta\right]\right]=0.
\end{equation}
Equation (\ref{eq:comi5}) shows that our form $\kappa$ does not 
depend on the choice of $\alpha$.

We can now choose $\alpha\in\mathcal{D}_{\mathrm{e}}$ as in Subsection \ref{subsec:siid}. Clearly,
	\begin{equation}\label{eq:clev7}
\mathbf{d}\kappa=\Trs\left[\mathbf{d}\alpha\mathbf{d}\delta\right].
\end{equation}
Using (\ref{eq:idab}), (\ref{eq:clev7})  and the fact that supertraces vanish on 
supercommutators, we get
\begin{equation}\label{eq:clev7a1}
\mathbf{d}\kappa=0,
\end{equation}
i.e., $\kappa$ is closed. Take $f$ as indicated in the proposition. 
Then
	\begin{equation}\label{eq:zim1}
\alpha\left[f,\delta\right]=-\left[\delta,\alpha f\right]+f.
\end{equation}
Since $\Trs$ vanishes on supercommutators,  by (\ref{eq:zim1}), we 
get (\ref{eq:gap1}). The proof of our proposition is complete. 
\end{proof}
\subsection{A simple identity}%
\label{subsec:sipid}
Let $\delta\in \mathcal{D}_{\mathrm{e}}^{-}$. Let $\alpha$ be a homotopy with 
respect to $\delta$. 
If $\beta$ is such that $\left[\delta,\beta\right]=0$, then
\begin{equation}\label{eq:sat1a1}
\beta=\left[\delta,\alpha\beta\right].
\end{equation}

Recall that $N$ is the number operator of $E$.
Clearly,
	\begin{equation}\label{eq:ida1}
\left[\delta,N\right]=\delta.
\end{equation}
In particular,
\begin{equation}\label{eq:ida1a}
\left[\delta,N-\alpha\delta\right]=0.
\end{equation}
By (\ref{eq:sat1a1}), (\ref{eq:ida1a}), we get
\begin{equation}\label{eq:ida1z1}
N-\alpha\delta=\left[\delta,\alpha\left(N-\alpha\delta\right)\right].
\end{equation}
\begin{proposition}\label{prop:ida}
The following identity holds:
\begin{equation}\label{eq:ida2}
N-\alpha\delta=\left[\delta,\alpha N\right].	
\end{equation}
\end{proposition}
\begin{proof}
	This is a consequence of (\ref{eq:comi1}), (\ref{eq:ida1}).
\end{proof}
\begin{remark}\label{rem:ident}
	If $\alpha^{2}=0$, i.e., if $\alpha\in \mathcal{D}_{\mathrm{e}}$, then 
	(\ref{eq:ida1z1}) and (\ref{eq:ida2}) are equivalent. In the case where $\alpha$ is just a homotopy, 
	\begin{equation}\label{eq:ida2a1}
\left[\delta,\alpha^{2}\delta\right]=0,
\end{equation}
which reflects the fact that
\begin{equation}\label{eq:ida2a2}
\alpha^{2}\delta=\left[\delta,\alpha^{3}\delta\right].
\end{equation}
From  (\ref{eq:ida1z1}), (\ref{eq:ida2a1}), we recover (\ref{eq:ida2}).
\end{remark}
\subsection{Determinants}%
\label{subsec:det}
Here, we follow the conventions of Knudsen-Mumford 
\cite{KnudsenMumford}. In particular we will consider the category 
$\mathcal{L}$ of 
complex lines equipped with a sign in $\Z_{2}$, the morphisms being  even isomorphisms. This category is 
equipped with a tensor product functor, so that if 
$\left(\lambda,\epsilon_{\lambda}\right),\left(\mu,\epsilon_{\mu}\right)$ are objects in $\mathcal{L}$, the canonical morphism mapping $\lambda \otimes \mu$ to $\mu\otimes \lambda$ being $\left(-1\right)^{\epsilon_{\lambda}
\epsilon_{\mu}}$.
This will be written in the form
\begin{equation}\label{eq:zebr2}
\lambda \otimes 
\mu=\left(-1\right)^{\epsilon_{\lambda}\epsilon_{\mu}}\mu \otimes \lambda.
\end{equation}

If $\lambda$ is an objet in $\mathcal{L}$, $\lambda^{-1}$ denotes its 
right inverse in the category $\mathcal{L}$, so that $\lambda \otimes \lambda^{-1}=\C$. More precisely, if $\left(\lambda,\epsilon_{\lambda}\right)$ is the full object in $\mathcal{L}$, its inverse is given by $\left(\lambda^{*},\epsilon_{\lambda}\right)$.  If $s\in \lambda$ is nonzero, $s^{-1}\in \lambda^{-1}$ denotes the right inverse of $s$, so that $s \otimes s^{-1}=1$. 

Let $\mathcal{K}$ be the category of complex vector spaces. The 
determinant functor $\det:\mathcal{K}\to\mathcal{L}$ is defined to be
\begin{equation}\label{eq:zebr3}
\det F=\left(\Lambda^{\max}F,\left[\dim F\right]\right),
\end{equation}
where $\left[\dim F\right]$ is the image of $\dim F$ in $\Z_{2}$.

If $F,F'$ are objects in $\mathcal{K}$, we have the canonical 
isomorphism,
\begin{equation}\label{eq:zebr4}
\det\left(F \oplus F'\right)=\det F \otimes \det F'.
\end{equation}
More generally, if
\begin{equation}\label{eq:zebr4a1}
0\to F\to F'\to F''\to 0
\end{equation}
 is an exact sequence in $\mathcal{K}$, we have a canonical 
 isomorphism,
 \begin{equation}\label{eq:zebr4a2}
\det F' =\det F \otimes \det F''.
\end{equation}

\subsection{The canonical section $\tau\left(d\right)$}%
\label{subsec:exad}
Now we use the notation of Subsection \ref{subsec:twdif}.
Put
\begin{align}\label{eq:dif2a1}
&\chi=\sum_{i=p}^{q}\left(-1\right)^{i}\dim 
E^{i},&\chi'=\sum_{i=p}^{q}\left(-1\right)^{i}i \dim E^{i}.
\end{align}

 Observe that 
\begin{align}\label{eq:chao}
&\chi=\dim E\,\mathrm{mod}\, 2,
& \chi'=\frac{1}{2}\left(\chi-\dim E\right)\,\mathrm{mod}\, 2.
\end{align}

Set
\begin{equation}\label{eq:zebra0}
\det E=\bigotimes_{i=p}^{q}\left(\det 
E^{i}\right)^{\left(-1\right)^{i}}.
\end{equation}
The product in the right-hand side is an ordered product from left to right, that starts 
at $i=p$.

If $d\in \mathcal{D}_{\mathrm{e}}$, then
\begin{equation}\label{eq:rong1-a0}
\chi=0,
\end{equation}
so that $\dim E$ is even, and $\det E$ is even in the sense of \cite{KnudsenMumford}.
There is an associated canonical nonzero 
section $\tau\left(d\right)\in \det E$. Let us recall its 
construction. For $p\le i\le q$, let $\sigma^{i}$ be a nonzero 
section of $\Lambda^{\max}\left(E^{i}/dE^{i-1}\right)$.  Then
\begin{equation}\label{eq:dif2a2}
\tau\left(d\right)=\bigotimes_{i=p}^{q}\left(d\sigma^{i-1}\we\sigma^{i}\right)^{\left(-1\right)^{i}}.
\end{equation}
As the notation indicates, the right-hand side does not depend on the choices.\footnote{\label{foot:order} When $i$ is odd, the proper interpretation of $\left(d\sigma^{i-1}\we \sigma^{i} \right) ^{-1}$ is as follows. If $L$ is an object in $\mathcal{L}$, and if $L^{-1}$ is its right inverse, when $s\in L$ is nonzero, let $s^{-1}_{\mathrm{r}},s^{-1}_{\mathrm{l}}$ denote the right and left inverses of $s$, i.e.,   the elements of $L^{-1}$ such that $ss^{-1}_{\mathrm{r}}=1,s^{-1}_{\mathrm{l}}s=1$. By convention, in (\ref{eq:dif2a2}), $\left(d\sigma^{i-1}\we \sigma^{i}\right)^{-1}=\left(d\sigma^{i-1}\right)^{-1}_{\mathrm{r}}\we \left( \sigma^{i} \right) ^{-1}_{\mathrm{l}}$.}

	Let $E'=\oplus_{i=p}^{q}E^{\prime i}$ be another $\Z$-graded complex. 
	Put $E''=E \oplus E'$. Then
	\begin{equation}\label{eq:rada0}
\det E''=\det E \otimes \det E'.
\end{equation}

By the existence and uniqueness of the determinant functor established by Knudsen-Mumford \cite[Theorem 1]{KnudsenMumford}, the sections 
$\tau\left(d\right)$ are multiplicative with respect to exact 
sequences of exact complexes.  
More precisely, assume that  $\left(E',d'\right)$ 
is exact. Then $E''=E 
\oplus E'$ is equipped with the exact differential $d''=d+d'$.
Also  $\det E''= \det E \otimes \det E'$ 
	is so defined that 
	\begin{equation}\label{eq:rada1}
\tau\left(d''\right)= \tau\left(d\right) \otimes 
\tau\left(d'\right).
\end{equation}

\subsection{The action of $\mathrm{Aut}^{0}\left(E\right)$ on $\det E$}%
\label{subsec:acau}
If $g\in \Aut^{0}\left(E\right)$, put
\begin{equation}\label{eq:co1x1}
	\det g\vert_{E}=\prod_{i=p}^{q}\det g\vert_{E^{i}}^{\left(-1\right)^{i}}.
\end{equation}
Then the action of $g$ on $\det E$ is given by
\begin{equation}\label{eq:colx1a1}
g\vert_{\det E}=\det g\vert_{E}.
\end{equation}
Also observe that $\Aut^{0}\left(E\right)$ acts on $\mathcal{D}_{\mathrm{e}}$ by conjugation.
\begin{proposition}\label{prop:pcoa1}
	If $d\in \mathcal{D}_{\mathrm{e}}$, then
	\begin{equation}\label{eq:coa2}
\tau\left(g. d\right)=\det g\vert_{E}\tau\left(d\right).
\end{equation}
In particular, if $g\in \mathrm{Aut}^{0}_{d}\left(E\right)$, then $\det g\vert_{E}=1$.

If $a\in \C^{*}$, then
\begin{equation}\label{eq:dif2a3}
\tau\left(ad\right)=a^{\chi'}\tau\left(d\right).
\end{equation}
\end{proposition}
\begin{proof}
	The proof of the first part is left to the reader. If $a\in 
	\C^{*}$, set $g=a^{N}$. Then $g.d=ad$. Equation 
	(\ref{eq:dif2a3}) is  a consequence of (\ref{eq:coa2}). The proof of our proposition is complete. 
\end{proof}
\begin{proposition}\label{prop:code}
If $d\in \mathcal{D}$,  if $g\in \mathrm{Aut}^{0}_{d}\left(E\right)$, then $g$ acts on $H\left(E,d\right)$. Moreover,
\begin{equation}\label{eq:code1}
\det g\vert_{E}=\det g\vert_{H\left(E,d\right)}.
\end{equation}
\end{proposition}
\begin{proof}
The first part of our proposition is obvious. Also we have the exact sequences,
\begin{align}\label{eq:code2}
&0\to H^{p}\left(E,d\right)\to E^{p}\to v\left(E^{p}\right)\to 0,   \notag \\
&0\to v\left(E^{p}\right)\to \ker v\vert_{E^{p+1}}\to H^{p+1}\left(E,d\right)\to  0,  \\
&0\to \ker v\vert_{E^{p+1}}\to E^{p+1}\to v\left(E^{p+1}\right)\to 0 ,\ldots \notag 
\end{align}
Moreover, $g$ acts on as an automorphism of the exact complexes in (\ref{eq:code2}). Using Proposition  \ref{prop:pcoa1}, our proposition follows.
\end{proof}
\subsection{The case where $d,\delta$ are exact}%
\label{subsec:exap}
Let $d\in \mathcal{D}, \delta\in \mathcal{D}^{-}$. Then 
$\left[d,\delta\right]\in\End^{0}\left(E\right)$, and
\begin{align}\label{eq:dif3}
&\left[d,\left[d,\delta\right]\right]=0,&\left[\delta,\left[d,\delta\right]\right]=0.
\end{align}
\begin{proposition}\label{prop:inv1}
The endomorphism $\left[d,\delta\right]$ is invertible if and only if 
$d\in \mathcal{D}_{\mathrm{e}},\delta\in \mathcal{D}^{-}_{\mathrm{e}}$, and $\mathrm{Im} d\cap 
\mathrm{Im} \delta=0$. Under these conditions, 
$E=\mathrm{Im}d \oplus \mathrm{Im}\delta$.
\end{proposition}
\begin{proof}
	Assume that $\left[d,\delta\right]$ is invertible. Put
	\begin{equation}\label{eq:dif4}
T=\left[d,\delta\right]^{-1}\delta.
\end{equation}
Then $T$ is a homotopy for $d$, i.e., (\ref{eq:dif2}) holds, so that 
$d$ is exact. The same argument also shows that $\delta$ is exact.
	 Under the same condition, 
	$\mathrm{Im}d +\mathrm{Im}\delta=E,\mathrm{Im}d\cap 
	\mathrm{Im}\delta=0$, so that $E=\mathrm{Im}d \oplus \mathrm{Im}\delta$.
	
	Let us now prove our theorem in the opposite direction, i.e., we 
	assume that the conclusions of our theorem do hold. If $f\in 
	E,\left[d,\delta\right]f=0$, then
	\begin{align}\label{eq:dif5}
&d\delta f=0,&\delta df=0.
\end{align}
Since $d\in \mathcal{D}_{\mathrm{e}},\delta\in \mathcal{D}^{-}_{\mathrm{e}}$, from 
(\ref{eq:dif5}), we get
\begin{align}\label{eq:dif6}
&df=0,&\delta f=0.
\end{align}
The same argument  shows that $f\in \mathrm{Im} d\cap 
\mathrm{Im}\delta$, and so $f=0$. Therefore $\left[d,\delta\right]$ 
is invertible, which completes the proof of our 
proposition.
\end{proof}

Under the  conditions of Proposition  \ref{prop:inv1}, $d$ restricts to a one to one map 
 $\mathrm{Im}\delta\to\mathrm{Im}d$, and  $\delta$ restricts to a one 
 to one map 
 $\Im d\to \mathrm{Im}\delta$.
 
 We will now use the previous construction to define 
 $\tau\left(\delta\right)\in\det E$. Let $\underline{E}$ be the 
 $\Z$-graded vector space $E$ in which the degree in $E^{i}$ is now 
 $-i$, so that $\delta$ acts as a differential of degree $1$ on 
 $\underline{E}$. Let $\underline{\chi}, \underline{\chi}'$ denote the analogues of $\chi, \chi'$ for $\underline{E}$. Then
 \begin{align}\label{eq:tsar1}
 &\underline{\chi}=\chi,&\underline{\chi}'=-\chi'.
 \end{align}

 Observe that
 \begin{equation}\label{eq:zebr0}
\det \underline{E}=\bigotimes_{i=q}^{p}\left(\det 
E^{i}\right)^{\left(-1\right)^{i}},
\end{equation}
the order being now decreasing from $q$ to $p$.

We claim that there  is a natural canonical anti-isomorphism $a$ from $\det 
\underline{E}$ to $\det E$.
First,  we consider  the case where $F$ is an object in 
$\mathcal{K}$. If $m=\dim F$, the action of $a$ on $\det F$ is just 
multiplication  by $\left(-1\right)^{\left(m-1\right)m/2}$. Then $a$ 
extends to an anti-automorphism from  $\det \underline{E}$ to $\det E$, so that $a$ should be thought of as just reversing the order globally.
 
 If $\delta\in\mathcal{D}^{-}_{\mathrm{e}}$, we denote by  $\underline{\tau}\left(\delta\right)$ 
 the nonzero section of $\det \underline{E}$ which is the analogue of $\tau\left(d\right)\in \det E$. Put
 \begin{equation}\label{eq:zebr-1}
\tau\left(\delta\right)=\left(-1\right)^{\dim 
E/2}\underline{\tau}\left(\delta\right).
\end{equation}
Equivalently,
\begin{equation}\label{eq:zebr-2}
\tau\left(\delta\right)=a\underline{\tau}\left(\delta\right).
\end{equation}

 Using obvious notation, we get
 \begin{equation}\label{eq:zebr-3}
\tau\left(\delta\right)=\bigotimes_{i=p}^{q}\left(\rho^{i}\we\delta \rho^{i+1}\right)^{\left(-1\right)^{i}}.
\end{equation}
In (\ref{eq:zebr-3}), the order is now increasing order from $p$ to 
$q$.\footnote{Using the notation in Footnote \ref{foot:order}, when $i$ is odd, the proper interpretation of $\left(\rho^{i}\we \delta
\rho^{i}\right)^{-1}$ is $\left(\rho^{i}\right)^{-1}_{\mathrm{r}}\we \left(\delta\rho^{i+1}\right)^{-1}_{\mathrm{l}}$. Equation (\ref{eq:zebr-3}) is ultimately true because $a$ exchanges right and left inverses.}

Since $\det E$ is an even line, right inverses  and left inverses of sections of $\det E$ are identical.
\begin{theorem}\label{thm:tidfo}
	Assume that $\left[d,\delta\right]$ is invertible. The following identities hold:
	\begin{align}\label{eq:dif7}
&\prod_{i=p}^{q}\left( \det\left[d,\delta\right]\vert_{E^{i}}\right) 
^{\left(-1\right)^{i}}=1,\\
&\tau\left(d\right) 
\otimes \tau\left(\delta\right)^{-1}=\prod_{i=p}^{q}\left( \det\left[d,\delta\right]\vert_{E^{i}}\right) 
^{\left(-1\right)^{i}i}. \nonumber 
\end{align}
In particular, if $\left[d,\delta\right]=1$, then
\begin{equation}\label{eq:dif7a1}
\tau\left(d\right) \otimes 
\tau\left(\delta\right)^{-1}=1.
\end{equation}
\end{theorem}
\begin{proof}
	If $g=\left[d,\delta\right]$, then $g\in \Aut^{0}\left(E\right)$, 
	and $g.d=d$. Using Proposition \ref{prop:pcoa1}, we get the first 
	identity  in (\ref{eq:dif7}).  

	In (\ref{eq:zebr-3}), by Proposition \ref{prop:inv1}, we can take $\rho^{i}$ nonzero in 
	$\Lambda^{\max}\left( dE^{i-1} \right) $. In  
	(\ref{eq:dif2a2}), we take
	\begin{equation}\label{eq:zebr-7}
\sigma^{i}=\delta \rho^{i+1}.
\end{equation}
By  (\ref{eq:dif2a2}), (\ref{eq:zebr-3}), we obtain
\begin{equation}\label{eq:zebr-8}
\tau\left(d\right) \otimes 
\tau\left(\delta\right)^{-1}=\prod_{i=p+1}^{q}\left[\det
d\delta\vert_{dE^{i-1}}\right]^{\left(-1\right)^{i}}.
\end{equation}
One verifies easily that
\begin{equation}\label{eq:zebr-9}
\prod_{i=p}^{q}\det\left[d,\delta\right]\vert_{E^{i}}^{\left(-1\right)^{i}i}=\prod_{i=p+1}^{q}\left[\det d\delta\vert_{dE^{i-1}}\right]^{\left(-1\right)^{i}}.
\end{equation}
By (\ref{eq:zebr-8}), (\ref{eq:zebr-9}), we get the second identity  in (\ref{eq:dif7}). The proof of our theorem is complete.
\end{proof}

\subsection{The canonical connection on $\det E$ and the form 
$\kappa$}%
\label{subsec:canco}
In the sequel, we view $\det E$ as a trivial line bundle over 
$\mathcal{D}^{-}_{\mathrm{e}}$. We still denote by $\mathbf{d}$ the 
corresponding trivial connection on $\det E$. 
\begin{proposition}\label{prop:conn}
	The following identity holds over $\mathcal{D}^{-}_{\mathrm{e}}$:
	\begin{equation}\label{eq:clev15}
\mathbf{d}\tau\left(\delta\right)=\kappa\tau\left(\delta\right).
\end{equation}
Equivalently, $\kappa$ is the connection form for $\mathbf{d}$ 
associated with the section $\tau\left(\delta\right)$.
Also $\mathrm{Re}\kappa$ is an exact $1$-form, and the periods of the 
closed $1$-form
$\mathrm{Im}\kappa/2\pi$ are integral. 
\end{proposition}
\begin{proof}
	We fix $\delta_{0}\in \mathcal{D}^{-}_{\mathrm{e}}$. Let $\alpha_{0}\in 
	\mathcal{D}_{\mathrm{e}}$, such that 
	$\left[\alpha_{0},\delta_{0}\right]=1$. For $\delta$ close to 
	$\delta_{0}$, $\left[\alpha_{0},\delta\right]$ is invertible. Using 
	(\ref{eq:dif7}), we get
	\begin{equation}\label{eq:clev16}
\tau\left(\alpha_{0}\right) \otimes 
\tau\left(\delta\right)^{-1}=\prod_{i=p}^{q}\left(\det\left[\alpha_{0},\delta\right]\vert_{E^{i}}\right)^{\left(-1\right)^{i}i}.
\end{equation}
By (\ref{eq:clev16}), we get
\begin{equation}\label{eq:clev17}
-\frac{\mathbf{d}\tau\left(\delta\right)}{\tau\left(\delta\right)}=-\Trs\left[N\left[\alpha_{0},\delta\right]^{-1}\left[\alpha_{0},\mathbf{d}\delta\right]\right].
\end{equation}
The minus sign in the right-hand side of  (\ref{eq:clev17}) comes 
from the fact that $\alpha_{0}$ is viewed as odd. By (\ref{eq:clev17}), we 
deduce that
\begin{equation}\label{eq:dev18}
\frac{\mathbf{d}\tau\left(\delta\right)}{\tau\left(\delta\right)}\vert_{\delta_{0}}=\Trs\left[N\left[\alpha_{0},\mathbf{d}\delta\right]\right]\vert_{\delta_{0}}.
\end{equation}
Since $\left[N,\alpha_{0}\right]=\alpha_{0}$, using the fact that 
supertraces vanish on supercommutators, from (\ref{eq:dev18}), we get
\begin{equation}\label{eq:dev19}
\frac{\mathbf{d}\tau\left(\delta\right)}{\tau\left(\delta\right)}\vert_{\delta_{0}}=\Trs\left[\alpha_{0}\mathbf{d}\delta\right]\vert_{\delta_{0}},
\end{equation}
from which the first part of our proposition follows. Since $\mathbf{d}$ is a flat connection on $\det E$, we 
recover the fact that $\kappa$ is a closed $1$-form. 

Let $h^{E}=\bigoplus_{i=p}^{q}h^{E^{i}}$ be a 
Hermitian metric on $E$, and let $\left\Vert  \,\right\Vert_{\det E}$ 
be the corresponding metric on $\det E$. Then
\begin{equation}\label{eq:dev20}
d\log\left\Vert  
\tau\left(\delta\right)\right\Vert^{2}=2\mathrm{Re}\kappa,
\end{equation}
which shows that $\mathrm{Re}\kappa$ is exact. By 
(\ref{eq:clev15}), we also obtain the last statement of our 
proposition, of which the proof is complete.\end{proof}
\begin{remark}\label{rem:clap}
	We use the notation of Proposition \ref{prop:clo}. The connection 
	$\n^{E}=\mathbf{d}+f$ on $E$ preserves the 
	$\Z_{2}$-grading of $E$. It induces a corresponding connection 
	$\mathbf{\n}^{\det E}$ on $\det E$, and $\n^{\End\left(E\right)}$ on $\End\left(E\right)$. By Proposition 
	\ref{prop:clo}, we deduce that
	\begin{equation}\label{eq:clev14a2}
\mathbf{\n}^{\det 
E}\tau\left(\delta\right)=\left(\mathbf{\kappa}+\Trs\left[f\right]\right)\tau\left(\delta\right).
\end{equation}
By (\ref{eq:clev6}), (\ref{eq:gap1}),  we can rewrite (\ref{eq:clev14a2}) in the form
\begin{equation}\label{eq:clev14a3}
\mathbf{\n}^{\det 
E}\tau\left(\delta\right)=\Trs\left[\alpha\mathbf{\n}^{\End\left(E\right)}\delta\right]\tau\left(\delta\right).
\end{equation}
\end{remark}

\subsection{Hermitian metrics: the case where $d,\delta$ are exact}%
\label{subsec:hefo}
In this Subsection, we assume that $d\in 
\mathcal{D}_{\mathrm{e}},\delta\in\mathcal{D}^{-}_{\mathrm{e}}$.

Let $h^{E}=\bigoplus_{i=p}^{q}h^{E^{i}}$ be a Hermitian form on 
$E=\bigoplus_{i=p}^{q}E^{i}$. Let $d^{*}$ denote the adjoint of $d$ with respect to 
$h^{E}$. Then $d^{*}\in \mathcal{D}^{-}_{\mathrm{e}}$. 

Let $T$ be a homotopy for $d$. 
Then $T^{*}$ is a homotopy for
 $d^{*}$.
 
Now we assume that $h^{E}$ 
 is positive. Then $\det E$ inherits a metric $\left\Vert  
 \,\right\Vert_{\det E}$. Also $\left[d,d^{*}\right]$ is invertible, and  the 
 conclusions of Theorem \ref{thm:tidfo} hold.
 \begin{proposition}\label{prop:nor}
 	The following identity holds:
	\begin{multline}\label{eq:dif13}
\tau\left(d\right) \otimes 
\tau\left(d^{*}\right)^{-1}=\left\Vert  \tau\left(d\right)\right\Vert^{2}_{\det E}=\left\Vert  
\tau\left(d^{*}\right)^{-1}\right\Vert^{2}_{\left( \det 
E \right)  ^{-1}}\\
=
\prod_{i=p}^{q}\left(\det\left[d,d^{*}\right]\vert_{E^{i}}\right)^{\left(-1\right)^{i}i}.
\end{multline}
If $g\in\Aut^{0}\left(E\right)$, 
then 
\begin{equation}\label{eq:dif13a1}
\prod_{i=p}^{q}\left(\det\left[g.d,\left( g.d \right) 
^{*}\right]\vert_{E^{i}}\right)^{\left(-1\right)^{i}i}=
\left\vert  \det g\vert _{E}\right\vert^{2}\prod_{i=p}^{q}\left(\det\left[d,d^{*}\right]\vert_{E^{i}}\right)^{\left(-1\right)^{i}i}.
\end{equation}
 \end{proposition}
 \begin{proof}
 	The first identity is a  trivial consequence of Theorem 
	\ref{thm:tidfo}. Combining (\ref{eq:coa2}) and (\ref{eq:dif13}), 
	we get (\ref{eq:dif13a1}). The proof of our proposition is complete. 
 \end{proof}
 
\begin{theorem}\label{thm:pidb}
	Assume that $\left[d,\delta\right]$ is invertible. The following identities hold:
	\begin{align}\label{eq:dif15}
&\tau\left(d\right) \otimes 
\tau\left(\delta\right)^{-1}=\prod_{i=p}^{q}\left(\det\left[d,\delta\right]\vert_{E^{i}}\right)^{\left(-1\right)^{i}i},\\
&\frac{\left\Vert  \tau\left(d\right)\right\Vert_{\det E}}{\left\Vert 
\tau\left(\delta\right) \right\Vert_{\det 
E}}=\left\vert  
\prod_{i=p}^{q}\left(\det\left[d,\delta\right]\vert_{E^{i}}\right)^{\left(-1\right)^{i}i}\right\vert. \nonumber 
\end{align}
In particular,
\begin{multline}\label{eq:dif16}
\prod_{i=p}^{q}\left(\det  \left[d,d^{*}\right]\vert_{E^{i}} \right) ^{\left(-1\right)^{i}i}
\prod_{i=p}^{q}\left(  \det\left[\delta,\delta^{*}\right]\vert_{E^{i}} \right) ^{\left(-1\right)^{i}i}\\
=\left\vert\prod_{i=p}^{q}\left(\det\left[d,\delta\right]\vert_{E^{i}}\right)^{\left(-1\right)^{i}i}\right\vert^{2}.
\end{multline}
\end{theorem}
\begin{proof}
	The first equation in (\ref{eq:dif15}) was already established in 
	Theorem \ref{thm:tidfo}. The second equation is an obvious 
	consequence of the first one. By the first identity in 
	(\ref{eq:dif13}), we get
	\begin{equation}\label{eq:zebr-9a1}
\left\Vert  \tau\left(\delta\right)\right\Vert^{2}_{\det 
E}=\prod_{i=-q}^{-p}\left( \det\left[\delta,\delta^{*}\right]\vert_{\underline{E}^{i}} \right) ^{\left(-1\right)^{i}i}.
\end{equation}
Equation (\ref{eq:zebr-9a1}) can be rewritten in  the form,
\begin{equation}\label{eq:zebr-10}
\left\Vert  \tau\left(\delta\right)\right\Vert^{2}_{\det 
E}=\prod_{i=p}^{q}\left( \det\left[\delta,\delta^{*}\right]\vert_{E^{i}} \right) ^{\left(-1\right)^{i-1}i},
\end{equation}
so that
\begin{equation}\label{eq:zebr-11}
\left\Vert  \tau\left(\delta\right)^{-1}\right\Vert^{2}_{ \left(\det 
E\right)^{-1}
}=\prod_{i=p}^{q}\left( \det\left[\delta,\delta^{*}\right]\vert_{E_{i}} \right) ^{\left(-1\right)^{i}i}.
\end{equation}

	Using (\ref{eq:dif13}), (\ref{eq:dif15}), and (\ref{eq:zebr-11}),  we get (\ref{eq:dif16}). The proof of our theorem is complete. 
\end{proof}
\subsection{The case where only $\delta$ is exact}%
\label{subsec:deex}
Here, we only assume that $\delta\in \mathcal{D}^{-}_{\mathrm{e}}$. In particular, 
$\left[d,\delta\right]$ is no longer assumed to be invertible. Also 
$\tau\left(\delta\right)$ is a nonzero section 
of $\det E$. 

Recall that we have  the canonical isomorphism,
\begin{equation}\label{eq:lab1}
\rho:\det H\left(E,d\right) \simeq \det E.
\end{equation}
When $H\left(E,d\right)=0$, this only says that $\tau\left(d\right)$ is a 
canonical nonzero section of $\det E$.

The nonzero section 
$\tau\left(\delta\right)$ of $\det E$ induces a 
corresponding nonzero section of $\det H\left(E,d\right)$.

For $a>0$, set
\index[not]{Da@$D_{a}$}%
\index[not]{Sa@$S_{a}$}%
\begin{align}\label{eq:disk1}
&D_{a}=\left\{z\in\C,\left\vert z\right\vert\le a\right\},&S_{a}=\left\{z\in\C,\left\vert z\right\vert=a\right\}.
\end{align}

Let $a>0$ be such that  no eigenvalue $\lambda$ of 
$\left[d,\delta\right]$ is such that $\lambda\in S_{a}$. Let $E_{<a},E_{>a}$ be the direct sum of the 
characteristic subspaces of $\left[d,\delta\right]$  associated with eigenvalues $\lambda$ such 
that $\left\vert \lambda  \right\vert<a,\left\vert  
\lambda\right\vert>a$. Then
\begin{equation}\label{eq:lab2}
E=E_{<a} \oplus E_{>a}.
\end{equation}
The splitting (\ref{eq:lab2}) is preserved by $d,\delta$.

By (\ref{eq:lab2}), we get
\begin{equation}\label{eq:lab3}
\det E=\det E_{<a}\otimes \det E_{>a}.
\end{equation}

Since $\left[d,\delta\right]$ is invertible on $E_{>a}$, the complex 
$\left(E_{>a},\delta\vert _{E>a}\right)$ is exact, so that the 
section 
$\tau\left(\delta\vert_{E_{>a}}\right)\in\det E_{>a}$ is well 
defined.
\begin{proposition}\label{prop:exa}
	The complex $\left(E_{<a},\delta\vert_{E_{<a}}\right)$ is exact.
\end{proposition}
\begin{proof}
	 Let $k$ 
	be a homotopy for $\delta$, so that $\left[k,\delta\right]=1$. Let 
	$P_{<a}$ denote the spectral projection on $E_{<a}$ with respect 
	to $\left[d,\delta\right]$. Then 
	$P_{<a}$ commutes with $\delta$. Put
	\begin{equation}\label{eq:zomb1}
k_{<a}=P_{<a}kP_{<a}.
\end{equation}
Then $k_{<a}$ acts as a morphism of degree $1$ on $E_{<a}$. Moreover,
\begin{equation}\label{eq:zomb2}
\left[\delta\vert_{E<a},k_{<a}\right]=1\vert_{E_{<a}}.
\end{equation}
Therefore $k_{<a}$ is a homotopy for $\delta\vert_{E_{<a}}$, which 
gives our proposition.
	\end{proof}

By the above, the  nonzero section 
$\tau\left(\delta\vert_{E<a}\right)\in\det E_{<a}$ is also 
well-defined. Moreover, we have 
\begin{equation}\label{eq:lab4}
\tau\left(\delta\right)=\tau\left(\delta\vert_{E<a}\right) \otimes 
\tau\left(\delta\vert_{E>a}\right).
\end{equation}

Since $\left(E_{>a},d_{>a}\right)$ is exact, we   have the canonical 
isomorphism,
\begin{equation}\label{eq:lab5}
\rho_{<a}:\det H\left(E,d\right) \simeq \det E_{<a}.
\end{equation}

From the above, it follows that the canonical isomorphism $\rho$ in 
(\ref{eq:lab1}) is given by
\begin{equation}\label{eq:lab6}
\rho=\rho_{<a} \otimes \tau\left(d\vert_{E_{>a}}\right).
\end{equation}
By (\ref{eq:lab6}), if $h\in \det H\left(E,d\right)$, then
\begin{equation}\label{eq:lab7}
\rho h=\rho_{<a}h \otimes \tau\left(d\vert_{E_{>a}}\right).
\end{equation}
By (\ref{eq:lab4}), (\ref{eq:lab7}), we have the identity of  complex numbers,
\begin{equation}\label{eq:lab8}
\frac{\rho h}{\tau\left(\delta\right)}
=\frac{\rho_{<a} 
h}{\tau\left(\delta\vert_{E_{<a}}\right)}  
\frac{\tau\left(d\vert_{E_{>a}}\right)}{\tau\left(\delta\vert_{E_{>a}}\right)}.
\end{equation}

Using Theorem \ref{thm:tidfo}, we can rewrite (\ref{eq:lab8}) in the 
form
\begin{equation}\label{eq:lab9}
\frac{\rho h}{\tau\left(\delta\right)}
=\frac{\rho_{<a} 
h}{\tau\left(\delta\vert_{E_{<a}}\right)} 
\prod_{i=p}^{q}\left( \det\left[d,\delta\right]\vert_{E^{i}_{>a}} \right) ^{\left(-1\right)^{i}i}.
\end{equation}
Needless to say, it can be shown directly that the right-hand side of 
(\ref{eq:lab9}) does not depend on $a$. In particular, if we need, we 
can take $a$ small enough so that $0$ is the only possible eigenvalue 
such that $\left\vert  \lambda\right\vert<a$.

Let $[ \,]_{\det E,\delta}$  be the linear 
form
   on $\det E$, such that  
\begin{equation}\label{eq:norm1}
[ \tau\left(\delta\right)]_{\det E,\delta}=  1.
\end{equation}
Similarly, we can define the linear form $[  \,]_{\det 
E_{<a},\delta}$ on $\det E_{<a}$. Let $\left\Vert  
\,\right\Vert_{\det E,\delta}$ be the norm on $\det E$ such that 
$\left\Vert  \tau\left(\delta\right)\right\Vert=1$. The norm 
$\left\Vert  \,\right\Vert_{\det E_{<a},\delta}$ on $\det E_{<a}$ is 
defined in the same way.
\begin{proposition}\label{prop:nori}
	If $h\in \det H\left(E,d\right)$, then
	\begin{equation}\label{eq:norm3}
[ \rho h]_{\det E,\delta}=[ 
\rho_{<a}h]_{\det E_{<a},\delta}  
\prod_{i=p}^{q}\left( \det\left[d,\delta\right]\vert_{E^{i}_{>a}} \right) ^{\left(-1\right)^{i}i}.
\end{equation}
\end{proposition}
\begin{proof}
	This is a trivial consequence of (\ref{eq:lab9}).
\end{proof}

Let $h^{E}$ be a Hermitian metric on $E$, and let $\left\Vert  
\,\right\Vert_{\det E}$ be the corresponding norm on $\det E$.  

By (\ref{eq:zebr-10}),  we obtain
\begin{equation}\label{eq:norm5}
\frac{\left\Vert  \,\right\Vert_{\det E}}{\left\Vert  
\,\right\Vert_{\det 
E,\delta}}=\left[\prod_{i=p}^{q}\left( \det\left[\delta,\delta^{*}\right]\vert_{E^{i}} \right)^{\left(-1\right)^{i}i}\right]^{-1/2}.
\end{equation}
\begin{theorem}\label{thm:idx}
	If $h\in \det H\left(E,d\right)$, the following identity holds:
	\begin{multline}\label{eq:norm6}
\left\Vert  \rho h\right\Vert_{\det E}=\left\Vert  
\rho_{<a}h\right\Vert_{\det E_{<a,\delta}}\left\vert 
\prod_{i=p}^{q}\left( \det\left[d,\delta\right]\vert_{E^{i}_{>a}}\right) ^{\left(-1\right)^{i}i} \right\vert\\
\left[\prod_{i=p}^{q}\left( \det\left[\delta,\delta^{*}\right]\vert_{E^{i}} \right) ^{\left(-1\right)^{i}i}\right]^{-1/2}.
\end{multline}
When $\left[d,\delta\right]$ is invertible, equation (\ref{eq:norm6}) reduces to 
(\ref{eq:dif16}).
\end{theorem}
\begin{proof}
	Equation (\ref{eq:norm6}) is a consequence of (\ref{eq:norm3}), (\ref{eq:norm5}). The last statement is trivial, which completes 
	the proof of our theorem.
\end{proof}
 \subsection{Duality}%
\label{subsec:dual}
We will now use the notation $\mathcal{D}_{E},\mathcal{D}^{-}_{E}$ 
instead of $\mathcal{D},\mathcal{D}^{-}$.

If $F$ is a vector space, let $F^{*}$ denote its dual. Then 
\begin{equation}\label{eq:sta1}
\det F^{*}=\left(\det F\right)^{*}.
\end{equation}

For $-q\le i\le -p$, put
\begin{equation}\label{eq:sta1a1}
E^{*i}=E^{-i*}.
\end{equation}
If $E^{*}=\bigoplus_{i=-q}^{-p}E^{*i}$, using the same notation as in (\ref{eq:zebra0}), we have
\begin{equation}\label{eq:sta1a2}
\det E^{*}=\bigotimes_{i=-q}^{-p}\left( \det E^{*i} \right) ^{\left(-1\right)^{i}}.
\end{equation}

As we saw in Subsection \ref{subsec:det}, $ \left( \det E^{i}\right) ^{\left(-1\right)^{i}} \otimes \left( \det E^{*-i} \right) ^{\left(-1\right)^{i}}$  has a canonical nonzero section so that $ \left(  \det E^{i} \right) ^{\left(-1\right)^{i}} \otimes \left( \det E^{*-i} \right) ^{\left(-1\right)^{i}} \simeq \C$. Then  we have canonical identifications,
\begin{align}\label{eq:sta4}
&\det E \otimes  \det E^{*} \simeq \C,&\det E^{*} \otimes \det E \simeq \C.
\end{align}

Let $\widetilde{d}$ be the transpose of $d$. Then $\widetilde{d}\in \mathcal{D}_{E^{*}}$. If $d\in \mathcal{D}_{E,\mathrm{e}}$, then $\widetilde{d}\in \mathcal{D}_{E^{*},\mathrm{e}}$.
 In that case,  $\tau\left(d\right),\tau\left(\widetilde{d}\right)$ are nonzero sections of $\det E,\det E^{*}$.  In this case, $\det E,\det E^{*}$ are even lines. 
\begin{proposition}\label{prop:dua}
	If $d\in \mathcal{D}_{E,\mathrm{e}}$, then
	\begin{align}\label{eq:dua2}
&\tau\left(d\right) \otimes \tau\left(\widetilde{d}\right)
\simeq 1,&\tau\left(\widetilde{d}\right) \otimes \tau\left(d\right)
\simeq 1.
\end{align}
\end{proposition}
\begin{proof}
	It is enough to consider the case where the complex $\left(E,d\right)$ splits, in which case the proposition is trivial. 
\end{proof}
\subsection{Determinants and shifts}%
\label{subsec:eveca}
Let $\left(E',d'\right)$ be a complex similar to $\left(E,d\right)$. Let $f:\left(E,d\right)\to \left(E',d'\right)$ be an isomorphism of complexes.  Then $f$ is an invertible morphism of degree $0$, that exchanges $d$ and $d'$. Since $\det$ is a functor, there is an induced isomorphism of lines $\det f:\det E\to \det E'$. 

If $\left(E,d\right),\left(E',d'\right)$ are exact, then 
\begin{equation}\label{eq:sta5}
\det f\tau\left(d\right)=\tau\left(d'\right).
\end{equation}
By (\ref{eq:sta5}), we can view $\det f$ as a section of $\left(\det E\right)^{-1} \otimes \det E'$. Also observe that since the two complexes are exact, $\det E$ and $\det E'$ are even lines.  By (\ref{eq:sta5}), we have
\begin{equation}\label{eq:sta6x1}
\det f=\tau \left(d\right)^{-1} \tau \left(d'\right)\,\mathrm{in}\,\left(\det E\right)^{-1} \otimes \det E'.
\end{equation}

Given $r\in \Z$, let 
\index[not]{Er@$E_{r}$}%
$\left(E_{r},d\right)$ be the shifted complex by $r$, so that $E_{r}^{\Ou}=E^{\Ou+r}$. If $r$ is even, there is an obvious canonical isomorphism $\det E \simeq \det E_{r}$. Let $f:\left(E,d\right)\to \left(E'_{r},d\right)$ be an isomorphism of complexes. When $r$ is even, there is still a determinant morphism $\det f:\det E\to \det E'$. 

Let us now consider the case $r=1$. Note that
\begin{equation}\label{eq:stax1}
\det E_{1}=\bigotimes _{i=p}^{q}\left( \det E^{i} \right) ^{\left(-1\right)^{i-1}}.
\end{equation}
There is an obvious canonical isomorphism,
\begin{equation}\label{eq:stax2}
\det E \otimes \det E_{1} \simeq \C.
\end{equation}

It will be convenient to equip $E_{1}$ with the differential $-d$. 

When $d\in \mathcal{D}_{E,\mathrm{e}}$, then $-d\in \mathcal{D}_{E_{1},\mathrm{e}}$. In this case, we denote by $\tau_{E}\left(d\right)\in \det E,\tau_{E_{1}}\left(-d\right)\in \det E_{1}$ the canonical sections of $\det E,\det E_{1}$ associated with the corresponding differentials.
\begin{proposition}\label{prop:sta}
	If $d\in \mathcal{D}_{E,\mathrm{e}}$, then
	\begin{equation}\label{eq:stax3}
\tau_{E}\left(d\right)\tau_{E_{1}}\left(-d\right)=1.
	\end{equation}
\end{proposition}
\begin{proof}
First, we consider the case $p=0,q=1$. Using the notation in Footnote \ref{foot:order}, we get
\begin{align}\label{eq:stax4}
&\tau_{E}\left(d\right)=s_{0}\left(ds_{0}\right)^{-1}_{\mathrm{r}}, 
&\tau_{E_{1}}\left(-d\right)= \left( s_{0} \right)_{\mathrm{l}} ^{-1}\left(-ds_{0}\right).
\end{align}
Using (\ref{eq:stax4}), we get (\ref{eq:stax3}). In the general case, we exploit the fact expressions similar to the ones in (\ref{eq:stax4}) are even, and also the arguments above. The proof of our proposition is complete.  
\end{proof}

Let $f:\left(E,d\right)\to \left(E'_{1},d\right)$ be an isomorphism of complexes. In particular $f$ maps $E^{i}$ in $E^{ \prime i+1}$. Also $\det f$ maps $\det E$ in $\det E'_{1}$. If $d\in \mathcal{D}_{E,\mathrm{e}}$, then
\begin{equation}\label{eq:stax5}
	\tau_{E'_{1}}\left(d\right)= \det f\tau_{E}\left(d\right).
\end{equation}
As in (\ref{eq:sta6x1}), by (\ref{eq:stax5}), we can view $\det f$ as a section of $\left(\det E\right)^{-1}  \otimes \det E'_{1}$. 

The above immediately extends to the case where $r$ is odd, and $f$ is an isomorphism $\left(E,d\right)\to \left(E'_{r},d\right)$.
\subsection{The determinant of a linear map}%
\label{subsec:fix}
Let $F,F'$ be finite dimensional complex vector spaces of the same dimension.  Then $\det F,\det F'$ are lines that are either even, or odd. Let $ \gamma$ be a morphism from $F$ into $F'$. Then 
\index[not]{detg@$\det \gamma$}%
$\det \gamma  $ is a  section of $ \left(  \det F \right) ^{-1} \otimes \det F'$. Using the notation of Footnote \ref{foot:order}, if $ \sigma\in \det F $ is non zero, then
\begin{equation}\label{eq:fix1}
\det \gamma= \left( \sigma \right)  ^{-1}_{\mathrm{r}} \gamma \sigma =\gamma  \sigma  \left( \sigma  \right)^{-1}_{\mathrm{l}}. 
\end{equation}
Observe that $\det \gamma $ is a version of $ \tau\left(d\right) $ that were considered before.
With the above definition, one verifies that $\det$ is multiplicative with respect to the composition of the $
 \gamma$. Moreover, if $F,F'$ are replaced by $F \oplus G, F' \oplus G'$, and using obvious notation, 
\begin{equation}\label{eq:fix2z1}
\det \gamma_{F \oplus G}=\det \gamma_{F}  \det \gamma_{G}.
\end{equation}

Following \cite[p. 30]{KnudsenMumford}, if $\gamma $ is invertible,  we define 
$\det \gamma \in \left(\det F'\right)^{-1}\otimes \det F $, so that if $\sigma$ is nonzero in $\det F$, then
\begin{equation}\label{eq:fix2z2}
	\det \gamma =\left( \gamma \sigma \right)^{-1}_{\mathrm{r}}\sigma =\sigma \left( \gamma \sigma \right)^{-1}_{\mathrm{l}}. 
\end{equation}
Observe that (\ref{eq:fix1}) and (\ref{eq:fix2z2}) can be deduced from each other by replacing $ \gamma  $ by $ \gamma^{-1}$.

Let $f:\left(E,d\right)\to \left(E',d'\right)$ be an isomorphism of complexes. 
 We define 
 \index[not]{detf@$\det f$}%
 $\det f\in \left(\det E\right)^{-1} \otimes  \det E'$ to be the product of the sections $\det f\vert_{E^{i}}$ defined as before, by using the definition (\ref{eq:fix1}) for $i$ even, and (\ref{eq:fix2z2}) for $i$ odd. Since each section is even, the order in the product is irrelevant.

 Our definition of $\det f$ makes equation (\ref{eq:sta5}) tautological.

\subsection{A special case}%
\label{subsec:spe}
Here, we assume that $q-p$ is odd. Let $d,\delta$ be taken as before, and such that $\left[d,\delta\right]=1$, so that by Proposition \ref{prop:inv1}, $E=\Im\, d \oplus \Im\, \delta$. In particular, the complexes  $\left(E,d\right)$ and $\left(E,\delta\right)$ are exact. Put
\begin{align}\label{eq:chada1}
&A=\bigoplus_{p\le i\le \frac{p+q-3}{2}}E^{i} \bigoplus dE^{\frac{p+q-3}{2}}, \notag \\
&B=\delta E^{\frac{p+q+1}{2}} \oplus dE^{\frac{p+q-1}{2}},\\
&C= \delta E^{\frac{p+q+3}{2}} \oplus \bigoplus _{\frac{p+q+3}{2}\le i\le q}E^{i}. \notag 
\end{align}
Then $A,B,C$ are stable by $d,\delta$, the restrictions of $d,\delta$ to $A,B,C$ define exact subcomplexes. Finally,
\begin{equation}\label{eq:chada2}
E=A \oplus B \oplus C ,
\end{equation}
and the direct sum in (\ref{eq:chada2}) is  a direct sum of complexes with respect to  $d$ and also $\delta$. We can then define the canonical sections $\tau_{A}\left(d\right),\tau_{B}\left(d\right),\tau_{C}\left(d\right)$ of $\det A,\det B,\det C$, and corresponding sections where $d$ is replaced by $\delta$. 
\begin{proposition}\label{prop:The following identity holds}
	The following identities hold:
	\begin{align}\label{eq:chada3}
	& \tau\left(d\right)=\tau_{A}\left(d\right)\tau_{B}\left(d\right)\tau_{C}\left(d\right), 
	&\tau\left(\delta\right)=\tau_{A}\left(\delta\right)\tau_{B}\left(\delta\right)\tau_{C}\left(\delta\right).
	\end{align}
\end{proposition}
\begin{proof}
	This follows from  the multiplicativity of $\tau$.
\end{proof}

Let $\Gamma$ be an odd  isomorphism  that maps  $E^{i}$ into $E^{p+q-i}$ that is such that $\Gamma^{2}=1$. Also we assume that $ \Gamma$ induces an isomorphism of complexes $\left(E^{\Ou},d\right)\to \left(\underline{E}^{\Ou-p-q},\delta\right)$. This exactly says that
\begin{equation}\label{eq:say1}
\delta\Gamma=\Gamma d,
\end{equation}
which is the same as 
\begin{equation}\label{eq:say2}
\Gamma\delta=d\Gamma.
\end{equation}

Then $\Gamma$ induces  isomorphisms of complexes $$\left(A^{\Ou},d\right) \to \left(\underline{C}^{\Ou-p-q},\delta\right),\left(B^{\Ou},
d\right) \to \left(\underline{B}^{\Ou-p-q},\delta\right),\left(C^{\Ou},d\right) \to \left(\underline{A}^{\Ou-p-q},\delta\right).$$

We make the assumption that
\begin{equation}\label{eq:chada4}
\Gamma\vert_{B}=d+\delta.
\end{equation}
Since $\left[d,\delta\right]=1$,  (\ref{eq:chada4}) is compatible with the fact that $\Gamma^{2}=1$.


Put
\begin{align}\label{eq:chada5a}
&E_{-}=\bigoplus_{p\le i\le \frac{p+q-1}{2}}E^{i}, &E_{+}=\bigoplus_{\frac{p+q+1}{2}\le i\le q}E^{i}.
\end{align}
Let $\Gamma_{-+}$ be the restriction of $\Gamma$ to  $E_{-}$. Then $\Gamma_{-+}$ is a morphism of degree $0$ from  $ \left( \underline{E}^{\Ou-p-q}_{-},\delta \right) $ to $ \left( E^{\Ou}_{+},d \right) $.

Now we follow the conventions of Subsection \ref{subsec:fix}.
\begin{definition}\label{def:defg}
	Let
	\index[not]{rG@$\rho_{\Gamma}$}%
	$ \rho _{\Gamma}\in \det E$ be given by 
	\begin{equation}\label{eq:chada5b1}
	\rho_{\Gamma}= \left( \det \Gamma _{-+} \right) ^{-1}.
	\end{equation}
\end{definition}

Let $\Gamma_{\underline{A},C}$ be the restriction of $\Gamma$ to $\underline{A}$. Then $\Gamma_{\underline{A},C}$ is a morphism of degree $0$ from  $\left(\underline{A}^{\Ou-p-q},\delta\right)$ to $\left(C^{\Ou},d\right)$. Let $\det\Gamma\vert_{\delta E^{\frac{p+q+1}{2}},dE^{\frac{p+q-1}{2}}}$ be the determinant of the restriction of $\Gamma$ to $\delta E^{\frac{p+q+1}{2}}$ taken in the sense of (\ref{eq:fix1}), (\ref{eq:fix2z2}).


Using (\ref{eq:fix2z1}), we get 
\begin{equation}\label{eq:chada5}
\rho_{\Gamma}= \left[\det \Gamma_{\underline{A},C}  \otimes
 \det\Gamma\vert_{\delta E^{\frac{p+q+1}{2}},dE^{\frac{p+q-1}{2}}}\right]^{-1}.
\end{equation}

Since $\left(A,d\right)$ is exact, $\dim A$ is even, so that when changing $\Gamma$ to $-\Gamma$, the first factor in the right-hand side of (\ref{eq:chada5}) is unchanged. As to the second factor, it is  fixed once and for in (\ref{eq:chada4}).
\begin{theorem}\label{thm:idt}
The following identity holds:
\begin{equation}\label{eq:chada6}
\tau\left(d\right)=\tau\left(\delta\right)=\rho_{\Gamma} \,\mathrm{in}\, \det E.
\end{equation}
\end{theorem}
\begin{proof}
Since $\left[d,\delta\right]=1$, the first equation was already established in Theorem \ref{thm:tidfo}. 

Using (\ref{eq:stax5}), we get
\begin{equation}\label{eq:chada4a}
	\tau_{C^{\Ou}}\left(d\right)=\det \Gamma_{\underline{A},C}\tau_{\underline{A}^{\Ou-p-q}}\left(\delta\right).
\end{equation}
By (\ref{eq:stax3}), we can identify $\tau_{\underline{A}^{\Ou-p-q}}\left(\delta\right)$ with $\tau_{\underline{A}}\left(-\delta\right)^{-1}$, so that (\ref{eq:chada4a}) can be rewritten in the form
\begin{equation}\label{eq:chada4b}
	\tau_{C}\left(d\right)=\det\Gamma_{\underline{A},C} \tau_{\underline{A}}\left(-\delta\right)^{-1}.
\end{equation}

By combining (\ref{eq:dif2a3}),  (\ref{eq:zebr-1}), (\ref{eq:dif7a1}), and (\ref{eq:chada4b}), we get
\begin{equation}\label{eq:chada4c}
	\tau_{C}\left(d\right)=\det\Gamma_{\underline{A},C} \left(-1\right)^{\dim A/2+\chi'\left(A\right)}\tau_{A}\left(d\right)^{-1}.
\end{equation}
Using (\ref{eq:chao}), we can rewrite (\ref{eq:chada4c}) in the form
\begin{equation}\label{eq:chada4d}
	\tau_{C}\left(d\right)=\det\Gamma_{\underline{A},C} \tau_{A}\left(d\right)^{-1}.
\end{equation}

By (\ref{eq:chada4}), we get
\begin{equation}\label{eq:chada12}
\tau_{B}\left(d\right)=\left( \det\Gamma\vert_{\delta E^{\frac{p+q+1}{2}},dE^{\frac{p+q-1}{2}}} \right) ^{-1}.
\end{equation}
By   (\ref{eq:chada3}), (\ref{eq:chada5}), (\ref{eq:chada4d}),   and (\ref{eq:chada12}), we get (\ref{eq:chada6}). The proof of our theorem is complete. 
\end{proof}
\section{Anosov vector fields and the Fried zeta function}%
\label{sec:frizet}
In this Section,  we consider an Anosov vector field on  a compact 
connected manifold $Y$, also equipped with a flat vector bundle $F$. 
We introduce the associated Fried zeta function 
$R_{Z,F}\left(\sigma\right)$ as an infinite product over primitive 
 closed trajectories of $-Z$. Also, we derive the main 
properties of that function, in connection with Poincaré duality.

This Section is organized as follows. In Subsection 
\ref{subsec:anvf}, we introduce the Anosov vector field $Z$, and  the associated
unstable and stable vector bundles.

In Subsection \ref{subsec:ruze}, we introduce the Ruelle zeta function 
$R_{Z,F}\left(\sigma\right),\sigma\in\C,\mathrm{Re}\,\sigma\gg 1$. The 
orientation line bundle $\mathrm{o}\left(T_{s}Y\right)$ appears 
explicitly in the definition.

In Subsection \ref{subsec:fripo}, we prove that the Fried zeta 
function is compatible with a version of Poincaré duality.

Finally, in Subsection \ref{subsec:invo}, we consider the case where 
an involution $\iota$ of $Y$ maps $Z$ to $-Z$.
\subsection{An Anosov vector field}%
\label{subsec:anvf}
Let $Y$ be a  smooth compact connected manifold of dimension $n$. Let $Z\in C^{\infty}(Y, TY)$ be 
an Anosov vector field  on $Y$.  Let 
\index[not]{ft@$\varphi_{t}$}%
$\varphi_{t}\vert_{t\in\R}$ be the group of diffeomorphisms associated 
with $Z$, which is obtained by integration of the differential 
equation
\begin{align}\label{eq:flo1}
&\dot y=Z,&y_{0}=y.
\end{align}
The vector field $Z$ does not vanish.

Let 
\index[not]{ToY@$T_{0}Y$}%
$T_{0}Y\subset TY$ be the real line 
 spanned by $Z$. Let 
 \index[not]{TuY@$T_{u}Y$}%
 \index[not]{TsY@$T_{s}Y$}%
	\begin{align}\label{eq:duaz0}
		TY=T_{0}Y\oplus T_{u}Y\oplus T_{s}Y
	\end{align} 
be the   Anosov decomposition of $TY$. 
In (\ref{eq:duaz0}),  $T_{u}Y,T_{s}Y$ are Hölder vector subbundles of 
$TY$. They are called the
\index[terms]{unstable}%
\index[terms]{stable}%
unstable and stable vector bundles respectively. The splitting 
(\ref{eq:duaz0}) is preserved by the flow $\varphi_{\cdot}$.

Set
\index[not]{nu@$n_{u}$}%
\index[not]{ns@$n_{s}$}%
\begin{align}\label{eq:dormi1}
	&n_{u}=\dim T_{u}Y,&n_{s}=\dim T_{s}Y.
\end{align}
Then
\begin{equation}\label{eq:dormi2}
n=n_{u}+n_{s}+1.
\end{equation}

If $U\in T_{u}Y\setminus 0,V\in 
T_{s}Y\setminus 0$, as $t\to + 
\infty $, 
\begin{align}\label{eq:duaz3b1}
&\left\vert  \varphi_{t*}U\right\vert\to + \infty,
&\left\vert \varphi_{t*}V\right\vert\to  0.
\end{align}

Let $T^{*}_{0}Y,T^{*}_{u}Y,T^{*}_{s}Y$ be dual to $T_{0}Y,T_{u}Y,T_{s}Y$.\footnote{In
\cite{FaureSjostrand11,DyatlovZworski16}, the notation for 
$T^{*}_{u}Y,T^{*}_{s}Y$ are interchanged.} By (\ref{eq:duaz0}), we have a corresponding  splitting of $T^{*}Y$,
	\begin{align}\label{eq:dua2z1}
		T^{*}Y=T^{*}_{0}Y\oplus T^{*}_{u}Y\oplus T^{*}_{s}Y.
	\end{align}

Let $T_{0}^{\perp}Y \subset 
T^{*}Y$ be the orthogonal bundle to $T_{0}Y$. Then $T_{0}^{\perp}Y$ 
is a smooth vector subbundle of $T^{*}Y$. By (\ref{eq:duaz0}), 
(\ref{eq:dua2z1}), 
\begin{equation}\label{eq:dua2z1x1}
T_{0}^{\perp}Y=T_{u}^{*}Y\oplus T^{*}_{s}Y.
\end{equation}
so that 
$T_{u}^{*}Y,T_{s}^{*}Y$ can be viewed as Hölder subbundles of 
$T_{0}^{\perp}Y$.
The splitting in (\ref{eq:dua2z1}) is also preserved by 
$\varphi_{\cdot}$. 

Let $\alpha$ be the  section of $ T^{*}Y$ such that 
\begin{align}\label{eq:duaz2}
&i_{Z}\alpha=1,&\alpha\vert_{T_{u}Y \oplus T_{s}Y}=0.
\end{align}
Then $\alpha$ can be viewed as the canonical section of $T^{*}_{0}Y$ 
that is dual to $Z$. As a section of $T^{*}Y$, $\alpha$ is Hölder. The section $\alpha$ is preserved by the flow $\varphi_{\cdot}$.

In the sequel, we will use the notation,
\index[not]{M@$M$}%
\begin{equation}\label{eq:duaz3a1}
M=T_{u}Y \oplus T_{s}Y.
\end{equation}
Then $M$ is a Hölder subbundle of $TY$, and it is isomorphic to the 
smooth bundles $TY/T_{0}Y$ and  $\left( T_{0}^{\perp}Y \right) 
^{*}$.

By (\ref{eq:duaz3a1}), we get
\begin{equation}\label{eq:gonf1}
M^{*}=T^{*}_{u}Y \oplus T^{*}_{s}Y.
\end{equation}
Then $M^{*}$ is isomorphic to the smooth vector bundle 
$T^{\perp}_{0}Y$ in (\ref{eq:dua2z1x1}).

Let $\mathrm{o}\left(TY\right)$ be the orientation bundle of $TY$. 
Then $\mathrm{o}\left(TY\right)$ is a $\Z_{2}$-bundle 
equipped with a canonical flat connection. Similarly, although 
$T_{u}Y,T_{s}Y$ are only Hölder vector bundles, their orientation  bundles 
$\mathrm{o}\left(T_{u}Y\right),\mathrm{o}\left(T_{s}Y\right)$ are 
$\Z_{2}$-bundles, and so they are equipped with canonical flat 
connections. By (\ref{eq:duaz0}), we get\footnote{These line bundles should also be viewed as $\Z_2$-graded. In particular, exchanging the order in (\ref{eq:gonf1}) introduces a sign similar to the sign in (\ref{eq:zebr2}).}
\begin{equation}\label{eq:dormi3}
\mathrm{o}\left(TY\right)=\mathrm{o}\left(T_{u}Y\right) \otimes 
\mathrm{o}\left(T_{s}Y\right).
\end{equation}
 After tensoring with $\R$ or $\C$, we obtain this way 
smooth flat line bundles. 

If $Y$ is orientable, then $\mathrm{o}\left(TY\right)$ has two global 
sections. In this case, $\varphi_{t}\vert_{t\in\R}$ acts 
trivially on $\mathrm{o}\left(TY\right)$, i.e., it preserves these 
global sections. Indeed this is obvious for small $\left\vert  t\right\vert$ and extends 
immediately to arbitrary $t$.


Let $\left(F,\n^{F}\right)$ be a  complex flat vector bundle on $Y$.  
Let $d^{Y}$ be the corresponding de Rham operator acting on  $\Omega
\left(Y,F\right)$, where
\begin{equation}\label{eq:tsar2}
\Omega\left(Y,F\right)=C^{\infty}\left(Y,\Lambda\left(T^{*}Y\right)\otimes_{\R}F\right).
\end{equation}

Set
	\begin{align}\label{eq:normi4}
		L_{Z}=\[d^{Y}, i_{Z}\]. 
	\end{align}
Then, $L_{Z}$ is a first order differential operator acting on 
$ \Omega \left( Y, F \right)$, and commuting with $d^{Y},i_{Z}$.  The 
group of diffeomorphisms 
$\varphi_{t}\vert_{t\in \R}$ of $Y$ associated with the vector field 
$Z$ lifts to a corresponding flow also acting on the vector bundle $F$ 
and preserving the flat connection $\n^{F}$. 

The above constructions can also be used when replacing $F$ by $F 
\otimes _{\Z_{2}}\mathrm{o}\left(TY\right),F 
\otimes _{\Z_{2}}\mathrm{o}\left(T_{u}Y\right),F 
\otimes _{\Z_{2}}\mathrm{o}\left(T_{s}Y\right)$. 
Since $L_{Z}$ acts on $C^{\infty 
}\left(Y,TY\right)$, there is an induced action on 
$\mathrm{o}\left(TY\right)$, which coincides with the definition 
coming from (\ref{eq:normi4}). The same is true for 
$\mathrm{o}\left(T_{u}Y\right),\mathrm{o}\left(T_{s}Y\right)$.
\subsection{The Fried zeta function}%
\label{subsec:ruze}
We will say that $y\in Y$ is a non-trivial periodic point  of 
$\varphi_{-t}\vert_{t>0}$ if there exists 
$t>0$ such that
\begin{equation}\label{eq:flo2}
\varphi_{-t}\left(y\right)=y.
\end{equation}
Let 
\index[not]{F@$\mathbf{F}$}%
$\mathbf{F}$ be the set of non-trivial periodic points. Then $\mathbf{F}$ is preserved 
by the flow $\varphi_{\cdot}$. Given $y\in \mathbf{F}$, the set of 
$t\in\R$ such that $\varphi_{-t}\left(y\right)=y$ forms a discrete 
subgroup $\Gamma_{y}$ of $\R$.  Let 
\index[not]{ty@$t_{y}$}%
$t_{y}>0$ be the generator of 
this discrete group. If $y\in\mathbf{F}$, we obtain this way a closed non 
self-intersecting trajectory of period $t_{y}$. Clearly, $t_{y}$ is 
invariant under the action of $\varphi_{\cdot}$.

If $y\in \mathbf{F},t\in \Gamma_{y}, t>0$, $\varphi_{-t*}=d\varphi_{-t}$ acts on 
$M_{y}$, and because of (\ref{eq:duaz3b1}), $1$ is not an eigenvalue 
of $\varphi_{-t*}$.
\begin{definition}\label{def:flo3}
	If $y\in\mathbf{F}, t\in \Gamma_{y}, t>0$, put
	\begin{align}\label{eq:flo5}
&\epsilon_{y,t}=\mathrm{sgn}\det\left(1-\varphi_{-t*}\right) 
\vert_{M_{y}},\nonumber \\
&\epsilon_{u,y,t}=\mathrm{sgn}\det\left(1-\varphi_{-t*}\right) 
\vert_{T_{u,y}Y},\\
&\epsilon_{s,y,t}=\mathrm{sgn}\det\left(1-\varphi_{-t*}\right) 
\vert_{T_{s,y}Y}. \nonumber 
\end{align}
\end{definition}
\begin{proposition}\label{prop:sign1}
	The following identity holds:
	\begin{align}\label{eq:flo3a1}
&\epsilon_{u,y,t}=1, 
&\epsilon_{y,t}=\epsilon_{s,y,t}=\left(-1\right)^{n_{s}}\mathrm{sgn}\det 
\varphi_{-t*}\vert_{T_{s,y}Y}.
\end{align}
\end{proposition}
\begin{proof}
	By (\ref{eq:duaz3b1}),  if $\lambda$ is an eigenvalue 
	of $\varphi_{-t*}$ on $T_{u,y}Y$, then $\left\vert  
	\lambda\right\vert<1$, so that the first equation in 
	(\ref{eq:flo3a1}) holds. Similarly, when replacing $T_{u,y}Y$ by 
	$T_{s,y}Y$, then $\left\vert  \lambda\right\vert>1$. Only the real 
	eigenvalues contribute to the sign $\epsilon_{s,y,t}$. If 
	$\lambda\in\R$, and $\left\vert  \lambda\right\vert>1$, then 
	\begin{equation}\label{eq:flip1}
\mathrm{sgn}\left(1-\lambda\right)=-\mathrm{sgn}\lambda.
\end{equation}
By (\ref{eq:flip1}), we get the second equation in (\ref{eq:flo3a1}). 
The proof of our proposition is complete. 
\end{proof}

In the sequel, we identify $y,y'\in \mathbf{F}$ if there is $t\in \R$ such 
that $y'=\varphi_{t}\left(y\right)$. We denote by
\index[not]{F@$\underline{\mathbf{F}}$}%
$\underline{\mathbf{F}}$ the 
corresponding quotient of $\mathbf{F}$. The corresponding periodic trajectory for $-Z$ of period $t_y$ is said to be a 
\index[terms]{primitive}%
primitive closed trajectory. If $\gamma$ is a primitive trajectory for $-Z$, let $\gamma^{-1}$ be the corresponding primitive trajectory for $Z$.

We recall a result given in  Dyatlov-Zworski \cite[Lemma 2.2]{DyatlovZworski16}, which is a mild form of a theorem of Margulis \cite[Theorem 1.1, p.78]{Margulis04}. The proof is based on simple estimates on  the growth of volumes.
\begin{proposition}\label{prop:gro}
Given $T>0$, the set $\left\{y\in\underline{\mathbf{F}},t_{y}\le T\right\}$ is finite. Moreover, there exist $c>0,C>0$ such that for $T>0$,
\begin{equation}\label{eq:gro0}
	\left\vert \left\{y\in\underline{\mathbf{F}},t_{y}\le T\right\} \right\vert\le Ce^{cT}.
\end{equation}
\end{proposition}
 
To define our dynamical zeta function,  we  follow the approach of Fried \cite[Section 5]{Fried87}.
\begin{definition}\label{def:friz1}
	For $\sigma\in\C,\mathrm{Re}\,\sigma\gg 1$, we define the 
	\index[terms]{Fried zeta function}%
	Fried zeta function,
	\index[not]{RZF@$R_{Z,F}\left(\sigma\right)$}%
	\begin{equation}\label{eq:flip2}
R_{Z,F}\left(\sigma\right)=\prod_{y\in\underline{\mathbf{F}}}^{}\left[\det\left(1-e^{-t_{y}\sigma}\varphi_{t_{y}*}\right)\vert_{ \left( F \otimes_{\Z_{2}}\mathrm{o}\left(T_{s}Y\right) \right) _{y}}\right]^{\left(-1\right)^{n_{s}+1}}.
\end{equation}
\end{definition}

The convergence of the infinite product in (\ref{eq:flip2}) for $\Re\,\sigma\gg 1$ follows immediately from Proposition \ref{prop:gro} and  from the fact there exist $c'>0,C'>0$ such that 
\begin{equation}\label{eq:gro2}
\left\vert \varphi_{t*}\vert_{F}\right\vert \le C'e^{c't}.
\end{equation}
Also for $\Re\,\sigma\gg 1$, $R_{Z,F}\left(\sigma\right)$ is a holomorphic function of $\sigma$. 

If the connection $\n^{F}$ is unitary, we get
\begin{equation}\label{eq:flip2a1}
R_{Z,F^{*}}\left(\sigma\right)=\overline{R_{Z,F}\left(\overline{\sigma}\right)}.
\end{equation}
\subsection{The Fried zeta function and Poincaré duality}%
\label{subsec:fripo}
When replacing $Z$ by $-Z$, $\underline{\mathbf{F}}$ is unchanged. 
\begin{proposition}\label{prop:psimp}
	For $\sigma\in\C, \Re\,\sigma\gg 1$, the following identity holds:
	\begin{equation}\label{eq:flip3a-1}
R_{-Z,F^{*} \otimes 
_{\Z_{2}}\mathrm{o}\left(TY\right)}\left(\sigma\right)=\prod_{y\in\underline{\mathbf{F}}}^{}\left[\det\left(1-e^{-t_{y}\sigma}\varphi_{-t_{y}*}\right)\vert_{ \left( F ^{*}\otimes_{\Z_{2}} \mathrm{o}\left(T_{s}Y\right) \right) _{y}}\right]^{\left(-1\right)^{n_{u}+1}}.
\end{equation}
\end{proposition}
\begin{proof}
We just use equation (\ref{eq:flip2}) for 
$R_{Z,F}\left(\sigma\right)$, while replacing  $F$ by $F^{*} \otimes 
_{\Z_{2}}\mathrm{o}\left(TY\right)$, and $\varphi_{\cdot}$ by 
$\varphi_{-\cdot}$, while noting that $\underline{\mathbf{F}}$ and 
the $t_{y}>0$ are unchanged. In principle 
$\mathrm{o}\left(T_{u}Y\right)$ should appear in the right-hand side, 
but the fact that originally we had introduced 
$\mathrm{o}\left(TY\right)$ makes that ultimately
$\mathrm{o}\left(T_{s}Y\right)$ appears in the right-hand side  
of (\ref{eq:flip3a-1}). The proof of our proposition is complete. 	
\end{proof}
\begin{theorem}\label{thm:poin}
	For $\sigma\in\C,\mathrm{Re}\,\sigma\gg 1$, the following identity 
	holds:
	\begin{equation}\label{eq:flip3}
R_{Z,F}\left(\sigma\right)=\left[R_{-Z,F^{*} \otimes 
_{\Z_{2}}\mathrm{o}\left(TY\right)}\left(\sigma\right)\right]^{\left(-1\right)^{n-1}}.
\end{equation}
\end{theorem}
\begin{proof}
	Clearly,
	\begin{equation}\label{eq:flip4}
\det\left(1-e^{-t_{y}\sigma}\varphi_{-t_{y}*}\right)\vert_{ \left( F 
^{*}\otimes_{\Z_{2}} \mathrm{o}\left(T_{s}Y\right) \right) _{y}}=\det\left(1-e^{-t_{y}\sigma}\varphi_{t_{y}*}\right)\vert_{ \left( F\otimes_{\Z_{2}} \mathrm{o}\left(T_{s}Y\right) \right) _{y}}
\end{equation}
Using (\ref{eq:dormi2}), (\ref{eq:flip2}), (\ref{eq:flip3a-1}), and  (\ref{eq:flip4}), we get 
(\ref{eq:flip3}). The proof of our theorem is complete.
\end{proof}
\begin{remark}\label{rem:myst}
	The title of this Subsection contained the mysterious `Poincaré 
	duality'. This title will be fully explained in Subsection 
	\ref{subsec:poinsec}. 
\end{remark}
\subsection{The case where an involution acts on $Y$}%
\label{subsec:invo}
Assume that 
\index[not]{i@$\iota$}%
$\iota$ is an involution of $Y$  that lifts to $\left(F,\n^{F}\right)$, that is such that 
\begin{equation}\label{eq:flip5}
\iota_{*}Z=-Z.
\end{equation}
In particular $\iota$ exchanges $T_{u}Y$ and $T_{s}Y$, so that
\begin{equation}\label{eq:flip6}
n_{u}=n_{s},
\end{equation}
 and $n$ is odd. By (\ref{eq:flip5}), $\iota$ does not have fixed 
 points.
 
 Let $\gamma$ be a primitive closed trajectory of $-Z$. Then 
 $\iota\gamma$ is  a primitive closed trajectory of $Z$.
 \begin{proposition}\label{prop:inte1}
 	The trajectories $\gamma$ and $\iota\gamma$ do not intersect.
 \end{proposition}
 \begin{proof}
 If $\gamma$ and 
 $\iota\gamma$ were intersecting, at the intersection point, up to 
 sign, their speeds would coincide.  This means that the unoriented 
 $\gamma$ is preserved by $\iota$. Since $\iota$ changes $Z$ to $-Z$, 
 this forces $\iota$ to have a fixed point on $\gamma$, which is 
 impossible. The proof of our proposition is complete. 
 \end{proof}
 
 Therefore, the primitive closed trajectories of $-Z$ always come by 
 distinct pairs $\gamma,\left[\iota\gamma\right]^{-1}$. 
 \begin{proposition}\label{prop:flip}
 	If $Y$ is orientable and $\n^{F}$ is unitary, for 
	$\sigma\in\C,\mathrm{Re}\,\sigma\gg 1$, then
	\begin{equation}\label{eq:flip7}
R_{Z,F}\left(\sigma\right)=\overline{R_{Z,F}\left(\overline{\sigma}\right)}.
\end{equation}
Under the same assumptions, if   $\sigma\in\R,\sigma\gg 1$, then
\begin{equation}\label{eq:flip7a1}
R_{Z,F}\left(\sigma\right)>0.
\end{equation}
 \end{proposition}
 \begin{proof}
 	Clearly,
	\begin{equation}\label{eq:flip8}
R_{-Z,F^{*}}\left(\sigma\right)=R_{Z,F^{*}}\left(\sigma\right).
\end{equation}
By combining (\ref{eq:flip2a1}), (\ref{eq:flip3}), and 
(\ref{eq:flip8}), we get (\ref{eq:flip7}). Since 
$\mathrm{o}\left(TY\right)=1$,  then
$\mathrm{o}\left(T_{u}Y\right)=\mathrm{o}\left(T_{s}Y\right)$. If 
$\gamma$ is a primitive closed trajectory of $-Z$, and
\index{tg@$t_{\gamma}$}
$t_{\gamma}$ is its period,  then
\begin{equation}\label{eq:flip8a1}
\overline{\det\left(1-e^{-t_{\gamma}\overline{\sigma}}\varphi_{t_{\gamma}}\right)\vert_{F \otimes 
_{\Z_{2}}\mathrm{o}\left(T_{s}Y\right)}}
=
\det\left(1-e^{-t_{\left(\iota\gamma \right) 
^{-1}}\sigma}\varphi_{t_{\left(\iota\gamma\right)^{-1}}}\right)\vert_{F \otimes 
_{\Z_{2}}\mathrm{o}\left(T_{s}Y\right)}.
\end{equation}
By (\ref{eq:flip8a1}), for $\sigma\in\R$, the 
contributions to $R_{Z,F}\left(\sigma\right)$ of $\gamma$ and 
$\left(\iota\gamma\right)^{-1}$ are conjugate, from which 
(\ref{eq:flip7a1}) follows. The proof of our proposition is complete. 
\end{proof}
\section{Operators, currents, and intersection supertraces}%
\label{sec:opcur}
The purpose of this section is to develop the theory of currents on 
$Y\times Y$, of their corresponding intersection supertraces, and 
of the corresponding Poincaré duality. Here our currents are 
$\Z$-graded by their degree in $\Lambda\left(T^{*}Y\right)$. 

This section is organized as follows. In Subsection 
\ref{subsec:poin}, we construct a properly normalized version of 
Poincaré duality.

In Subsection \ref{subsec:semi}, given a current $T$ on $Y\times Y$, we 
define a corresponding operator $M_{T}$ acting on forms on $Y$ 
with kernel $M_T$. A companion operator $N_T$ is also defined.

In Subsection \ref{subsec:sclatr}, when $T$ has a well-defined 
restriction  to the diagonal, we define the intersection supertrace of the operators $M_T,N_T$.
 We show that this is
 a supertrace on a suitable algebra, i.e., it 
vanishes on supercommutators.

Finally, in Subsection \ref{subsec:adjo}, we construct the proper 
adjoint of the operators $M_{T},N_T$.
\subsection{Poincaré duality and integration along the fiber}%
\label{subsec:poin}
Let $Y$ be a compact manifold.
It will be useful to normalize Poincaré duality in the proper 
 way. Here, we follow \cite[Section 6.6]{Bismut10b} and 
 \cite[Sections 4.4 and 
 7.1]{BismutShenWei23}. If $f\in T^{*}Y$, set
 \index[not]{f@$\widetilde f$}%
 \begin{equation}\label{eq:normi1}
\widetilde f=-f.
\end{equation}
We denote by\, 
\index[not]{t@$\widetilde{}$}%
$\,\widetilde{}\,$\, the corresponding anti-automorphism of 
$\Lambda\left(T^{*}Y\right)$. Observe that if $\alpha\in 
\Lambda^{p}\left(T^{*}Y\right)$, then
\begin{equation}\label{eq:norm2a1}
\widetilde \alpha=\left(-1\right)^{p\left(p+1\right)/2}\alpha.
\end{equation}

 Let 
 \index[not]{DY@$\mathcal{D}\left(Y,F\right)$}%
 $\mathcal{D}\left(Y,F\right)$ be the vector space of 
 currents with coefficients in $F$. This is a $\Z_{2}$-graded vector space. Also $\Omega\left(Y,F\right)$ is dense in $\mathcal{D}\left(Y,F\right)$.   
 
 If 
$\alpha\in\Omega\left(Y,F\right),\beta\in\Omega\left(Y,F^{*} \otimes_{\Z_{2}}\mathrm{o}\left(TY\right)\right)$, set
\begin{equation}\label{eq:normi3}
\left\langle \beta,\alpha\right\rangle=\int_{Y}^{} 
\left\langle \widetilde \beta\we \alpha \right\rangle.
\end{equation}
Equation (\ref{eq:normi3}) extends to the case where $\beta\in\mathcal{D}\left(Y,F^{*} \otimes_{\Z_{2}}\mathrm{o}\left(TY\right)\right) $. No choice of a volume form or of an orientation of $Y$ has been made in (\ref{eq:normi3}). 

Then (\ref{eq:normi3}) is a non-degenerate bilinear pairing so that 
$\mathcal{D}\left(Y,F^{*} \otimes _{\Z_{2}}\mathrm{o}\left(TY\right)\right)$ can be identified  
with the topological dual to $\Omega\left(Y,F\right)$. Up to sign, 
it coincides with the traditional Poincaré duality, where $\widetilde 
\beta$ is replaced by $\beta$. When interchanging $F$ 
and $\Ft$, we get
\begin{equation}\label{eq:normio3a1}
\left\langle  
\beta,\alpha\right\rangle=\left(-1\right)^{n\left(n+1\right)/2}\left\langle  \alpha,\beta\right\rangle.
\end{equation}
The above bilinear pairing is not symmetric. However, the fact that 
the lack of symmetry is just expressed via a constant makes that the 
definition of transposes does not depend on the choice of one or the 
other pairings.

If $A$ is a bounded operator acting on $\Omega\left(Y,F\right)$, we 
denote by
\index[not]{As@$A^{\dag}$}%
$A^{\dag}$ the transpose of $A$ acting on 
$\mathcal{D}\left(Y,F^{*} \otimes 
_{\Z_{2}}\mathrm{o}\left(TY\right)\right)$.
If $f,U$ are smooth sections of 
$T^{*}Y, TY$, then
\begin{align}\label{eq:normi4a}
&f\we^{\dag}=-f\we, &i_{U}^{\dag}=i_{-U}.
\end{align}

If $\tau:M\to S$ is a proper submersion with compact oriented fiber 
$X$, if $\alpha\in \Omega\left(M,\R\right)$, the integral along 
the fiber $\pi_{*}\alpha$ is defined to be such that if 
$\beta\in\Omega\left(S,\R\right)$, then
\begin{equation}\label{eq:normi5}
\tau_{*}\left(\tau^{*}\beta\alpha\right)=\beta\tau_{*}\alpha.
\end{equation}
If $X$ is not oriented, we should take $\alpha\in\Omega\left(M,\R \otimes_{\Z_{2}}\mathrm{o}\left(TX\right)\right)$.
\subsection{Operators and currents}%
\label{subsec:semi}
Let 
\index[not]{DY@$\Delta^{Y}$}%
$\Delta^{Y}$ be the diagonal in $Y\times Y$, and let
\index[not]{i@$i$}%
$i:\Delta^{Y}\to Y\times Y$ be the obvious embedding.

Let $\pi_{1},\pi_{2}$ denote the projections from $Y\times Y$ to $Y$. 
If $0\le k\le 2n$, let $T\in \mathcal{D}^{k}\left(Y\times Y,\pi_{1}^{*}F \otimes 
\pi_{2}^{*}\left(F^{*} \otimes _{\Z_{2}}
\mathrm{o}\left(TY\right)\right)\right)$ be a current of 
degree $k$. Let
\index[not]{MT@$M_{T}$}%
$M_{T}:\Omega\left(Y,F\right)\to 
\mathcal{D}\left(Y,F\right)$ be such that if 
$s\in\Omega\left(Y,F\right)$, then 
\begin{equation}\label{eq:rec1}
M_{T}s=\pi_{1*}\left[ T\we \pi_{2}^{*}s\right].
\end{equation}
Implicitly, $F^{*}$ in $T$ is contracted with $F$ in $s$. The 
integral in (\ref{eq:rec1}) is unambiguously defined because of the 
presence of $\mathrm{o}\left(TY\right)$.
Then $M_{T}$ is an operator of degree $k-n$. In particular if $k=n$, $M_{T}$ 
is an operator of total degree $0$. 

 We 
have the obvious identity,
\begin{equation}\label{eq:rec1a0}
T\left(Y\times Y\right)=\pi_{1}^{*}TY \oplus \pi_{2}^{*}TY.
\end{equation}
Let $Z^{Y\times Y}$ be the vector field on $Y\times Y$  given by
\begin{equation}\label{eq:rec1a-1}
Z^{Y\times Y}=\pi_{1}^{*}Z \oplus \pi_{2}^{*}Z.
\end{equation}
\begin{proposition}\label{prop:comid}
	The following identities hold:
	\begin{align}\label{eq:comid1}
	&d^{Y}M_{T}=M_{d^{Y\times 
	Y}T}+\left(-1\right)^{\mathrm{deg}T}M_{T}d^{Y},\\
	&i_{Z}M_{T}=M_{i_{Z^{Y\times 
	Y}}T}+\left(-1\right)^{\mathrm{deg}T}M_{T}i_{Z}. \nonumber 
\end{align}
\end{proposition}
\begin{proof}
	The proof of these identities is elementary and is left to the 
	reader.
\end{proof}
 
Let 
\index[not]{N@$N$}%
$N$ be the number operator of $\Lambda\left(T^{*}Y\right)$. Then $N$  acts on $\Omega\left(Y,F\right)$ and on $\mathcal{D}\left(Y,F\right)$.

Put
\index[not]{NT@$N_{T}$}%
\begin{equation}\label{eq:prig1}
N_{T}=M_{T}\left(-1\right)^{nN}.
\end{equation}
Note that $\left(-1\right)^{nN}$ also acts on the current $T$ on the left 
and on the right. If  $T$ has degree $p,q$ on 
the first and second copy of $Y$, the action of 
$\left(-1\right)^{nN}$ on $T$ on the left is just 
$\left(-1\right)^{pn}$ and on the right is $\left(-1\right)^{qn}$.  
Equation (\ref{eq:prig1}) can be rewritten as 
\begin{equation}\label{eq:prig1aa}
N_{T}=M_{T\left(-1\right)^{n\left(n-N\right)}}.
\end{equation}
Therefore, the classes of the operators $M_{T}$ and $N_{T}$ are the same.

In the sequel, $\left[\ \right]$ still denotes the supercommutator.
\begin{proposition}\label{prop:scm}
	The following identities hold:
\begin{align}\label{eq:orig2}
&\left[d^{Y},N_{T}\right]=N_{d^{Y\times Y}T},&\left[i_{Z},N_{T}\right]=N_{i_{Z^{Y\times Y}}T}.
\end{align}	
\end{proposition}
\begin{proof}
By (\ref{eq:comid1}), we get
\begin{equation}\label{eq:prig3}
d^{Y}N_{T}=N_{d^{Y\times Y}T}+\left(-1\right)^{\deg T-n}N_{T}d^{Y}.
\end{equation}
Moreover, 
\begin{equation}\label{eq:prig4}
\deg N_{T}=\deg\, T-n.
\end{equation}
By (\ref{eq:prig3}), (\ref{eq:prig4}), we obtain the first identity 
in (\ref{eq:orig2}). The proof of the second identity is similar.
\end{proof}
\begin{remark}\label{rem:cur}
Equation (\ref{eq:orig2}) allows us to pass from the formalism of operators to the formalism of currents. The fact that $N_{T}$ supercommutes with $d^{Y}$ just says that the current $T$  is closed.
\end{remark}

If $T$ is taken as before, let 
\index[not]{Tpq@$T^{p,q}$}%
$T^{p,q}$ denote the component of $T$ which is of partial degree $p,q$ as a current in the first and second copies of $Y$, so that
\begin{equation}\label{eq:prig4a1}
T=\sum_{0\le p,q\le n}^{}T^{p,q}.
\end{equation}
 
If  $T\in \Omega^{n}\left(Y\times Y,\pi_{1}^{*}F \otimes 
\pi_{2}^{*}\left(F^{*} \otimes_{\Z_{2}} 
\mathrm{o}\left(TY\right)\right) \right) $, then $M_{T}$ has 
a smooth kernel, and so it is trace class. Also $M_{T}$ preserves the grading in
$\Omega\left(Y,F\right)$. Note that 
$\int_{\Delta^{Y}}^{}\Tr^{F}\left[i^{*}T\right]$ is  well-defined.
\begin{proposition}\label{prop:trp}
	For $0\le p\le n$, the following identities hold:
	\begin{align}\label{eq:prig5}
&\Tr^{\Omega^{p}\left(Y,F\right)}\left[M_{T}\right]=\left(-1\right)^{p\left(n-p\right)}\int_{\Delta^{Y}}^{}\Tr^{F}\left[i^{*}T^{p,n-p}\right],\\
&\Tr^{\Omega^{p}\left(Y,F\right)}\left[N_{T}\right]=\left(-1\right)^{p}\int_{\Delta^{Y}}^{}\Tr^{F}\left[i^{*}T^{p,n-p}\right]. \notag
\end{align}
In particular,
\begin{equation}\label{eq:prig6}
\Trs\left[N_{T}\right]=\int_{\Delta^{Y}}^{}\Tr^{F}\left[i^{*}T\right].
\end{equation}
\end{proposition}
\begin{proof}
	The proof of the first identity (\ref{eq:prig5}) follows from an easy computation 
	left to the reader. The second identity is a consequence of the first one. The proof of our proposition is complete. 
\end{proof}
\begin{remark}\label{rem:nogr}
	The second identity in (\ref{eq:prig5}) can be viewed as consequence of (\ref{eq:prig6}), by assuming that $T$  is only concentrated in degree $p,n-p$. If $\deg\, T\neq n$, both sides of (\ref{eq:prig6}) 
	vanish identically, so that (\ref{eq:prig6}) holds for $T$ of 
	arbitrary degree.
\end{remark}

Clearly,   $\Delta^{Y}$ lies in $\mathcal{D}^{n}\left(Y\times Y,\pi^{*}_{1}F \otimes \pi^{*}_{2} \left(\Ft\right) \right)$, and so $d^{Y}\Delta^{Y}$ is of degree $n+1$, and $i_Z\Delta^{Y}$ is of degree $n-1$. Also $L_Z\Delta^{Y}$ is of degree $n$.
\begin{proposition}\label{prop:diag}
	The following identities hold: 
	\begin{align}\label{eq:sacr0}
&M_{\Delta^{Y}}=\left(-1\right)^{nN}, &N_{\Delta^{Y}}=1.
\end{align}
Moreover, 
\begin{align}\label{eq:sacr0a1}
&d^{Y}=N_{d^{Y}\Delta^{Y}}=\left(-1\right)^{n}N_{\Delta^{Y}d^{Y}},\notag\\
&i_{Z}=N_{i_{Z}\Delta^{Y}}=\left(-1\right)^{n}N_{\Delta^{Y}i_{Z}},\\
& L_Z=N_{L_Z\Delta^{Y}}=N_{\Delta^{Y}L_{Z}}.\notag
\end{align}
\end{proposition}
\begin{proof}
	If $s\in \Omega\left(Y,F\right)$, since $\Delta^{Y}$ is a 
	current of degree $n$, then
	\begin{equation}\label{eq:sacr0a2}
\Delta^{Y}\we \pi^{*}_{2}s=\left(-1\right)^{nN}\pi^{*}_{2}s\we\Delta^{Y}.
\end{equation}
By integrating (\ref{eq:sacr0a2}), we get (\ref{eq:sacr0}). From (\ref{eq:sacr0}), we obtain (\ref{eq:sacr0a1}). The proof of our proposition is complete.
\end{proof} 
\subsection{Intersection traces and the composition of currents}%
\label{subsec:sclatr}
Let $T,T'$ be taken as before so that $T \otimes T'$ is a well defined 
current on $Y^{4}$. Then $Y\times\Delta^{Y}\times 
Y$ is a current  on  $Y^{4}$, whose 
\index[terms]{wave front set}%
wave front set  coincides with 
$\pi_{2,3}N^{*}_{\Delta^{Y}/Y\times Y}\setminus Y^{4}$. Let
\index[not]{WF@$\mathrm{WF}$}%
$\mathrm{WF}\left(T \otimes T'\right)$ be the
\index[terms]{wave front set}%
wave front set of $T \otimes T'$. Assume that 
$\mathrm{WF}\left(T \otimes 
T'\right)+\pi_{2,3}N^{*}_{\Delta^{Y}/Y\times Y}\setminus Y$ does not 
contain $0$. Then the restriction of $T \otimes T'$ to $Y\times 
\Delta^{Y}\times Y$ is well defined.

Set
\begin{equation}\label{eq:sv1}
T\circ T'=\pi_{1,4*}\left[\left( T \otimes T' \right) \delta_{Y\times 
\Delta^{Y}\times Y}\right].
\end{equation}
Then $T\circ T'$ is a current similar to $T,T'$ with degree $\mathrm{deg}\,T+\mathrm{deg}\,T'-n$.  If $s\in 
\Omega\left(Y,F\right)$, $M_{T}M_{T'}s$ is well-defined.
\begin{proposition}\label{prop:compo}
	The operators $M_{T\circ T'},N_{T\circ T'}$ are  well defined. 
	Moreover, 	
	\begin{align}\label{eq:sv2}
		&M_{T}M_{T'}=N_{T\circ T'},
&N_{T}N_{T'}=\left(-1\right)^{n\left(\deg\, T'-n\right)}N_{T\circ T'}.
\end{align}
\end{proposition}
\begin{proof}
	Using an obvious notation, we get
	\begin{equation}\label{eq:prigb1}
M_{T}M_{T'}s=\pi_{1*}\left[T\we\pi_{2}^{*}\pi_{2*}\left[T'\we 
\pi_{3}^{*}s\right]\right].
\end{equation}
We claim that the first identity in (\ref{eq:sv2}) follows from 
(\ref{eq:prigb1}). This is because of our rules of integration along 
the fiber\footnote{In principle, an extra factor $\left(-1\right)^{n}$ should appear in the right-hand side. This extra term disappears because the two factors $\mathrm{o}\left(TY\right)$ are also interchanged, and this also introduces the same factor $\left(-1\right)^{n}$.}. The second identity is a consequence of the first one.
\end{proof}
\begin{remark}\label{rem:asso}
	As a consequence of associativity, if $T,T',T''$ verify the 
	proper wave front set conditions, we get the identity
	\begin{equation}\label{eq:asso1}
\left(T\circ T'\right)\circ 
T''=\left(-1\right)^{n\deg\, T''-n}T\circ\left(T'\circ T''\right),
\end{equation}
which can be proved directly.

Also note that
\begin{align}\label{eq:sacr-1}
&T\circ \Delta^{Y}=T,&\Delta^{Y}\circ T=\left(-1\right)^{n\left(\deg\,T-n\right)}T.
\end{align}
Equation  (\ref{eq:prig1}) is a consequence of  (\ref{eq:sacr0}), of the first equation in (\ref{eq:sv2}) with $T'=\Delta^{Y}$, and of (\ref{eq:sacr-1}), and the second equation in 
(\ref{eq:sv2}) can be understood in the same way.  Finally, observe 
that the identity $L_{Z}=\left[d^{Y},i_{Z}\right]$ can be obtained as 
a consequence of (\ref{eq:sacr0a1}) and (\ref{eq:sv2}).
\end{remark}

Let 
\index[not]{Y@$\mathcal{Y}^{*}$}%
$\pi:\mathcal{Y}^{*}\to Y$ be the total space of the cotangent bundle of $Y$, with fiber $T^{*}Y$. Let 
\index[not]{q@$q$}%
$q$ be the tautological section of $\pi^{*}T^{*}Y$ on $\mathcal{Y}^{*}$. Then $Y$ embeds in $\mathcal{Y}^{*}$ as the zero set of $q$.

If $s\in \mathcal{D}\left(Y,\C\right)$, let $\mathrm{WF}s\in 
\mathcal{Y}^{*}\setminus Y$ be the wave front set of $s$, which is a 
conical subset of $\mathcal{Y}^{*}\setminus Y$. If $s'\in 
\mathcal{D}\left(Y,\C\right)$ and if $\mathrm{WF}s +
\mathrm{WF}s'$ does not contain $0$, by \cite[Theorem 8.2.10]{Hormander83a}, the product 
of currents $s\we s'$ is unambiguously defined, and 
\begin{equation}\label{eq:sanfa2a1}
\mathrm{WF}\left( s\we s' \right)  \subset \left( \mathrm{WF}s+\mathrm{WF}s' \right) \cup 
\mathrm{WF}s\cup\mathrm{WF}s'.
\end{equation}
The same considerations apply when $\C$ is replaced by $F$ or $\Ft$.

Let $T$ be taken as above. Let $s\in \mathcal{D}\left(Y,F\right)$ be 
such that $0\notin\mathrm{WF}T+\pi_{2}^{*}\mathrm{WF}s$.
The product $T\we\pi_{2}^{*}s$ is still well 
 defined, and its wave front set is contained in  
 $ \left(\mathrm{WF}T+\pi_{2}^{*}\mathrm{WF}s \right) \cup \mathrm{WF}T\cup \pi^{*}_{2}\mathrm{WF}s$.
 We can  define the current $M_{T}s$ as in (\ref{eq:rec1}), and 
$\mathrm{WF}M_{T}s$ is contained in the $q\in 
 T^{*}Y\setminus 0$ such that $\pi_{1}^{*}q$ lies in 
 $\mathrm{WF}\left[ T\we\pi_{2}^{*}s\right]$. 
 Assume that $0\notin
\mathrm{WF}\left(T\right) + N^{*}_{\Delta^{Y}/Y\times Y}$. Suppose 
that $T$ is of degree $k$. 
 Then $T\Delta^{Y}$ is a well-defined
 current of degree $k+n$ on $Y\times Y$. Equivalently, $i^{*}T$ is a well-defined 
element of $\mathcal{D}^{k}\left(Y,F \otimes F^{*} \otimes 
_{\Z_{2}}\mathrm{o}\left(TY\right)\right)$.

Here, we inspire ourselves from Proposition \ref{prop:trp}.
\begin{definition}\label{def:sctra}
	If $T$ is of degree $n$, let $\mathrm{Tr}^{\Omega^{p}\left(Y,F\right)}_{\mathrm{int}}\left[M_{T}\right]$ be the 
	intersection trace\footnote{In \cite{AtiyahBott64}, Atiyah and Bott introduced  the idea of graph trace. In dynamical systems, this became instead the flat trace. We prefer to use the terminology `intersection trace' to make clear that we are dealing with the standard intersection of a current with the diagonal.} of $M_{T}$ in degree $p$, given by
	\begin{equation}\label{eq:sacr1}
		\mathrm{Tr}^{\Omega^{p}\left(Y,F\right)}_{\mathrm{int}}\left[M_{T}\right]=\left(-1\right)^{p\left(n-p\right)}\int_{\Delta^{Y}}^{}\Tr^{F}\left[i^{*}T^{p,n-p}\right].
\end{equation}
\end{definition}

The corresponding intersection trace of $N_{T}$ is given by
\begin{equation}\label{eq:scar1a1}
	\mathrm{Tr}^{\Omega^{p}\left(Y,F\right)}_{\mathrm{int}}\left[N_{T}\right]=\left(-1\right)^{p}\int_{\Delta^{Y}}^{}\Tr^{F}\left[i^{*}T^{p,n-p}\right].
\end{equation}
By (\ref{eq:scar1a1}), we may use the notation
	\begin{equation}\label{eq:sct1}
\Tr_{\mathrm{s,int}}\left[N_{T}\right]=\int_{Y}^{}\Tr^{F}\left[i^{*}T\right].
\end{equation}
 As in Remark \ref{rem:nogr}, we find that (\ref{eq:scar1a1}) is a consequence of (\ref{eq:sct1}). The integral in the right-hand side of (\ref{eq:sct1}) vanishes except when $\deg\, T=n$. 
 Equation (\ref{eq:sct1}) is also valid when $T$ is of arbitrary degree.

As in Subsection \ref{subsec:sclatr}, we assume that $T \otimes T'$ has a well defined restriction to $Y\times\Delta^{Y}\times Y$, which will be also denoted $\Delta^{Y}_{2,3}$. 
 Let $\Delta^{Y}_{1,4}$ be the submanifold of $Y^{4}$ in whose first and last components in $Y$ coincide. Assume  also that $T \otimes T'$ has a well defined restriction to 
$\Delta_{1,4}^{Y}$. Then $T \circ T'$ and $T'\circ T$ are 
both well-defined.  We also assume that $T \otimes T'$ has a well defined restriction 
to $\Delta_{1,4}^{Y}\cap\Delta_{2,3}^{Y}$. Then $N_{T}N_{T'}$ and $N_{T'}N_{T}$ have well-defined restrictions to $\Delta^{Y}$.
\begin{proposition}\label{prop:vano1}
	The following identity holds:
	\begin{equation}\label{eq:comid6}
\mathrm{Tr_{s,int}}\left[\left[N_{T},N_{T'}\right]\right]=0.
\end{equation}
\end{proposition}
\begin{proof}
Put $k=\deg\, T,k'=\deg\,T'$. We have the 
identity,
\begin{equation}\label{eq:comid6b1}
\left[N_{T},N_{T'}\right]=N_{T}N_{T'}-\left(-1\right)^{\left(k-n\right)\left(k'-n\right)}N_{T'}N_{T}.
\end{equation}
By (\ref{eq:sv2}), (\ref{eq:comid6b1}), we obtain
\begin{equation}\label{eq:comid6b2}
\left[N_{T},N_{T'}\right]=N_{\left(-1\right)^{n\left(k'-n\right)}T\circ T'-\left(-1\right)^{\left(k-n\right)\left(k'-n\right)+n\left(k-n\right)}T'\circ T}.
\end{equation}
Equation (\ref{eq:comid6b2}) can be rewritten in the form
\begin{equation}\label{eq:comid6b3}
\left[N_{T},N_{T'}\right]=\left(-1\right)^{n\left(k'-1\right)}N_{T\circ T'-\left(-1\right)^{kk'+n}T'\circ T}.
\end{equation}

By (\ref{eq:sct1}), (\ref{eq:comid6b3}), we obtain
\begin{equation}\label{eq:prigb4}
\mathrm{Tr_{s,int}}\left[\left[N_{T},N_{T'}\right]\right]=\left(-1\right)^{n\left(k'-1\right)}\int_{Y}^{}\Tr^{F}\left[i^{*}\left(T\circ T'-\left(-1\right)^{kk'+n}T'\circ T\right)\right].
\end{equation}
We claim that (\ref{eq:prigb4}) vanishes. Indeed let 
\index[not]{j@$j$}%
$j$ be the 
obvious involution of $Y\times Y$. Then
\begin{align}\label{eq:prigb5}
&\int_{Y}^{}\Tr^{F}\left[i^{*}\left(T \circ 
T'\right)\right]=\int_{Y\times Y}^{}\Tr^{F}\left[T\we 
j^{*}T'\right],\\
&\int_{Y}^{}\Tr^{F}\left[i^{*}\left(T' \circ 
T\right)\right]=\int_{Y\times Y}^{}\Tr^{F}\left[T'\we j^{*}T\right]. 
\nonumber 
\end{align}
Also,
\begin{equation}\label{eq:prigb6}
\int_{Y\times Y}^{}\Tr^{F}\left[T\we 
j^{*}T'\right]=\int_{Y\times Y}^{}j^{*}\Tr^{F}\left[T\we 
j^{*}T'\right].
\end{equation}
Equivalently,
\begin{equation}\label{eq:prigb7}
	\int_{Y\times Y}^{}\Tr^{F}\left[T\we 
j^{*}T'\right]=\int_{Y\times Y}^{}\Tr^{F}\left[j^{*}T\we 
T'\right].
\end{equation}
We claim that
\begin{equation}\label{eq:prigb8}
\Tr^{F}\left[j^{*}T\we 
T'\right]=\left(-1\right)^{kk'+n}\Tr^{F}\left[T'\we j^{*}T\right].
\end{equation}
This is because $T,T'$ have degrees $k,k'$, and also because of the 
fact that $\mathrm{o}\left(TY\right)$ is a line bundle of 
$\Z_{2}$-degree $\left(-1\right)^{n}$. By (\ref{eq:prigb4})--(\ref{eq:prigb8}), we get (\ref{eq:comid6}). The proof of our proposition is complete.  
\end{proof}
\begin{remark}\label{rem:comp}
When $T$ is a smooth current, by Proposition \ref{prop:trp}, $\mathrm{Tr}_{\mathrm{s,int}}\left[N_{T}\right]$ is just the classical supertrace $\Trs\left[N_{T}\right]$, which itself vanishes on supercommutators. In this case, when $T,T'$ are smooth, Proposition \ref{prop:vano1} is trivial. In general, $T,T'$ can be suitably approximated by smooth currents, so that Proposition \ref{prop:vano1} can be obtained as a consequence of Proposition \ref{prop:trp}. This is exactly what was done in Chaubet-Dang \cite[Proposition 4.1]{ChaubetDang24}, in which they gave a version of Proposition \ref{prop:vano1}.
\end{remark}
\subsection{Adjoints}%
\label{subsec:adjo}
We will now show how to take adjoints properly. We use the 
conventions in (\ref{eq:normi3}), (\ref{eq:normi4a}). Let 
\index[j]{j@$j$}%
$j$ be the obvious involution of $Y\times Y$. The anti-automorphism \,$\,\widetilde {}\,$\, also acts naturally on 
$\pi_{1}^{*}\Lambda\left(T^{*}Y\right) \otimes 
\pi_{2}^{*}\Lambda\left(T^{*}Y\right) $. We make this 
anti-automorphism act on $T$ by coupling it with the involution 
$j^{*}$, and we denote the ultimate result by 
$\underline{j}^{*}T$.
\begin{proposition}\label{prop:adj}
	If $\alpha,\beta$ are taken as in (\ref{eq:normio3a1}), then
	\begin{equation}\label{eq:sv5}
\left\langle \beta, 
M_{T}\alpha\right\rangle=\left(-1\right)^{n\left(n+1\right)/2}\left\langle 
M_{\underline{j}^{*}T}\beta,\alpha\right\rangle.
\end{equation}
\end{proposition}
\begin{proof}
	Clearly,
	\begin{equation}\label{eq:sv5a1}
\left\langle  \beta,M_{T}\alpha\right\rangle=\int_{Y\times Y}^{}\left\langle  
\pi_{1}^{*} \widetilde \beta\we T\we\pi_{2}^{*}\alpha\right\rangle.
\end{equation}
Using (\ref{eq:normio3a1}) with $n$ replaced by $2n$, and also the 
fact that two copies of $\mathrm{o}\left(TY\right)$ are interchanged, 
(\ref{eq:sv5a1}) can be written in the form,
\begin{equation}\label{eq:sv5a2}
\left\langle  \beta,M_{T}\alpha\right\rangle=\int_{Y\times 
Y}^{}\left\langle  \pi_{1}^{*}\widetilde \alpha\we \underline{j}^{*}T\we 
\pi_{2}^{*}\beta\right\rangle.
\end{equation}
Equation (\ref{eq:sv5a2}) can be rewritten as
\begin{equation}\label{eq:sv5a3}
\left\langle  \beta,M_{T}\alpha\right\rangle=\left\langle  
\alpha,M_{\underline{j}^{*}T}\beta\right\rangle.
\end{equation}
Using (\ref{eq:normio3a1}) again and (\ref{eq:sv5a3}), we get 
(\ref{eq:sv5}). The proof of our proposition is complete. 
\end{proof}

In the sequel, we denote by $M_{T}^{\dag}$ the transpose of $M_{T}$ 
with respect to $\left\langle  \,\right\rangle$. Equation 
(\ref{eq:sv5}) can be rewritten in the form,
\begin{equation}\label{eq:sv5a4}
M_{T}^{\dag}=\left(-1\right)^{n\left(n+1\right)/2}M_{\underline{j}^{*}T}.
\end{equation}
By (\ref{eq:sv5a4}), we deduce that
\begin{equation}\label{eq:sv5a5}
N_{T}^{\dag}=\left(-1\right)^{n\left(n+1\right)/2}\left(-1\right)^{n\left(n-N\right)}M_{\underline{j}^{*}T}.
\end{equation}
By (\ref{eq:sv5a4}), (\ref{eq:sv5a5}), we get
\begin{equation}\label{eq:sv5a6}
N_{T}^{\dag}=\left(-1\right)^{n\left(n+1\right)/2+n\deg\,T}N_{\underline{j}^{*}T}.
\end{equation}

Let us prove that (\ref{eq:sct1}) and (\ref{eq:sv5a6}) are 
compatible. We only need to consider the case where $\deg\,T=n$, in which 
case (\ref{eq:sv5a6}) can be written in the form
\begin{equation}\label{eq:sv5a7}
N_{T}^{\dag}=\left(-1\right)^{n\left(n+1\right)/2+n}N_{\underline{j}^{*}T}.
\end{equation}

Using (\ref{eq:sct1}) and (\ref{eq:sv5a7}), we get
\begin{equation}\label{eq:sv5a8}
\mathrm{Tr_{s,int}}\left[N^{\dag}_{T}\right]=\left(-1\right)^{n}\mathrm{Tr_{s,int}}\left[N_{T}\right].
\end{equation}
Equation (\ref{eq:sv5a8}) is indeed compatible with the construction 
of a supertrace and with Poincaré duality. 
 \section{Anosov vector fields  and anisotropic Sobolev spaces}%
\label{sec:avf}
The purpose of this section is to recall basic results on the 
anisotropic Sobolev spaces associated with the Anosov vector field 
$Z$, and to construct the resolvent of $L_{Z}$ as a meromorphic 
family of Fredholm operators. We rely entirely on former results of Dyatlov-Zworski \cite{DyatlovZworski16} and Faure-Sjöstrand \cite{FaureSjostrand11} and Faure-Roy-Sjöstrand \cite{FaureRoySjostrand08}. We have provided the above results in great detail because we will need to reference them very precisely later.

This section is organized as follows. In Subsection \ref{subsec:res}, 
we define the resolvent $\left(L_{Z}-z\right)^{-1}$ when $\Re\,z\ll 
-1$.

In Subsection \ref{subsec:esca}, we define the escape functions and 
the order functions on the total space $\mathcal{Y}^{*}$ of $T^{*}Y$.

In Subsection \ref{subsec:anis}, we recall the construction of the 
anisotropic Sobolev spaces.

In Subsection \ref{subsec:opasp}, we prove that   $L_{Z}$ acts as an 
unbounded operator on the anisotropic Sobolev spaces.

Finally, in Subsection \ref{subsec:relz}, we show that in the proper 
sense, the resolvent $\left(L_{Z}-z\right)^{-1}$ is a meromorphic 
family of Fredholm operators. 

In this section, we make the same assumptions as in Section \ref{sec:frizet}.
\subsection{The resolvent of $L_{Z}$}%
\label{subsec:res}
Let  $g^{TY}$ be a Riemannian metric on $TY$, let $g^{F}$ be a 
Hermitian  metric on $F$. These metrics induce a corresponding 
$L_{2}$ Hermitian product $\left\langle \ \right\rangle$ on $\Omega\left(Y,F\right)$.  Let $\Omega_{2}\left(X,F\right)$ be the Hilbert space of 
square-integrable sections of $\Lambda\left(T^{*}Y\right) \otimes 
_{\R}F$. Note that $\Omega_{2}\left(Y,F\right)$ does not depend on 
the choice of $g^{TY},g^{F}$.

Since $L_{Z}$ 
is a first order differential operator, there exists $c_{0}>0$ such 
that for  $u\in \Omega(Y, F)$, 
\begin{align}\label{eq:eqc0}
\Re	\<L_{Z}u, u\>\g -c_{0}\<u, u\>. 
\end{align}
In the sequel, we will consider $L_{Z}$ as an unbounded operator acting on 
$\Omega_{2}\left(Y,F\right)$. 

In  \cite[Lemma A.1]{FaureSjostrand11}, Faure and Sjöstrand showed 
that if $Q$ is any scalar pseudo-differential operator of order $1$, 
viewed as an unbounded operator acting on $\Omega_{2}\left(Y,F\right)$,  the minimal 
and maximal extensions of $Q$  coincide. Since $L_{Z}$ is a first order operator with scalar principal symbol, their result extends 
trivially to the operator $L_{Z}$ acting on 
$\Omega_{2}\left(Y,F\right)$. 
By 
\eqref{eq:eqc0}, the spectrum of this extension is contained in $\{z: 
\Re{z}\g -c_{0}\}$, so that if $z\in \bC$ with $\Re z< -c_{0}$,  
\begin{align}\label{eq:song0}
	\(L_{Z}-z\)^{-1}: \Omega_{2}\left(Y,F\right)\to 
	\Omega_{2}\left(Y,F\right)
\end{align}
is a holomorphic family of bounded operators. The resolvent takes in 
fact its values in the smaller domain of $L_{Z}$.

By restriction, we get a holomorphic family of bounded 
operators\footnote{This means that for 
$u\in\Omega\left(Y,F\right),v\in\Omega\left(Y,F \otimes 
_{\Z_{2}}\mathrm{o}\left(TY\right)\right)$, the function  
$z\in\C,\mathrm{Re}z<-c_{0}\to\left\langle  
\left(L_{Z}-z\right)^{-1}u,v\right\rangle$ is holomorphic.}
\begin{equation}\label{eq:eqhlo1}
\left(L_{Z}-z\right)^{-1}: 
\Omega\left(Y,F\right)\to\mathcal{D}\left(Y,F\right), 
z\in\C,\mathrm{Re}z<-c_{0}.
\end{equation}
\subsection{Escape functions and order functions}%
\label{subsec:esca}
 Note that 
$\mathcal{Y}^{*}$ is a symplectic manifold, whose symplectic form is 
denoted 
\index[not]{om@$\omega^{\mathcal{Y}^{*}}$}%
$\omega^{\mathcal{Y}^{*}}$. Also 
\index[not]{q@$q$}%
$q$ still denotes  the canonical section of $\pi^{*}T^{*}Y$ on 
$\mathcal{Y}^{*}$. 

The vector field $Z$ lifts to a vector field $\mathcal{Z}$ on 
$\mathcal{Y}^{*}$.   If 
$\left(y,q\right)$ is the generic point of $\mathcal{Y}^{*}$ with 
$y\in Y,q\in T^{*}_{y}Y$, then $\mathcal{Z}$ is a Hamiltonian vector field  
whose Hamiltonian $\mathcal{H}$ is given by 
\begin{equation}\label{eq:duaz4}
\mathcal{H}\left(y,q\right)=\left\langle  q,Z_{y}\right\rangle.
\end{equation}

The flow $\varphi_{t}\vert_{t\in 
\R}$ lifts to the corresponding flow $\varphi_{t*}\vert_{t\in \R}$ on 
$\mathcal{Y}^{*}$ associated with $\mathcal{Z}$. Note that 
$\varphi_{t*}\vert_{t\in\R}$ maps fibers into fibers.

Let
\index[not]{R@$R$}%
$R$ be the Euler vector field along the fibers of $T^{*}Y$. Then 
$R$ can be viewed as a vector field on $\mathcal{Y}^{*}$, and
\begin{equation}\label{eq:duaz5}
\left[\mathcal{Z},R\right]=0.
\end{equation}
The vector field $R$ generates the group of dilations along the fibers $T^{*}Y$. If $a\in \R^{*}$, let 
\index[not]{r@$r_{a}$}%
$r_{a}$ be the dilation $\left(y,q\right)\to \left(y,aq\right)$.  

Observe that  
$\varphi_{t*}\vert_{t\in \R}$ preserves  $\mathcal{Y}^{*}\setminus 
Y$. The vector field $\mathcal{Z}$ restricts to a vector field on 
$\mathcal{Y}^{*}\setminus Y$. On $\mathcal{Y}^{*}\setminus Y$, $R$ 
is the generator of  a vector subbundle 
$\mathcal{R}$ of $T\mathcal{Y}^{*}$.

Given $y\in Y$, if $C^{*}_{y}\in T^{*}Y\setminus 0$ is a conical set, we will say that 
$\mathcal{C}^{*}=\left\{y,C^{*}_{y}\right\}$ is a conical set in 
$\mathcal{Y}^{*}\setminus Y$.

Let 
\index[not]{Y@$\mathcal{Y}^{*}_{0}$}%
\index[not]{Y@$\mathcal{Y}^{*}_{u}$}%
\index[not]{Y@$\mathcal{Y}^{*}_{s}$}%
$\mathcal{Y}^{*}_{0},\mathcal{Y}^{*}_{u},\mathcal{Y}^{*}_{s}$ be the total spaces of 
$T^{*}_{0}Y,T^{*}_{u}Y,T^{*}_{s}Y$. Then $\mathcal{Y}^{*}_{0}\setminus Y,\mathcal{Y}^{*}_{u}\setminus 
Y,\mathcal{Y}^{*}_{s}\setminus Y$ are conical subsets in 
$\mathcal{Y}^{*}\setminus Y$.
\begin{definition}\label{def:cono}
		Let $\mathcal{C}^{*}_{u}, \mathcal{C}^{*}_{s}\subset 
		\mathcal{Y}^{*}\setminus Y$ be  closed non-intersecting conical 
		neighborhoods of $\mathcal{Y}^{*}_{u}\setminus Y, 
		\mathcal{Y}^{*}_{s}\setminus Y$.  Given $\beta>0$,  
		a smooth function $f: \mathcal{Y}^{*}\backslash Y\to 
		\bR_{+}^{*}$  is 
		called an 
		\index[terms]{adapted function}%
		adapted function for $\mathcal{Z}$ with respect to $\(\mathcal{C}^{*}_{u}, 
		\mathcal{C}^{*}_{s}, \beta\)$ if $f$	 is a positively $1$-homogenous 
		function\footnote{For $h\in\R$, a function $g$ is said to be positively $h$-homogeneous if  for $a>0$, $g\left(y,aq\right)=a^{h}g\left(y,q\right)$.}, and
		\begin{align}\label{eq:tann1}
	\begin{aligned}
&\mathcal{Z} f \le- \beta f	\,\text{on a conical neighborhood of} \ \mathcal{C}^{*}_{u}, 
\\
&\mathcal{Z}f\ge  \beta f	\,\text{on a conical neighborhood of}\ \mathcal{C}^{*}_{s}. 
\end{aligned}
\end{align} 
	\end{definition} 
	\begin{remark}\label{rem:neg}
	If $f$ is an adapted function for the flow $\mathcal{Z}$, then $f$ is 
	also an adapted function for the flow $-\mathcal{Z}$. 
\end{remark}

We recall a result in \cite[Section 2, eq. (2.8)]{FaureSjostrand11}, \cite[Appendix 
A.2]{DangShen20}. 
	\begin{proposition}\label{prop:pco}
If  $\mathcal{C}^{*}_{u}, \mathcal{C}^{*}_{s}$  
are small enough, there exist $\beta>0$ and  an adapted 
	function $f$ with respect to  $\(\mathcal{C}^{*}_{u}, \mathcal{C}^{*}_{s}, \beta\)$. 
\end{proposition} 
\begin{definition}\label{def:ord}
		If 
		$\underline{\mathcal{C}}^{*}_{0},\underline{\mathcal{C}}_{u}^{*}, \underline{\mathcal{C}}^{*}_{s}\subset \mathcal{Y}^{*}\setminus Y$ are closed  non-intersecting conical neighborhoods 
		of	$\mathcal{Y}^{*}_{0}\setminus Y, 
		\mathcal{Y}^{*}_{u}\setminus Y, \mathcal{Y}^{*}_{s}\setminus 
		Y$, and if $\gamma>0$, a smooth  function $m: 
		\mathcal{Y}^{*}\setminus Y\to [-2,2]$ 
		is called an 
		\index[terms]{order function}%
		order function for $\mathcal{Z}$ with respect to 	$ 
		\left( 
		\underline{\mathcal{C}}^{*}_{0},\underline{\mathcal{C}}_{u}^{*}, \underline{\mathcal{C}}^{*}_{s},\gamma \right)$, 	if 
			\begin{enumerate}
		\item $m$ is positively $0$-homogenous. 
		\item  We have the inequality:
		\begin{align}
	\mathcal{Z}  m\l 0, 
\end{align} 
\item Also,
\begin{align}
	\mathcal{Z} m \vert_{\left(\mathcal{Y}^{*}\setminus 
	Y\right)\setminus 
	\left(\underline{\mathcal{C}}^{*}_{0}\cup\underline{\mathcal{C}}_{u}^{*}\cup \underline{\mathcal{C}}^{*}_{s}\right)} \l -\gamma.
\end{align} 
		\item Finally,  
		\begin{align}\label{eq:eqmz}
	m(z)&	\g 1\ \mathrm{on}\ \underline{\mathcal{C}}^{*}_{u},\\
	&\l -1\ \mathrm{on} \ \underline{\mathcal{C}}^{*}_{s}. \nonumber 
\end{align} 
	\end{enumerate} 
	\end{definition} 
	
\begin{remark}\label{rem:chs}
	If $m$ is an order function for $\mathcal{Z}$, then $-m$ is an order 
	function for $-\mathcal{Z}$. 
\end{remark}

\begin{theorem}\label{thm:thm1}
	If 
	\index[not]{C@$\mathcal{C}^{*}_{0}$}%
	\index[not]{Cu@$\mathcal{C}^{*}_{u}$}%
	\index[not]{Cs@$\mathcal{C}^{*}_{s}$}%
	$\mathcal{C}^{*}_{0}, 
	\mathcal{C}^{*}_{u},\mathcal{C}^{*}_{s}$ are closed 
	non-intersecting conical neighborhoods of 
	$\mathcal{Y}^{*}_{0}\setminus 
	Y,\mathcal{Y}^{*}_{u}\setminus 
	Y,\mathcal{Y}^{*}_{s}\setminus Y$, there exist $\gamma>0$, closed conical 
	neighborhoods
	\index[not]{C@$\underline{\mathcal{C}}^{*}_{0}$}%
	\index[not]{Cu@$\underline{\mathcal{C}}^{*}_{u}$}%
	\index[not]{Cs@$\underline{\mathcal{C}}^{*}_{s}$}%
	$\underline{\mathcal{C}}^{*}_{0}\subset 
	\mathcal{C}^{*}_{0},\underline{\mathcal{C}}^{*}_{u} \subset 
	\mathcal{C}^{*}_{u},\underline{\mathcal{C}}^{*}_{s} \subset 
	\mathcal{C}^{*}_{s}$ of $\mathcal{Y}^{*}_{0}\setminus 
	Y,\mathcal{Y}^{*}_{u}\setminus Y,\mathcal{Y}^{*}_{s}\setminus Y$, and an order function $m: 
	\mathcal{Y}^{*}\setminus Y\to [-2,2]$ with 
	respect to 
	$\left(\underline{\mathcal{C}}^{*}_{0},\underline{\mathcal{C}}^{*}_{u},\underline{\mathcal{C}}^{*}_{s},\gamma\right)$.
\end{theorem} 
\begin{proof}
We refer to  \cite[Section 2]{FaureSjostrand11} and  to \cite[Appendix A.1 and eq. (A.10)]{DangShen20}. 
\end{proof} 
\begin{proposition}\label{prop:pextra}
	With the same notation and assumptions as in Theorem 
	\ref{thm:thm1},   $m$ can be taken such that
	\begin{align}\label{eq:4}
m	\g 1/4 \ \mathrm{on} \ \left(\mathcal{Y}^{*}\setminus 
Y\right)\setminus \mathcal{C}^{*}_{s}.
\end{align} 
\end{proposition}
\begin{proof}
	This statement follows from the same references as in the proof 
	of Theorem \ref{thm:thm1}.
\end{proof}

We  fix small closed non-intersecting conical neighborhoods 
$\mathcal{C}^{*}_{0},\mathcal{C}^{*}_{u},\mathcal{C}^{*}_{s}$ of 
$\mathcal{Y}^{*}_{0}\setminus Y,\mathcal{Y}^{*}_{u}\setminus 
Y,\mathcal{Y}^{*}_{s}\setminus Y$.  Let $f$ be an adapted function 
with respect to $ \mathcal{C}_{u}^{*}, \mathcal{C}_{s}^{*}$ and some 
$\beta>0$, and let $m$ be an order function with  respect to  smaller 
closed conical neighborhood 
$\underline{\mathcal{C}}^{*}_{0},\underline{\mathrm{C}}^{*}_{u},\underline{\mathcal{C}}^{*}_{s}$ and some $\gamma>0$.

Let $U_{1}, U_{2}\subset \mathcal{Y}^{*}$ be two relatively compact 
open neighborhoods of $Y$ such that $\ol U_{1}\subset U_{2}$. 
	\begin{definition}\label{def:esc}
Let $g:\mathcal{Y}^{*}\to \R$ be a smooth function such that 
\begin{align}\label{eq:res1}
		g(z)=&
m(z)\ln\(1+ f(z)\),& z\notin U_{2}, \\
&0, &z\in U_{1}. \nonumber 
\end{align} 		
	\end{definition} 
Clearly, there is $C>0$ such that 
\begin{equation}\label{eq:res2}
\mathcal{Z} g\l C. 
\end{equation}

\begin{theorem}\label{thm:thmrec}
There exists 
\index[not]{d@$\delta>0$}%
$\delta>0$ and there is  a  relatively compact open 
	neighborhood $U_{3}$ of $Y$ with  $\ol{U}_{2}  \subset U_{3}$ such that 
  \begin{align}\label{eq:eqZgde}
	\mathcal{Z} g\vert_{\mathcal{Y}^{*}\setminus \left(U_{3}\cup 
	\underline{\mathcal{C}}^{*}_{0}\right)}\l -\delta.  
\end{align} 
\end{theorem} 
\begin{proof}
Our theorem is a consequence of Definitions \ref{def:cono} and  
\ref{def:ord},
 and of Theorem \ref{thm:thm1}. 
\end{proof} 

The function $g$ is called an 
\index[terms]{escape function}%
escape function for the vector field 
$\mathcal{Z}$.
\begin{remark}\label{rem:esc}
	 If $g$ is an escape 
	function for $\mathcal{Z}$, then $-g$ is an escape function for $-\mathcal{Z}$. 
\end{remark}

\subsection{The anisotropic Sobolev spaces}%
\label{subsec:anis}
For $a\in \bR$, let  
\index[not]{SaY@$\mathcal{S}^{a}\left(\mathcal{Y}^{*}\right)$}%
$\mathcal{S}^{a}\left( \mathcal{Y}^{*} \right)$ be the 
space of scalar  symbols of order $a$ in the sense of Hörmander 
\cite[Definition 18.1.1]{Hormander85a}. Put
\begin{equation}\label{eq:sym1}
\mathcal{S}^{a^{+}}\(\mathcal{Y}^{*}\)=\bigcap_{a'>a}\mathcal{S}^{a'}\(\mathcal{Y}^{*}\). 
\end{equation}
Then $g\in \mathcal{S}^{0^{+}}\(\mathcal{Y}^{*}\)$. 

Let
\index[not]{SaY@$\mathcal{S}^{a}\left(\mathcal{Y}^{*},\End\left(\Lambda\left(T^{*}Y\right) \otimes _{\R}F\right)\right)$}%
$\mathcal{S}^{a}\left(\mathcal{Y}^{*},\End\left(\Lambda\left(T^{*}Y\right) \otimes _{\R}F\right)\right)$ be the corresponding space of  symbols of order $a$ with coefficients in $\End\left(\Lambda\left(T^{*}Y\right) \otimes _{\R}F\right)$.

We define $\mathcal{S}^{a^{+}}\left(\mathcal{Y}^{*},\End\left(\Lambda\left(T^{*}Y\right) \otimes _{\R}F\right)\right)$ as in (\ref{eq:sym1}). 

For $a\in \R$, let 
\index[not]{Ps@$\Psi^{a}(Y, \Lambda(T^{*}Y)\otimes_{\R} F)$}%
$\Psi^{a}(Y, \Lambda(T^{*}Y)\otimes_{\R} F)$ be the 
 space  of pseudo-differential operators of order $a$ acting on 
$C^{\infty }\left(Y,\Lambda\left(T^{*}Y\right) \otimes 
_{\R}F\right)$. These are  the pseudo-differential operators 
with symbol in  
$\mathcal{S}^{a}\left(\mathcal{Y}^{*},\End\left(\Lambda\left(T^{*}Y\right) \otimes _{\R}F\right)\right)$. 
There is an associated algebra of pseudo-differential operators, that is denoted  
\index[not]{Ps@$\Psi \left( Y,\End\left(\Lambda\left(T^{*}Y\right) \otimes _{\R}F\right)\right)$}%
$\Psi \left( Y,\Lambda\left(T^{*}Y\right) \otimes _{\R}F\right)$.

Set
\index[not]{Ps@$\Psi^{a^{+}}(Y, \Lambda(T^{*}Y)\otimes_{\R}F)$}%
\begin{equation}\label{eq:sym2}
\Psi^{a^{+}}(Y, \Lambda(T^{*}Y)\otimes_{\R} 
F)=\bigcap_{a'>a}\Psi^{a'}(Y, \Lambda(T^{*}Y)\otimes_{\R} F).
\end{equation}

In \cite{FaureRoySjostrand08}, Faure, Roy and Sjöstrand introduced scalar symbols with variable 
	order. Given a smooth function $a\left(y,q\right)$ of  positive 
	homogeneity $0$ on $\mathcal{Y}^{*}\setminus Y$, 
	$\mathcal{S}^{a}\left(\mathcal{Y}^{*}\right)$ consists of the 
	scalar functions $b$ defined on $\mathcal{Y}^{*}$ 
that verify  the same uniform conditions as the usual 
symbols with a fixed $a\left(y,q\right)=a$. Instead of 
(\ref{eq:sym1}), put
\index[not]{SaY@$\mathcal{S}^{a^{+}}\(\mathcal{Y}^{*}\)$}%
\begin{equation}\label{eq:rada3}
\mathcal{S}^{a^{+}}\left(\mathcal{Y}^{*}\right)=\bigcap_{\epsilon>0}\mathcal{S}^{a+\epsilon}\left(\mathcal{Y}^{*}\right).
\end{equation}
Faure, Roy and Sjöstrand \cite{FaureRoySjostrand08} showed that the above considerations on the 
case of constant $a$ extend to the case of a variable $a$. 

In the case of variable $a$, we will use the same notation as in  the case of constant $a$. We define 
\index[not]{Ps@$\Psi^{a}(Y, \Lambda(T^{*}Y)\otimes_{\R} F)$}%
$\Psi^{a}\left(Y,\Lambda\left(T^{*}Y\right) \otimes_{\R}F\right)$ as before. 
As shown by the above authors, there is a corresponding algebra of 
pseudo-differential operators. We define 
\index[not]{Ps@$\Psi^{a^{+}}(Y, \Lambda(T^{*}Y)\otimes_{\R}F)$}%
$\Psi^{a^{+}}\left(Y,\Lambda\left(T^{*}Y\right) \otimes 
_{\R}F\right)$ as in (\ref{eq:sym2}).

Let 
\index[not]{G@$G$}%
$G\in \Psi^{0^{+}}(Y, \Lambda(T^{*}Y)\otimes_{\R} F)$ be a
 pseudo-differential operator whose principal symbol 
$\sigma\left(G\right)$ is such that
\begin{align}\label{eq:normi0}
	\sigma\(G\)=g \in 
	\mathcal{S}^{0^{+}}/\mathcal{S}^{0^{+}-1}\(\mathcal{Y}^{*}, 
	\End\(\Lambda(T^{*}Y)\otimes _{\R}F\)\).
	\end{align}

By (\ref{eq:res1}), if $p\in \bR$, 	 
\begin{equation}\label{eq:sym2a1}
\exp(pg)=(1+f)^{pm} \ \mathrm{on}\ \mathcal{Y}^{*}\setminus U_{2}. 
\end{equation}
For $p\in\R$,  $\exp\(pg\)\in \mathcal{S}^{(pm)^{+}}\left( \mathcal{Y}^{*}, 
\End\(\Lambda(T^{*}Y)\otimes _{\R}F\) \right) $ is a symbol with variable order 
$(pm)^{+}$ 
in the sense of \cite[Appendix]{FaureRoySjostrand08}. 
By 
\cite[Section 3.1]{DyatlovZworski16},  $\exp\(pG\)$ is a
  pseudo-differential  operator in 
$\Psi^{(pm)^{+}}(Y, \Lambda(T^{*}Y)\otimes_{\R} F)$ with variable order 
$(pm)^{+}$. More precisely, by \cite[Theorem 8.6]{Zworski12}, we get
\begin{equation}\label{eq:sym2a1x}
\sigma\left(\exp\left(pG\right)\right)=\exp\left(pg\right).
\end{equation}

\begin{definition}\label{def:wfs}
	Let
	\index[not]{DYu@$\mathcal{D}_{\mathcal{Y}^{*}_{u}}\left(Y,F\right)$}%
	\index[not]{DYs@$\mathcal{D}_{\mathcal{Y}^{*}_{s}}\left(Y,F\right)$}%
	$\mathcal{D}_{\mathcal{Y}^{*}_{u}}\left(Y,F\right),\mathcal{D}_{\mathcal{Y}^{*}_{s}}\left(Y,F\right)$ be the vector spaces of currents on $Y$ with values in $F$ whose wave front set is included in $\mathcal{Y}^{*}_{u},\mathcal{Y}^{*}_{s}$.

	Similarly, let 
	\index[not]{DCu@$\cD_{\underline{\mathcal{C}}^{*}_{u}}(Y,F)$}%
	\index[not]{DCs@$\cD_{\underline{\mathcal{C}}^{*}_{s}}(Y,F)$}%
	$\cD_{\underline{\mathcal{C}}^{*}_{u}}(Y,F),\cD_{\underline{\mathcal{C}}^{*}_{s}}(Y,F)$ be the vector spaces of currents on $Y$ with values in $F$ whose wave front set is included in $\underline{\mathcal{C}}^{*}_{u}, \underline{\mathcal{C}}^{*}_{s}$.
\end{definition}

If $q\in \bR$, let 
\index[not]{HqX@$\mathcal{H}^{q}\(X,\Lambda(T^{*}Y)\otimes_{\R} F\)$}%
$\mathcal{H}^{q}\(X,\Lambda(T^{*}Y)\otimes_{\R} F\)$ be 
the classical Sobolev space of sections of $\Lambda\left(T^{*}Y \right) 
\otimes _{\R}F$ of order $q$.

\begin{definition}\label{def:sob}
For $p, q\in \bR$, set
\index[not]{Hpq@$\cH_{pG}^{q}\(Y, \Lambda(T^{*}Y)\otimes_{\R} F\)$}%
\begin{align}\label{eq:sobobo}
		\cH_{pG}^{q}\(Y, \Lambda(T^{*}Y)\otimes_{\R} F\)= \exp\(-pG\)\cH^{q}\(Y, 	\Lambda(T^{*}Y)\otimes_{\R} F\). 
	\end{align}
The vector space $\cH_{pG}^{q}\(Y, \Lambda(T^{*}Y)\otimes_{\R} F\)$ is called an 
\index[terms]{anisotropic Sobolev space}%
anisotropic Sobolev space. 
For $q=0$, we will write instead $\cH_{pG}\(Y, \Lambda(T^{*}Y)\otimes_{\R} F\)$. 
\end{definition}

The vector space $\cH_{pG}^{q}\(Y, \Lambda(T^{*}Y)\otimes_{\R} F\)$ inherits a corresponding norm
from the norm on $\mathcal{H}^{q}\left(Y,\Lambda(T^{*}Y)\otimes_{\R} 
F\right)$.

Observe that
\begin{equation}\label{eq:soba1}
\mathcal{H}^{0}\left(Y,\Lambda\left(T^{*}Y\right) \otimes 
_{\R}F\right)=\Omega_{2}\left(Y,F\right).
\end{equation}
The anisotropic Sobolev spaces neither increase nor decrease with $p$.
We have
\begin{align}\label{eq:tubi0}
	&\cD_{\underline{\mathcal{C}}^{*}_{s}}(Y,F)\subset \bigcup_{p\ge 0}\cH_{pG}(Y, 
			\Lambda(T^{*}Y)\otimes_{\R} F),\\
	&\cD_{\underline{\mathcal{C}}^{*}_{u}}(Y, F)\subset \bigcup_{p\le 0}\cH_{pG}(Y, 
			\Lambda(T^{*}Y)\otimes_{\R}F).
			 \nonumber 
\end{align}
 Moreover,   by  (\ref{eq:4}), we have 
	\begin{align}\label{eq:eqDs}
		&\bigcap_{p\ge 0, g} \cH_{pG}(Y, \Lambda(T^{*}Y)\otimes_{\R} F) \subset 
		\cD_{\mathcal{Y}^{*}_{s}}(Y,  F),\\
		&\bigcap_{p\le 0, g} \cH_{pG}(Y, \Lambda(T^{*}Y)\otimes_{\R} F) \subset 
		\cD_{\mathcal{Y}^{*}_{u}}(Y, F). \nonumber 
	\end{align} 
\subsection{The operator $L_{Z}$ and  the anisotropic 
Sobolev spaces}%
\label{subsec:opasp}
In the sequel, we fix $g$.

Assume  that $p\g0$. We will view $L_{Z}$ as an unbounded 
operator acting on $\cH_{pG}\left(Y,\Lambda\left(T^{*}Y 
\right)\otimes _{\R}F\right)$. 	By  \cite[Lemma A.1]{FaureSjostrand11}, the 
minimal and maximal extensions of $\(L_{Z},\cH_{pG}\left(Y,\Lambda\left(T^{*}Y 
\right)\otimes _{\R}F\right)\)$  
coincide, and will be denoted 
\index[not]{LZp@$L_{Z,p}$}%
$L_{Z,p}$. The corresponding  domain is given by
\begin{align}
		{\rm dom}(L_{Z,p})= \left\{u\in 
		\cH_{pG}\left(Y,\Lambda(T^{*}Y)\otimes_{\R} F\right), L_{Z}u\in \cH_{pG}\left(Y,\Lambda(T^{*}Y)\otimes_{\R} F\right)\right \}. 
	\end{align} 
	
Equivalently, we  consider the conjugate operator 
\index[not]{MZp@$M_{Z,p}$}%
\begin{align}\label{eq:eqLZG}
M_{Z,p}=	e^{pG} L_{Z} e^{-pG}.
\end{align} 
Then
\begin{equation}\label{eq:eqLZG1}
M_{Z,p}= L_{Z} + \[e^{pG}, L_{Z}\] e^{-pG}.
\end{equation}
The operator $M_{Z,p}$ will be viewed  as an unbounded operator
acting on the Hilbert space 
$\Omega_{2}\left(Y,\Lambda(T^{*}Y)\otimes_{\R} F\right)$. Then 
$L_{Z,p}$ and $M_{Z,p}$ are unitarily equivalent. We still denote by 
$M_{Z,p}$ its minimal and maximal extension.

Note that 
$e^{\pm pG}$  are 
pseudo-differential operators (possibly with variable order). Since $L_{Z}$ is a differential 
operator of order $1$ with scalar principal symbol,  $\[e^{pG}, L_{Z}\] 
e^{-pG}\in \Psi^{0^{+}}\left(Y,\Lambda\left(T^{*}Y\right) \otimes _{\R}F\right)$, and its principal symbol is given by
\begin{align}\label{eq:eqELZ}
	\sigma\(\[e^{pG}, L_{Z}\] e^{-pG}\)= -p \mathcal{Z} g. 
\end{align} 
	
Recall that $c_{0}$ is defined in \eqref{eq:eqc0}, and that $\delta>0$ 
was defined in Theorem \ref{thm:thmrec}. We have the 
fundamental result of Faure-Sjöstrand \cite[Theorem 1.4]{FaureSjostrand11}. 

 \begin{theorem}\label{thm:FaureSjostrand}
	 Given $p\ge 0$, there exists $c_{p}\ge 
	 c_{0}$ such  that the spectrum of $L_{Z,p}$ is  contained in $\{z\in \mathbf{C}: \Re(z)\ge -c_{p}\}$. Moreover, on  
	 $\{z\in \bC: {\rm Re}(z)< -c_{0}+\delta p\}$, $L_{Z,p}-z$ is a 
	 holomorphic family of Fredholm operators. In particular, 
	 $L_{Z,p}$ has only discrete  spectrum in this region,  and  the corresponding characteristic 
	 subspaces are finite dimensional. 
 \end{theorem} 
\begin{proof}
 We will recall the essential steps of the proof of the second 
 statement. Since $L_{Z,p}$ and $M_{Z,p}$ are unitarily equivalent, we may as 
 well establish our result when replacing $L_{Z,p}$ by $M_{Z,p}$.  
 
 We fix $\e>0$. It is 
 enough to show that if  $z\in \bC,\Re(z)\l 
 -c_{0}+\delta p -\e$,  $	L_{p, Z}-z $ is a Fredholm operator from
 $\Dom(L_{p, Z})$ into $\Omega_{2}\left(Y,F\right)$.
 
 By \eqref{eq:eqZgde} and 
 \eqref{eq:eqELZ}, outside $U_{3}\cup \ul{\mathcal{C}}^{*}_{0}$, we have 
	\begin{align}
	\sigma\(\[e^{p G}, L_{Z}\] e^{-p G}\) \g  \delta p. 
\end{align}

Let $A\in \Psi^{0}\left(Y,\Lambda\left(T^{*}Y\right)  \otimes_{\R}F\right)$ be such that its  principal symbol  is scalar and supported in a 
neighborhood of $U_{3}\cup \ul{\mathcal{C}}^{*}_{0}$  and is equal to a  
constant $\kappa>0$ on
$U_{3}\cup \ul{\mathcal{C}}^{*}_{0}$. Using (\ref{eq:res2}), (\ref{eq:eqZgde}), 
by taking $\kappa$ large enough, for any $p\ge 0$, 
\begin{align}\label{eq:coma}
	\sigma\(\[e^{p G}, L_{Z}\] e^{-p G}+\delta p A\) \ge \delta p. 
\end{align} 
By (\ref{eq:eqc0}), (\ref{eq:eqLZG1}), and (\ref{eq:coma}), and by the sharp
Gårding inequality, there is $a_{\delta, 
p}>0$ such that for $\Re(z)\le -c_{0}+\delta p -\e$, and $u\in \Omega(Y, F)$, we have 
\begin{align}\label{eq:doma2}
\Re \<\(	M_{p, Z}+\delta p A-z\) u, u \> \ge 
\e \|u\|_{\Omega_{2}\left(Y,F\right)}^{2} - a_{\delta, 
p}\|u\|_{\cH^{-1/2}\left(Y,\Lambda\left(T^{*}Y\right) \otimes _{\R}F\right)}^{2}. 
\end{align}
The last term in the right-hand side of (\ref{eq:doma2}) comes from 
the fact that the difference between the full symbol of $\left[e^{pG},L_{Z}\right]e^{-pG}$ and its 
principal symbol is an operator of order $-1$.
Therefore, there is $A_{\delta, p}\in \Psi^{0}\left(Y,\Lambda\left(T^{*}Y\right) \otimes _{\R}F\right)$ with compactly supported principal symbol such that 
\begin{align}\label{eq:dom3}
	\Re \<\(	M_{Z,p}+\delta p A+ A_{\delta, p}-z\) u, u \> \ge 
\frac{\e}{2} \|u\|_{\Omega_{2}\left(Y,F\right)}^{2}.  
\end{align} 

Put
\begin{equation}\label{eq:dom3z-1}
B_{\delta, p}=\delta p A+ A_{\delta, p}.
\end{equation}
The operator 
$B_{\delta,p}$ is a bounded operator on $\Omega_{2}\left(Y,F\right)$. 
It induces a corresponding bounded operator from $\mathrm{dom} M_{Z,p}$ into 
$\Omega_{2}\left(Y,F\right)$. Put
\begin{equation}\label{eq:dom0}
C_{\delta,p}=M_{Z,p}+B_{\delta,p}.
\end{equation}
The operator
$C_{\delta,p}$ is a bounded operator from 
$\mathrm{dom}M_{Z,p}$ into $\Omega_{2}\left(Y,F\right)$. By 
(\ref{eq:dom3}), $C_{\delta,p}-z$ is injective and its image is closed. Also (\ref{eq:dom3}) 
can be rewritten in the form,
\begin{equation}\label{eq:dom3x1}
\Re \<\(	C_{\delta,p}-z\) u, u \> \ge 
\frac{\e}{2} \|u\|_{\Omega_{2}\left(Y,F\right)}^{2}.  
\end{equation}

The operator $C_{\delta,p}$ is a pseudo-differential operator of 
order $1$ with scalar principal symbol. Therefore it has an unambiguously defined formal $L_{2}$-adjoint 
$C^{*}_{\delta,p}$.  Inspection of (\ref{eq:eqLZG1}) shows that 
the structure of $C^{*}_{\delta,p}$ is essentially the same as the 
structure of $C_{\delta,p}$. This
argument shows that $C^{*}_{\delta,p}-\overline{z}$ is also injective. Therefore 
the image of $C_{\delta,p}-z$ coincides with
$\Omega_{2}\left(Y,F\right)$. By (\ref{eq:dom3x1}), we conclude that 
$C_{\delta,p}-z$ is a bounded invertible operator from 
$\mathrm{dom}\,M_{Z,p}$ into $\Omega_{2}\left(Y,F\right)$.
Therefore $\(C_{\delta,p}-z\)^{-1}$ is a 
holomorphic family of bounded operators parametrized by $z\in \bC$ 
with $\Re z\le -c_{0}+\delta p-\e$. 

The family of operators $M_{Z,p}-z$ depends holomorphically on $z$, 
and moreover, for $\mathrm{Re}\,z \ll -1$, it is invertible. We will 
now prove that for $\mathrm{Re}\,z\le -c_{0}+\delta p-\epsilon$, it is a 
family of Fredholm operators.
Observe that
\begin{align}\label{eq:song1}
	M_{Z,p}-z= \(1- B_{\delta, p}\left( C_{\delta,p}-z \right)^{-1}\) \(C_{\delta,p}-z\). 
\end{align} 
It remains to show $B_{\delta, p}\(C_{\delta,p}-z\)^{-1}$ is a 
compact operator from 
	$\Omega_{2}\left(Y,F\right)$ into itself. By (\ref{eq:dom3z-1}), we can write 
	$B_{\delta,p}$ in the form 
	\begin{equation}\label{eq:doma4}
B_{\delta, p}= B_{\delta, p,1}+  B_{\delta, p,2}, 
\end{equation}
	where 
the principal symbol of $B_{\delta, p,1}$ is compactly supported, 
and the principal symbol of $B_{\delta, p,2}$ is supported away from 
$0$ and near $\underline{\mathcal{C}}_{0}^{*}$. 

The operator  $B_{\delta, p,1}$ is a 
regularizing 
	operator, and so it is  compact. 
	We claim that $B_{\delta, 
	p,2}\(C_{\delta,p}-z\)^{-1}$ is a pseudo-differential 
	operator of order $-1$. Indeed $B_{\delta,p,2}$ is a 
	pseudo-differential operator of order $0$. Moreover, on the 
	support of the principal symbol of $B_{\delta,p,2}$, 
	$\left(C_{\delta,p}-z\right)^{-1}$ is a classical elliptic
	pseudo-differential operator of order $-1$, which proves our claim. In 
	particular, $B_{\delta, 
	p,2}\(C_{\delta,p}-z\)^{-1}$ is also a compact operator, which 
	concludes the proof of the fact that 
	$B_{\delta,p}\left(C_{\delta,p}-z\right)^{-1}$ is compact, and 
	completes the proof of our theorem.
 \end{proof} 
\subsection{The resolvent of $L_{Z}$}%
\label{subsec:relz}
If $p\ge 0$, by Theorem \ref{thm:FaureSjostrand},  the 
operator 
\begin{align}\label{eq:eqsLz}
	(L_{Z, p}-z)^{-1}: \cH_{pG}\left(Y,\Lambda(T^{*}Y)\otimes_{\R} F\right)\to \cH_{pG}\left(Y,\Lambda(T^{*}Y)\otimes_{\R} F\right)
\end{align} 
is holomorphic in the variable $\{z\in \mathbf{C}: \Re z<-c_{p}\}$. By 
analytic Fredholm theory \cite[Appendix C, Theorem 
C.8]{DyatlovZworski19}, $\left(L_{Z,p}-z\right)^{-1}$ has a meromorphic extension by Fredholm operators to $\{z\in \mathbf{C}: 
\Re z < -c_{0} + p\delta\}$. That means that  for any $z_{0}\in 
\bC$ with $\Re z_{0} <-c_{0}+p\delta$, near $z_{0}$, there are
finite ranks operators $B_{a}$ on 
$\cH_{pG}\left(Y,\Lambda\left(T^{*}Y\right) \otimes _{\R}F \right) $ 
with $1\le a\le A_{Z,p,z_{0}}$ such that 
\begin{align}\label{eq:eqAfred}
	(L_{Z, p}-z)^{-1}-\sum_{a=1}^{A_{Z,p,z_{0}}}(z-z_{0})^{-a} 	B_{a}. 
\end{align} 
is holomorphic near $z_{0}$. Formula can be made more precise. Assume 
that $z_{0}$ is a pole for $\left(L_{Z,p}-z\right)^{-1}$. Let 
\index[not]{DZp@$\cD_{Z,p,z_{0}}(Y,F)$}%
$\cD_{Z,p,z_{0}}(Y,F)$ be the 
finite dimensional characteristic subspace for the eigenvalue 
$z_{0}$. Put
\index[not]{eZp@$e_{Z,p,z_{0}}$}%
\begin{equation}\label{eq:doma6}
e_{Z,p,z_{0}}=\dim \cD_{Z,p,z_{0}}(Y, 
F)
\end{equation}
Then
\begin{equation}\label{eq:doma8}
\left(L_{Z,p}-z_{0}\right)^{e_{Z,p,z_{0}}}\vert_{\cD_{Z,p,z_{0}}(Y, 
F)}=0.
\end{equation}

 Let 
 \index[not]{PZp@$P_{Z,p,z_{0}}$}%
 $P_{Z,p,z_{0}}$ be the spectral projector on $\cD_{Z,p,z_{0}}(Y, 
F)$. Then 
\begin{equation}\label{eq:doma7}
\left(L_{Z,p}-z\right)^{-1}+P_{Z,p,z_{0}}\sum_{0}^{e_{Z,p,z_{0}}}\left(L_{Z,p}-z_{0}\right)^{i}\left(z-z_{0}\right)^{-(i+1)}
\end{equation}
is holomorphic near $z_{0}$.

By restriction, we get a  
meromorphic family of continuous operators 
\begin{align}\label{eq:eqsLz2}
		(L_{Z,p}-z)^{-1}: \Omega(Y, F)\to \cD(Y, F), z\in \mathbf{C}: 
\Re z<-c_{0}+p\delta.
\end{align}

We have the result of Faure-Sjöstrand \cite[Theorem 1.5]{FaureSjostrand11}. 
\begin{proposition}\label{prop:propextReso}
	Given $ p\ge p'\ge 0$, we have an identity of meromorphic families of continuous  
	operators from $\Omega(Y, F)$ into 
$\cD(Y, F)$ parametrized by $\{z\in  \mathbf{C}: 
\Re z < -c_{0}+p'\delta\}$, 
\begin{align}\label{eq:eqLppp}
(L_{Z,p}-z)^{-1}=(L_{Z,p'}-z)^{-1}.	
\end{align}
\end{proposition} 
\begin{proof}
By the uniqueness of meromorphic extensions, it is enough to establish 
\eqref{eq:eqLppp} when $p'=0$ and  $\Re z\ll -1$. 

Recall that $m$ takes its values in $[-2, 2]$. For $p\ge 0$, then 
\begin{align}\label{eq:eq2p1}
&	\cH^{2p+1}_{pG}\left(Y,\Lambda(T^{*}Y)\otimes_{\R} F\right)\subset 
\Omega_{2}\left(Y,F\right), \\
& \cH^{2p+1}_{pG}\left(Y,\Lambda(T^{*}Y)\otimes_{\R} F\right)\subset 
\cH_{pG}\left(Y,\Lambda(T^{*}Y)\otimes_{\R} F\right). \nonumber 
\end{align} 

Let $\( L_{Z,p,2p+1},\mathrm{dom}(L_{Z,p,2p+1})\)$ be the unbounded 
operator $L_{Z}$ acting on  $\cH^{2p+1}_{pG}$.  For $\Re z\ll -1$, we 
claim that the resolvent  $\( L_{Z,p,2p+1}-z\)^{-1}$ is  
well-defined. Indeed, we 
already did this when replacing $\mathcal{H}^{2p+1}_{pG}$ by 
$\mathcal{H}_{pG}$. When adding the extra Sobolev weight $2p+1$, the 
argument is exactly the same.

 If 
$u\in \Omega(Y, F)$, then  
\begin{align}\label{eq:ziga1}
&u\in \cH^{2p+1}_{pG}\left(Y,\Lambda\left(T^{*}Y\right) \otimes _{\R}F\right),\\
&	\(L_{Z,p,2p+1}-z\)^{-1}u\in \cH^{2p+1}_{pG}\left(Y,\Lambda\left(T^{*}Y\right) \otimes _{\R}F\right). \notag 
\end{align} 
By \eqref{eq:eq2p1}, we know that for $\Re z\ll -1$, 
\begin{align}
	\( L_{Z,p,2p+1}-z\)^{-1}u=\(
	L_{Z, p}-z\)^{-1}u=\(L_{Z,p'}-z\)^{-1}u,
\end{align} 
which implies our proposition. 
\end{proof} 

By Proposition \ref{prop:propextReso},  the family of 
operators in  
\eqref{eq:eqhlo1} extends to a meromorphic family on $\bC$. 
This extension does not depend on the choice of the escape function  $g$ 
or on the pseudo-differential operator  $G$. Recall that by 
(\ref{eq:eqZgde}), $\delta$ depends on $g$. Using  \eqref{eq:eqDs}, 
given $a>0$, we have a meromorphic family of operators,
\begin{multline}\label{eq:pox}
	(L_{Z}-z)^{-1}: \Omega(Y, F)\to\bigcap_{g} \bigcap_{p\ge \left(a+c_{0}+1\right)/\delta} \cH_{pG}(Y, 
	\Lambda(T^{*}Y)\otimes_{\R} F) \\
	\subset \cD_{\mathcal{Y}^{*}_{s}}(Y, F), z\in\C, \Re z<a.
\end{multline}

In the sequel, we will denote by 
\index[not]{SpL@$\mathrm{Sp}\, L_{Z}$}%
$\mathrm{Sp}\, L_{Z}$ the union in 
$p\ge 0$ of the discrete spectrums of the $L_{Z,p}$. 
\section{Poincaré duality, spectral 
projectors, and the truncated Fried zeta function}%
\label{sec:frize}
The purpose of this section is to explore various aspects of Poincaré 
duality associated with the operator $L_{Z}$. We study in detail the 
wave front set properties of its resolvent, and of the associated 
spectral projectors, and also the truncated spectral 
projectors associated with eigenvalues contained in the disk of 
radius $a>0$. We also construct the truncated Fried zeta function 
$R_{Z,F,>a}\left(\sigma\right)$. Our constructions rely on previous results by Dyatlov-Zworski \cite{DyatlovZworski16}. Some of our results were previously obtained by Chaubet-Dang \cite{ChaubetDang24} in the context of contact manifolds.

This section is organized as follows. In Subsection 
\ref{subsec:podan}, we connect the construction of the anisotropic 
Sobolev spaces with  Poincaré duality.

In Subsection \ref{subsec:micr}, by following Dyatlov-Zworski \cite{DyatlovZworski16}, we describe the wave front set of  the resolvent $\left(L_{Z}-z\right)^{-1}$.

In Subsection \ref{subsec:sibi}, for $\sigma\in\C,\Re\,\sigma\gg 1$, we relate the intersection supertrace of $N\left(L_{Z}+\sigma\right)^{-1}$ with  the logarithmic derivative of $R_{Z,F}\left(\sigma\right)$.

In Subsection \ref{subsec:resd}, we show that the resolvent commutes 
with the operators $d^{Y},i_{Z}$.

In Subsection \ref{subsec:sptru}, if $z_{0}\in \mathrm{Sp}L_{Z}$, we 
consider the associated characteristic space 
\index[not]{DZz@$\mathcal{D}_{Z,z_{0}}\left(Y,F\right)$}%
$\mathcal{D}_{Z,z_{0}}\left(Y,F\right)$ and 
 the corresponding spectral projector 
 \index[not]{PZz@$P_{Z,z_{0}}$}%
 $P_{Z,z_{0}}$. Also we study 
 the wave front set of $P_{Z,z_{0}}$.

 In Subsection \ref{subsec:sptrdu}, we study the effect of Poincaré duality on the spectral projectors.
 
 In Subsection \ref{subsec:vspa}, we introduce the truncated direct 
 sum of characteristic subspaces
 \index[not]{DZa@$\mathcal{D}_{Z,<a}\left(Y,F\right)$}%
 $\mathcal{D}_{Z,<a}\left(Y,F\right)$, and the associated projectors 
 \index[not]{PZa@$P_{Z,<a}$}%
 $P_{Z,<a}$.
 
 In Subsection \ref{subsec:trqi}, we prove  that the 
 embedding 
 $$\left(\Omega\left(Y,F\right),d^{Y}\right)\to \left( 
 \mathcal{D}_{\mathcal{Y}^{*}_{s}}\left(Y,F\right),d^{Y} \right)$$ is 
 a quasi-isomorphism.
 
 In Subsection \ref{subsec:truqi}, we show that $P_{Z,<a}$ is a 
 quasi-isomorphism.  Also we study the intersection supertraces of 
 operators involving the resolvent $\left(L_{Z}-z\right)^{-1}$.
 
 In Subsection \ref{subsec:trstr}, we study the truncated 
 intersection supertraces associated with the resolvent and the 
 operator $N$.

In Subsection \ref{subsec:mero}, following Dyatlov-Sworski \cite[p. 561]{DyatlovZworski16}, we reprove the fact that 
 $R_{Z,F}\left(\sigma\right)$ extends to a meromorphic function.
 
 Finally, in Subsection \ref{subsec:truz}, we define the truncated Fried zeta 
 function $R_{Z,F,>a}\left(\sigma\right)$.

 In this section, we make the same assumptions as in Sections \ref{sec:frizet}, \ref{sec:avf}, and we use the notation of Section \ref{sec:opcur}.
\subsection{Poincaré duality and anisotropic Sobolev spaces}%
\label{subsec:podan}
The operator $d^{Y},i_{Z}, L_{Z}$ acts on $\Omega\left(Y,F\right)$ and on 
$\mathcal{D}\left(Y,\Ft\right)$. By (\ref{eq:normi4}), we deduce 
easily that
\begin{align}\label{eq:normi5a1}
&d^{Y\dag}=d^{Y}, &i_{Z}^{\dag}=i_{-Z}, \qquad L_{Z}^{\dag}=L_{-Z}.
\end{align}
The last identity in (\ref{eq:normi5a1}) should be expected, because 
$L_{Z}$ preserves the degree. It is also a consequence of the first 
two identities.

Let 
\index[not]{G@$G^{\dag}$}%
$G^{\dag}\in \Psi^{0^{+}}\left(Y,\Lambda\left(T^{*}Y\right) 
\otimes_{\R}\Ft \right) $ be the transpose 
of $G$ with respect to  Poincaré duality in (\ref{eq:normi3}). Let
\index[not]{g@$g^{\dag}$}%
$g^{\dag}:\mathcal{Y}^{*}\to\R$  be given by
\begin{equation}\label{eq:vac1}
g^{\dag}\left(y,q\right)=g\left(y,-q\right).
\end{equation}

By 
(\ref{eq:normi0}), we get
\begin{equation}\label{eq:normi6}
\sigma\left(- G^{\dag}\right)=-g^{\dag}.
\end{equation}

Recall that our cones in $\mathcal{Y}^{*}\setminus Y$ are not assumed to be symmetric with respect to the involution $\left(y,q\right)\to \left(y,-q\right)$. We will denote with an extra $\dag$ the images of our cones by this involution.

By Remarks \ref{rem:neg} and \ref{rem:chs}, $-G^{\dag}$ verifies the same 
assumptions with respect to  $-Z$ and the dashed cones   as 
$G$ with respect to $Z$. The same considerations apply to 
$\exp\left(-p G^{\dag}\right)\in 
\Psi^{\left(-pm\right)^{+}}\left(Y,\Lambda\left(T^{*}Y\right) 
\otimes _{\R}\Ft \right) $.
There is an induced Poincaré   duality of  Hilbert spaces 
\begin{multline}
	\cH_{p G}(Y,\Lambda^{\bullet}(T^{*}Y)\otimes_{\R} F)\\
	\times \cH_{-p G^{\dag}}(Y, 
	\Lambda^{n-\bullet}(T^{*}Y)\otimes_{\R} F^{*} \otimes 
	_{\Z_{2}}\mathrm{o}\left(TY\right))\to \mathbf{C}. 
\end{multline} 
The  operator $L_{-Z}$\footnote{It will be useful to distinguish 
$-L_{Z}$ and $L_{-Z}$. These operators  are made here to act on different spaces.} acts on $\Omega\left(Y,F^{*}\otimes 
_{\Z_{2}}\mathrm{o}\left(TY\right)\right)$ and   $\mathcal{D}\left(Y,F^{*}\otimes 
_{\Z_{2}}\mathrm{o}\left(TY\right)\right)$. 
The unbounded operator $L_{-Z}\vert_{\mathcal{H}_{-p 
G^{\dag}}\left(Y,\Lambda^{n-\Ou}\left(T^{*}Y\right) \otimes 
_{\R}\Ft\right)}$ is  the transpose of the unbounded operator 
$L_{Z}\vert_{\mathcal{H}_{pG}\left(Y,\Lambda^{\Ou}\left(T^{*}Y \right)  \otimes_{\R} 
F\right)}$.
Therefore the operators $L_{Z,p}$ 
and $L_{-Z,p}$ are isospectral. 

Recall that by (\ref{eq:eqZgde}), $\delta$ depends on $g$. The analog of (\ref{eq:pox}) says that given $a>0$, for $z\in \C,\Re z<a$, we 
have a meromorphic family
\begin{multline}\label{eq:pox1}
	(L_{-Z}-z)^{-1}: \Omega(Y, \Ft)\\
	\to \bigcap_{g}\bigcap_{p\ge 
	\left(a+c_{0}+1\right)/\delta} \cH_{-p G^{\dag}}(Y, 
	\Lambda(T^{*}Y)\otimes_{\R} \Ft) \\
	\subset \cD_{\mathcal{Y}^{*}_{u}}(Y, \Ft).
\end{multline}

As we saw before, the operators $L_{Z}$ and $L_{-Z}$ are isospectral. The meromorphic 
family $\left(L_{-Z}-z\right)^{-1}$ is the transpose of the family 
$\left(L_{Z}-z\right)^{-1}$.

By (\ref{eq:tubi0}), if 
$\varphi\in\mathcal{D}_{\underline{\mathcal{C}}^{*}_{s}}\left(Y,F\right), 
\psi\in\mathcal{D}_{\underline{\mathcal{C}}^{*\dag}_{u}}\left(Y,\Ft\right)$,  for 
$p\in \R$ large enough,
\begin{align}\label{eq:casa1}
&\varphi\in \mathcal{H}_{pG}\left(Y,\Lambda\left(T^{*}Y\right) 
\otimes _{\R}F\right),
&\psi\in\mathcal{H}_{-p 
G^{\dag}}\left(Y,\Lambda\left(T^{*}Y\right)\otimes _{\R}\Ft\right) .
\end{align}
Using Theorem \ref{thm:FaureSjostrand}, we find that 
$\left(L_{Z}-z\right)^{-1}\varphi$ is a meromorphic function with 
values in $\mathcal{H}_{pG}\left(Y,\Lambda\left(T^{*}Y\right) \otimes 
_{\R}F\right)$. In particular $\left\langle  
\left(L_{Z}-z\right)^{-1}\varphi,\psi\right\rangle$ is a meromorphic 
function, and its poles are included in $\mathrm{Sp}L_{Z}$. Here, the 
dependence on $p$ has disappeared.
\subsection{The microlocal properties of the resolvent}%
\label{subsec:micr}
Recall that $\mathcal{Y}^{*}$ is equipped with the symplectic form 
$\omega^{\mathcal{Y}^{*}}$. Put
\index[not]{oY@$\omega^{\mathcal{Y}^{*}\times\mathcal{Y}^{*}}$}%
\begin{equation}\label{eq:gonf0}
	\omega^{\mathcal{Y}^{*}\times\mathcal{Y}^{*}}=\pi_{1}^{*}\omega^{\mathcal{Y}^{*}}+\pi_{2}^{*}\omega^{\mathcal{Y}^{*}}.
\end{equation}
Then $\omega^{\mathcal{Y}^{*}\times\mathcal{Y}^{*}}$ is a symplectic form on $\mathcal{Y}^{¨*}\times\mathcal{Y}^{*}$.

Recall that $M,M^{*}$ were defined in (\ref{eq:duaz3a1}), 
(\ref{eq:gonf1}).
Let 
\index[not]{M@$\mathcal{M}^{*}$}%
$\mathcal{M}^{*} \subset \mathcal{Y}^{*}$ be the total space of 
$M^{*}$.

As we already explained, the flow $\varphi_{t}\vert_{t\in\R}$ on $Y$ associated with the vector field $Z$ lifts to the flow $\varphi_{t*}\vert_{t\in\R}$  on $\mathcal{Y}^{*}$ associated with the vector field $\mathcal{Z}$. 

Recall that the flow $\varphi_{t*}\vert_{t\in\R}$  preserves 
$T^{*}_{u}Y$ and $T^{*}_{s}Y$, and so it preserves $\mathcal{M}^{*}$. 
 
Given $t\in \R$, let 
\index[not]{Grf@$\mathrm{Gr}\varphi_{t}$}%
$\mathrm{Gr}\varphi_{t} \subset Y\times Y$ be the graph of 
$\varphi_{t}$. Let $N^{*}_{\mathrm{Gr}\varphi_{t}/\left(Y\times Y\right)} $ be its conormal 
bundle.  Its total space 
\index[not]{NGr@$\mathcal{N}^{*}_{\mathrm{Gr}\varphi_{t}/\left( Y\times Y\right)}$}%
$\mathcal{N}^{*}_{\mathrm{Gr}\varphi_{t}/\left( Y\times Y\right)}$ is a Lagrangian submanifold of 
$\mathcal{Y}^{*}\times \mathcal{Y}^{*}$. More precisely, using the definition of $r_{-1}$ after (\ref{eq:duaz5}), we get
\begin{equation}\label{eq:gagari1}
	\mathcal{N}^{*}_{\mathrm{Gr}\varphi_{t}/\left( 
		Y\times Y\right)}=\mathrm{Gr}\left[r_{-1}\varphi_{t*}\right].
\end{equation}

We view $\mathrm{Gr}\varphi_{t}$ as a distribution in $\mathcal{D}\left(Y\times 
Y ,\pi_{1}^{*}F \otimes\pi_{2}^{*}\left(F^{*} \otimes 
_{\Z_{2}}\mathrm{o}\left(TY\right)\right)\right)$. Then
\begin{equation}\label{eq:gaga0}
\mathrm{WF}\left(\mathrm{Gr}\varphi_{t}\right)=\mathcal{N}^{*}_{\mathrm{Gr}\varphi_{t}/ 
Y\times Y}\setminus \left(Y\times Y\right).
\end{equation}
\begin{definition}\label{def:omr}
	For $r\in \R$, put
	\index[not]{Or@$\Omega_{\ge r}$}%
	\begin{equation}\label{eq:gaga1}
\Omega_{\ge r}=\left\{\left(m,m'\right)\in \mathcal{M}^{*} \times 
\mathcal{M}^{*}, \mathrm{there\,exists}\, t\ge r,\varphi_{-t*}m+m'=0\right\}.
\end{equation}
\end{definition}
If $m=\left(y,q\right), m'=\left(y',q'\right)$, with $q\in M^{*}_{y}, q'\in M^{*}_{y'}$,  the condition in (\ref{eq:gaga1}) says that
\begin{align}\label{eq:gaga1z1}
&y'=\varphi_{-t}y, &q'+ \varphi_{-t*}q=0.
\end{align}

Let $\overline{\Omega_{\ge r}\setminus \left(Y\times Y\right)}$ be the 
closure of $\Omega_{\ge r}\setminus \left(Y\times Y\right)$ in 
$\left( \mathcal{M}^{*}\times\mathcal{M}^{*}\right) \setminus 
\left(Y\times Y\right)$.
\begin{proposition}\label{prop:b1}
	The set $\Omega_{\ge r}\setminus\left(Y\times Y\right)\cup 
	\left(\mathcal{Y}^{*}_{s}\times 
	\mathcal{Y}^{*}_{u} \right) \setminus \left(Y\times Y\right)$ is closed in $\left(\mathcal{Y}^{*}\times\mathcal{Y}^{*}\right)\setminus Y\times Y$. In particular, 
	\begin{equation}\label{eq:gaga2}
\overline{\Omega_{\ge r}\setminus \left(Y\times 
Y\right)} \subset \Omega_{\ge r}\setminus\left(Y\times Y\right)\cup 
	\left(\mathcal{Y}^{*}_{s}\times
	\mathcal{Y}^{*}_{u} \right) \setminus \left(Y\times Y\right).
\end{equation}
 Moreover, there exists $C>0$ such that if $y\in \mathbf{F}$, if $t>0$ is such that $\varphi_{-t}y=y$, then
\begin{equation}\label{eq:singa4}
\left\vert \det\left(1-\varphi_{-t*}\right)\vert_{M_{y}}\right\vert\ge C.
\end{equation}
\end{proposition}
\begin{proof}
	We fix $t\ge r$. If $q\in M^{*}$, we can write $q$ in the form
	\begin{align}\label{eq:gaga3}
&q=q_{u} + q_{s}, &q_{u}\in M^{*}_{u},\qquad q_{s}\in M^{*}_{s}.
\end{align}
If $q,q'$ verify the conditions in (\ref{eq:gaga1z1}), then
\begin{align}\label{eq:guaga5}
&q=-\varphi_{t*}q'_{u}+q_{s},&q'=q'_{u}-\varphi_{-t*}q_{s}.
\end{align}

As $t\to + \infty $, the norm of $\varphi_{t*}$ tends 
uniformly to $0$ on $T^{*}_{u}Y$, and the norm of 
$\varphi_{-t*}$ tends uniformly to $0$ on 
$T^{*}_{s}Y$. Using (\ref{eq:guaga5}), we obtain the first part of our proposition.	

To establish the second part, it is enough to show that given $T>0$ large enough, there is $C>0$ such that if $t\ge T$, for $f\in M_{y}$, 
\begin{equation}\label{eq:singa5}
\left\vert \left( \varphi_{-t*}-1\right)f\right\vert\ge C\left\vert f\right\vert.
\end{equation}
Since $\varphi_{-t*}$ preserves the splitting (\ref{eq:duaz3a1}) of $M$, we may as well assume that $f\in T_{u,y}Y$ or $f\in T_{s,y}Y$. In any of these two cases, since $\varphi_{-t*}$ contracts or expands uniformly, equation (\ref{eq:singa5}) is obvious. The proof of our theorem is complete.
\end{proof}

Let $\mathcal{N}^{*}_{\Delta^{Y}/\left(Y\times Y\right)}$ be the 
total space of $N^{*}_{\Delta^{Y}/\left(Y\times Y\right)}$.
\begin{definition}\label{def:enla}
	For $r\in \R$, put
	\index[not]{Tr@$\Theta_{\ge r}$}%
	\begin{equation}\label{eq:sig2}
\Theta_{\ge r}=\mathcal{N}^{*}_{\Delta^{Y}/\left(Y\times Y\right)}\cup \Omega_{\ge 
r}\cup \left( 
\mathcal{Y}^{*}_{s}\times\mathcal{Y}^{*}_{u} 
\right) .
\end{equation}
\end{definition}
By Proposition \ref{prop:b1}, $\Theta_{\ge 
r}\setminus\left(Y\times Y\right)$ is 
closed in $\left(\mathcal{Y}^{*}\times 
\mathcal{Y}^{*}\right)\setminus\left(Y\times Y\right)$.
\begin{definition}\label{def:alg}
	Let 
	\index[not]{Dri@$\mathcal{D}_{r,i}$}%
	$\mathcal{D}_{r,1},\mathcal{D}_{r,2}$ be the spaces of   
	currents in 
	$$\mathcal{D}^{n}\left(Y\times Y,\pi_{1}^{*}F \otimes 
	\pi_{2}^{*}\left( \Ft \right) \right)$$
	whose wave front set 
	is included in $\Omega_{\ge r}\setminus\left(Y\times Y\right)\cup 
	\left(\mathcal{Y}^{*}_{s}\times 
	\mathcal{Y}^{*}_{u} \right) \setminus \left(Y\times Y\right),\Theta_{\ge r}\setminus\left(Y\times Y\right)$.
\end{definition}

We claim that if $T\in \mathcal{D}_{r,i}, i=1,2$, and if $u\in 
\mathcal{D}_{\mathcal{Y}^{*}_{s}}\left(Y,F\right)$, the product of currents 
$T\we\pi_{2}^{*}u$ is well-defined. This is essentially because 
$T_{u}^{*}Y\cap T^{*}_{s}Y=0$, and because the flow 
$\varphi_{t*}\vert_{t\in\R}$ 
preserves $T^{*}_{u}Y,T^{*}_{s}Y$.  We can then define $M_{T}u$ as in 
(\ref{eq:rec1}).
\begin{theorem}\label{thm:cur}
	Given $r\in\R$, if $T\in \mathcal{D}_{r,i},i=1,2$, if $u\in 
	\mathcal{D}_{\mathcal{Y}^{*}_{s}}\left(Y,F\right)$, then $M_{T}u\in 
	\mathcal{D}_{\mathcal{Y}^{*}_{s}}\left(Y,F\right)$. The operator 
	$M_{T}$ is an 
	operator acting on 
	$\mathcal{D}_{\mathcal{Y}^{*}_{s}}\left(Y,F\right)$ and preserving the degree. Moreover, if
	$r_{1}\ge 0,r_{2}\ge 0$, if $T_{1}\in 
	\mathcal{D}_{r_{1},i},T_{2}\in \mathcal{D}_{r_{2},i}$, for $i=1$,  
	$T_{1}\circ T_{2}\in \mathcal{D}_{r_{1}+r_{2},1}$, and for $i=2$, 
	$T_{1}\circ T_{2}\in 
	\mathcal{D}_{\inf\left(r_{1},r_{2}\right),2}$. In particular the 
	currents in $\mathcal{D}_{0,i}\vert_{i=1,2}$ form an algebra.
\end{theorem}
\begin{proof}
	The fact that $M_{T}u\in 
	\mathcal{D}_{\mathcal{Y}^{*}_{s}}\left(Y,F\right)$ is an obvious 
	consequence of (\ref{eq:gaga2}), (\ref{eq:sig2}). By the 
	considerations that follow (\ref{eq:rec1}), $M_{T}$ preserves the 
	degree. The last statements can be proved easily.
\end{proof}

By proceeding as in the proof of (\ref{eq:sacr0}), we have the identity of operators,
\begin{equation}\label{eq:bas1}
e^{-rL_{Z}}=N_{\mathrm{Gr}\varphi_{-r}}.
\end{equation}

If $T\in \mathcal{D}_{0,2}$,  if $ r\ge 0$, the composition 
$\mathrm{Gr}\varphi_{-r}\circ T$ is a 
well-defined element of $\mathcal{D}^{n}\left(Y\times Y,\pi_{1}^{*}F \otimes\pi_{2}^{*}\left(F^{*} \otimes 
_{\Z_{2}}\mathrm{o}\left(TY\right)\right)\right)$. Also
\begin{equation}\label{eq:pdz2}
\mathrm{WF}\left( \mathrm{Gr}\varphi_{-r}\circ T \right) \subset 
\left(\mathcal{N}^{*}_{\mathrm{Gr}\varphi_{-r}/Y\times Y}\cup\Omega_{\ge r}\cup 
	\mathcal{Y}^{*}_{s}\times 
	\mathcal{Y}^{*}_{u} \right) \setminus \left(Y\times Y\right). 
\end{equation}

Finally, by Proposition \ref{prop:compo} and by (\ref{eq:bas1}), we get
\begin{equation}\label{eq:bas2}
e^{-rL_{Z}}N_{T}=N_{\mathrm{Gr}\varphi_{-r}\circ T}.
\end{equation}

Let 
\index[not]{r@$r_{0}$}%
$r_{0}\in \R_{+}^{*}$ be given by
\begin{equation}\label{eq:bas2a1}
r_{0}=\inf_{y\in\mathbf{F}}t_{y}.
\end{equation}

\begin{proposition}\label{prop:wff}
	For $0<r<r_{0}$, if $T\in \mathcal{D}_{0,2}$,  
	\begin{equation}\label{eq:int1}
\mathrm{WF} \left( \mathrm{Gr}\varphi_{-r}\circ T \right) \cap 
\mathcal{N}^{*}_{\Delta^{Y}/\left(Y\times Y\right)}\setminus 
\left(Y\times Y\right)=\emptyset.
\end{equation}
\end{proposition}
\begin{proof}
	We will use equation (\ref{eq:pdz2}). If $r$ is taken as indicated, $\varphi_{-r}$ 
	has no fixed point, which takes care of the first term in the 
	right-hand side of (\ref{eq:pdz2}). Also for $t>0$, even if 
	$\varphi_{-t}$ may have fixed points, $\varphi_{-t}^{*}$ has no 
	fixed point in $\mathcal{M}^{*}\setminus Y$. 
	This takes care of the second term in the right-hand side of 
	(\ref{eq:pdz2}). The proof of our proposition is complete. 
\end{proof}

By (\ref{eq:pox}), $\left(L_{Z}-z\right)^{-1}$ defines a meromorphic family of currents on $Y\times Y$. More precisely, if $z\notin\mathrm{Sp}L_{Z}$, there is a current $T_{z}$ such that we have the identity of operators,
\begin{equation}\label{eq:bost1}
	\left(L_{Z}-z\right)^{-1}=N_{T_{z}}.
\end{equation}
In the sequel, we will not distinguish between $\left(L_{Z}-z\right)^{-1}$ and $T_{z}$, so that $\left(L_{Z}-z\right)^{-1}$ is viewed both as a current and an operator.

Now we state 
 a fundamental result of Dyatlov-Zworski \cite[Proposition 
3.3]{DyatlovZworski19}.
\begin{theorem}\label{thm:pdz1}
If $z\in\C\setminus \mathrm{Sp}L_{Z}$, then 
\begin{equation}\label{eq:sing}
\left(L_{Z}-z\right)^{-1}\in \mathcal{D}_{0,2}.
\end{equation}
\end{theorem}


By (\ref{eq:pox1}),  we have a meromorphic family 
of operators 
\begin{equation}\label{eq:pox1a}
\left(L_{-Z}-z\right)^{-1}:\Omega\left(Y,\Ft\right)\to \mathcal{D}_{\mathcal{Y}^{*}_{u}}\left(Y,\Ft\right), z\in\C.
\end{equation}

Theorem \ref{thm:pdz1} is 
compatible to Poincaré duality.

Since $\left(L_{-Z}-z\right)^{-1}$ is just 
the transpose of $\left(L_{Z}-z\right)^{-1}$, the corresponding 
currents $\left(L_{-Z}-z\right)^{-1}$ are obtained from the currents 
for $\left(L_{Z}-z\right)^{-1}$ as in (\ref{eq:sv5a6}) with $k=n$. More 
precisely, we have the identity of currents, 
\begin{equation}\label{eq:comer0a}
\left(L_{-Z}-z\right)^{-1}=\left(-1\right)^{n\left(n-1\right)/2}\underline{j}^{*}\left(L_{Z}-z\right)^{-1}.
\end{equation}
\subsection{The Fried zeta function for $\Re\,\sigma\gg 1$}%
\label{subsec:sibi}
For $p\in \R$, let 
$$H^{p}\left(Y\times Y, \pi_{1}^{*}\left( \Lambda\left(T^{*}Y\right) \otimes_{\R}F \right) \otimes _{\R}\pi_{2}^{*}\left(\Lambda\left(T^{*}Y\right) \otimes _{\R}\Ft\right)\right)$$
 be the obvious Sobolev space of index $p$. For $p>n$, this Sobolev space embeds in the corresponding vector space of continuous sections. Moreover, for $r\in \R$, if $s\in \Omega\left(Y,F\right)$, 
\begin{equation}\label{eq:singa1}
e^{-rL_{Z}}s=\varphi^{*}_{-r}s.
\end{equation}
By (\ref{eq:singa1}),  for $p<-n$, $e^{-rL_{Z}}$ defines a distribution on $Y\times Y$ that lies in $H^{p}\left(Y\times Y,\pi_{1}^{*}\left( \Lambda\left(T^{*}Y\right) \otimes _{\R}F \right) \otimes _{\R}\pi_{2}^{*}\left(\Lambda\left(T^{*}Y\right) \otimes _{\R}\Ft\right)\right)$. Let $\left\Vert \,\right\Vert_{p}$ denote a norm on this last vector space. 
\begin{proposition}\label{prop:singa1}
Given $p<-n$, there exist $C_{p}>0,c_{p}>0$ such that
\begin{equation}\label{eq:singa2}
\left\Vert e^{-rL_{Z}}\right\Vert_{p}\le C_{p}e^{c_{p}r}.
\end{equation}
\end{proposition}
\begin{proof}
Let $g^{TY}$ be a Riemannian metric on $TY$, and let $g^{F}$ be a Hermitian metric on $F$. Let $\n^{TY}$ be the Levi-Civita connection on $TY$, and let $\n^{F,u}$ be a unitary connection on $F$. Let $\n^{\Lambda\left(T^{*}Y\right) \otimes _{\R}F,u}$ be the induced connection on $\Lambda\left(T^{*}Y\right) \otimes _{\R}F$. There is a smooth section $A$ of $\End\left(\Lambda\left(T^{*}Y \otimes _{\R}F\right)\right)$ such that
\begin{equation}\label{eq:singa6}
L_{Z}=\n^{\Lambda\left(T^{*}Y\right) \otimes _{\R}F,u}_{Z}+A.
\end{equation}
Let 
\index[not]{CY@$\mathcal{C}\left(Y,\Lambda\left(T^{*}Y\right) \otimes _{\R}F\right)$}%
$\mathcal{C}\left(Y,\Lambda\left(T^{*}Y\right) \otimes _{\R}F\right)$ be the  space of continuous sections of $\Lambda\left(T^{*}Y\right) \otimes _{\R}F$.
By (\ref{eq:singa6}), we conclude that if $s\in \mathcal{C}\left(Y,\Lambda\left(T^{*}Y\right) \otimes _{\R}F\right)$, 
\begin{equation}\label{eq:singa7}
\left\Vert e^{-rL_{Z}}s\right\Vert_{\mathcal{C}\left(Y,\Lambda\left(T^{*}Y\right) \otimes _{\R}F\right)}\le Ce^{cr}\left\Vert s\right\Vert_{\mathcal{C}\left(Y,\Lambda\left(T^{*}Y\right) \otimes _{\R}F\right)}.
\end{equation}
By (\ref{eq:singa7}), we get (\ref{eq:singa2}). The proof of our proposition is complete. 
\end{proof}
\begin{proposition}\label{prop:sobo}
Given $p<-n$, for $\sigma\in\C,\Re\,\sigma\gg 1$, there is $C_{p}>0$ such that 
\begin{equation}\label{eq:singa8}
\left\Vert \left(L_{Z}+\sigma\right)^{-1}\right\Vert_{p}\le C_{p}.
\end{equation}
 \end{proposition}
\begin{proof}
	 By (\ref{eq:song0}), for $\sigma\in\C,\Re\,\sigma\gg 1$, if $s\in\Omega\left(Y,F\right)$, then $\left(L_{Z}+\sigma\right)^{-1}s$ lies in $\Omega_{2}\left(Y,F\right)$. 
We claim that  for $\sigma\in\C,\Re\,\sigma\gg 1$, we have the identity,
\begin{equation}\label{eq:singa9}
\left(L_{Z}+\sigma\right)^{-1}s=\int_{0}^{+ \infty}e^{-u\left(L_{Z}+\sigma\right)}s du.
\end{equation}
To prove (\ref{eq:singa9}), note that by (\ref{eq:singa7}), for $\Re\,\sigma\gg 1$, the right-hand side of (\ref{eq:singa9}) lies $\mathcal{C}\left(Y,\Lambda\left(T^{*}Y\right) \otimes _{\R}F\right)$, and therefore it lies  in $\Omega_{2}\left(Y,F\right)$.  Moreover, we have the identity of distributions, 
\begin{equation}\label{eq:singa9aa1}
\left(L_{Z}+\sigma\right)\int_{0}^{+ \infty}e^{-u\left(L_{Z}+\sigma\right)}s du=s.
\end{equation}
By the definition of $\left(L_{Z}+\sigma\right)^{-1}$, this forces the equality in (\ref{eq:singa9}).
Our proposition now follows from Proposition \ref{prop:singa1} and from (\ref{eq:singa9}).
\end{proof}

  Let $\chi\in C^{\infty}_{\mathrm{c}}\left(\R,\R\right)$ be equal to $1$ near $0$. For $\alpha>0,\sigma\in\C$, set
  \index[not]{Ua@$U_{\alpha,\sigma}$}%
\begin{equation}\label{eq:singa10}
	U_{\alpha,\sigma}=\int_{\R_{+}}^{}\chi\left(\alpha u\right)e^{-u \left( L_{Z}+\sigma \right)}du.
\end{equation}
By construction, 
\begin{equation}\label{eq:singa11}
\mathrm{WF}\left(U_{\alpha,\sigma}\right)  \subset  \mathcal{N}^{*}_{\Delta^{Y}\setminus\left(Y\times Y\right)}\cup \Omega_{\ge 0}.
\end{equation}
By Proposition \ref{prop:compo} and by (\ref{eq:bas1}), for $r\ge 0$, 
\begin{equation}\label{eq:singa12}
e^{-rL_{Z}}U_{\alpha,\sigma}=N_{\mathrm{Gr}_{\varphi_{-r}}\circ U_{\alpha,\sigma}}.
\end{equation}
By proceeding as in (\ref{eq:pdz2}) and using (\ref{eq:singa12}), we get
\begin{equation}\label{eq:singa13}
\mathrm{WF}\left(e^{-rL_{Z}}U_{\alpha,\sigma}\right) \subset 
\mathcal{N}^{*}_{\mathrm{Gr}\varphi_{-r}} \cup\Omega_{\ge r}.
\end{equation}

By  Proposition \ref{prop:wff} and using (\ref{eq:singa13}),  for $0<r< r_{0}$,
\begin{equation}\label{eq:singa14}
\mathrm{WF}\left(e^{-rL_{Z}}U_{\alpha,\sigma}\right)\cap\mathcal{N}^{*}_{\Delta^{Y}/\left(Y\times Y\right)}\setminus\left(Y\times Y\right)=\emptyset.
\end{equation}
By (\ref{eq:singa14}), under the same conditions, the product $e^{-rL_{Z}}U_{\alpha,\sigma}\Delta^{Y}$ is well defined.

We now state a result which is the consequence of an important theorem  by Chaubet and Guedes Bonthonneau \cite[Theorem 2]{Chaubet:2022aa}.
\begin{theorem}\label{thm:chge}
For $\sigma\in\C,\Re\,\sigma\gg 1$, for $0<r<r_{0}$, as $\alpha\to 0$, we have the convergence of $2n$ currents on $Y\times Y$, 
\begin{equation}\label{eq:singa15}
e^{-rL_{Z}}U_{\alpha,\sigma}\Delta^{Y}\to e^{-rL_{Z}}\left(L_{Z}+\sigma\right)^{-1}\Delta^{Y}. 
\end{equation}
\end{theorem}

By integration of (\ref{eq:singa15}), we deduce that for $0<r<r_{0}$,  $\sigma\in\C,\Re\,\sigma\gg 1$, as $\alpha\to 0$, 
\begin{equation}\label{eq:singa16}
\mathrm{Tr_{s,int}}\left[Ne^{-rL_{Z}}U_{\alpha,\sigma}\right]\to
\mathrm{Tr_{s,int}}\left[Ne^{-rL_{Z}}\left(L_{Z}+\sigma\right)^{-1}\right].
\end{equation}

Recall that for $\Re\,\sigma\gg 1$, the Fried zeta function $R_{Z,F}\left(\sigma\right)$ was 
defined in Definition \ref{def:friz1}. By the results of Subsection \ref{subsec:ruze}, for $\Re\,\sigma\gg 1$, 
$R_{Z,F}\left(\sigma\right)$ is a holomorphic function of $\sigma$, 
and  as $\sigma\to + \infty $,
\begin{equation}\label{eq:toto6}
R_{Z,F}\left(\sigma\right)\to 1.
\end{equation}

Now we recall a result of Dyatlov-Zworski \cite[p.~561]{DyatlovZworski16}.
\begin{theorem}\label{thm:delo}
	For $0<r<r_{0}$ and $\mathrm{Re}\,\sigma\gg 1$, 
	\begin{equation}\label{eq:toto5}
\Tr_{\mathrm{s,int}}\left[Ne^{-r \left(L_{Z}+\sigma \right)} \left(L_{Z}+\sigma\right)^{-1}\right]=\frac{d}{d\sigma}\log R_{Z,F}\left(\sigma\right).
\end{equation}
\end{theorem}
\begin{proof}
	Using the Guillemin trace formula \cite[Theorem 8] {Guillemin77}, we get
	\begin{multline}\label{eq:singa16a}
		\mathrm{Tr_{s,int}}\left[Ne^{-r\left(L_{Z}+\sigma\right)}U_{\alpha,\sigma}\right]=\sum_{y\in\underline{\mathbf{F}}}^{}\sum_{k\ge 1}^{}
	t_{y}\frac{1}{\left\vert \det\left(1-\varphi_{-kt_{y}*}\right)\vert_{M}\right\vert}\\
	\Tr_{\mathrm{s}}^{\Lambda\left(T^{*}Y\right) \otimes _{\R}F}\left[N\varphi^{*}_{-kt_{y}}\right]
	\chi\left(\alpha \left(k  t_{y}-r\right) \right) e^{-kt_{y}\sigma}.
	\end{multline}
	Given $\alpha>0$, the sum in  the right-hand side of (\ref{eq:singa16a}) is finite.	 Clearly, 
	\begin{equation}\label{eq:singa17}
		\Tr_{\mathrm{s}}^{\Lambda\left(T^{*}Y\right) \otimes _{\R}F}\left[N\varphi^{*}_{-kt_{y}}\right]=\Tr_{\mathrm{s}}^{\Lambda\left(T^{*}Y\right) }\left[N\varphi^{*}_{-kt_{y}}\right]
		\Tr^{F}\left[\varphi_{kt_{y}*}\right].
	\end{equation}
	Moreover, 
	\begin{equation}\label{eq:singa18}
		\Tr_{\mathrm{s}}^{\Lambda\left(T^{*}Y\right) }\left[N\varphi^{*}_{-kt_{y}}\right]=-\Tr_{\mathrm{s}}^{\Lambda\left(M^{*}\right) }\left[\varphi^{*}_{-kt_{y}}\right].	
	\end{equation}
	Finally,
	\begin{equation}\label{eq:singa19}
		\Tr_{\mathrm{s}}^{\Lambda\left(M^{*}\right) }\left[\varphi^{*}_{-kt_{y}}\right]=\det\left(1-\varphi_{-kt_{y}*}\right) \vert_{M}.
	\end{equation}

	By Proposition \ref{prop:sign1} and by (\ref{eq:singa16a})--(\ref{eq:singa19}), we obtain
	\begin{multline}\label{eq:singa20}
		\mathrm{Tr_{s,int}}\left[Ne^{-r\left(L_{Z}+\sigma\right)}U_{\alpha,\sigma}\right]\\
		=\left(-1\right)^{n_{s}+1}\sum_{y\in\underline{\mathbf{F}}}^{}\sum_{k\ge 1}^{}
		t_{y}\Tr^{F \otimes _{\Z_{2}}\mathrm{o}\left(T_{s}Y\right)}\left[\varphi_{kt_{y}*}\right]\chi\left(\alpha \left( kt_{y}-r \right) \right)e^{-kt_{y}\sigma}.
		\end{multline}
		Using (\ref{eq:singa16}) and  making $\alpha\to 0$ in (\ref{eq:singa20}), we get
\begin{multline}\label{eq:ny3}
	\Tr_{\mathrm{s,int}}\left[Ne^{-r \left( L_{Z}+\sigma \right) }\left(L_{Z}+\sigma\right)^{-1}\right]\\
	=\left(-1\right)^{n_{s}+1}\sum_{y\in\underline{\mathbf{F}}}^{}\sum_{k\ge 1}^{}
	t_{y}\Tr^{F \otimes_{\Z_{2}} \mathrm{o}\left(T_{s}Y\right)}\left[\varphi^{k}_{t_{y}*}\right]e^{-kt_{y}\sigma}.
\end{multline}

Also by (\ref{eq:flip2}), for $\sigma\in \C, \Re\sigma\gg 1$, we obtain
\begin{equation}\label{eq:ny4}
\log R_{Z,F}\left(\sigma\right)
=\left(-1\right)^{n_{s}}\sum_{y\in\underline{\mathbf{F}}}^{}\sum_{k\ge 1}^{}\frac{1}{k}\Tr^{F \otimes  _{\Z_{2}}\mathrm{o}\left(T_{s}Y\right)}\left[\varphi^{k}_{t_{y}*}\right]e^{-kt_{y}\sigma}.	
\end{equation}
By comparing (\ref{eq:ny3}) and (\ref{eq:ny4}), we get (\ref{eq:toto5}). The proof of our theorem is complete. 
\end{proof}
\subsection{The resolvent and the operators $d^{Y},i_{Z}$}%
\label{subsec:resd}
\begin{proposition}\label{prop:comd}
	For $z\in \C, z\notin \mathrm{Sp}L_{Z}$, we have the identities,
	\begin{align}\label{eq:rond1a1-}
&\left[d^{Y},\left(L_{Z}-z\right)^{-1}\right]=0,&\left[i_{Z},\left(L_{Z}-z\right)^{-1}\right]=0.
\end{align}
In particular, if $\left(L_{Z}-z\right)^{-1}$ is viewed as a current, then
\begin{align}\label{eq:sabr1}
&d^{Y\times Y}\left(L_{Z}-z\right)^{-1}=0, 
&i_{Z^{Y\times Y}}\left(L_{Z}-z\right)^{-1}=0.
\end{align}
\end{proposition}
\begin{proof}
We have the obvious identity of operators acting on 
	$\mathcal{D}\left(Y,F\right)$, 
	\begin{equation}\label{eq:doma17-}
\left(L_{Z}-z\right)d^{Y}=d^{Y}\left(L_{Z}-z\right).
\end{equation}

Assume that $z\in\C,\Re z\ll -1$.
Take $u\in \Omega\left(Y,F\right)$. Then 
$\left(L_{Z}-z\right)^{-1}u\in \Omega_{2}\left(Y,F\right)$, and 
\begin{equation}\label{eq:doma19-}
\left(L_{Z}-z\right)\left(L_{Z}-z\right)^{-1}u=u.
\end{equation}

By (\ref{eq:doma17-}), (\ref{eq:doma19-}), we obtain
\begin{equation}\label{eq:doma20-}
\left(L_{Z}-z\right)d^{Y}\left(L_{Z}-z\right)^{-1}u=d^{Y}u.
\end{equation}
Note that 
\begin{equation}\label{eq:som1-}
d^{Y}\left(L_{Z}-z\right)^{-1}u\in  \mathcal{H}^{-1}\left(Y,\Lambda\left(T^{*}Y\right) \otimes 
_{\R}F\right)
\end{equation}
Since $d^{Y}u\in \Omega\left(Y,F\right)$, by (\ref{eq:som1-}), we 
deduce that $d^{Y}\left(L_{Z}-z\right)^{-1}u$ lies in the domain of 
$L_{Z}$ on $\mathcal{H}^{-1}\left(Y,\Lambda\left(T^{*}Y\right) 
\otimes _{\R}F\right)$. As we saw before, by a result of Faure and 
Sjöstrand \cite{FaureSjostrand11}, the minimal and maximal domain of $L_{Z}$ on $\mathcal{H}^{-1}\left(Y,\Lambda\left(T^{*}Y\right) 
\otimes _{\R}F\right)$ coincide. By (\ref{eq:doma20-}), we ultimately 
find that 
\begin{equation}\label{eq:doma21-}
d^{Y}\left(L_{Z}-z\right)^{-1}u=\left(L_{Z}-z\right)^{-1}d^{Y}u.
\end{equation}
By (\ref{eq:doma21-}), and also using the fact that 
$\left(L_{Z}-z\right)^{-1}$ is meromorphic, we get the first equation in 
(\ref{eq:rond1a1-}). The proof of the second equation is similar. Using Proposition \ref{prop:scm} and (\ref{eq:rond1a1-}), we get (\ref{eq:sabr1}). The proof of our proposition is complete. 
\end{proof}

\subsection{The spectral projectors $P_{Z,z_{0}}$}%
\label{subsec:sptru}
Let $z_{0}\in \bC$. Take $p>0$ such that $\Re\left(z_{0}\right)+2\le 
-c_{0}+\delta p$.  By Theorem \ref{thm:FaureSjostrand}, $L_{Z,p}$ has discrete spectrum on $\{z\in 
\mathbf{C}: \Re(z)\le \Re(z_{0})+1\}$, and  on  the same domain in $\C$, 
$(L_{Z,p}-z)^{-1} $
is a meromorphic family of bounded operators on $\cH_{pG}\left(Y,\Lambda\left(T^{*}Y \right)  \otimes _{\R}F\right)$.

If $z_{0}\in\mathrm{Sp}L_{Z}$, for 
$\epsilon>0$ small enough, $z_{0}$ is the unique element in the 
spectrum that lies in the closed disk centered at $z_{0}$ and of radius 
$\epsilon$.
 The associated spectral projector
 \index[not]{PZp@$P_{Z,p,z_{0}}$}%
 $P_{Z,p,z_{0}}$ is given by 
\begin{align}\label{eq:duc}
	P_{Z,p,z_{0}}=\frac{1}{2i\pi}\int_{\substack{z\in\C\\
\left\vert  
z-z_{0}\right\vert=\epsilon}}-(L_{Z,p}-z)^{-1}dz. 
\end{align} 
By \eqref{eq:eqAfred},  the operator $P_{Z,p,z_{0}}$  is bounded and has finite rank.  
Recall that
\index[not]{DZp@$\cD_{Z,p,z_{0}}(Y,F)$}%
$\cD_{Z,p,z_{0}}(Y, 
F)\subset \cH_{pG}\left(Y,\Lambda\left(T^{*}Y\right) \otimes 
_{\R}F\right)$ is the image of  $P_{Z,p,z_{0}}$.  

Since $\Omega(Y, 
F)$ is dense in $\cH_{pG}\left(Y,\Lambda\left(T^{*}Y\right) \otimes 
_{\R}F\right)$,  and since $P_{Z,p,z_{0}}\Omega(Y, F)$ is closed in 
$\cD_{Z,p,z_{0}}(Y, F)$,
\begin{align}\label{eq:doma9}
\cD_{Z,p,z_{0}}(Y, F)=P_{Z,p,z_{0}}\Omega(Y, F). 
\end{align} 
By Proposition \ref{prop:propextReso} and by (\ref{eq:doma9}), 
$\cD_{Z,p,z_{0}}\left(Y,F\right)$ does not depend on $p$. In the sequel, it will be instead
denoted 
\index[not]{DZz@$\cD_{Z,z_{0}}\left(Y,F\right)$}%
$\cD_{Z,z_{0}}\left(Y,F\right)$. By (\ref{eq:pox}), we get
\begin{equation}\label{eq:doma11}
\cD_{Z,z_{0}}\left(Y,F\right) \subset \bigcap_{p\g 0, g} \cH_{pG}(Y, 
\Lambda(T^{*}Y)\otimes_{\R} F)  \subset  \mathcal{D}_{\mathcal{Y}^{*}_{s}}\left(Y,F\right).
\end{equation}
By Proposition \ref{prop:propextReso}, the restriction on 
$P_{Z,p,z_{0}}$ to $\Omega\left(Y,F\right)$ does not depend on $p$. 
It will instead be denoted 
\index[not]{PZz@$P_{Z,z_{0}}$}%
$P_{Z,z_{0}}$.

Similarly, when replacing $Z$ by $-Z$ and $F$ by $F^{*} \otimes 
_{\Z_{2}}\mathrm{o}\left(TY\right)$, we obtain a corresponding projection 
operator $P_{-Z,z_{0}}$, which is naturally the transpose of 
$P_{Z,z_{0}}$, and projects on $\mathcal{D}_{-Z,z_{0}}\left(Y,F^{*} \otimes _{\Z_{2}}\mathrm{o}\left(TY\right)\right)$. 
This new projector will be denoted 
\index[not]{PZz@$P_{-Z,z_{0}}$}%
$P_{-Z,z_{0}}$.  The following result is implicitly proved in Faure-Sjöstrand \cite[Section 4]{FaureSjostrand11} and Dyatlov-Zworski \cite[Lemma 2.2]{DyatlovZworski17}.
\begin{proposition}\label{prop:spec}
Given $z_{0}\in\C$, then
\begin{align}\label{eq:eqDm}
	\cD_{Z,z_{0}}(Y, F)=\{u\in \cD_{\mathcal{Y}^{*}_{s}}(Y, F) : 
	(L_{Z}-z_{0})^{e_{Z,z_{0}}}u=0\}. 
\end{align} 
\end{proposition} 
\begin{proof}
In (\ref{eq:eqDm}), we already know that the left-hand side is included in the right-hand 
side. Take $u\in \cD_{\mathcal{Y}^{*}_{s}}(Y, F)$ be such  that for a given $m\in 
\N$,   $(L_{Z}-z_{0})^{m}u=0$. For 
$p>0$ big enough, 
$$u, L_{Z}u,\ldots \(L_{Z}\)^{m}u\in 
\cH_{pG}\left(Y,\Lambda\left(T^{*}Y\right) \otimes _{\R}F\right),$$ so 
that $u$ is in the domain of $L_{Z,p}$, and $z_{0}$ 
lies in $\mathrm{Sp}L_{Z}$, and moreover, $u\in \cD_{Z,p,z_{0}}(Y, 
F)$, from which we find that $u\in 
\mathcal{D}_{Z,z_{0}}\left(Y,F\right)$. The proof of our proposition is complete. 
\end{proof}

	By (\ref{eq:duc}), we have the identity of currents on $Y\times Y$,
	\begin{equation}\label{eq:doma14a1}
P_{Z,z_{0}}=\frac{1}{2i\pi}\int_{\substack{z\in\C\\
\left\vert  
z-z_{0}\right\vert=\epsilon}}^{}-\left(L_{Z}-z\right)^{-1}dz.
\end{equation}
\begin{theorem}\label{thm:wfp}
	The following identities hold:
	\begin{align}\label{eq:comer2}
&\left[d^{Y},P_{Z,z_{0}}\right]=0,&\left[i_{Z},P_{Z,z_{0}}\right]=0.
\end{align}
Equivalently,
\begin{align}\label{eq:comer2a1}
&d^{Y\times Y}P_{Z,z_{0}}=0,&i_{Z^{Y\times Y}}P_{Z,z_{0}}=0.
\end{align}
Moreover, 
	\begin{equation}\label{eq:wfpz}
\mathrm{WF}\left( P_{Z,z_{0}}\right) \subset  \left(\mathcal{Y}^{*}_{s} 
\times\mathcal{Y}^{*}_{u}\right)\setminus\left(Y\times Y\right).
\end{equation}
Finally, 
\begin{equation}\label{eq:rago1}
\mathcal{D}_{Z,z_{0}}\left(Y,F\right)=P_{Z,z_{0}}\Omega\left(Y,F\right).
\end{equation}
\end{theorem}
\begin{proof}
	By Proposition \ref{prop:comd} and by  (\ref{eq:doma14a1}), we 
	get (\ref{eq:comer2}) and (\ref{eq:comer2a1}).  By Theorem \ref{thm:pdz1} and by (\ref{eq:doma14a1}), we deduce that 
$P_{Z,z_{0}}\in\mathcal{D}_{0,2}$, so that
\begin{equation}\label{eq:coume1}
\mathrm{WF}\left(P_{Z,z_{0}}\right)\subset 
\left(\mathcal{N}^{*}_{\Delta^{Y}/Y\times Y}\cup\Omega_{\ge 0}\cup \left( 
	\mathcal{Y}^{*}_{s}\times 
	\mathcal{Y}^{*}_{u} \right)  \right) \setminus \left(Y\times Y\right). 
\end{equation}

We can rewrite 
(\ref{eq:doma14a1}) in the form,
\begin{equation}\label{eq:doma14a2}
P_{Z,z_{0}}=e^{rL_{Z}}\int_{\substack{z\in\C\\
\left\vert  
z-z_{0}\right\vert=\epsilon}}^{}-e^{-rL_{Z}}\left(L_{Z}-z\right)^{-1}dz.
\end{equation}
The first term $e^{rL_{Z}}$ in the right-hand side can be viewed 
as a matrix acting on $\mathcal{D}_{Z,z_{0}}\left(Y,F\right)$,  and has no effect on 
the wave front set. To deal with  the integrand, we  combine 
(\ref{eq:pdz2}) and (\ref{eq:sing}), which guarantee that for any 
$r>0$, 
\begin{equation}\label{eq:doma14a3}
\mathrm{WF}\left(P_{Z,z_{0}}\right) \subset \left(\mathcal{N}^{*}_{\mathrm{Gr}\varphi_{-r}/Y\times Y}\cup\Omega_{\ge r}\cup 
	\left( \mathcal{Y}^{*}_{s}\times 
	\mathcal{Y}^{*}_{u}\right)  \right) \setminus \left(Y\times Y\right). 
\end{equation}
Observe that for $0<r<r_{0}$, 
\begin{equation}\label{eq:coume2}
\mathcal{N}^{*}_{\Delta^{Y}/Y\times 
Y}\cap\mathcal{N}^{*}_{\mathrm{Gr}\varphi_{-r}/Y\times Y}=\emptyset.
\end{equation}
By (\ref{eq:doma14a3}), (\ref{eq:coume2}), we conclude that for $r\ge 
0$,
\begin{equation}\label{eq:coume3}
\mathrm{WF}\left(P_{Z,z_{0}}\right) \subset \Omega_{\ge r}\cup 
	\left(\mathcal{Y}^{*}_{s}\times 
	\mathcal{Y}^{*}_{u} \right) \setminus \left(Y\times Y\right). 
\end{equation}
By proceeding as in the proof of Proposition \ref{prop:b1} and using 
(\ref{eq:coume3}), we get (\ref{eq:wfpz}). 
Equation (\ref{eq:rago1}) is a consequence of (\ref{eq:doma9}). The proof of our theorem is complete. 
\end{proof}

By Theorem \ref{thm:wfp}, we conclude that $P_{Z,z_{0}}$ induces an 
endomorphism of 
$\mathcal{D}_{\mathcal{Y}^{*}_{s}}\left(Y,F\right)$, which is a 
projection. Set
\index[not]{QZz@$Q_{Z,z_{0}}$}%
\begin{equation}\label{eq:comer4}
Q_{Z,z_{0}}=1-P_{Z,z_{0}}.
\end{equation}
By Proposition \ref{prop:diag}, the current corresponding 
to $1$ is just $\Delta^{Y}$. In particular,
\begin{equation}\label{eq:comer5}
\mathrm{WF}\left( Q_{Z,z_{0}} \right)  \subset \left(\mathcal{N}^{*}_{\Delta^{Y}/\left(Y\times Y\right)} \cup \left( \mathcal{Y}^{*}_{s} 
\times\mathcal{Y}^{*}_{u}\right) \right) 
\setminus\left(Y\times Y\right).
\end{equation}

Set
\index[not]{DZz@$\mathcal{D}_{Z,z_{0}}^{\perp}\left(Y,F\right)$}%
\begin{equation}\label{eq:comer6}
\mathcal{D}_{Z,z_{0}}^{\perp}\left(Y,F\right)=Q_{Z,z_{0}}
\mathcal{D}_{\mathcal{Y}^{*}_{s}}\left(Y,F\right)
\end{equation}
Then we have the splitting
\begin{equation}\label{eq:comer7}
\mathcal{D}_{\mathcal{Y}^{*}_{s}}\left(Y,F\right)=\mathcal{D}_{Z,z_{0}}\left(Y,F\right) \oplus \mathcal{D}_{Z,z_{0}}^{\perp}\left(Y,F\right).
\end{equation}
\subsection{Spectral projectors and Poincaré duality}%
\label{subsec:sptrdu}
 Observe that
\begin{equation}\label{eq:vega0}
P_{Z,z_{0}}^{\dag}=P_{-Z,z_{0}}.
\end{equation}
By (\ref{eq:sv5a6}), we can rewrite (\ref{eq:vega0}) in the form
\begin{equation}\label{eq:comer3}
P_{-Z,z_{0}}=\left(-1\right)^{n\left(n-1\right)/2}\underline{j}^{*}P_{Z,z_{0}}.
\end{equation}
Since $\underline{j}^{*}$  exchanges the two wave front sets 
in the right-hand side  of (\ref{eq:wfpz}),
equations (\ref{eq:wfpz}) and (\ref{eq:comer3}) are compatible.
\begin{theorem}\label{thm:podua1}
	For  $0\le p\le n$, the bilinear map $\left(\phi,\psi\right)\in 
	\mathcal{D}^{p}_{Z,z_{0}}\left(Y,F\right) 
	\times\mathcal{D}^{n-p}_{-Z,z_{0}}\left(Y,\Ft\right)\to \left\langle  
	\phi,\psi\right\rangle$ is non-degenerate. Also 
	$\mathcal{D}_{Z,z_{0}}^{\perp}\left(Y,F\right)$ and 
	$\mathcal{D}_{-Z,z_{0}}\left(Y,\Ft\right)$ are orthogonal.
\end{theorem}
\begin{proof}
	Take  $\varphi\in 
\Omega\left(Y,F\right),\psi\in\Omega\left(Y,\Ft\right)$. Observe that
$P_{Z,z_{0}}\varphi\in\mathcal{D}_{\mathcal{Y}^{*}_{s}}\left(Y,F\right),
P_{-Z,z_{0}}\psi\in\mathcal{D}_{\mathcal{Y}^{*}_{u}}\left(Y,\Ft\right)$. Using (\ref{eq:vega0}), we get 
\begin{equation}\label{eq:vega1}
\left\langle  P_{Z,z_{0}}\varphi,\psi\right\rangle=\left\langle  
\varphi,P_{-Z,z_{0}}\psi\right\rangle=\left\langle  
P_{Z,z_{0}}\varphi,P_{-Z,z_{0}}\psi\right\rangle.
\end{equation}
The bilinear map in the left-hand side of (\ref{eq:vega1}) is 
 non-degenerate. By (\ref{eq:rago1}), $P_{Z,z_{0}}$ 
maps $\Omega\left(Y,F\right)$ onto 
$\mathcal{D}_{Z,z_{0}}\left(Y,F\right)$. Using  (\ref{eq:vega1}), we obtain the first part of our 
theorem. 
Since $Q_{-Z,z_{0}}$ vanishes on 
$\mathcal{D}_{-Z,z_{0}}\left(Y,\Ft\right)$, the second statement is 
trivial.  The proof of our theorem is complete. 
\end{proof}

Given $p,0\le p\le n$, let $\left(\varphi^{p}_{i}\right)\vert_{1\le i\le 
e_{Z,z_{0}}^{p}}$ be a basis of $\mathcal{D}^{p}_{Z,z_{0}}\left(Y,F\right)$, let 
$\left(\psi^{n-p}_{i}\right)\vert_{1\le i\le e^{p}_{Z,z_{0}}}$ be the associated 
dual basis of $\mathcal{D}^{n-p}_{-Z,z_{0}}\left(Y,\Ft\right)$, so that 
\begin{equation}\label{eq:sacr1b1}
\left\langle  \psi^{p}_{i},\varphi^{p}_{j}\right\rangle=\delta_{i,j}.
\end{equation}
\begin{proposition}\label{prop:proja}
	The following identity of operators acting on $\Omega^{p}\left(Y,F\right)$ holds:
	\begin{equation}\label{eq:vega2}
P_{Z,z_{0}}\vert_{\Omega^{p}\left(Y,F\right)}=M_{\sum_{1}^{e^{p}_{Z,z_{0}}}\pi_{1}^{*}\varphi^{p}_{i}\we  \pi_{2}^{*}\widetilde\psi^{n-p}_{i}}.
\end{equation}
\end{proposition}
\begin{proof}
	This follows from an easy computation that is left to the reader.
\end{proof}

Recall that we have identified $P_{Z,z_{0}}$ with a current in 
$$\mathcal{D}_{\pi_{1}^{*}\mathcal{Y}^{*}_{s}\times\pi_{2}^{*}\mathcal{Y}^{*}_{u}}\left(Y\times Y,\pi_{1}^{*}F \otimes \pi_{2}^{*}\left(F^{*} \otimes _{\Z_{2}}\left(TY\right)\right)\right).$$
\begin{proposition}\label{prop:pidproj}
The following identity of currents holds:
\begin{equation}\label{eq:proxy1}
P_{Z,z_0}=\sum_{p=0}^{n}\left(-1\right)^{np}
\sum_{1}^{e^{p}_{Z,z_{0}}}\pi_{1}^{*}\varphi^{p}_{i}\we  \pi_{2}^{*}\widetilde\psi^{n-p}_{i}.
\end{equation}
We have the identities:
\begin{align}\label{eq:proxy2}
&\mathrm{Tr}^{\Omega^{p}\left(Y,F\right)}_{\mathrm{int}}\left[P_{Z,z_{0}}\right]=\dim \mathcal{D}^{p}_{Z,z_{0}}\left(Y,F\right),\\
&\mathrm{\Tr_{s,int}}\left[P_{Z,z_{0}}\right]=
\sum_{p=0}^{n}\left(-1\right)^{p}\dim \mathcal{D}^{p}_{Z,z_{0}}\left(Y,F\right)=0. \notag 
\end{align}
\end{proposition}
\begin{proof}
We use the fact that the correspondence between currents and operators is via $N$, and also we use Proposition \ref{prop:proja}, and we get (\ref{eq:proxy1}). Also we use (\ref{eq:scar1a1}) with $T=P_{Z,z_0}$, together with (\ref{eq:sacr1b1}) and (\ref{eq:proxy1}), and we obtain the first identity in (\ref{eq:proxy2}). The first part of the second identity is now a consequence of the first identity. The vector field $Z^{Y\times Y}$ restricts to a nonzero vector field on $\Delta_{Y}$. By the second identity in (\ref{eq:comer2a1}), we obtain the last identity in (\ref{eq:proxy2}). The proof of our proposition is complete.  
\end{proof}
\begin{remark}\label{rem:dual}
Using (\ref{eq:normio3a1}), one verifies easily that (\ref{eq:comer3}) and (\ref{eq:proxy1}) are compatible.
\end{remark}

Let $z_{0}\in \mathrm{Sp}\,L_{Z}$ be a pole of $\left(L_{Z}-z\right)^{-1}$. We have 
the obvious identity
\begin{equation}\label{eq:dz1}
\left(L_{Z}-z\right)^{-1}=\left(L_{Z}-z\right)^{-1}P_{Z,z_{0}}+\left(L_{Z}-z\right)^{-1}Q_{Z,z_{0}}.
\end{equation}
Since $\left(L_{Z}-z\right)^{-1}$ is a meromorphic family, we have the expansion of the first 
term near $z=z_{0}$,
\begin{equation}\label{eq:dz2}
\left(L_{Z}-z\right)^{-1}P_{Z,z_{0}}=-\sum_{a=0}^{e_{Z,z_{0}}}\left(L_{Z}-z_{0}\right)^{a}P_{Z,z_{0}}\left(z-z_{0}\right)^{-a-1}.
\end{equation}
The second term in the right-hand side of (\ref{eq:dz1}) is also a meromorphic family of operators that is 
holomorphic at $z_{0}$. A result of Dyatlov-Zworski 
\cite[Proposition 3.3]{DyatlovZworski16} says that the family 
$\left(L_{Z}-z\right)^{-1}Q_{Z,z_{0}}$ has the same wave front set 
properties as the original $\left(L_{Z}-z\right)^{-1}$, except 
that $z_{0}$ is also included in the domain of definition.

Let $z_{0},z'_{0}\in \mathrm{Sp}L_{Z}$. The composition 
 of currents $P_{Z,z_{0}}\circ P_{Z,z'_{0}},P_{Z,z'_{0}}\circ P_{Z,z_{0}}$ 
are well-defined, and that the corresponding wave front sets are the 
same as in (\ref{eq:wfpz}). By Proposition \ref{prop:compo}, the corresponding operators are just $P_{Z,z_{0}}P_{Z,z'_{0}},P_{Z,z'_{0}}P_{Z,z_{0}}$.
\begin{proposition}\label{prop:dif}
	If $z_{0}\neq z'_{0}$, then
	\begin{align}\label{eq:comer7z1}
&\mathcal{D}_{Z,z'_{0}}\left(Y,F\right) \subset 
\mathcal{D}^{\perp}_{Z,z_{0}}\left(Y,F\right), &\mathcal{D}_{Z,z_{0}}\left(Y,F\right) \subset 
\mathcal{D}^{\perp}_{Z,z'_{0}}\left(Y,F\right),\\
&P_{Z,z_{0}}P_{Z,z'_{0}}=0,&P_{Z,z'_{0}}P_{Z,z_{0}}=0. \nonumber 
\end{align}
Under the above conditions, $\mathcal{D}_{Z,z_{0}}\left(Y,F\right)$ 
and $\mathcal{D}_{-Z,z'_{0}}\left(Y,\Ft\right)$ are orthogonal with 
respect to the duality defined in (\ref{eq:normi3}).
\end{proposition}
\begin{proof}
We know that $\left(L_{Z}-z\right)^{-1}$ 	 preserves 
$\mathcal{D}^{\perp}_{Z,z_{0}}\left(Y,F\right)$ and that this 
restriction has the same properties as the original resolvent 
$\left(L_{Z}-z\right)^{-1}$ except that the eigenvalue $z_{0}$ has 
been eliminated, so that (\ref{eq:comer7}) holds. The proof of our proposition is complete. 
\end{proof}

\subsection{The truncated vector spaces 
$\mathcal{D}_{Z,<a}\left(Y,F\right)$}%
\label{subsec:vspa}
We will now put together our characteristic subspaces in a way similar to Chaubet-Dang \cite[Section 6.3]{ChaubetDang24}.
\begin{definition}\label{def:sptru}
If $a>0$ is such that no eigenvalue of $L_{Z}$ is contained in the 
circle $S_{a}$, put
\index[not]{DZa@$\mathcal{D}_{Z,<a}\left(Y,F\right)$}%
\begin{align}\label{eq:comer0}
	\cD_{Z,<a}(Y, F)=\bigoplus_{\substack{z_{0}\in\mathrm{Sp}L_{Z}\\
	|z_{0}|<a}}\cD_{Z,z_{0}}(Y, F). 
\end{align}
By (\ref{eq:eqDm}), we get
\begin{equation}\label{eq:doma15}
\cD_{Z,<a}(Y, F) \subset \mathcal{D}_{\mathcal{Y}^{*}_{s}}\left(Y,F\right).
\end{equation}

Let 
\index[not]{PZa@$P_{Z,<a}$}%
$P_{Z,<a}$ be the  spectral projector from 
$\mathcal{D}_{\mathcal{Y}^{*}_{s}}\left(Y,F\right)$ on $\cD_{Z,<a}(Y, F)$, 
i.e.,
\begin{equation}\label{eq:doma16}
P_{Z,<a}=\sum_{\substack{z_{0}\in\mathrm{Sp}L_{Z}\\\left\vert  
z_{0}\right\vert<a}}^{}P_{Z,z_{0}}.
\end{equation}
\end{definition}
By (\ref{eq:rago1}), we get
\begin{equation}\label{eq:doma16a1}
\mathcal{D}_{Z,<a}\left(Y,F\right)=P_{Z,<a}\Omega\left(Y,F\right).
\end{equation}

There is a  splitting corresponding to (\ref{eq:comer0}),
\begin{equation}\label{eq:comer8}
\mathcal{D}_{-Z,<a}\left(Y,\Ft\right)=\bigoplus_{\substack{z_{0}\in\mathrm{Sp}L_{-Z}\\
	|z_{0}|<a}}\cD_{-Z,z_{0}}(Y, \Ft). 
\end{equation}
By Theorem \ref{thm:podua1} and by Proposition  \ref{prop:dif}, $\left\langle  \,\right\rangle$ induces  a 
duality between $\mathcal{D}_{Z,<a}\left(Y,F\right)$, and 
$\mathcal{D}_{-Z,<a}\left(Y,\Ft\right)$ and the components in 
(\ref{eq:comer0}), (\ref{eq:comer8}) for distinct $z_{0},z'_{0}$ are 
mutually orthogonal.

Using (\ref{eq:duc}), we can rewrite 
(\ref{eq:doma16}) as an equality of operators 
$\Omega\left(X,F\right)\to \mathcal{D}_{\mathcal{Y}^{*}_{s}}\left(Y,F\right)$,
\begin{equation}\label{eq:doma16a2}
P_{Z,<a}=\frac{1}{2i\pi}\int_{z\in S_{a}}^{}-\left(L_{Z}-z\right)^{-1}dz.
\end{equation}
By (\ref{eq:wfpz}), (\ref{eq:doma16}), we get
\begin{equation}\label{eq:doma14a4}
\mathrm{WF}\left(P_{Z,<a}\right) \subset \left(\mathcal{Y}^{*}_{s} \times 
\mathcal{Y}_{u}^{*}\right)\setminus\left(Y\times Y\right).
\end{equation}

Put 
\index[not]{PZa@$P_{Z,>a}$}%
\begin{align}
	P_{Z,>a}=1-P_{Z,<a}. 
\end{align}
Then $P_{Z,>a}$ is also a projector acting on  
$\mathcal{D}_{\mathcal{Y}^{*}_{s}}\left(Y,F\right)$. By (\ref{eq:doma14a4}), we 
get
\begin{equation}\label{eq:doma14a4z}
\mathrm{WF}\left(P_{Z,>a}\right) \subset \left( 
\mathcal{N}^{*}_{\Delta^{Y}/\left(Y\times Y\right)}\cup
\left(\mathcal{Y}^{*}_{s}\times\mathcal{Y}^{*}_{u}\right)\right) \setminus\left(Y\times Y\right).
\end{equation}

By (\ref{eq:dz2}), (\ref{eq:doma14a4}), and 
(\ref{eq:doma14a4z}), if $z\in \C\setminus\mathrm{Sp}L_{Z}$, we get
\begin{equation}\label{eq:runa1}
\mathrm{WF}\left( \left(L_{Z}-z\right)^{-1}P_{Z,<a}\right) \subset 
\left( \mathcal{Y}^{*}_{s}\times\mathcal{Y}^{*}_{u} \right) \setminus\left(Y\times Y\right).
\end{equation}
For $0<r<r_{0}$, 
$\mathrm{WF}\left( \left(L_{Z}-z\right)^{-1}P_{Z,>a}\right) , 
\mathrm{WF}\left( e^{-rZ}\left(L_{Z}-z\right)^{-1}P_{Z,>a}\right) $ have  
the same properties as 
$\mathrm{WF}\left( \left(L_{Z}-z\right)^{-1}\right) ,\mathrm{WF}\left( e^{-rL_{Z}}\left(L_{Z}-z\right)^{-1}\right)$.

Put 
\index[not]{DZa@$\cD_{Z,>a}(Y, F)$}%
\begin{align}
	\cD_{Z,>a}(Y, F)= P_{Z,>a}\cD_{\mathcal{Y}^{*}_{s}}(Y, F). 
\end{align} 
Then $ P_{Z,>a}$ is a projector from 
$\mathcal{D}_{\mathcal{Y}^{*}_{s}}\left(Y,F\right)$ on $\cD_{Z,>a}(Y, F)$. Moreover,
\begin{align}\label{eq:spli}
	\cD_{\mathcal{Y}^{*}_{s}}(Y, F)= \cD_{Z,<a}(Y, F) \oplus \cD_{Z,>a}(Y, F).
\end{align} 
By (\ref{eq:doma16a1}), $\Omega\left(Y,F\right)$ is transverse to 
$\mathcal{D}_{Z,>a}\left(Y,F\right)$ in 
$\mathcal{D}_{\mathcal{Y}^{*}_{s}}\left(Y,F\right)$, i.e., 
\begin{equation}\label{eq:spli1}
\mathcal{D}_{\mathcal{Y}^{*}_{s}}\left(Y,F\right)=\Omega\left(Y,F\right)+\cD_{Z,>a}(Y, F).
\end{equation}
Also $P_{Z,<a}$ commutes with $L_{Z}$, so that the splitting 
(\ref{eq:spli}) is also preserved by $L_{Z}$.

As in (\ref{eq:dz1}), we have the identity
\begin{equation}\label{eq:romb0}
\left(L_{Z}-z\right)^{-1}=\left(L_{Z}-z\right)^{-1}P_{Z,<a}+\left(L_{Z}-z\right)^{-1}P_{Z,>a}.
\end{equation}
The first term in (\ref{eq:romb0}) only has poles in 
$\left\{z\in\C,\left\vert  z\right\vert<a\right\}$. The second term 
has the same wave front set properties as 
$\left(L_{Z}-z\right)^{-1}$, except that by a result of 
Dyatlov-Zworski \cite[Proposition 3.3]{DyatlovZworski16}, these properties also extend to $z\in \C,\left\vert  
z\right\vert<a$. 
\subsection{Two quasi-isomorphisms}%
\label{subsec:trqi}
\begin{proposition}\label{prop:embco}
	The embeddings 
$\left(\Omega\left(Y,F\right),d^{Y}\right)\to\left(\mathcal{D}\left(Y,F\right),d^{Y}\right)$ and
$\left(\Omega\left(Y,F\right),d^{Y}\right)\to\left(\mathcal{D}_{\mathcal{Y}^{*}_{s}}\left(Y,F\right),d^{Y}\right)$ are quasi-isomorphisms. Also the  complexes $\left(\Omega\left(Y,F\right),i_{Z}\right)$, $\left(\mathcal{D}\left(Y,F\right),i_{Z}\right)$, and $\left(\mathcal{D}_{\mathcal{Y}^{*}_{s}}\left(Y,F\right),i_{Z}\right)$ are exact.
\end{proposition}
\begin{proof}
 The first part of our proposition is well-known. For the sake of completeness, we recall its proof. Let $d^{Y*}$ denote the formal adjoint of $d^{Y}$ with respect to 
 the Hermitian product on $\Omega_{2}\left(Y,F\right)$, and let 
 $\square^{Y}$ denote the corresponding Hodge Laplacian, which is a 
 second order elliptic operator. For $t>0$, we have
 \begin{equation}\label{eq:domer1}
\exp\left(-t\square^{Y}\right)-1=-\int_{0}^{t}\square^{Y}\exp\left(-s\square^{Y}\right)ds.
\end{equation}
Put
\begin{equation}\label{eq:domer3}
T_{t}=-d^{Y^{*}}\int_{0}^{t}\exp\left(-s\square^{Y}\right)ds.
\end{equation}
Equation (\ref{eq:domer1}) can be rewritten in the form
\begin{equation}\label{eq:domer2}
\exp\left(-t\square^{Y}\right)-1=\left[d^{Y},T_{t}\right].
\end{equation}

Observe that
\begin{equation}\label{eq:dome3a1}
T_{t}=d^{Y*}\frac{\exp\left(-t\square^{Y}\right)-1}{\square^{Y}}.
\end{equation}
In particular $T_{t}$ is a pseudo-differential operator of order $-1$.

For $t>0$, the operator $\exp\left(-t\square^{Y}\right)$ maps 
$\mathcal{D}\left(Y,F\right)$ into $\Omega\left(Y,F\right)$, and the 
operator 
$T_{t}$ acts on 
$\Omega\left(Y,F\right)$ and on $\mathcal{D}\left(Y,F\right)$.
The first part of our proposition  is a consequence of the above.
Also $\mathcal{D}_{\mathcal{Y}^{*}_{s }}\left(Y,F\right)$ is preserved by the 
 pseudo-differential operators. The above arguments can then 
be used to complete the proof of the second part our proposition.

Recall that $Z$ is a non-vanishing vector field. Therefore there is a 
smooth $1$-form $\alpha$ on $Y$ such that $i_{Z}\alpha=1$. This 
relation can also be written in the form
\begin{equation}\label{eq:rond1}
\left[i_{Z},\alpha\we\right]=1.
\end{equation}
The operator $\alpha\we$ acts on the complexes mentioned in the 
proposition as a homotopy with respect to $i_{Z}$, which shows they 
are exact. The proof of our proposition is complete. 
\end{proof}
\subsection{A truncated quasi-isomorphism}%
\label{subsec:truqi}
\begin{theorem}\label{thm:presa}
	We have the identities,
	\begin{align}\label{eq:rond1a2}
&\left[d^{Y},P_{Z,<a}\right]=0,&\left[i_{Z},P_{Z,<a}\right]=0,\\
&d^{Y\times Y}P_{Z,<a}=0,&i_{Z^{Y\times Y}}P_{Z,<a}=0. \notag
\end{align}
The vector spaces
	$\mathcal{D}_{Z,<a}\left(Y,F\right)$ and $\mathcal{D}_{Z,>a}\left(Y,F\right)$ are preserved by 
	$d^{Y}$ and $i_{Z}$.   
	
	The map 
$P_{Z,<a}:\left(\Omega\left(Y,F\right),d^{Y}\right)\to 
\left(\mathcal{D}_{Z,<a}\left(Y,F\right),d^{Y}\right)$ is a surjective quasi-isomorphism,  and 
the complex $\left(\mathcal{D}_{Z,>a}\left(Y,F\right),d^{Y}\right)$ is exact. 

Finally, the complexes $\left(\mathcal{D}_{Z,<a}\left(Y,F\right),i_{Z}\right)$ and 
$\left(\mathcal{D}_{Z,>a}\left(Y,F\right),i_{Z}\right)$ are exact.
\end{theorem}
\begin{proof}
By Proposition \ref{prop:comd} and by (\ref{eq:doma16a2}), we get 
	(\ref{eq:rond1a2}).  By \cite[Proposition 3.3]{DyatlovZworski16}, the spectrum of 
$L_{Z}\vert_{\mathcal{D}_{Z,>a}\left(Y,F\right)}$ does not contain $0$, 
so that $L_{Z}^{-1}$ maps $\mathcal{D}_{Z,>a}\left(Y,F\right)$ into 
itself. Also $L_{Z}^{-1}$ and $i_{Z}$ commute on 
$\mathcal{D}_{Z,>a}\left(Y,F\right)$, and we have the identity of 
operators acting on $\mathcal{D}_{Z,>a}\left(Y,F\right)$,
\begin{equation}\label{eq:doma22}
\left[d^{Y},L_{Z}^{-1}i_{Z}\right]=1,
\end{equation}
so that $L_{Z}^{-1}i_{Z}$ is a homotopy for $d^{Y}$ on 
$\mathcal{D}_{Z,>a}\left(Y,F\right)$, which implies that
$\left(D_{Z,>a}\left(Y,F\right),d^{Y}\right)$ is exact. Using (\ref{eq:doma9}), (\ref{eq:spli}) and Proposition \ref{prop:embco}, we find that
$P_{Z,<a}:\left( \Omega\left(Y,F\right),d^{Y}\right) \to 
\left( \mathcal{D}_{Z,<a}\left(Y,F\right),d^{Y} \right) $ is a surjective quasi-isomorphism.

Let $\alpha$ be a smooth $1$-form on $Y$ taken as in the proof of 
Proposition \ref{prop:embco}, so that (\ref{eq:rond1}) holds. We will 
now proceed as in the proof of Proposition \ref{prop:exa}. Put
\begin{align}\label{eq:rond2}
&\alpha_{<a}=P_{Z,<a}\alpha \we 
P_{Z,<a},&\alpha_{>a}=P_{Z,>a}\alpha\we P_{Z,>a}.
\end{align}
Then $\alpha_{<a},\alpha_{>a}$ act as morphisms of degree $-1$ 
on 
$\mathcal{D}_{Z,<a}\left(Y,F\right),\mathcal{D}_{Z,>a}\left(Y,F\right)$. Using (\ref{eq:rond1}) and the second equation in (\ref{eq:rond1a2}), we get
\begin{align}\label{eq:rond3}
&\left[i_{Z},\alpha_{<a}\right]=1\,\mathrm{on}\,\mathcal{D}_{Z,<a}\left(Y,F\right),&\left[i_{Z},\alpha_{>a}\right]=1\,\mathrm{on}\,\mathcal{D}_{Z,>a}\left(Y,F\right).
\end{align}
By (\ref{eq:rond3}), $\alpha_{<a}$ and $\alpha_{>a}$ are homotopies 
of the complexes 
$\left(\mathcal{D}_{Z,<a}\left(Y,F\right),i_{Z}\right)$ and 
$\left(\mathcal{D}_{Z,>a}\left(Y,F\right),i_{Z}\right)$, so that these 
complexes are exact. The proof of our theorem is complete. 
\end{proof}
\begin{remark}\label{rem:dari}
In \cite{DangRiviere20}, Dang and Rivière had already proved that $P_{Z,<a}$ is a quasi-isomorphism with respect to $d^{Y}$.
\end{remark}

We take $a>0$ as in Definition \ref{def:sptru}. By Proposition \ref{prop:wff} and by Theorem \ref{thm:pdz1}, for $0<r<r_{0}$, if $z\notin \mathrm{Sp}L_{Z}$, 
$e^{-rZ}\left(L_{Z}-z\right)^{-1}$ and $Ne^{-rZ}\left(L_{Z}-z\right)^{-1}$ have a intersection supertrace, that 
depends meromorphically on $z$.\footnote{Let $R_{t}\vert_{t>0}$ be a family of smooth kernels that approximate the identity. Then $\Trs\left[e^{-rL_{Z}}\left(L_{Z}-z\right)^{-1}R_{t}\right]$ is holomorphic on the complement of $\mathrm{Sp}\,L_{Z}$. By making $t\to 0$, since the limit of bounded holomorphic functions is still holomorphic, $\mathrm{Tr}_{\mathrm{s,int}}\left[e^{-rL_{Z}}\left(L_{Z}-z\right)^{-1}\right]$ is holomorphic on the complement of $\mathrm{Sp}\,L_{Z}$. Near $z_{0}\in\mathrm{Sp}\,L_{Z}$, a similar argument  can be used to prove that $\mathrm{Tr}_{\mathrm{s,int}}\left[e^{-rL_{Z}}\left(L_{Z}-z\right)^{-1}\right]$ is meromorphic at $z_{0}$.}  In the above statements, we may  
replace $\left(L_{Z}-z\right)^{-1}$ by 
$\left(L_{Z}-z\right)^{-1}P_{Z,>a}$. Similarly by
(\ref{eq:doma14a4}), $P_{Z,<a}$ and $\left(L_{Z}-z\right)^{-1}P_{Z,<a}$ have well-defined 
intersection supertraces.

Note that $P_{Z,<a}$ is a projector from 
$\mathcal{D}_{\mathcal{Y}^{*}_{s}}\left(Y,F\right)$ on  the finite 
dimensional vector space 
$\mathcal{D}_{Z,<a}\left(Y,F\right)$. In particular, 
$P_{Z,z_{0}}\vert_{\mathcal{D}^{p}_{Z,z_{0}}\left(Y,F\right)}$ has a 
well-defined standard trace, which is just 
$\dim\mathcal{D}^{p}_{Z,<a}\left(Y,F\right)$. Now we prove a result 
which gives an analogue of Proposition \ref{prop:pidproj}.
\begin{proposition}\label{prop:trtrs}
	The following identities holds:
	\begin{align}\label{eq:tra1}
&\mathrm{Tr}^{\Omega_{p}\left(Y,F\right)}_{\mathrm{int}}\left[P_{Z,<a}\right]=\dim \mathcal{D}^{p}_{Z,<a}\left(Y,F\right),\\
&\mathrm{\Tr_{s,int}}\left[P_{Z,<a}\right]=\sum_{0}^{n}\left(-1\right)^{p}\dim\mathcal{D}^{p}_{Z,<a}\left(Y,F\right)=0. \notag 
\end{align}
Moreover, if $z\notin \mathrm{Sp}L_{Z}$, then
\begin{equation}\label{eq:trax1}
	\mathrm{\Tr_{s,int}}\left[P_{Z,<a}\left(L_{Z}-z\right)^{-1}\right]=0.
\end{equation}
\end{proposition}
\begin{proof}
	The proof is the same as the proof of Proposition \ref{prop:pidproj}. 
\end{proof}
\begin{proposition}\label{prop:scmat}
For $0<r<r_{0}$, if $z \notin\mathrm{Sp}L_{Z}$, then
\begin{align}\label{eq:chaeu}
&\Tr_{\mathrm{s,int}}\left[e^{-rL_{Z}}\left(L_{Z}-z\right)^{-1}\right]=0,
&\Tr_{\mathrm{s,int}}\left[e^{-rL_{Z}}P_{Z,>a}\left(L_{Z}-z\right)^{-1}\right]=0.
\nonumber 
\end{align}
\end{proposition}
\begin{proof}
Using the second identity in (\ref{eq:sabr1}), the last identity in (\ref{eq:rond1a2}), and proceeding as in the proofs of Propositions \ref{prop:pidproj} and  \ref{prop:trtrs}, our proposition follows. 
\end{proof}
\subsection{Truncated supertraces}%
\label{subsec:trstr}
For $z\notin \mathrm{Sp}L_{Z}$, given $a>0$, by (\ref{eq:romb0}),
\begin{equation}\label{eq:romb1}
N\left(L_{Z}-z\right)^{-1}=N\left(L_{Z}-z\right)^{-1}P_{Z,<a}+N\left(L_{Z}-z\right)^{-1}P_{Z,>a}.
\end{equation}
By the arguments that follow (\ref{eq:romb0}), the second term in the right-hand side of (\ref{eq:romb1}) 
extends to a meromorphic function of $z\in \C$ which is also
holomorphic on the disk $D_{a}=\left\{z\in \C,\left\vert  
z\right\vert<a\right\}$ with the same wave front set properties as 
before.

For $z\notin \mathrm{Sp}L_{Z}$ and $0<r<r_{0}$,  let us 
multiply (\ref{eq:romb1}) by $e^{-r \left(L_{Z}-z \right)}$ and take the 
corresponding intersection supertrace. We obtain
\begin{multline}\label{eq:romb2}
\Tr_{\mathrm{s,int}}\left[Ne^{-r\left(L_{Z}-z\right)}\left(L_{Z}-z\right)^{-1}\right]=
\Tr_{\mathrm{s,int}}\left[Ne^{-r\left(L_{Z}-z\right)}\left(L_{Z}-z\right)^{-1}P_{Z,<a}\right ]\\
+\Tr_{\mathrm{s,int}}\left[Ne^{-r\left(L_{Z}-z\right)}\left(L_{Z}-z\right)^{-1}P_{Z,>a}\right].
\end{multline}
\begin{proposition}\label{prop:inde}
	For $0<r<r_{0}, z\notin \mathrm{Sp}L_{Z}$,
	$$\Tr_{\mathrm{s,int}}\left[Ne^{-r \left( L_{Z}-z \right)}\left(L_{Z}-z\right)^{-1}\right]$$ 
	does not depend on $r$.
\end{proposition}
\begin{proof}
For $z\in\C, \Re\,z\ll -1$, this follows from Theorem \ref{thm:delo}. Our proposition follows by analytic continuation.
\end{proof}

 We denote the above common value as the 
regularized supertrace
\index[not]{Trsre@$\Tr_{\mathrm{s,int}}^{\mathrm{reg}}$}%
$$\Tr_{\mathrm{s,int}}^{\mathrm{reg}}\left[N\left(L_{Z}-z\right)^{-1}\right].$$
As $r\to 0$, 
\begin{equation}\label{eq:romb3}
\Tr_{\mathrm{s,int}}\left[Ne^{-r \left(L_{Z}-z\right)}\left(L_{Z}-z\right)^{-1}P_{Z,<a}\right 
]\to\Tr_{\mathrm{s,int}}\left[N\left(L_{Z}-z\right)^{-1}P_{Z,<a}\right ].
\end{equation}

Using (\ref{eq:romb2}), (\ref{eq:romb3}), we define the regularized supertrace 
$$\Tr_{\mathrm{s,int}}^{\mathrm{reg}}\left[N\left(L_{Z}-z\right)^{-1}P_{Z,>a}\right]$$
 by the formula
\begin{multline}\label{eq:romb4}
\Tr^{\mathrm{reg}}_{\mathrm{s,int}}\left[N\left(L_{Z}-z\right)^{-1}\right]=
\Tr_{\mathrm{s,int}}\left[N\left(L_{Z}-z\right)^{-1}P_{Z,<a}\right ]\\
+\Tr_{\mathrm{s,int}}^{\mathrm{reg}}\left[N\left(L_{Z}-z\right)^{-1}P_{Z,>a}\right].
\end{multline}
\begin{remark}\label{rem:repla}
	As in  (\ref{eq:ida2}), we get
\begin{equation}\label{eq:erda1}
N-\alpha i_{Z}=\left[i_{Z},\alpha N\right].
\end{equation}
Similarly, we have the identity of operators acting on  $\mathcal{D}_{Z,<a}\left(Y,F\right)$ or on 
$\mathcal{D}_{Z,>a}\left(Y,F\right)$, 
\begin{align}\label{eq:erda2}
&N-\alpha_{<a}i_{Z}=\left[i_{Z},\alpha_{<a}N\right],
&N-\alpha_{>a}i_{Z}=\left[i_{Z},\alpha_{>a}N\right].
\end{align}
It follows that in   the above supertraces or regularized supertraces, we can as well replace 
$N$ by $\alpha i_{Z},\alpha_{<a}i_{Z}$, or $\alpha_{>a}i_{Z}$. More 
generally, this will be the case for other analogous expressions.
\end{remark}
\begin{proposition}\label{prop:exten}
	The function 
	$\mathrm{Tr^{reg}_{s,int}}\left[N\left(L_{Z}-z\right)^{-1}P_{Z,>a}\right]$ 
	extends to a holomorphic function on $D_{a}$.
\end{proposition}
\begin{proof}
By (\ref{eq:romb2}), (\ref{eq:romb4}), for $r>0$ small 
enough, we get
\begin{multline}\label{eq:romb5}
\Tr_{\mathrm{s,int}}^{\mathrm{reg}}\left[N\left(L_{Z}-z\right)^{-1}P_{Z,>a}\right]
=\Tr_{\mathrm{s,int}}\left[Ne^{-r \left( L_{Z}-z \right) }\left(L_{Z}-z\right)^{-1}P_{Z,>a}\right]\\
+\Tr_{\mathrm{s,int}}\left[N \left( e^{-r \left( L_{Z}-z \right) }-1 \right) \left(L_{Z}-z\right)^{-1}P_{Z,<a}\right].
\end{multline}
The above equality is only valid for $z\notin\mathrm{Sp}L_{Z}$. 
However, the first term in the right-hand side extends 
holomorphically to $D_{a}$. 
We claim that for $r>0$, the second term also extends holomorphically on $D_{a}$. This is because the matrix 
$$\left(e^{-r\left(L_{Z}-z\right)}-1\right)\left(L_{Z}-z\right)^{-1}$$
extends to a holomorphic function of $\left(L_{Z}-z\right)$.
	
Combining (\ref{eq:romb5}) and the above concludes the proof.
\end{proof}

\subsection{The Fried zeta function as a meromorphic function}%
\label{subsec:mero}
Now we state a fundamental result of Giulietti-Liverani-Pollicott \cite[Corollary 2.2]{GiuliettiLiveraniPollicott13} and  Dyatlov-Zworski \cite[Theorem p.\,543]{DyatlovZworski16}. Here, we  follow the approach of Dyatlov-Zworski.
\begin{theorem}\label{thm:exten}
	The function $\frac{d}{d\sigma}\log R _{Z,F}\left( \mathrm{\sigma}\right) $ extends to a 
	meromorphic function of $\sigma$ with simple poles in 
	$-\mathrm{Sp}L_{Z}$ and integral 
	residues. More precisely, if $z_{0}\in 
	\mathrm{Sp}L_{Z}$, the residue at $-z_{0}$ is given 
	by $\Trs\left[NP_{Z,z_{0}}\right]$.

	The function $R_{Z,F}\left(\sigma\right)$ extends to a meromorphic 
	function on $\C$, whose zeros and poles are included in 
	$-\mathrm{Sp}L_{Z}$. If 
	$z_{0}\in\mathrm{Sp}L_{Z}$, the order of  $R_{Z,F}$ at $-z_{0}$ is given by 
	$\Trs\left[NP_{Z,z_{0}}\right]$.
\end{theorem}
\begin{proof}
	Take $a>0$ as before.  We use equation (\ref{eq:romb4}) in the form
\begin{multline}\label{eq:romb6}
\Tr^{\mathrm{reg}}_{\mathrm{s,int}}\left[N\left(L_{Z}+\sigma\right)^{-1}\right]=
\Tr_{\mathrm{s,int}}\left[N\left(L_{Z}+\sigma\right)^{-1}P_{Z,<a}\right ]\\
+\Tr_{\mathrm{s,int}}^{\mathrm{reg}}\left[N\left(L_{Z}+\sigma\right)^{-1}P_{Z,>a}\right].
\end{multline}
By Proposition \ref{prop:exten}, the second term in the right-hand 
side of (\ref{eq:romb6}) extends to a holomorphic function on  
 $\left\{\sigma\in \C,\left\vert  \sigma\right\vert<a\right\}$. The first 
term in the right-hand side of (\ref{eq:romb6}) is meromorphic on the same disk with poles 
of order $1$ and integral residues. More precisely, if $z_{0}\in 
\C,\left\vert  z_{0}\right\vert<a$ is in the spectrum of $L_{Z}$, the 
residue at $-z_{0}$ is given by  the integer 
$\Trs\left[NP_{Z,z_{0}}\right]$. 

Let $\pa$ be the classical Dolbeault operator on $\C$. By the above, if $\gamma$ is the meromorphic 
	$1$-form with integral residues $\pa\log R_{Z,F}\left(\sigma\right)$ on $\C$, for 
	$\sigma_{0}\in \C, \mathrm{Re}\,\sigma_{0}\gg 1$, we can define 
	$R_{Z,F}\left(\sigma\right)$ by the formula
	\begin{equation}\label{eq:romb7}
R_{Z,F}\left(\sigma\right)=R_{Z,F}\left(\sigma_{0}\right)\exp\left(\int_{\sigma_{0}}^{\sigma}\gamma\right),
\end{equation}
which completes the proof of our theorem.
\end{proof}

\subsection{The truncated Fried zeta function}%
\label{subsec:truz}
We fix $a>0$.
\begin{definition}\label{def:trulo}
For $\sigma\in \C$, put
\begin{equation}\label{eq:romb8}
R_{Z,F,<a}\left(\sigma\right)=\prod_{0}^{n}\det\left(\left(L_{Z}+\sigma\right)\vert_{\mathcal{D}^{i}_{Z,<a}\left(Y,F\right)}\right)^{\left(-1\right)^{i}i}.
\end{equation}
\end{definition}

The function $R_{Z,F,<a}\left(\sigma\right)$ is a meromorphic function 
whose zeroes and poles are included in $\mathrm{Sp}L_{Z}\cap 
\left\{z\in \C,\left\vert  z\right\vert<a\right\}$. If $z_{0}\in 
\C,\left\vert  z_{0}\right\vert<a$, the order of 
$R_{Z,F,<a}\left(\sigma\right)$ at $-z_{0}$ is just 
$\Trs\left[NP_{Z,z_{0}}\right]$. Moreover,
\begin{equation}\label{eq:romb9}
\frac{d}{d\sigma}\log R_{Z,F,<a}\left(\sigma\right)=\Trs\left[N\left(L_{Z}+\sigma\right)^{-1}P_{Z,<a}\right].
\end{equation}
\begin{definition}\label{def:rtru}
	For $\sigma\in \C,\Re\sigma\gg 1$, let $R_{Z,F,>a}\left(\sigma\right)$ be defined by 
	\begin{equation}\label{eq:romb9z1}
		R_{Z,F}\left(\sigma\right)=R_{Z,F,<a}\left(\sigma\right)R_{Z,F,>a}\left(\sigma\right).
	\end{equation}
\end{definition}
Then $R_{Z,F,>a}\left(\sigma\right)$ is a holomorphic function of 
$\sigma\in\C,\mathrm{Re}\,\sigma\gg 1$.
\begin{proposition}\label{prop:trho}
	The function $R_{Z,F,>a}\left(\sigma\right)$ extends meromorphically to a
	function of $\sigma\in\C$ with zeros and poles are the ones of 
	$R_{Z,F}\left(\sigma\right)$ that are contained in 
	$\left\{z\in\C,\left\vert  z\right\vert>a\right\}$. Moreover, we 
	have the identity,
	\begin{equation}\label{eq:romb10}
\frac{d}{d\sigma}\log R_{Z,F,>a}\left(\sigma\right)=\mathrm{Tr_{s,int}^{reg}}\left[N\left(L_{Z}+\sigma\right)^{-1}P_{Z,>a}\right].
\end{equation}
\end{proposition}
\begin{proof}
	The first part of our proposition is a consequence of Theorem 
	\ref{thm:exten} and of the definition of 
	$R_{Z,F,<a}\left(\sigma\right)$. By (\ref{eq:romb6}), 
	(\ref{eq:romb9}), we get (\ref{eq:romb10}). The proof of our proposition is complete. 
\end{proof}
\section{The determinant line and the section 
$\tau_{\nu}\left(i_{Z}\right)$}%
\label{sec:detli}
The purpose of this section is to construct a canonical nonzero 
section $\tau_{\nu}\left(i_{Z}\right)$ of 
the determinant line $\nu=\det H\left(Y,F\right)$.  This section is 
well-defined, even when $H\left(Y,F\right)\neq 0$. It will be obtained 
using the truncations of Section \ref{sec:frize}, and also the value 
at $0$ of the truncated Fried zeta function 
$R_{Z,F,>a}\left(\sigma\right)$.

This section is organized as follows. In Subsection 
\ref{subsec:frime}, we construct two complex lines $\lambda,\mu$ 
using the truncated determinants 
$\vartheta_{<a}=\det\mathcal{D}_{<a}\left(Z,F\right)$, and also the 
nonzero section 
$\tau_{\lambda}\left(i_{Z}\right)$ of $\lambda$.

In Subsection \ref{subsec:setn}, we construct the nonzero section $\tau_{\nu}\left(i_{Z}\right)$ of $\nu$.

In Subsection \ref{subsec:poinsec}, we show that the above section 
behaves properly under Poincaré duality.

In Subsection \ref{subsec:orun}, we specialize our results 
to the case considered in Subsection \ref{subsec:invo} where an 
involution $\iota$ acts on $Y$ that maps $Z$ to $-Z$, and where $F$ is 
unitary.

Finally, in Subsection \ref{subsec:chada}, when $Y$ is a contact manifold and when $Z$ is the associated Reeb vector field, we show that the section $\tau_{\nu} ^{\mathrm{CD}}$ of $\nu $ constructed by Chaubet-Dang \cite{ChaubetDang24} coincides with our section $\tau_{\nu } \left(i_{Z}\right)$.

We make the same assumptions as in Section \ref{sec:frize} and we use 
the corresponding notation.
\subsection{Two complex lines}%
\label{subsec:frime}
If $0<a<b<+ \infty $ are such that there is no eigenvalue $z_{0}$ of $L_{Z}$ 
such that $\left\vert  z_{0}\right\vert=a$ or $\left\vert  
z_{0}\right\vert=b$, set
\index[not]{DZab@$\mathcal{D}_{Z,(a,b)}\left(Y,F\right)$}%
\begin{equation}\label{eq:zomb16}
\mathcal{D}_{Z,(a,b)}\left(Y,F\right)=\bigoplus_{\substack{z_{0}\in\mathrm{Sp}L_{Z}
\\
a<\left\vert z_{0} \right\vert<b}}\mathcal{D}_{Z,z_{0}}\left(Y,F\right).
\end{equation}
Then
\begin{equation}\label{eq:zomb17}
\mathcal{D}_{Z,<b}\left(Y,F\right)=\mathcal{D}_{Z,<a}\left(Y,F\right) 
\oplus \mathcal{D}_{Z,(a,b)}\left(Y,F\right).
\end{equation}
\begin{definition}\label{def:libu}
	Set	
\index[not]{th@$\vartheta_{<a}$}%
\index[not]{thab@$\vartheta_{(a,b)}$}%
	\begin{align}\label{eq:bal11}
&\vartheta_{<a}=\det \mathcal{D}_{Z,<a}\left(Y,F\right),
&\vartheta_{(a,b)}=\det\mathcal{D}_{Z,(a,b)}\left(Y,F\right).
\end{align}
\end{definition}

Then we have the canonical isomorphism
\begin{equation}\label{eq:bal2}
\vartheta_{<b}=  \vartheta_{<a} \otimes \vartheta_{(a,b)}.
\end{equation}

	By Theorem \ref{thm:presa}, 
	$P_{Z,<a}:\left( \Omega\left(Y,F\right),d^{Y}\right) \to \left( 
	\mathcal{D}_{Z,<a}\left(Y,F\right),d^{Y} \right)$ is a 
	quasi-isomorphism, so that we have  the canonical isomorphism,
\begin{equation}\label{eq:zomb18}
\vartheta_{<a} \simeq \det H\left(Y,F\right).
\end{equation}
In particular, if $H\left(Y,F\right)=0$, equation just says that 
$\vartheta_{<a}$ is equipped with the canonical nonzero section
\index[not]{tad@$\tau_{<a}\left(d^{Y}\right)$}%
$\tau_{<a}\left(d^{Y}\right)$.

By Theorem \ref{thm:presa},  the complex 
$\left(\mathcal{D}_{Z,<a}\left(Y,F\right),i_{Z}\right)$ is exact, and so 
$\vartheta_{<a}$ has a canonical nonzero section 
\index[not]{taZ@$\tau_{<a}\left(i_{Z}\right)$}%
$\tau_{<a}\left(i_{Z}\right)$.
 Similarly $\theta_{(a,b)}$ is equipped with two 
canonical nonzero sections $\tau_{(a,b)}\left(d^{Y}\right)$ and 
$\tau_{(a,b)}\left(i_{Z}\right)$. Moreover,
\begin{equation}\label{eq:rex1}
\tau_{<b}\left(i_{Z}\right)=\tau_{<a}\left(i_{Z}\right) \otimes 
\tau_{(a,b)}\left(i_{Z}\right).
\end{equation}

By proceeding as in Definition \ref{def:trulo}, we  define the 
function $R_{Z,F,(a,b)}\left(\sigma\right)$. Then
\begin{equation}\label{eq:zomb19}
R_{Z,F,>a}\left(\sigma\right)=R_{Z,F,(a,b)}\left(\sigma\right)R_{Z,F,>b}\left(\sigma\right).
\end{equation}
The functions that appear in (\ref{eq:zomb19}) are holomorphic at 
$0$, where they take nonzero values.

By  (\ref{eq:dif15}), we get
	\begin{equation}\label{eq:zomb20}
\tau_{(a,b)}\left(d^{Y}\right)=R_{Z,F,(a,b)}\left(0\right)\tau_{(a,b)}\left(i_{Z}\right).
\end{equation}

\begin{definition}\label{def:libua1}
	Let 
	\index[not]{l@$\lambda$}%
	$\lambda$ be the complex line  which is canonically 
	identified with $\vartheta_{<a}$, $\vartheta_{<a}$ being 
	identified with $\vartheta_{<b}$ via the map $ s\in\vartheta_{<a}\to s \otimes 
	\tau_{(a,b)}\left(d^{Y}\right)\in \vartheta_{<b}$.
	
	Let 
	\index[not]{mu@$\mu$}%
	$\mu$ be the complex line  which is canonically 
	identified with $\vartheta_{<a}$, $\vartheta_{<a}$ being
	identified with $\vartheta_{<b}$ via the map $s\in\vartheta_{<a}\to s \otimes 
	\tau_{(a,b)}\left(i_{Z}\right)\in\vartheta_{<b}$.
\end{definition}
\begin{theorem}\label{thm:liba}
	The  line $\lambda$ is canonically isomorphic to $\det 
	H\left(Y,F\right)$.
	
	The nonzero sections $\tau_{<a}\left(i_{Z}\right)$ of $\vartheta_{<a}$ induce 
	a corresponding nonzero section 
	\index[not]{tmiZ@$\tau_{\mu}\left(i_{Z}\right)$}%
	$\tau_{\mu}\left(i_{Z}\right)$ of 
	$\mu$.
	
The map
	$s\in \vartheta_{<a}\to R^{-1}_{Z,F,>a}\left(0\right)s\in \vartheta_{<a}$
	induces a canonical isomorphism $\mu \simeq \lambda$.
\end{theorem}
\begin{proof}
	The identification (\ref{eq:zomb18}) is compatible with the 
	identification $\vartheta_{<a}$ via $\vartheta_{<b}$ that is 
	given in the definition of $\lambda$, which implies the first 
	part of our theorem.
	
	By (\ref{eq:rex1}), we get the corresponding result for $\mu$, 
	and we obtain this way the existence of the nonzero section 
	$\tau_{\mu}\left(i_{Z}\right)$.
	
	Finally, equations (\ref{eq:zomb19}),  (\ref{eq:zomb20}) show that  we have 
	indeed an isomorphism $\mu \simeq \lambda$. The proof of our theorem is complete. 
\end{proof}

Once $\lambda$ and $\mu$ have been canonically identified, the 
section $\tau_{\mu}\left(i_{Z}\right)$ can be viewed as a canonical 
nonzero section of $\lambda$.
\begin{definition}\label{def:cansec}
	In the sequel, we will denote by 
	\index[not]{tliZ@$\tau_{\lambda}\left(i_{Z}\right)$}%
	$\tau_{\lambda}\left(i_{Z}\right)$ the nonzero canonical 
	section of $\lambda$ that was constructed before.
\end{definition}
\begin{remark}\label{rem:expli}
	If $\mathcal{D}_{Z,0}\left(Y,F\right)=0$, the sections 
	$\tau_{<a}\left(d^{Y}\right)$ define a corresponding nonzero 
	section $\tau_{\lambda}\left(d^{Y}\right)$ of $\lambda$, so that
	  $\lambda \simeq  \C$. Also $R_{Z,F}\left(0\right)$ is a 
	well defined nonzero complex number. Moreover, we have the 
	identity of nonzero sections of $\lambda$,
	\begin{equation}\label{eq:co1}
\tau_{\lambda}\left(i_{Z}\right)=R^{-1}_{Z,F}\left(0\right)\tau_{\lambda}\left(d^{Y}\right).
\end{equation}

More generally, if $H\left(Y,F\right)=0$,  then 
$\tau_{\lambda}\left(d^{Y}\right),\tau_{\lambda}\left(i_{Z}\right)$ are still nonzero sections of $\lambda$.  Then $\mathcal{D}_{0}\left(Y,F\right)$ may 
not vanish, but the complex 
$\left(D_{0}\left(Y,F\right),d^{Y}\right)$ is still exact. Equation 
(\ref{eq:co1}) is replaced by
\begin{equation}\label{eq:co2}
\frac{\tau_{\lambda}\left(d^{Y}\right)}{\tau_{\lambda}\left(i_{Z}\right)}=R_{Z,F,>0}\left(0\right)
\frac{\tau\left(d^{Y}\vert_{\mathcal{D}_{0}\left(Y,F\right)} \right) 
}{\tau\left(i_{Z}\vert_{\mathcal{D}_{0}\left(Y,F\right)} \right) }.
\end{equation}

In the general case, let $\gamma\in \det H\left(Y,F\right)$. Then 
$\gamma$ can be viewed as a section $\gamma_{\lambda}$ of $\lambda$, and a 
section $\gamma_{\vartheta_{0}}$ of $\vartheta_{0}=\det 
\mathcal{D}_{0}\left(Y,F\right)$. Moreover, 
\begin{equation}\label{eq:co3}
\frac{\gamma_{\lambda}}{\tau_{\lambda}\left(i_{Z}\right)}=R_{Z,F,>0}\left(0\right)\frac{\gamma_{\vartheta_{0}}}{\tau\left(i_{Z}\vert_{\mathcal{D}_{0}\left(Y,F\right)}\right)}.
\end{equation}
\end{remark}
\subsection{The section $\tau_{\nu}\left(i_{Z}\right)$}%
\label{subsec:setn}
Put
\index[not]{nu@$\nu$}%
\begin{equation}\label{eq:nua1}
\nu=\det H\left(Y,F\right).
\end{equation}
 By Theorem \ref{thm:liba}, we have the canonical isomorphism
 \begin{equation}\label{eq:nua2}
 \lambda \simeq  \nu.
 \end{equation}
 \begin{definition}\label{def:tn}
 We denote by 
 \index[not]{tniZ@$\tau_{\nu}\left(i_{Z}\right)$}%
 $\tau_{\nu}\left(i_{Z}\right)$ the section of $\nu$ that corresponds to $\tau_{\lambda}\left(i_{Z}\right)$ via the isomorphism in (\ref{eq:nua2}).
 \end{definition}
\subsection{The section 
$\tau_{\lambda}\left(i_{Z}\right)$ and Poincaré 
duality}%
\label{subsec:poinsec}
In this subsection, for greater clarity, the dependence of the above 
 objects on $Z,F$ will be denoted explicitly with the subscript $Z,F$. In particular,
 \index[not]{nF@$\nu_{F}$}%
\begin{equation}\label{eq:co3a1}
\nu_{F}=\det H\left(Y,F\right).
\end{equation}
Also 
\index[not]{tlZ@$\tau_{\lambda_{Z,F}}\left(i_{Z}\right)$}%
$\tau_{\lambda_{Z,F}}\left(i_{Z}\right)$ is 
 a nonzero section of $\lambda_{Z,F}$.
 
Using Poincaré duality, we  have the canonical isomorphism,
\begin{equation}\label{eq:co3a3}
\nu_{F^{*} \otimes 
_{\Z_{2}}\mathrm{o}\left(TY\right)}=\nu_{F}^{\left(-1\right)^{n-1}}.
\end{equation}
By Theorem \ref{thm:podua1}, we have the canonical isomorphism,
 \begin{equation}\label{eq:zeph1}
\lambda_{-Z,F^{*} \otimes 
_{\Z_{2}}\mathrm{o}\left(TY\right)}=\lambda_{Z,F} 
^{\left(-1\right)^{n-1}}.
\end{equation}
The identifications in (\ref{eq:co3a3}), (\ref{eq:zeph1}) are 
compatible.

Also $\tau_{\lambda_{-Z,\Ft}}\left(i_{-Z}\right)$ is a nonzero 
section of $\lambda_{-Z,\Ft}$. 
By (\ref{eq:zeph1}),  $\left[\tau_{\lambda_{-Z,\Ft}}\left(i_{-Z}\right)\right]^{\left(-1\right)^{n-1}}$ is a 
section of $\lambda_{Z,F}$.

\begin{theorem}\label{thm:podu}
	The following identities hold:
	\begin{align}\label{eq:zeph2}
&\left[\tau_{\lambda_{-Z,F^{*} \otimes 
_{\Z_{2}}\mathrm{o}\left(TY\right)}}\left(i_{-Z}\right)\right]^{\left(-1\right)^{n-1}}=\tau_{\lambda_{Z,F}} \left(i_{Z}\right)\, \mathrm{in}\, \lambda_{Z,F},\\
&\left[\tau_{\nu_{\Ft}}\left(i_{-Z}\right)\right]^{\left(-1\right)^{n-1}}=\tau_{\nu_{F}}\left(i_{Z}\right)
\,\mathrm{in}\,\nu_{F}. \notag 
\end{align}
\end{theorem}
\begin{proof}
	By (\ref{eq:normi4a}), we know that
	\begin{equation}\label{eq:zeph3}
i_{-Z}^{\dag}=i_{Z}.
\end{equation}
By (\ref{eq:dua2}), (\ref{eq:zeph3}), we obtain
\begin{equation}\label{eq:zeph4}
\left( \tau_{-Z,F^{*} \otimes 
_{\Z_{2}}\mathrm{o}\left(TY\right),<a}\left(i_{-Z}\right) 
\right)^{\left(-1\right)^{n-1}}=\tau_{Z,F,<a}\left(i_{Z}\right)\,\mathrm{in}\, \vartheta_{Z,F,<a}.
\end{equation}

Moreover, one has the easy identity,
\begin{equation}\label{eq:zeph5}
R_{-Z,F^{*} \otimes 
_{\Z_{2}}\mathrm{o}\left(TY\right),<a}\left(\sigma\right)=R_{Z,F,<a}^{\left(-1\right)^{n-1}}\left(\sigma\right).
\end{equation}
Using Theorem \ref{thm:poin} and (\ref{eq:zeph5}), for 
$\sigma\in\C,\mathrm{Re}\,\sigma\gg 1$, we get
\begin{equation}\label{eq:zeph6}
R_{-Z,F^{*} \otimes 
_{\Z_{2}}\mathrm{o}\left(TY\right),>a}\left(\sigma\right)=R_{Z,F,>a}^{\left(-1\right)^{n-1}}\left(\sigma\right).
\end{equation}
Since both sides of (\ref{eq:zeph6}) extend holomorphically at 
$\sigma=0$, we obtain
\begin{equation}\label{eq:zeph7}
R_{-Z,F^{*} \otimes 
_{\Z_{2}}\mathrm{o}\left(TY\right),>a}\left(0\right)=R_{Z,F,>a}^{\left(-1\right)^{n-1}}\left(0\right).
\end{equation}

By (\ref{eq:zeph4}), (\ref{eq:zeph7}), we deduce that
\begin{multline}\label{eq:zeph8}
\left[ R^{-1}_{-Z,F^{*} \otimes 
_{\Z_{2}}\mathrm{o}\left(TY\right),>a}\left(0\right)\tau_{-Z,F^{*} \otimes 
_{\Z_{2}}\mathrm{o}\left(TY\right),<a}\left(i_{-Z}\right)\right]^{\left(-1\right)^{n-1}}\\
=R_{Z,F,>a}^{-1}\left(0\right)\tau_{Z,F,<a}\left(i_{Z}\right)\,\mathrm{in}\,\vartheta_{Z,F,<a},
\end{multline}
which is just the first identity in  (\ref{eq:zeph2}). The second identity follows from the first one.  The proof of our theorem is complete. 
\end{proof}
\subsection{The case where $Y$ is orientable and $F$ is unitary}%
\label{subsec:orun}
In this subsection, we assume that the assumptions of Subsection 
\ref{subsec:invo} are verified. In particular, $n$ is odd. Also we 
assume that $Y$ is orientable, and that $\n^{F}$ is unitary with respect to a Hermitian metric $g^{F }$. 

Then $\mathrm{o}\left(TY\right)$ is a constant $\Z_{2}$-space with 
two canonical sections. In particular,
\begin{equation}\label{eq:fla1}
\det H^{p}\left(Y,F^{*} \otimes 
_{\Z_{2}}\mathrm{o}\left(TY\right)\right)=\det H^{p}\left(Y,F^{*}\right) 
\otimes _{\Z_{2}}\mathrm{o}^{\dim H^{p}\left(Y,F^{*}\right)}\left(TY\right).
\end{equation}
Let $\chi\left(Y,F^{*}\right)$ be the Euler characteristic of $\left(\Omega\left(Y,F^{*}\right),d^{Y}\right)$. By (\ref{eq:fla1}), we obtain 
\begin{equation}\label{eq:fla2}
\nu_{F^{*} \otimes 
_{\Z_{2}}\mathrm{o}\left(TY\right)}=\nu_{F^{*}} \otimes 
\mathrm{o}^{\chi\left(Y,F^{*}\right)}\left(TY\right).
\end{equation}
Since $n$ is odd,  $\chi\left(Y,F^{*}\right)=0$. From (\ref{eq:fla2}), we get
\begin{equation}\label{eq:fla3}
\nu_{F^{*} \otimes 
_{\Z_{2}}\mathrm{o}\left(TY\right)}=\nu_{F^{*}}.
\end{equation}
By (\ref{eq:co3a3}), (\ref{eq:fla3}), we obtain
\begin{equation}\label{eq:zeph8a-1}
\nu_{F^{*}} =\nu_{F}
\end{equation}

Since $\n^{F}$ is unitary, we have an 
identification of flat vector bundles,
\begin{equation}\label{eq:zeph8a0}
F^{*} \simeq \overline{F}.
\end{equation}
The identification in (\ref{eq:zeph8a0}) depends on the choice of 
$g^{F}$. 
By (\ref{eq:zeph8a0}), we deduce that
\begin{equation}\label{eq:zeph8a-1x}
\nu_{F^{*}} \simeq \nu_{\overline{F}}.
\end{equation}
By (\ref{eq:zeph8a-1}), (\ref{eq:zeph8a-1x}), we obtain
\begin{equation}\label{eq:zeph8a-2}
\nu_{F} \simeq \nu_{\overline{F}}.
\end{equation}
Also, 
\begin{equation}\label{eq:zepha-3}
\nu_{\overline{F}}=\overline{\nu_{F}}.
\end{equation}
By (\ref{eq:zeph8a-2}), (\ref{eq:zepha-3}), we have an isomorphism,
\begin{equation}\label{eq:zepha-4}
\nu_{F}=\overline{\nu_{F}}.
\end{equation}
By (\ref{eq:zepha-4}), we conclude that $\nu_{F}$ is  
the complexification of a real line, the real line consisting of 
the conjugation invariant elements of $\nu_{F}$, this real line not depending on the choice of $g^{F}$.

Observe that $\iota$ acts on $\nu_{F}$ like $\pm 1$.
\begin{proposition}\label{prop:ide}
	The action of $\iota$ on $\nu_{F}$ is given by
	\begin{equation}\label{eq:zepha-5}
\iota_{\nu_{F}}=\iota\vert_{\mathrm{o}\left(TY\right)}^{\sum_{p\le 
\left(n-1\right)/2}^{}\dim H^{p}\left(Y,F\right)}.
\end{equation}
In particular, if $\iota\vert_{\mathrm{o}\left(TY\right)}=1$, then 
$\iota\vert_{\nu_{F}}=1$.
\end{proposition}
\begin{proof}
	We may and we will assume that $g^{F}$ is an $\iota$-invariant 
	metric. Clearly $\iota$ acts on $H^{p}\left(Y,F\right)$ as an 
	involution, and 
	 the action on 
	$H^{p}\left(Y,\overline{F}\right)$ is just the conjugate action. By Poincaré duality, the 
	actions of $\iota$ on $\det H^{p}\left(Y,F\right)$ and 
	$\det H^{n-p}\left(Y,F^{*} \otimes_{\Z_{2}} 
	\mathrm{o}\left(TY\right)\right)$ are given by the same scalar.
By (\ref{eq:fla1}), we find that
\begin{equation}\label{eq:fla5}
\iota\vert_{\det H^{n-p}\left(Y,F ^{*}\otimes _{\Z_{2}}\mathrm{o}\left(TY\right)\right)}=\iota\vert_{\det  
H^{n-p}\left(Y,F^{*}\right)}\iota\vert_{\mathrm{o}\left(TY\right)}^{\dim H^{n-p}\left(Y,F^{*}\right)}.
\end{equation}

By the above, we deduce that
\begin{equation}\label{eq:fla6}
\iota_{\nu_{F}}=\iota\vert_{\mathrm{o}\left(TY\right)}^{\sum_{p\ge\left(n+1\right)/2}^{}\dim 
H^{n-p}\left(Y,F^{*}\right)}.
\end{equation}
In the right-hand side of (\ref{eq:fla6}), we may as well replace 
$F^{*}$ by $F$.
Since the Euler characteristic of $Y$ with coefficients in $F$ 
vanishes, (\ref{eq:zepha-5}) is an equivalent form of (\ref{eq:fla6}). 
The proof of our proposition is complete. 
\end{proof}

Put
\index[not]{nF@$\nu_{F}^{<n/2}$}%
\begin{equation}\label{eq:aujo1}
\nu_{F}^{<n/2} =\bigotimes _{p<n/2}\det 
H^{p}\left(Y,F\right)^{\left(-1\right)^{p}},
\end{equation}
Using Poincaré duality, we get
\begin{equation}\label{eq:aujo2}
\nu_{F}=\nu_{F}^{<n/2} \otimes \nu^{<n/2}_{F^{*} \otimes _{\Z_{2}} 
\mathrm{o}\left(TY\right)}.
\end{equation}
Equation (\ref{eq:aujo2}) can be rewritten in the form
\begin{equation}\label{eq:aujo3}
\nu_{F}=\nu^{<n/2}_{F} \otimes \nu^{<n/2}_{F^{*}} \otimes \mathrm{o}
^{\sum_{p<n/2}^{}\dim H^{p}\left(Y,F\right)}\left(TY\right).
\end{equation}

By choosing a flat metric $g^{F}$, equation (\ref{eq:aujo3}) can be written in the 
form,
\begin{equation}\label{eq:aujo4}
\nu_{F}=\nu^{<n/2}_{F} \otimes \overline{\nu^{<n/2}_{F}}
 \otimes \mathrm{o}
^{\sum_{p<n/2}^{}\dim H^{p}\left(Y,F\right)}\left(TY\right).
\end{equation}
By (\ref{eq:aujo4}), $\nu_{F}$ is a real line bundle. 
Since the action of $\iota$  on $\nu^{<n/2}_{F} \otimes 
\overline{\nu^{<n/2}_{F}}$ is trivial, from (\ref{eq:aujo4}), we 
recover Proposition \ref{prop:ide}.

Observe that $\iota$ maps 
$\lambda_{Z,F},\mu_{Z,F}$ into $\lambda_{-Z,F}, \mu_{-Z,F}$.   
By Proposition \ref{prop:ide}, the map $\iota$ is not the 
trivial one obtained via the identifications $\lambda_{Z,F} \simeq 
\nu_{F},\lambda_{-Z,F} \simeq \nu_{F}$. 

By Poincaré duality, we
get
\begin{equation}\label{eq:fla7-a1}
\vartheta_{Z,F,<a} =\vartheta_{-Z,\Ft,<a}.
\end{equation}

The Euler characteristic of the complex 
$\left( \mathcal{D}_{-Z,<a}\left(Y,F^{*}\otimes 
_{\Z_{2}}\mathrm{o}\left(TY\right)\right), i_{Z} \right) $ vanishes. As in 
(\ref{eq:fla3}), we get
\begin{equation}\label{eq:fla7}
\vartheta_{-Z,F^{*} \otimes 
_{\Z_{2}}\mathrm{o}\left(TY\right),<a}=\vartheta_{-Z,F^{*},<a}.
\end{equation}
In particular, the actions of $\iota$ on both sides of 
(\ref{eq:fla7}) are the same. Note that $\iota$ maps 
$\vartheta_{-Z,F^{*},<a}$ into $\vartheta_{Z,F^{*},<a}$.

Again, we have an isomorphism,
\begin{equation}\label{eq:fla8}
\vartheta_{Z,F^{*},<a} \simeq \vartheta_{Z,\overline{F},<a}.
\end{equation}
 
By (\ref{eq:fla7-a1})--(\ref{eq:fla8}), we ultimately get an isomorphism,
\begin{equation}\label{eq:fla10}
\iota_{Z,F,<a}: \vartheta_{Z,F,<a}\to \vartheta_{Z,\overline{F},<a}.
\end{equation}
One verifies easily that 
\begin{equation}\label{eq:fla11a1}
\iota_{Z,F,<a}\tau_{Z,F,<a}\left(i_{Z}\right)=\overline{\tau_{Z,F,<a}\left(i_{Z}\right)}.
\end{equation}
Similar considerations show that
\begin{equation}\label{eq:fla12}
\iota_{Z,F,<a}\tau_{Z,F,<a}\left(d^{Y}\right)=\overline{\tau_{Z,F,<a}\left(d^{Y}\right)}.
\end{equation}

Using the above, we find that
\begin{equation}\label{eq:fla13}
R_{Z,F,<a}\left(\sigma\right)=\overline{R_{Z,F,<a}\left(\overline{\sigma}\right)}.
\end{equation}
By combining Proposition \ref{prop:flip} and (\ref{eq:fla13}), for 
$\sigma\in\C,\mathrm{Re}\,\sigma\gg 1$, we get
\begin{equation}\label{eq:fla14}
R_{Z,F,>a}\left(\sigma\right)=\overline{R_{Z,F,>a}\left(\overline{\sigma} \right) }.
\end{equation}
In particular,
\begin{equation}\label{eq:fla15}
R_{Z,F,>a}\left(0\right)=\overline{R_{Z,F,>a}\left(0\right) }.
\end{equation}
\begin{theorem}\label{thm:conju}
	The following identities hold:
	\begin{align}\label{eq:fla16}
&\iota\tau_{\lambda_{Z,F}}\left(i_{Z}\right)=\tau_{\lambda_{-Z,F}}\left(i_{-Z}\right)=\overline{\tau_{\lambda_{Z,F}}\left(i_{Z}\right)}\,\text{in}\, \lambda_{Z,\overline{F}},\\
&\iota\tau_{\nu_{F}}\left(i_{Z}\right)=\tau_{\nu_{F}}\left(i_{-Z}\right)=\overline{\tau_{\nu_{F}}\left(i_{Z}\right)}\,\text{in}\,\nu_{\overline{F}}. \notag 
\end{align}
\end{theorem}
\begin{proof}
	The first identity  is an obvious consequence of (\ref{eq:fla11a1}) and 
	(\ref{eq:fla15}). From the first identity, we obtain the second identity, which completes the proof of our theorem.
\end{proof}

Set
\begin{equation}\label{eq:fla16a1}
\vartheta^{<n/2}_{Z,F}=\bigotimes_{p<n/2}\left[\det\mathcal{D}^{p}_{Z,<a}\left(Y,F\right)\right]^{\left(-1\right)^{p}}.
\end{equation}
Using Poincaré duality, we get
\begin{equation}\label{eq:fla17}
\vartheta_{Z,F,<a}=\vartheta^{<n/2}_{Z,F,<a} \otimes \vartheta^{<n/2}_{-Z,F^{*} 
\otimes _{\Z_{2}}\mathrm{o}\left(TY\right),<a}.
\end{equation}
\subsection{Comparison with the  section of Chaubet-Dang}%
\label{subsec:chada}
In this subsection, we follow Chaubet-Dang \cite{ChaubetDang24}. We assume that $Y$ is a contact manifold. In particular, $n$ is odd. Let $\alpha$ be the contact form, which is such that $ \alpha\we \left(d^{Y}\alpha 
\right)^{\left(n-1\right)/2}\neq 0$. We assume that $Z$ is the Reeb vector field associated with $\alpha$, so that
\begin{align}\label{eq:chab1}
&i_{Z}\alpha =1,&i_{Z}d \alpha  =0.
\end{align}
In particular,
\begin{align}\label{eq:chab2}
&\left[ \alpha ,i_{Z}\right]=1, &\left[L_{Z},\alpha\right] =0, && \left[L_{Z},i_{Z}\right]=0.
\end{align}

Using the above data, in \cite[Definition 6.1]{ChaubetDang24}, Chaubet-Dang defined a smooth odd morphism $ \Gamma $ acting on $\Lambda \left(T^{*}Y\right)$, that verifies pointwise the assumptions of Subsection   \ref{subsec:spe} with respect to $d=\alpha ,\delta =i_{Z}$, with $p=0,q=n$, and is such that
\begin{equation}\label{eq:chab3}
L_{Z}\Gamma =0.
\end{equation}

The action of $ \Gamma $ obviously extends to $\Omega\left(Y,F\right)$, and  commutes with $L_{Z}$. Using the notation of Subsection \ref{subsec:vspa}, $\Gamma $ also acts on $\mathcal{D}_{Z,<a}\left(Y,F\right)$. Let $\Gamma _{<a}$ be the restriction of $\Gamma $ to $\mathcal{D}_{Z,<a}\left(Y,F\right)$. 

Using Theorem \ref{thm:idt}, we get
\begin{equation}\label{eq:chab4}
\tau_{<a}\left(i_{Z}\right)=\rho _{\Gamma_{<a}}. 
\end{equation}

Let 
\index[not]{tCD@$\tau_{\nu }^{\mathrm{CD}}$}%
$\tau_{\nu }^{\mathrm{CD}}$ be the section of $\nu$ constructed by Chaubet-Dang \cite[Definition 6.3]{ChaubetDang24}.
\begin{theorem}\label{thm:chada}
	The following identity holds:
	\begin{equation}\label{eq:chab5}
\tau_{\nu}\left(i_{Z}\right)=\tau_{\nu } ^{\mathrm{CD}}\, \mathrm{in} \, \nu .
	\end{equation}
\end{theorem}
\begin{proof}
In the construction of their section $\tau_{\nu }^{\mathrm{CD}} $, instead of $\tau_{<a}\left(i_{Z}\right) $, Chaubet and Dang use  $\rho _{\Gamma _{<a}}$. However, by (\ref{eq:chab4}), these two sections coincide. 
Using the fact that the gluing construction in \cite[Definition 6.3 and Remark 6.4]{ChaubetDang24} for the section $\tau_{\nu } ^{\mathrm{CD}}$  is the same as ours, we get (\ref{eq:chab5}). The proof of our theorem is complete. 
\end{proof}

\section{The determinant line bundle}%
\label{sec:para}
The purpose of this section is to show that our previous 
constructions of the determinant line can also be done in  families, 
exactly like in the theory of determinant line bundles associated 
with relative de Rham or Dolbeault complex.  In particular, if $F$ 
is now flat on the total space of a fibration with fiber $Y$, we show 
that $\tau_{\nu}\left(i_{Z}\right)$ is 
flat with respect to  the Gauss-Manin connection on $\nu$.

This section is organized as follows. In Subsection \ref{subsec:stst}, we recall basic results of Llave, Marco and Moriyòn on the structural stability of Anosov vector fields. 

In Subsection 
\ref{subsec:faso}, if $\pi:M\to S$ is a submersion of smooth 
manifolds with compact fiber $Y$, we construct a family $Z$ of vector 
fields along the fibers $Y$ that are Anosov along the fibers. 

In Subsection \ref{subsec:defr}, we recall a result by \cite{DangShen20} asserting that for $\sigma\in\C,\Re\,\sigma\gg 1$, the fiberwise Fried zeta function $R_{Z,F}\left(\sigma\right)$ is also a $C^{1}$ function of $s\in S$.

In Subsection \ref{subsec:smotr}, we show that in the proper sense, the corresponding 
resolvents are $C^{1}$ in the parameter $s\in S$. 

In Subsection \ref{subsec:runpro}, we establish similar results for 
the truncation projectors $P_{Z,<a}$. Also we prove that for $a>0$, on the 
open set $V_{a}$ where no eigenvalue of  $L_{Z}$ is contained in $S_{a}=\left\{z\in\C,\left\vert z\right\vert=a\right\}$, 
the $D_{Z,<a}\left(Y,F\right)$ form a $C^{1}$ vector bundle. More precisely, near 
$s_{0}\in V_{a}$, we construct an explicit $C^{1}$ trivialization of 
$\mathcal{D}_{Z,<a}\left(Y,F\right)$. 

In Subsection \ref{subsec:derfl}, we compute the proper derivative of 
the action of $i_{Z}$ on $\mathcal{D}_{Z,<a}\left(Y,F\right)$ with respect to the 
above trivialization.

Finally, in Subsection \ref{subsec:flam}, we consider the case where 
$F$ is flat on $M$.
\subsection{Structural stability of Anosov vector fields}%
\label{subsec:stst}
We make the same assumptions as in Section \ref{sec:frizet}. In particular $Y$ is a compact manifold, and $Z$ is an Anosov vector field on $Y$. 

A fundamental result of Llave, Marco and Moriyón \cite[Apprendix, Theorem A.1]{Llave86} asserts that Anosov vector fields are structurally stable. This implies in particular that if $\mathcal{O} \subset \R^{m}$ is an open set, and if $Z_{s}\vert_{s\in\mathcal{O}}$ is a smooth family of Anosov vector fields on $Y$, if $s_{0}\in \mathcal{O}$, for $\mathcal{O}$ small enough,  there exist a function $\tau_{s}\left(y\right):\mathcal{O}\times Y\to \R_{+}^{*}$ and a map $h_{s}\left(y\right):\mathcal{O}\times Y\to Y$, which are continuous, which for fixed $y\in Y$ are $C^{ \infty }$ in $s\in\mathcal{O}$ (the derivatives being themselves jointly continuous), and for given $s\in \mathcal{O}$, $h_{s}$ is a homeomorphism of $Y$,   is $C^{1}$ along $Z_{s_{0}}$, and  such that $Z_{s}\left(h_{s}\left(y\right)\right)=h_{s*}\tau_{s}\left(y\right)Z_{s_{0}}\left(y\right)$.

Let $\underline{\mathbf{F}}_{s}$ denote the quotient fixed point set associated with $Z_{s}$. By the above, for $s\in \mathcal{O}$ close enough to $s_{0}$, 
\begin{equation}\label{eq:fixa1}
\underline{\mathbf{F}}_{s}=h_{s}\underline{\mathbf{F}}_{s_{0}}.
\end{equation}

Given $y\in \underline{\mathbf{F}}_{s}$, let $t_{s,y}$ be the analogue of $t_{y}$ for the vector field $Z_{s}$.
Also by proceeding as in \cite[Section 4.1, Remark 4]{DangShen20}, we deduce that if $\mathcal{O}$ is small enough, if $y\in \underline{\mathbf{F}}_{s_{0}}$, then
\begin{equation}\label{eq:fixa2}
\frac{1}{2}t_{s_{0},y}\le t_{s,h_{s}y}\le 2t_{s_{0},y}.
\end{equation}
By (\ref{eq:fixa2}),  for $\mathcal{O}$ small enough,   Proposition \ref{prop:gro} holds uniformly on $\mathcal{O}$. More precisely, if $\mathcal{O}$ is small enough, there exist $c>0,C>0$ such that for $s\in \mathcal{O}$, 
\begin{equation}\label{eq:fixa3}
\left\vert \left\{y\in\underline{\mathbf{F}}_{s}, t_{s,y}\le T\right\}\right\vert\le Ce^{cT}.
\end{equation}
 \subsection{A family of Anosov vector fields}%
\label{subsec:faso}
Let $\pi:M\to S$ be a proper submersion of smooth manifolds, with 
compact fiber $Y$ of dimension $n$. Let 
$\left(F,\n^{F}\right)$ be a complex  
vector bundle on $M$ with  connection which is flat along the fibers 
$Y$.\footnote{Let  
\index[not]{ga@$\gamma$}%
$\gamma:Y\to M$ be the  embedding of the generic fiber $Y$. Let $R^{F}$ be the curvature of $\n^{F}$. Our assumption just says that $\gamma^{*}R^{F}=0$.} We have the exact sequences of vector bundles on $M$,
\begin{align}\label{eq:co4}
&0\to TY\to TM\to\pi^{*}TS\to 0,\\
&0\to \pi^{*}T^{*}S\to T^{*}M\to T^{*}Y\to 0. \nonumber 
\end{align}

Let $\left(\Omega\left(Y,F\right),d^{Y}\right)$ denote the de Rham 
complex along the fibers $Y$, so that
\begin{equation}\label{eq:co5}
C^{\infty } \left(S,\Omega\left(Y,F\right) \right)=C^{\infty }\left(M,\Lambda\left(T^{*}Y\right) 
\otimes _{\R}F\right),
\end{equation}
and $d^{Y}$ is the fiberwise de Rham operator.

Let $\nu$ denote the  line bundle on $S$ which is the 
determinant of the cohomology of $\left(\Omega\left(Y,F\right),d^{Y}\right)$ along $Y$. The theory of 
determinant bundles \cite{Quillen85b,BismutFreed86a,BismutFreed86b} asserts that $\nu$ is a smooth line bundle on 
$S$, and moreover for any $s\in S$, we have a canonical isomorphism 
\begin{equation}\label{eq:co6}
\nu_{s} \simeq \det H\left(Y_{s},F\vert_{Y_{s}}\right).
\end{equation}

Let $Z$ be a smooth section of $TY$.
In the sequel, we assume that $Z$ is a fiberwise  Anosov vector 
field. By \cite[Chapter 2, \S 7]{Anosov69}, \cite[Appendix 1]{DangShen20}, $T_{u}Y,T_{s}Y$ are continuous 
		vector subbundles of $TY$.

By the results of Section \ref{sec:detli}, for each fiber $Y_{s}$, given the Anosov vector field 
$Z\vert_{Y_{s}}$, we can construct the complex line $\lambda_{s}$, 
and its canonical nonzero section 
$\tau_{\lambda_{s}}\left(i_{Z\vert_{Y_{s}}}\right)$. By Theorem 
\ref{thm:liba}, 
\begin{equation}\label{eq:co7}
\lambda_{s} \simeq \det H\left(Y_{s},F\vert_{Y_{s}}\right).
\end{equation}
Using (\ref{eq:co6}), (\ref{eq:co7}),  we have a canonical 
isomorphism,
\begin{equation}\label{eq:smo1a1}
\lambda_{s} \simeq \nu_{s}.
\end{equation}

Similarly, for  $s\in S$, the complex line $\mu_{s}$ is 
well-defined, and $\tau_{\mu_{s}}\left(i_{Z\vert_{Y_{s}}}\right)$ is a 
section of $\mu_{s}$.
\subsection{The $C^{1}$ dependence of the Fried zeta function}%
\label{subsec:defr}
For any $s\in S$, we can define the fiberwise zeta function $R_{Z\vert_{Y_{s}},F\vert_{Y_{s}}}\left(\sigma\right)$ for $\sigma\in\C,\Re\,\sigma\gg 1$. 

Let us now recall a basic result of Dang, Guillarmou, Rivière, and Shen \cite[Lemma 4.3]{DangShen20}. Its proof is based on the results of Llave, Marco and Moriyón \cite{Llave86} that were described in Subsection \ref{subsec:stst}.
\begin{theorem}\label{thm:derv}
If $s_{0}\in S$, and if $\mathcal{O} \subset S $ is a small open neighborhood of $s_{0}$, there exists $C>0$ such that for $\sigma\in \C,\Re\,\sigma\ge C$, the function $R_{Z,F}\left(\sigma\right)$ is both  holomorphic in $\sigma$ and $C^{1}$ in $s\in\mathcal{O}$, with uniform bounds on $R_{Z,F}\left(\sigma\right), \frac{\pa}{\pa s}R_{Z,F}\left(\sigma\right)$.
\end{theorem}
\subsection{Truncations and $C^{1}$ dependence of the resolvents}%
\label{subsec:smotr}
Given $a>0$, let $V_{a} \subset S$ be  such that if 
	$s\in V_{a}$, there is no eigenvalue $\lambda$ of 
	$L_{Z\vert_{Y_{s}}}$  such that $\left\vert  
	\lambda\right\vert=a$. If $s\in V_{a}$, by the results of 
	Subsection \ref{subsec:sptru}, we can construct the vector space 
	$\mathcal{D}_{Z\vert_{Y_{s}},<a}\left(Y_{s},F\vert_{Y_{s}}\right)$, which is a  
	vector subspace of 
	$\mathcal{D}_{\mathcal{Y}^{*}_{s}}\left(Y_{s},F\vert_{Y_{s}}\right)$. Let $P_{Z\vert_{Y_{s}},<a}$ denote the associated projector.
	By  (\ref{eq:doma16a1}),  we get
	\begin{equation}\label{eq:song-3}
\mathcal{D}_{Z\vert_{Y_{s}},<a}\left(Y,F\right)=P_{Z\vert_{Y_{s}},<a}\Omega\left(Y_{s},F\vert_{Y_{s}}\right).
\end{equation}

We use the notation of Subsection \ref{subsec:esca}, and we assume that $s_{0}\in V_{a}$. We   may 
		replace $S$ by a small open neighborhood  $\mathcal{O}$ of $s_{0}\in S$, such that $\pi^{-1}\mathcal{O} 
		\simeq \mathcal{O}\times 
		Y_{s_{0}}$, and  $F$ is just 
	the pull back of its restriction to $Y_{s_{0}}$.  If $F$ is flat, 
	we may assume that the flat connection $\n^{F}$ is the pull-back 
	of the
	restriction of the flat connection of $F\vert_{Y_{s_{0}}}$. The vector field $Z$ 
		is now a smooth vector field along the fixed $Y_{s_{0}}$, that also 
		depends on $s\in S$. In the sequel, we will  be led to take a 
		smaller $\mathcal{O}$ depending on what is needed. Moreover, we will use the notation $Y=Y_{s_{0}}, 
		F=F\vert_{Y_{s_{0}}}$.
		
		By \cite[Appendix 1]{DangShen20}, if $\mathcal{O}$ 
		is small enough, we can choose the closed conical neighborhoods
		$\mathcal{C}^{*}_{0},\mathcal{C}^{*}_{u},\mathcal{C}^{*}_{s}$
		 independent of $s\in S$, and  find 
		$f,\beta,\gamma$ and		
		$\underline{\mathcal{C}}^{*}_{0},\underline{\mathcal{C}}^{*}_{u},\underline{\mathcal{C}}^{*}_{s}$,  not depending on $s\in \mathcal{O}$, and a smooth function $m_{s}$ which is also smooth in $s$  that have all the properties described in Subsection \ref{subsec:esca} for every $s\in \mathcal{O}$.   In Definition \ref{def:esc}, we can choose $U_{1},U_{2}$ not depending on $s$, and we define $g_{s}\left(z\right)$ as in (\ref{eq:res1}) by a partition of unity not depending on $s$, so that  (\ref{eq:res2}) holds with $C>0$ not depending on $s$. Then there exist $\delta>0$ and $U_{3}$ not depending on $s$ that are such that   Theorem \ref{thm:thmrec} is valid for $s\in \mathcal{O}$.	 
		
		By (\ref{eq:doma15}), for $s\in \mathcal{O}$, 
		\begin{align}\label{eq:rumi1}
&\mathcal{D}_{Z\vert_{Y_{s}},<a}\left(Y,F\right) \subset 
\mathcal{D}_{\underline{\mathcal{C}}^{*}_{s}}\left(Y,F\right),\\
&\mathcal{D}_{-Z\vert_{Y_{s}},<a}\left(Y,\Ft\right) \subset 
\mathcal{D}_{\underline{\mathcal{C}}^{*\dag}_{u}}\left(Y,\Ft\right). \nonumber 
\end{align}

Let $\Gamma_{u}\subset \underline{\mathcal{C}}^{*\dag}_{u},\Gamma_{s} \subset \underline{\mathcal{C}}^{*}_{s}$ be closed conical neighborhoods of $\mathcal{Y}^{*}_{u},\mathcal{Y}^{*}_{s}$. Let $\varphi\in 
\mathcal{D}_{\Gamma_{s}}\left(Y,F\right),\psi\in 
\mathcal{D}_{\Gamma_{u}}\left(Y,\Ft \right)$. As we saw in 
Subsection \ref{subsec:podan}, for $s\in\mathcal{O}$, $\left\langle  
\left(L_{Z\vert_{Y_{s}}}-z\right)^{-1}\varphi,\psi\right\rangle$ is a 
meromorphic function of $z$, and its poles are included in 
$\mathrm{Sp}L_{Z\vert_{Y_{s}}}$.

By the results of Subsection \ref{subsec:podan}, if $z\in S_{a}$, 
 we have the following equality of well-defined functions 
on $\mathcal{O}\cap V_{a}$, 
\begin{equation}\label{eq:rumi2x1}
\left\langle  
\left(L_{Z\vert_{Y_{s}}}-z\right)^{-1}\varphi,\psi\right\rangle=
\left\langle  \varphi,
\left(L_{-Z\vert_{Y_{s}}}-z\right)^{-1}\psi\right\rangle.
\end{equation}

A version of the following result was already established by Chaubet-Dang \cite[Proposition 7.2 (2)]{ChaubetDang24} in the case where $\n^{F}$ does not depend on $s\in \mathcal{O}$. Their method uses the anisotropic Sobolev spaces of Guedes Bonthonneau \cite{Bonthonneau20}.
\begin{theorem}\label{thm:smo2-a}
	The sets  $V_{a}$ form an open cover of $S$. If $\Gamma_{u},\Gamma_{s}$ are small enough,  the 
		functions in (\ref{eq:rumi2x1}) are  $C^{1}$ functions on 
		$\mathcal{O}\cap V_{a}$.
\end{theorem}
\begin{proof}
		In  (\ref{eq:eqc0}),  for $\mathcal{O}$ small enough, we can take $c_{0}>0$ 
		independent of $s\in \mathcal{O}$. By  
		(\ref{eq:eqhlo1}), for $s\in \mathcal{O},\mathrm{Re}z< 
		-c_{0}$, we get a holomorphic family of operators,
		\begin{equation}\label{eq:song-1}
\left(L_{Z\vert_{Y_{s}}}-z\right)^{-1}:\Omega\left(Y,F\right)\to\mathcal{D}_{\mathcal{Y}^{*}_{s}}\left(Y,F\right).
\end{equation}

		The anisotropic Sobolev spaces of Subsection 
		\ref{subsec:anis} now depend on $s\in\mathcal{O}$. In the 
		proof of Theorem \ref{thm:FaureSjostrand}, we can take 
		$p\ge 0$ as large as necessary for every $s\in \mathcal{O}$, with
		$A,\epsilon,A_{\delta,p},B_{\delta,p}$ not depending on 
		$s\in \mathcal{O}$. The operator $M_{Z,p}$ in 
		(\ref{eq:eqLZG1}) depends on $s$ through $Z\vert_{Y_{s}}$ and 
		also because $G$ depends on $s$. The operator $C_{\delta,p,s}$ in 
		(\ref{eq:dom0})  depends 
		on $s\in \mathcal{O}$, because $M_{Z,p}$ depends on $s$. These operators act on the fixed vector 
		space $\Omega_{2}\left(Y,F\right)$.

		In (\ref{eq:song1}),  for $\mathrm{Re}z\le 
		-c_{0}+\delta p-\epsilon$, the operator $C_{\delta,p,s}-z$ is 
		a bounded invertible operator from $\mathrm{dom}\left( M_{Z\vert_{Y_{s}},p}\right)$ 
		into $\Omega_{2}\left(Y,F\right)$ 
		with a uniform bounded norm on 
		$\left(C_{\delta,p,s}-z\right)^{-1}$ with respect to the $L_{2}$ norm 
		on $\Omega_{2}\left(Y,F\right)$. 
		
		To
		obtain a uniform control of the resolvent of $M_{Z\vert_{Y_{s}},p}$ under 
		the stated conditions on $z$, we only need to control the 
		inverse of 
		$1-B_{\delta,p}\left(C_{\delta,p,s}-z\right)^{-1}$.   For $s\in 
		\mathcal{O}$, the operator 
		$B_{\delta,p}\left(C_{\delta,p,s}-z\right)^{-1}$ is compact for 
		$s\in \mathcal{O}$, the arguments in the proof of Theorem 
		\ref{thm:FaureSjostrand}	giving the required uniformity.

		Let $a>0$ be such that there is no pole for the resolvent 
		$ \left( M_{Z\vert_{Y_{s_{0}}},p}-z \right) ^{-1}$ in  
		$S_{a}$. By the 
		above considerations, if $\mathcal{O}$ is small enough, this 
		is also the case for 
		$\left(M_{Z\vert_{Y_{s}},p}-z\right)^{-1}$,  the norm
		of this resolvent with respect to 
		$\Omega_{2}\left(Y_{s_{0}},F\vert_{Y_{s_{0}}}\right)$ being 
		also uniformly controlled.  We
		may as well replace 
		$\left(M_{Z\vert_{Y_{s}},p}-z\right)^{-1}$ by 
		$\left(L_{Z\vert_{Y_{s}},p}-z\right)^{-1}$, the  corresponding norms 
		being  associated with the norm on 
		$\mathcal{H}_{pG_{s}}\left(Y,\Lambda\left(T^{*}Y 
		\right) \otimes_{\R} F\right)$. This shows that $V_{a}$ is 
		open in $S$. Therefore the $V_{a}$ form an open cover of $S$.
		
		We take $p$ large enough so that
		\begin{equation}\label{eq:tubi1}
a<-c_{0}+\delta p.
\end{equation}
By the results of Subsection \ref{subsec:relz}, for $z\in 
S_{a}$, the operator 
$\left(L_{Z\vert_{Y_{s}}}-z\right)^{-1}$ maps  
$\mathcal{H}_{pG_{s}}\left(Y,\Lambda\left(T^{*}Y\right) \otimes_{\R} 
F\right)$ into itself, with a norm that is bounded 
independently of $s$.  If 
$$\varphi\in 
\mathcal{H}_{pG_{s}}\left(Y,\Lambda\left(T^{*}Y\right) \otimes 
_{\R}F\right),\psi\in 
\mathcal{H}_{-pG^{\dag}_{s}}\left(Y,\Lambda\left(T^{*}Y\right)\otimes_{\R}\Ft\right),$$
 then
		\begin{multline}\label{eq:rong1}
\left\vert\left\langle  
\left(L_{Z\vert_{Y_{s}}}-z\right)^{-1}\varphi, 
\psi\right\rangle\right\vert\\
\le C_{p}\left\Vert  \varphi\right\Vert_{\mathcal{H}_{pG_{s}}\left(Y,\Lambda\left(T^{*}Y\right) \otimes 
F\right)}\left\Vert  
\psi\right\Vert_{\mathcal{H}_{-pG^{\dag}_{s}}\left(Y,\Lambda\left(T^{*}Y\right) \otimes_{\R}
 \Ft \right)}.
\end{multline}

By (\ref{eq:tubi0}), if 
$\varphi\in\mathcal{D}_{\Gamma_{s}}\left(Y,F\right),\psi\in\mathcal{D}_{\Gamma_{u}}\left(Y,F^{*} \otimes _{\Z_{2}}\mathrm{o}\left(TY\right)\right)$,
for $p\in\R_{+}$ large enough, for $s\in \mathcal{O}$, 
$$\varphi\in 
\mathcal{H}_{pG_{s}}\left(Y,\Lambda\left(T^{*}Y \right)\otimes 
_{\R}F\right),\psi\in\mathcal{H}_{-pG^{\dag}_{s}}\left(Y,\Lambda\left(T^{*}Y\right) \otimes _{\R}F^{*} \otimes _{\Z_{2}}\mathrm{o}\left(TY\right)\right).$$
 In particular, (\ref{eq:rong1}) still holds.

We will now establish a version of the resolvent identity. Clearly,
\begin{equation}\label{eq:tubi4}
L_{Z\vert_{Y_{s}}}-L_{Z\vert_{Y_{s_{0}}}}=
\left(L_{Z\vert_{Y_{s}}}-z\right)-\left(L_{Z\vert_{Y_{s_{0}}}}-z\right).
\end{equation}

Take $\varphi\in 
\mathcal{D}_{\Gamma_{s}}\left(Y,F\right)$.  Using the results of Subsection 
\ref{subsec:relz} and the above inclusions, we get
\begin{equation}\label{eq:tubi5}
\left(L_{Z\vert_{Y_{s_{0}}}}-z\right)\left(L_{Z\vert_{Y_{s_{0}}}}-z\right)^{-1}\varphi=\varphi.
\end{equation}
By (\ref{eq:tubi4}), (\ref{eq:tubi5}), we deduce that
\begin{equation}\label{eq:tubi6}
\left( L_{Z\vert_{Y_{s}}}-L_{Z\vert_{Y_{s_{0}}}}\right)
\left(L_{Z\vert_{Y_{s_{0}}}}-z\right)^{-1}\varphi=\left(L_{Z\vert_{Y_{s}}}-z\right)
\left(L_{Z\vert_{Y_{s_{0}}}}-z\right)^{-1}\varphi-\varphi.
\end{equation}

We claim that for $p\in \R_{+}$ large enough,  if  $s\in \mathcal{O}$,  
\begin{equation}\label{eq:tubi8}
\left(L_{Z\vert_{Y_{s_{0}}}}-z\right)^{-1}\varphi\in 
\mathcal{H}^{1}_{pG_{s}}\left(Y,\Lambda\left(T^{*}Y \right)  \otimes_{\R} 
F\right).
\end{equation}
Indeed, let $ m'_{s_{0 } } $ be another order function as in Theorem \ref{thm:thm1}.
By Proposition \ref{prop:pextra}, we can assume that on  $ \left\{m_{s}  \ge -1/4\right\}$, 
\begin{align}\label{eq:htro}
 m'_{s_{0 } } >   1/4.
\end{align}
Therefore,  $ \left\{m_{s} < - 1/4\right\}$,  $ \left\{m'_{s_{0 } } >   1/4\right\}$ is an open cover of $ \mathcal{Y}^{*} \backslash Y$.

Since $ m'_{s_{0 } } $ takes its values in $ [-2,2]$,  on  $\left\{m_{s} < - 1/4\right\} $, we have 
\begin{align}\label{eq:nupz}
 m_{s } < -1/4 \le m'_{s_{0 }  }/8 .
\end{align}
Similarly, on  $\left\{m'_{s_{0 } }  >   1/4\right\} $, we have 
\begin{align}\label{eq:nupz1}
 m_{s } \le 2 < 8m_{s_{0 } }'.
\end{align}
Let $G'_{s_{0}}$ be the analogue of $G_{s}$ that is associated with $f,m'_{s_{0}}$. By \eqref{eq:nupz}, \eqref{eq:nupz1}, we have 
\begin{multline}\label{eq:xhxt}
\left\|\left(L_{Z|_{Y_{s_{0 } } }  } -z\right)^{-1} \varphi\right\|_{_{\mathcal{H}^{1}_{pG_{s} } \left(Y,\Lambda \left(T^{*} Y\right)\otimes_{\mathbf{R} } F\right)} } \\
\le C\Biggl(\left\|\left(L_{Z|_{Y_{s_{0 } } }  } -z\right)^{-1} \varphi\right\|_{{\mathcal{H}^{1}_{\frac{1}{8}pG'_{s_{0 }} } \left(Y,\Lambda \left(T^{*} Y\right)\otimes_{\mathbf{R} } F\right)} }\\
+ \left\|\left(L_{Z|_{Y_{s_{0 } } }  } -z\right)^{-1} \varphi\right\|_{{\mathcal{H}^{1}_{8pG'_{s_{0 } } } \left(Y,\Lambda \left(T^{*} Y\right)\otimes_{\mathbf{R} } F\right)} } \Biggr).
\end{multline}
If $\Gamma_{s}$ is small enough and $p\in\R_+$ is large enough,  then 
$$\varphi\in \mathcal{H}^{1}_{\frac{1}{8}pG'_{s_{0 }} } \left(Y,\Lambda \left(T^{*} Y\right)\otimes_{\mathbf{R} } F\right)\cap \mathcal{H}^{1}_{8pG'_{s_{0 } } } \left(Y,\Lambda \left(T^{*} Y\right)\otimes_{\mathbf{R} } F\right).$$ 
For $ p\in \mathbf{R}_{+} $ large enough, by Theorem \ref{thm:FaureSjostrand} and by \eqref{eq:xhxt}, we get 
\begin{multline}\label{eq:pd23}
  \left\|\left(L_{Z|_{Y_{s_{0 } } }  } -z\right)^{-1} \varphi\right\|_{_{\mathcal{H}^{1}_{pG_{s} } \left(Y,\Lambda \left(T^{*} Y\right)\otimes_{\mathbf{R} } F\right)} } \\
  \le C\left(\left\|\varphi\right\|_{{\mathcal{H}^{1}_{\frac{1}{8}pG'_{s_{0 }} } \left(Y,\Lambda \left(T^{*} Y\right)\otimes_{\mathbf{R} } F\right)} } + \left\| \varphi\right\|_{{\mathcal{H}^{1}_{8pG'_{s_{0 } } } \left(Y,\Lambda \left(T^{*} Y\right)\otimes_{\mathbf{R} } F\right)} } \right).
\end{multline}
By \eqref{eq:pd23}, we get \eqref{eq:tubi8}.

By Theorem \ref{thm:FaureSjostrand} and  by (\ref{eq:tubi8}), we deduce that
\begin{equation}\label{eq:tubi9}
\left(L_{Z\vert_{Y_{s}}}-z\right)^{-1}\left(L_{Z\vert_{Y_{s}}}-z\right)
\left(L_{Z\vert_{Y_{s_{0}}}}-z\right)^{-1}\varphi=\left(L_{Z\vert_{Y_{s_{0}}}}-z\right)^{-1}\varphi.
\end{equation}

By (\ref{eq:tubi6}), (\ref{eq:tubi9}), we obtain a version of the 
resolvent identity,
\begin{multline}\label{eq:tubi10a-1}
\left(L_{Z\vert_{Y_{s}}}-z\right)^{-1}\varphi-\left(L_{Z\vert_{Y_{s_{0}}}}-z\right)^{-1}\varphi\\
=
\left(L_{Z\vert_{Y_{s}}}-z\right)^{-1}L_{Z\vert_{Y_{s_{0}}}-Z\vert_{Y_{s}}}\left(L_{Z\vert_{Y_{s_{0}}}}-z\right)^{-1}\varphi.
\end{multline}

By (\ref{eq:tubi10a-1}), we conclude that
\begin{multline}\label{eq:tubi11}
\left\langle \left( 
\left(L_{Z\vert_{Y_{s}}}-z\right)^{-1}-\left( 
L_{Z\vert_{Y_{s_{0}}}}-z\right)^{-1}\right) \varphi, \psi  \right\rangle\\
=\left\langle  
L_{Z\vert_{Y_{s_{0}}}-Z\vert_{Y_{s}}}\left(L_{Z\vert_{Y_{s_{0}}}}-z\right)^{-1}\varphi,
\left(L_{-Z\vert_{Y_{s}}}-z\right)^{-1}\psi\right\rangle.
\end{multline}
By (\ref{eq:tubi11}), we get
\begin{multline}\label{eq:tubi11a1}
\left\vert\left\langle \left( 
\left(L_{Z\vert_{Y_{s}}}-z\right)^{-1}-\left( 
L_{Z\vert_{Y_{s_{0}}}}-z\right)^{-1}\right) \varphi, \psi  \right\rangle\right\vert\\
\le 
\left\Vert L_{Z\vert_{Y_{s}}-Z\vert_{Y_{s_{0}}}}
\left(L_{Z\vert_{Y_{s_{0}}}}-z\right)^{-1}\varphi  
\right\Vert _{\mathcal{H}_{pG_{s_{0}}}\left(Y,\Lambda\left(T^{*}Y\right) 
\otimes _{\R}F\right)}\\
\left\Vert 
\left(L_{-Z\vert_{Y_{s}}}-z\right)^{-1}\psi\right\Vert_{\mathcal{H}_{-p 
G^{\dag}_{s_{0}}}
\left(Y,\Lambda\left(T^{*}Y\right) \otimes _{\R}\Ft\right)}.
\end{multline}

Let $\left\Vert  \,\right\Vert _{1, \infty }$ denote the $\sup$-norm 
on the vector space of $C^{1}$ vector fields on $Y_{s_{0}}$. Clearly,
\begin{multline}\label{eq:tubi11a2}
\left\Vert L_{Z\vert_{Y_{s}}-Z\vert_{Y_{s_{0}}}}\left(L_{Z\vert_{Y_{s_{0}}}}-z\right)^{-1}\varphi  
\right\Vert_{\mathcal{H}_{pG_{s_{0}}}\left(Y,\Lambda\left(T^{*}Y\right) \otimes _{\R}F\right)}\\
\le C\left\Vert 
Z\vert_{Y_{s}}-Z\vert_{Y_{s_{0}}}\right\Vert_{1,\infty }\left\Vert  \left(L_{Z\vert_{Y_{s_{0}}}}-z\right)^{-1}\varphi  
\right\Vert_{\mathcal{H}^{1}_{pG_{s_{0}}}\left(Y,\Lambda\left(T^{*}Y\right) \otimes _{\R}F\right)}.
\end{multline}
By proceeding as in the proof of Proposition \ref{prop:propextReso} and using (\ref{eq:tubi11a2}), we get
\begin{multline}\label{eq:tubi11a3}
\left\Vert L_{Z\vert_{Y_{s}}-Z\vert_{Y_{s_{0}}}}\left(L_{Z\vert_{Y_{s_{0}}}}-z\right)^{-1}\varphi  
\right\Vert_{\mathcal{H}_{pG_{s_{0}}}\left(Y,\Lambda\left(T^{*}Y\right) \otimes _{\R}F\right)}\\
\le C\left\Vert 
Z\vert_{Y_{s}}-Z\vert_{Y_{s_{0}}}\right\Vert_{1,\infty }\left\Vert  \varphi  
\right\Vert_{\mathcal{H}^{1}_{pG_{s_{0}}}\left(Y,\Lambda\left(T^{*}Y\right) \otimes _{\R}F\right)}.
\end{multline}

Since 
$\left\{{}^{c}\mathcal{C}^{*\dag}_{s},{}^{c}\mathcal{C}^{*\dag}_{u}\right\}$ is an open cover of 
$\mathcal{Y}^{*}\setminus Y$,  by \cite[Appendix, Proposition 
E.30]{DyatlovZworski19}, there exist two pseudo-differential operators 
$Q_{s},Q_{u}$ of order $0$ such that
\begin{align}\label{eq:tubi10}
&Q_{s}+Q_{u}=1,
&\mathrm{WF}\left(Q_{s}\right) \subset {}^{c}\mathcal{C}^{*\dag}_{s},\qquad 
\mathrm{WF}\left(Q_{u}\right) \subset {}^{c}\mathcal{C}^{*\dag}_{u}.
\end{align}
Here, $\mathrm{WF}\left(Q_{s}\right), \mathrm{WF}\left(Q_{u}\right)$ are closed conical sets  outside of which
the full symbols of $Q_{s},Q_{u}$ behave like $\left\vert  
\xi\right\vert^{- \infty }$.

The dual version of equation \eqref{eq:eqsLz} says that for $p_{0}\in\R_{+}$ large enough, if $s\in \mathcal{O}$,  
\begin{multline}\label{eq:tubi12a-1}
\left\Vert  
\left(L_{-Z\vert_{Y_{s}}}-z\right)^{-1}\psi\right\Vert_{\mathcal{H}_{-p_{0} G^{\dag}_{s}}\left(Y,\Lambda\left(T^{*}Y\right) \otimes _{\R}\Ft\right)}\\
\le C\left\Vert  \psi\right\Vert_{\mathcal{H}_{-p_{0} G^{\dag}_{s}}\left(Y,\Lambda\left(T^{*}Y\right) \otimes _{\R}\Ft\right)}.
\end{multline}
Also
\begin{multline}\label{eq:tubi12}
\left\Vert  
Q_{s}\left(L_{-Z\vert_{Y_{s}}}-z\right)^{-1}\psi\right\Vert_{\mathcal{H}_{-p G^{\dag}_{s_{0}}}\left(\Lambda\left(T^{*}Y\right) \otimes _{\R}\Ft\right)}\\
=\left\Vert  e^{p_{0} G^{\dag}_{s}}e^{-p G^{\dag}_{s_{0}}}
Q_{s}\left(L_{-Z\vert_{Y_{s}}}-z\right)^{-1}\psi\right\Vert_{\mathcal{H}_{-p_{0}G^{\dag}_{s}}\left(Y,\Lambda\left(T^{*}Y\right) \otimes _{\R}\Ft\right)}.
\end{multline}

By \cite[Theorem 8]{FaureRoySjostrand08}, on ${}^{c}U^{\dag}_{2}$, we have the identity,
\begin{equation}\label{eq:tubi13}
\sigma\left(e^{p_{0} G^{\dag}_{s}}e^{-p G^{\dag}_{s_{0}}}Q_{s}\right)=\left(1+f^{\dag}\right)^{p_{0}m^{\dag}_{s}-pm^{\dag}_{s_{0}}}\sigma\left(Q_{s}\right).
\end{equation}
Since $\sigma \left(  Q_{s} \right) $ is  supported in 
${}^{c}\mathcal{C}^{*\dag}_{s}$ and $m^{\dag}_{s_{0}}\ge 1/4$ on ${}^{c}\mathcal{C}^{*\dag}_{s}$, for 
$p\in \R_{+}$ large enough, for $s\in\mathcal{O}$, 
$e^{p_{0} G^{\dag}_{s}}e^{-p G^{\dag}_{s_{0}}}Q_{s}$ is a pseudo-differential 
operator of order $0$. By (\ref{eq:tubi12}), (\ref{eq:tubi13}), we 
conclude that for $p\in\R_{+}$ large enough, 
\begin{multline}\label{eq:tubi14}
\left\Vert  
Q_{s}\left(L_{-Z\vert_{Y_{s}}}-z\right)^{-1}\psi\right\Vert_{\mathcal{H}_{-p G^{\dag}_{s_{0}}}\left(Y,\Lambda\left(T^{*}Y\right) \otimes _{\R}\Ft\right)}\\
\le
C_{p}\left\Vert  \psi\right\Vert_{\mathcal{H}_{-p_{0} G^{\dag}_{s}}\left(Y,\Lambda\left(T^{*}Y\right) \otimes _{\R}\Ft\right)}.
\end{multline}

Now we fix $p\in\R_{+}$ as in (\ref{eq:tubi14}). Let $m^{\prime\dag}_{s}$ be an 
analogue of $m_{s}$ in Theorem \ref{thm:thm1} and in Proposition 
\ref{prop:pextra} with 
respect to $-Z\vert_{Y_{s}}$ and the dashed cones. In particular instead of (\ref{eq:4}), 
we have 
\begin{equation}\label{eq:tubi15}
m^{\prime\dag}_{s}\ge 1/4\,\mathrm{on}\,\left(\mathcal{Y}^{*}\setminus 
\mathcal{Y}\right) \setminus\mathcal{C}^{*\dag}_{u}.
\end{equation}
We take $f$ as before, and define $g^{\prime\dag}_{s}\left(z\right)$ as in 
(\ref{eq:res1}) while replacing $f,m$ by $f^{\dag}, m^{\prime\dag}_{s}$. The same is true for 
$G^{\prime\dag}_{s}$.

The obvious analogue of (\ref{eq:eqsLz}) says that for $p'\in\R_{+}$ 
large enough, 
\begin{multline}\label{eq:tubi16}
\left\Vert \left(L_{-Z\vert_{Y_{s}}}-z\right)^{-1}\psi 
\right\Vert_{\mathcal{H}_{p' G^{\dag}_{s}}\left(Y,\Lambda\left(T^{*}Y\right) \otimes _{\R}\Ft\right)}\\
\le 
C_{p'}\left\Vert  \psi\right\Vert_{\mathcal{H}_{p' G^{\dag}_{s}}\left(Y,\Lambda\left(T^{*}Y\right) \otimes _{\R}\Ft\right)}.
\end{multline}
Observe that
\begin{multline}\label{eq:tubi17}
\left\Vert  
Q_{u}\left(L_{-Z\vert_{Y_{s}}}-z\right)^{-1}\psi\right\Vert_{\mathcal{H}_{-p G^{\dag}_{s_{0}}}\left(Y,\Lambda\left(T^{*}Y\right) \otimes _{\R}\Ft\right)}\\
=
\left\Vert  
e^{-p'G^{\prime \dag}_{s}}e^{-p G^{\dag}_{s_{0}}}Q_{u}\left(L_{-Z\vert_{Y_{s}}}-z\right)^{-1}\psi\right\Vert_{\mathcal{H}_{p'  G^{\prime\dag}_{s}}\left(Y,\Lambda\left(T^{*}Y\right) \otimes _{\R}\Ft\right)}.
\end{multline}
Instead of (\ref{eq:tubi13}), on ${}^{c}U_{2}^{\dag}$, we get
\begin{equation}\label{eq:tubi18a1}
\sigma\left(e^{-p'G^{\prime\dag}_{s}}e^{-p G^{\dag}_{s_{0}}}Q_{u}\right)=\left(1+f^{\dag}\right)^{-\left(p'm^{\prime\dag}_{s}+pm^{\dag}_{s_{0}}\right)}\sigma\left(Q_{u}\right).
\end{equation}
By (\ref{eq:tubi10}), (\ref{eq:tubi15}),  for $p'\in\R_{+}$ large enough, 
for $s\in \mathcal{O}$, $e^{-p'G^{\prime\dag}_{s}}e^{-p G^{\dag}_{s_{0}}}Q_{u}$ is a pseudo-differential operator of 
order $0$. By (\ref{eq:tubi16}), (\ref{eq:tubi17}), we deduce that for $\Gamma_{u}$ small enough and $p'\in\R_{+}$ 
large enough, for $s\in\mathcal{O}$, we get
\begin{multline}\label{eq:tubi19}
\left\Vert  
Q_{u}\left(L_{-Z\vert_{Y_{s}}}-z\right)^{-1}\psi\right\Vert_{\mathcal{H}_{-p G^{\dag}_{s_{0}}}\left(Y,\Lambda\left(T^{*}Y\right) \otimes _{\R}\Ft\right)}\\
\le C_{p,p'}\left\Vert  \psi\right\Vert_{\mathcal{H}_{p'G^{\prime\dag}_{s}}\left(Y,\Lambda\left(T^{*}Y\right) \otimes _{\R}\Ft\right)}.
\end{multline}

By (\ref{eq:tubi10}), (\ref{eq:tubi14}), and (\ref{eq:tubi19}), for $p_{0}\in\R_{+},p\in\R_{+},p'\in\R_{+}$ large 
enough,
we obtain
\begin{multline}\label{eq:tubi19a1}
\left\Vert  
\left(L_{-Z\vert_{Y_{s}}}-z\right)^{-1}\psi\right\Vert_{\mathcal{H}_{-pG^{\dag}_{s_{0}}}\left(Y,\Lambda\left(T^{*}Y\right) \otimes _{\R}\Ft\right)}\\
\le C_{p}\left\Vert  
\psi\right\Vert_{\mathcal{H}_{-p_{0} G^{\dag}_{s}}\left(Y,\Lambda\left(T^{*}Y\right) \otimes _{\R}\Ft\right)}\\
+ C_{p,p'}\left\Vert  \psi\right\Vert_{\mathcal{H}_{p'G^{\prime\dag}_{s}}\left(Y,\Lambda\left(T^{*}Y\right) \otimes _{\R}\Ft\right)}.
\end{multline}

By (\ref{eq:tubi11a1}), (\ref{eq:tubi11a3}), and (\ref{eq:tubi19a1}), 
we get
\begin{multline}\label{eq:tubi20}
\left\vert\left\langle \left( 
\left(L_{Z\vert_{Y_{s_{0}}}}-z\right)^{-1}-\left(L_{Z\vert_{Y_{s}}}-z\right)^{-1}\right) \varphi, \psi  \right\rangle\right\vert\\
\le 
\left\Vert 
Z\vert_{Y_{s}}-Z\vert_{Y_{s_{0}}}\right\Vert_{1,\infty }
\left\Vert  \varphi  
\right\Vert_{\mathcal{H}^{1}_{pG_{s_{0}}}\left(Y,\Lambda\left(T^{*}Y\right) \otimes _{\R}F\right)}\\
\Biggl(C_{p}\left\Vert  
\psi\right\Vert_{\mathcal{H}_{-p_{0} G^{\dag}_{s}}\left(Y,\Lambda\left(T^{*}Y\right) \otimes _{\R}\Ft\right)}\\
+ 
C_{p,p'}\left\Vert  
\psi\right\Vert_{\mathcal{H}_{p'G^{\prime\dag}_{s}}\left(Y,\Lambda\left(T^{*}Y\right) \otimes _{\R}\Ft\right)}\Biggr).
\end{multline}

By (\ref{eq:tubi20}), we deduce that the function
\begin{equation}\label{eq:tubi20a1}
\left\langle  
\left(L_{Z\vert_{Y_{s}}}-z\right)^{-1}\varphi, \psi\right\rangle
\end{equation}
 is continuous on $\mathcal{O}\cap V_{a}$. 
 
 We will now prove it is a $C^{1}$ function. A 
 version of equation (\ref{eq:tubi10a-1}) asserts that
 \begin{multline}\label{eq:tubi20a2}
\left(L_{-Z\vert_{Y_{s}}}-z\right)^{-1}\psi-\left(L_{-Z\vert_{Y_{s_{0}}}}-z\right)^{-1}\psi\\
=
\left(L_{-Z\vert_{Y_{s}}}-z\right)^{-1}L_{-\left(Z\vert_{Y_{s_{0}}}-Z\vert_{Y_{s}}\right)}\left(L_{-Z\vert_{Y_{s_{0}}}}-z\right)^{-1}\psi.
\end{multline}

By (\ref{eq:tubi11}), (\ref{eq:tubi20a2}), we conclude that
\begin{multline}\label{eq:tubi20a3}
\left\langle \left( 
\left(L_{Z\vert_{Y_{s}}}-z\right)^{-1}-\left(L_{Z\vert_{Y_{s_{0}}}}-z\right)^{-1}\right) \varphi, \psi  \right\rangle\\
=\left\langle  L_{Z\vert_{Y_{s_{0}}}-Z\vert_{Y_{s}}}
\left(L_{Z\vert_{Y_{s_{0}}}}-z\right)^{-1}\varphi,
\left(L_{-Z\vert_{Y_{s_{0}}}}-z\right)^{-1}\psi\right\rangle\\
+\Bigg\langle  L_{Z\vert_{Y_{s_{0}}}-Z\vert_{Y_{s}}}
\left(L_{Z\vert_{Y_{s_{0}}}}-z\right)^{-1}\varphi,\\
\left(L_{-Z\vert_{Y_{s}}}-z\right)^{-1}L_{- \left( Z\vert_{Y_{s_{0}}}-Z\vert_{Y_{s}} \right) }\left(L_{-Z\vert_{Y_{s_{0}}}}-z\right)^{-1}\psi\Bigg\rangle.
\end{multline}

Using (\ref{eq:tubi20a3}), it is now easy to conclude 
that (\ref{eq:tubi20a1}) is a differentiable function in $s\in 
\mathcal{O}$. If 
\index[not]{dS@$d^{S}$}%
$d^{S}$ denotes differentiation on $S$, from the above, we get
\begin{multline}\label{eq:20a4}
d^{S}\left\langle  
\left(L_{Z\vert_{Y_{s}}}-z\right)^{-1}\varphi, 
\psi\right\rangle\\
=-
\left\langle  
L_{d^{S}Z\vert_{Y_{s}}}\left(L_{Z\vert_{Y_{s}}}-z\right)^{-1}\varphi,
\left(L_{-Z\vert_{Y_{s}}}-z\right)^{-1}\psi\right\rangle.
\end{multline}
An iteration procedure shows that (\ref{eq:20a4}) is also a 
continuous function of $s$. We have completed the proof that the 
functions in (\ref{eq:rumi2x1}) are $C^{1}$, which concludes the proof 
of our theorem.
\end{proof}
\subsection{The $C^{1}$ dependence of the truncation projectors}%
\label{subsec:runpro}
We take $\Gamma_{u},\Gamma_{s}$ as in Theorem \ref{thm:smo2-a}. By (\ref{eq:doma14a4}), for any $s\in \mathcal{O}$, the maps 
\begin{align}\label{eq:rumi2}
&\varphi\in 
\mathcal{D}_{\Gamma_{s}}\left(Y,F\right)\to 
P_{Z\vert_{Y_{s}},<a}\varphi\in\mathcal{D}_{Z\vert_{Y_{s}},<a}\left(Y,F\right),\\
&\psi\in 
\mathcal{D}_{\Gamma_{u}}\left(Y,\Ft\right)\to 
P_{-Z\vert_{Y_{s}},<a}\psi\in\mathcal{D}_{-Z\vert_{Y_{s}},<a}\left(Y,\Ft\right), \nonumber 
\end{align}
are well-defined and continuous. 

Then
\begin{equation}\label{eq:song-4}
\left\langle  P_{Z\vert_{Y_{s}},<a}\varphi, 
\psi\right\rangle=\left\langle  \varphi, 
P_{-Z\vert_{Y_{s}},<a}\psi\right\rangle=
\left\langle  P_{Z\vert_{Y_{s}},<a}\varphi, P_{-Z\vert_{Y_{s}},<a}
\psi \right\rangle
\end{equation}
is an identity of functions on $\mathcal{O}$. 
\begin{theorem}\label{thm:smo2}
			The functions in (\ref{eq:song-4}) are $C^{1}$ on $\mathcal{O}\cap V_{a}$. The vector spaces
		$\mathcal{D}_{Z\vert_{Y_{s}},<a}\left(Y,F\right)$ are the fibers of a 
		$C^{1}$ vector bundle 
		 on $V_{a}$. In particular, if $s\in \mathcal{O}\cap V_{a}$, the map $\varphi\in 
		 \mathcal{D}_{Z\vert_{Y_{s_{0}}},<a}\left(Y,F\right)\to 
		 P_{Z\vert_{Y_{s}},<a}\varphi\in \mathcal{D}_{Z\vert_{Y_{s}},<a}\left(Y,F\right)$ 
		 defines a $C^{1}$-trivialization of the vector bundle 
		 $\mathcal{D}_{Z,<a}\left(Y,F\right)$ on $\mathcal{O}\cap V_{a}$. Near 
		 $s_{0}$, a section $\Phi$ of $\mathcal{D}_{Z,<a}\left(Y,F\right)$ is 
		 $C^{1}$ if and only if $P_{Z\vert_{Y_{s_{0}}},<a}\Phi$ is 
		 a $C^{1}$ section of $\mathcal{D}_{Z\vert_{Y_{s_{0}},<a}}\left(Y,F\right)$.

		 The 
		 operators $d^{Y}, i_{Z},L_{Z}$ restrict to $C^{1}$ 
		 endomorphisms of $\mathcal{D}_{Z,<a}\left(Y,F\right)$.
		 The complex lines $\lambda,\mu$ inherit a corresponding $C^{1}$ 
		structure, and $\tau_{\mu}\left(i_{Z}\right)$ is a nonzero 
		$C^{1}$ section of $\mu$.  The line bundles $\lambda$ and 
		$\mu$ are $C^{1}$ isomorphic. In particular, the section 
		$\tau_{\lambda}\left(i_{Z}\right)$ is a $C^{1}$ nonzero 
		section of $\lambda$.
		
		On $V_{a}$, the map $P_{Z,<a}:\left( 
		\Omega\left(Y,F\right),d^{Y} \right) \to 
		\left( 
		\mathcal{D}_{Z,<a}\left(Y,F\right),d^{Y}\vert_{\mathcal{D}_{Z,<a}\left(Y,F\right)} 
		\right) $ is a $C^{1}$ family of quasi-isomorphisms on $V_{a}$. The line bundles $\lambda$ and $\nu$ are $C^{1}$  
		isomorphic, and this isomorphism induces the canonical 
		isomorphism of the fibers in (\ref{eq:smo1a1}).
		In particular, the section 
		$\tau_{\lambda}\left(i_{Z}\right)$ of $\lambda$ induces a 
		corresponding $C^{1}$ nonzero section
		\index[not]{tniZ@$\tau_{\nu}\left(i_{Z}\right)$}%
		$\tau_{\nu}\left(i_{Z}\right)$ of the smooth line bundle 
		$\nu$.
		\end{theorem}
\begin{proof}
	By (\ref{eq:doma16a2}), we get
\begin{equation}\label{eq:tubi21}
\left\langle  P_{Z\vert_{Y_{s}},<a}\varphi,\psi\right\rangle
=-\left\langle  \frac{1}{2i\pi}\int_{z\in S_{a}}^{}\left(L_{Z\vert_{Y_{s}}}-z\right)^{-1}\varphi,\psi\right\rangle.
\end{equation}
Using Theorem \ref{thm:smo2-a} and (\ref{eq:tubi21}), we find that 
the functions in (\ref{eq:song-4}) are $C^{1}$.

Take 
$\varphi\in\mathcal{D}_{Z\vert_{Y_{s_{0}}},<a}\left(Y,F\right),\psi\in\mathcal{D}_{-Z\vert_{Y_{s_{0}}},<a}\left(Y,\Ft\right)$. Clearly,
\begin{equation}\label{eq:rubi1}
\left\langle  
P_{Z\vert_{Y_{s_{0}}},<a}\varphi,\psi\right\rangle=
\left\langle  \varphi,\psi\right\rangle
\end{equation}
gives the   pairing of
$\mathcal{D}_{Z\vert_{Y_{s_{0}}},<a}\left(Y,F\right)$ with 
$\mathcal{D}_{-Z\vert_{Y_{s_{0}}},<a}\left(Y,\Ft\right)$.  The functions in 
(\ref{eq:song-4}) being $C^{1}$, for $s\in \mathcal{O}$ close enough to 
$s_{0}$, the pairing in (\ref{eq:song-4}) is still non-degenerate. 
Therefore, for  $s\in\mathcal{O}$ close to $s_{0}$, the map $\varphi\in 
		 \mathcal{D}_{Z\vert_{Y_{s_{0}},<a}}\left(Y,F\right)\to 
		 P_{Z\vert_{Y_{s}},<a}\varphi\in \mathcal{D}_{Z\vert_{Y_{s}},<a}\left(Y,F\right)$ is 
		 injective.   In particular, for $s\in \mathcal{O}$ close 
		 enough to $s_{0}$, 
		 \begin{equation}\label{eq:rubi2}
\dim \mathcal{D}_{Z\vert_{Y_{s_{0}}},<a}\left(Y,F\right)\le \dim 
\mathcal{D}_{Z\vert_{Y_{s}},<a}\left(Y,F\right).
\end{equation}

It is not possible to exchange the roles of $s$ and $s_{0}$ because 
of lack of uniformity. We will give a direct proof of the fact that 
for $s$ close to $s_{0}$, there is equality in (\ref{eq:rubi2}).

For $0\le p\le n$, let $N_{p}$ be the operator that acts like $1$ on 
$\Lambda^{p}\left(T^{*}Y\right)$ and like $0$ in every degree 
distinct of $p$. Then $N_{p}$ acts on our spaces of currents.

Let $\varphi_{s,t}\vert_{t\in\R}$ be the group of diffeomorphisms of $Y$ 
associated with $Z\vert_{Y_{s}}$. For $r>0$ small enough, if $\mathcal{O}$ is 
small enough, $\varphi_{s,r}$ does not have fixed points in $Y$. 
By combining Proposition \ref{prop:wff} and Theorem \ref{thm:pdz1}, 
given $s$ close to $s_{0}$, 
$\mathrm{Tr_{sc}}\left[N_{p}e^{-rL_{Z\vert_{Y_{s}}}}\left(L_{Z\vert_{Y_{s}}}-z\right)^{-1}\right]$  is a meromorphic function of $z$, whose poles are included of $\mathrm{Sp}\,L_{Z\vert_{Y_{s}}}$. In particular, if $s\in V_{a}$, this function is holomorphic near $S_{a}$. Also
\begin{multline}\label{eq:trog1}
\mathrm{Tr_{sc}}\left[N_{p}e^{-rL_{Z\vert_{Y_{s}}}}\left(L_{Z\vert_{Y_{s}}}-z\right)^{-1}\right]\\
=\mathrm{Tr_{sc}}\left[N_{p}e^{-rL_{Z\vert_{Y_{s}}}}\left(L_{Z\vert_{Y_{s}}}-z\right)^{-1}P_{Z\vert_{Y_{s}},>a}\right]\\
+\mathrm{Tr_{sc}}\left[N_{p}e^{-rL_{Z\vert_{Y_{s}}}}\left(L_{Z\vert_{Y_{s}}}-z\right)^{-1}P_{Z\vert_{Y_{s}},<a}\right].
\end{multline}
By (\ref{eq:trog1}), if $s$ is close enough to $s_{0}$,
\begin{multline}\label{eq:trog2}
-\frac{1}{2i\pi}\int_{S_{a}}^{}\mathrm{Tr_{sc}}\left[N_{p}e^{-r\left(L_{Z\vert_{Y_{s}}}-z\right)}\left(L_{Z\vert_{Y_{s}}}-z\right)^{-1}\right]dz\\
=-\frac{1}{2i\pi}\int_{S_{a}}^{}\mathrm{Tr_{sc}}\left[N_{p}e^{-r\left(L_{Z\vert_{Y_{s}}}-z\right)}\left(L_{Z\vert_{Y_{s}}}-z\right)^{-1}P_{Z\vert_{Y_{s}},>a}\right]\\
-\frac{1}{2i\pi}\int_{S_{a}}^{}\mathrm{Tr_{sc}}\left[N_{p}e^{-r\left(L_{Z\vert_{Y_{s}}}-z\right)}\left(L_{Z\vert_{Y_{s}}}-z\right)^{-1}P_{Z\vert_{Y_{s}},<a}\right]dz.
\end{multline}
By Proposition \ref{prop:exten}, the first term in the right-hand side is the integral on 
$S_{a}$ of a holomorphic function that has no pole in $D_{a}$, and so 
it vanishes. By (\ref{eq:trog2}), we get
\begin{multline}\label{eq:trog3}
-\frac{1}{2i\pi}\int_{S_{a}}^{}\mathrm{Tr_{sc}}\left[N_{p}e^{-r\left(L_{Z\vert_{Y_{s}}}-z\right)}\left(L_{Z\vert_{Y_{s}}}-z\right)^{-1}\right]dz\\
=-\frac{1}{2i\pi}\int_{S_{a}}^{}\mathrm{Tr_{sc}}\left[N_{p}e^{-r\left(L_{Z\vert_{Y_{s}}}-z\right)}\left(L_{Z\vert_{Y_{s}}}-z\right)^{-1}P_{Z\vert_{Y_{s}},<a}\right]dz.
\end{multline}
From (\ref{eq:trog3}), we deduce that
\begin{equation}\label{eq:trog4}
-\frac{1}{2i\pi}\int_{S_{a}}^{}\mathrm{Tr_{sc}}\left[N_{p}e^{-r\left(L_{Z\vert_{Y_{s}}}-z\right)}\left(L_{Z\vert_{Y_{s}}}-z\right)^{-1}\right]dz
=\Tr^{\mathcal{D}_{Z\vert_{Y_{s}},<a}\left(Y,F\right)}\left[N_{p}\right].
\end{equation}
Equation (\ref{eq:trog4}) can be rewritten in the form,
\begin{equation}\label{eq:trog5}
-\frac{1}{2i\pi}\int_{S_{a}}^{}\mathrm{Tr_{sc}}\left[N_{p}e^{-r\left(L_{Z\vert_{Y_{s}}}-z\right)}\left(L_{Z\vert_{Y_{s}}}-z\right)^{-1}\right]dz=
\dim \mathcal{D}^{p}_{Z\vert_{Y_{s}},<a}\left(Y,F\right).
\end{equation}

By \cite[Theorem 4]{DangShen20}, the left-hand side of (\ref{eq:trog5}) depends 
continuously on $s$ near $s_{0}$, so that $\dim 
\mathcal{D}^{p}_{Z\vert_{Y_{s}},<a}\left(Y,F\right)$ is a locally constant function near 
$s_{0}$. 

By the above, we deduce that for $s\in\mathcal{O}$ close enough to 
$s_{0}$, the map $\varphi\in\mathcal{D}_{Z\vert_{Y_{s_{0}}},<a}\left(Y,F\right)\to 
P_{Z\vert_{Y_{s}},<a}\varphi\in\mathcal{D}_{Z\vert_{Y_{s}},<a}\left(Y,F\right)$ is an isomorphism.
We have obtained this way a trivialization 
of $\mathcal{D}_{Z,<a}\left(Y,F\right)$ near $s_{0}$. 

Take $s'_{0}\in \mathcal{O}$ close to $s_{0}$. If  
$$\varphi\in \mathcal{D}_{Z\vert_{Y_{s_{0}}},<a}\left(Y,F\right),\psi
\in \mathcal{D}_{-Z\vert_{Y_{s'_{0}}},<a}\left(Y,\Ft\right),$$
 then
\begin{equation}\label{eq:rubi4}
\left\langle  P_{Z\vert_{Y_{s}},<a}\varphi,\psi\right\rangle=
\left\langle 
P_{Z\vert_{Y_{s'_{0}}},<a}P_{Z\vert_{Y_{s}},<a}P_{Z\vert_{Y_{s_{0}}},<a}\varphi,\psi 
\right\rangle.
\end{equation}
If $\mathcal{O}$ is small enough, (\ref{eq:rubi4}) is a $C^{1}$ 
function of $s\in\mathcal{O}$. Therefore 
$$P_{Z\vert_{Y_{s'_{0}}},<a}P_{Z\vert_{Y_{s}},<a}P_{Z\vert_{Y_{s_{0}}},<a}$$
 is a $C^{1}$ section 
of 
$\Hom\left(\mathcal{D}_{Z\vert_{Y_{s_{0}}},<a}\left(Y,F\right),\mathcal{D}_{Z\vert_{Y_{s'_{0}}},<a}\left(Y,F\right)\right)$. As we saw before, if $\mathcal{O}$ is small enough
$P_{Z\vert_{Y_{s}},<a}P_{Z\vert_{Y_{s_{0}}},<a}$ is an isomorphism. By duality, this is 
also the case for $P_{Z\vert_{Y_{s'_{0}}},<a}P_{Z\vert_{Y_{s}},<a}$. 
Therefore
$$P_{Z\vert_{Y_{s'_{0}}},<a}P_{Z\vert_{Y_{s}},<a}P_{Z\vert_{Y_{s_{0}}},<a}$$
 is also an isomorphism.

As we saw before, if  $\mathcal{O}$ is small enough, $P_{Z\vert_{Y_{s'_{0}}},<a}$ 
induces an isomorphism from 
$\mathcal{D}_{Z\vert_{Y_{s}},<a}\left(Y,F\right)$ into 
$\mathcal{D}_{Z\vert_{Y_{s'_{0}}},<a}\left(Y,F\right)$.  If $\mathcal{O}$ is small enough, 
\begin{multline}\label{eq:rubi5}
\left(P_{Z\vert_{Y_{s}},<a}P_{Z\vert_{Y_{s'_{0}}},<a}\right) 
^{-1}\left( P_{Z\vert_{Y_{s}},<a}P_{Z\vert_{Y_{s_{0}}},<a}\right)\\
=
\left(P_{Z\vert_{Y_{s'_{0}}},<a}P_{Z\vert_{Y_{s}},<a}P_{Z\vert_{Y_{s'_{0}}},<a}\right) 
^{-1}\left(P_{Z\vert_{Y_{s'_{0}}},<a} 
P_{Z\vert_{Y_{s}},<a}P_{Z\vert_{Y_{s_{0}}},<a}\right). 
\end{multline}
The right-hand side is the product of two $C^{1}$ functions. 
Therefore it is
 $C^{1}$. Equation (\ref{eq:rubi5}) just says that the transition 
 functions for $\mathcal{D}_{Z,<a}\left(Y,F\right)$ are $C^{1}$, so that 
 $\mathcal{D}_{Z,<a}\left(Y,F\right)$ is  a $C^{1}$ vector bundle on $V_{a}$.

Let us now prove that the above $C^{1}$ structure does not depend on 
 the choice of the trivialization of $M$ and $F$ over $\mathcal{O}$.  
 Equivalently, for $s\in S$, let $\psi_{s}$ be a smooth family of 
 diffeomorphisms  of $Y=Y_{s_{0}}$.  Let $\Psi$ be the corresponding 
 diffeomorphism of $S\times Y$ given by $\left(s,y\right)\to 
 \left(s,\psi_{s}\left(y\right)\right)$. Set
 \begin{equation}\label{eq:rubi16}
Z'_{s}=\psi^{*}_{s}Z\vert_{Y_{s}}.
\end{equation}
Then $Z'$ is another family of Anosov vector field along the fiber 
$Y=Y_{s_{0}}$. 
Let $F'=\Psi^{*}F$. Here $\psi_{s}^{*}$ identifies 
$\mathcal{D}_{Z\vert_{Y_{s}},<a}\left(Y,F\right)$ and $\mathcal{D}_{Z'_{s},<a}\left(Y,F'\right)$. 
We have to prove that this 
identification is $C^{1}$. Let $Q_{s_{0},s}: \mathcal{D}_{Z'_{s},<a}\left(Y,F'\right)\to \mathcal{D}_{Z'_{s_{0}},<a}\left(Y,F'\right)$ be given by the inverse to 
$P_{Z'_{s},<a}P_{Z'_{s_{0}},<a}:\mathcal{D}_{Z'_{s_{0}},<a}\left(Y,F'\right)\to \mathcal{D}_{Z'_{s},<a}\left(Y,F'\right)$. Tautologically, $Q_{s_{0},s}$ is $C^{1}$ in $s\in\mathcal{O}$.
 The transition map $T_{s}:\mathcal{D}_{Z\vert_{Y_{s_{0}},<a}}\left(Y,F\right)\to\mathcal{D}_{Z'_{s_{0}},<a}\left(Y,F'\right)$ is given by
\begin{equation}\label{eq:treb1}
T_{s}=	Q_{s_{0},s}\psi^{*}_{s}P_{Z\vert_{Y_{s},<a}}P_{Z\vert_{Y_{s_{0}},<a}},
\end{equation}
which is also $C^{1}$. This shows that the $C^{1}$ structure on  $\mathcal{D}_{Z,<a}\left(Y,F\right)$ does not depend on the trivialization.

We will now prove that $\varphi$ is a $C^{1}$ section of 
$\mathcal{D}_{Z,<a}\left(Y,F\right)$ near $s_{0}$ if and only if 
$P_{Z\vert_{Y_{s_{0}}},<a}\varphi$ is a $C^{1}$ section of  
$\mathcal{D}_{Z\vert_{Y_{s_{0}}},<a}\left(Y,F\right)$. First, assume that $\varphi$ is a $C^{1}$ section of 
$\mathcal{D}_{Z,<a}\left(Y,F\right)$ near $s_{0}$. For our purpose, we may even 
assume that $\varphi_{s}=P_{Z\vert_{Y_{s}},<a}\varphi_{s_{0}}$. Then
\begin{equation}\label{eq:rubi17}
P_{Z\vert_{Y_{s_{0}}},<a}\varphi_{s}=P_{Z\vert_{Y_{s_{0}}},<a}
P_{Z\vert_{Y_{s}},<a}P_{Z\vert_{Y_{s_{0}}},<a}\varphi_{s_{0}}.
\end{equation}
As we saw before, the right-hand side is a $C^{1}$ function of $s$.
Conversely, if $P_{Z\vert_{Y_{s_{0}}},<a}\varphi $ is $C^{1}$, note 
that
\begin{equation}\label{eq:tubi18}
\varphi_{s}=P_{Z\vert_{Y_{s}},<a}\left(P_{Z\vert_{Y_{s_{0}}},<a}
P_{Z\vert_{Y_{s}},<a}P_{Z\vert_{Y_{s_{0}}},<a}\right)^{-1}
P_{Z\vert_{Y_{s_{0}}},<a}\varphi_{s}.
\end{equation}
By (\ref{eq:tubi18}), we conclude that $\varphi$ is a $C^{1}$ 
section of $\mathcal{D}_{Z,<a}\left(Y,F\right)$.

The fact that $d^{Y},i_{Z},L_{Z}$ restrict to $C^{1}$ endomorphisms 
of $\mathcal{D}_{Z,<a}\left(Y,F\right)$ is trivial.

For $0<a<b$, the splitting (\ref{eq:zomb17}) is a splitting of 
$C^{1}$ vector bundles on $V_{a}\cap V_{b}$, from which we conclude that $\lambda,\mu$ are
$C^{1}$ line bundles on $S$, and that
$\tau_{\mu}\left(i_{Z}\right)$ is a $C^{1}$ section of $\mu$.

Let us  prove that the fiberwise canonical isomorphism $\mu_{s} \simeq 
\lambda_{s}$ that was obtained in Theorem \ref{thm:liba} is indeed 
$C^{1}$. It is enough to prove that $R_{Z,F,>a}\left(0\right)$ is a 
$C^{1}$ function on $V_{a}$. On $V_{a}$, since $R_{Z,F,>a}$ has no zeroes 
and no poles in the disk $D_{a}$, we have the identity,
\begin{equation}\label{eq:tubi19z1}
	R_{Z,F,>a}\left(0\right)=\frac{1}{2i\pi}\int_{S_{a}}^{}R_{Z,F,>a}\left(z\right)\frac{dz}{z}.
\end{equation}
By the results of \cite[Section 5.2]{DangShen20},  the restriction of $R_{Z,F}\left(z\right)$ 
to $S_{a}$ is also a $C^{1}$ function on $V_{a}$. Also since 
$D_{Z,<a}\left(Y,F\right)$ is a $C^{1}$ vector bundle 
on $V_{a}$, $R_{Z,F,<a}\left(z\right)$ is also $C^{1}$ on $S_{a}$. 
By (\ref{eq:romb9z1}), we deduce that $R_{Z,F,>a}\left(z\right)$ is 
$C^{1}$ on $V_{a}$. By (\ref{eq:tubi19z1}), we 
conclude that $R_{Z,F,>a}\left(0\right)$ is $C^{1}$ on $V_{a}$. Using 
Theorem \ref{thm:liba}, we find that the identification $\lambda 
\simeq \mu$ is indeed $C^{1}$. In particular, 
$\tau_{\lambda}\left(i_{Z}\right)$ is  a $C^{1}$ section of 
$\lambda$.

We will prove that 
$$P_{Z,<a}:\left( \Omega\left(Y,F\right),d^{Y} \right) \to 
\left( \mathcal{D}_{Z,<a}\left(Y,F\right), d^{Y}\vert_{\mathcal{D}_{Z,<a}\left(Y,F\right)} \right) $$
 is a $C^{1}$ family of quasi-isomorphisms. We 
only need to show that if $\varphi\in 
\Omega\left(Y_{s_{0}},F\right)$, then 
$P_{Z,<a}\varphi$ is a $C^{1}$ section of $\mathcal{D}_{Z,<a}\left(Y,F\right)$. If we 
take $\psi\in \mathcal{D}_{-Z\vert_{Y_{s_{0}}},<a}\left(Y,\Ft\right)$ in 
(\ref{eq:song-4}) which   is a $C^{1}$ function, we deduce that 
$P_{Z\vert_{Y_{s_{0}}},<a}P_{Z\vert_{Y_{s}},<a}\varphi$ is a $C^{1}$ section of 
$\mathcal{D}_{Z\vert_{Y_{s_{0}}},<a}\left(Y,F\right)$. By the first part 
of our theorem, we deduce that $P_{Z\vert_{Y_{s}},<a}\varphi$ 
is a $C^{1}$ section of $\mathcal{D}_{Z,<a}\left(Y,F\right)$. This proves that 
$P_{Z,<a}:\left( \Omega\left(Y,F\right),d^{Y} \right) \to \left( 
\mathcal{D}_{Z,<a}\left(Y,F\right),d^{Y} \right) $ is a $C^{1}$ 
family of morphisms of complexes. By Theorem \ref{thm:presa}, this is 
a fiberwise quasi-isomorphism.

Let us now prove that the complex lines $\lambda$ and $\nu$ are 
canonically isomorphic as $C^{1}$ line bundles. In (\ref{eq:smo1a1}), 
we already saw that their fibers are canonically isomorphic.  We are 
forced here to recall how the smooth structure on the line bundle 
$\nu$ is constructed. Let $g^{TY}$ be a Euclidean metric on $TY$, 
let $g^{F}$ be a Hermitian metric on $F$. Let $\square^{Y}$ denote 
the corresponding fiberwise Hodge Laplacian. Given $b>0$, let 
$\Omega_{<b}\left(Y,F\right)$ denote the direct sum of eigenspaces 
associated with eigenvalues $<b$. Put
\begin{equation}\label{eq:rubi20}
\nu_{<b}=\det\Omega_{<b}\left(Y,F\right).
\end{equation}
Let $U_{b} \subset S$ be the open set where $b$ is not an eigenvalue 
of $\square^{Y}$. The $U_{b}$ form an open cover of $S$. Then 
$\nu_{<b}$ is a smooth line bundle on $U_{b}$. The smooth line bundle 
$\nu$ is obtained by a formal construction similar to the 
construction of the line bundle $\lambda$ in Definition 
\ref{def:libu}. The line bundle $\nu$ is a smooth line bundle. The 
arguments of \cite{BismutGilletSoule88c} show that when we choose 
other metrics, we obtain another smooth line bundle which is 
canonically isomorphic to our given $\nu$.  So we may as well work 
with our given construction of $\nu$. 

It follows from the above that on $V_{a}\cap U_{b}$, $P_{Z,<a}$ 
restricts to a $C^{1}$ map 
$\Omega_{<b}\left(Y,F\right)\to\mathcal{D}_{Z,<a}\left(Y,F\right)$ which is pointwise 
a quasi-isomorphism. We can now form the cone 
$C^{\bullet}=C\left(\Omega_{<b}^{\bullet}\left(Y,F\right),\mathcal{D}_{Z,<a}^{\bullet}\left(Y,F\right)\right)$. The cone $C^{\bullet}$ is pointwise exact, so that its determinant has a canonical $C^{1}$ nonzero section. By proceeding as in \cite[proof of Theorem 2.8]{BismutGilletSoule88c},  this canonical section is the one that identifies the fibers of $\lambda$ and $\nu$. We have found that the complex lines $\lambda$ and $\nu$ are canonically $C^{1}$ isomorphic. In particular, the section $\tau_{\lambda}\left(i_{Z}\right)$ induces a corresponding $C^{1}$ section of the smooth line bundle $\nu$. The proof of our theorem is complete. 
\end{proof}
\begin{remark}\label{rem:rap}
	The first sentence in our theorem was already established in Chaubet-Dang \cite[Proposition 7.2 (3)]{ChaubetDang24}. When the fibers $Y$ are contact, in Theorem \ref{thm:chada}, we showed that pointwise, $\tau_{\nu}\left(i_{Z}\right)$ coincides with $\tau_{\nu}^{\mathrm{CD}} $. Our Theorem implies in particular that $\tau _{\nu }^{\mathrm{CD}}$ is $C^{1}$ with respect to the canonical smooth structure on $\nu$. This extends \cite[Proposition 8.1]{ChaubetDang24} to the case where $H\left(Y,F\right)$ does not vanish. 
\end{remark}
\subsection{The derivative of $i_{Z}$ with respect to a flat 
connection}%
\label{subsec:derfl}
We fix $s_{0}\in S$ as before. We have defined a trivialization of 
$\mathcal{D}_{Z,<a}\left(Y,F\right)$ near $s_{0}$. In particular, near $s_{0}$, 
$\mathcal{D}_{Z,<a}\left(Y,F\right)$ is equipped with a flat connection, which is 
denoted $\n^{\mathcal{D}_{Z,<a}\left(Y,F\right)}_{s_{0}}$.  
On $\mathcal{O}$, $i_{Z}\vert_{\mathcal{D}_{Z,<a}\left(Y,F\right)}$ is a $C^{1}$ section of 
$\mathcal{D}_{Z,<a}\left(Y,F\right)$. As before, 
\index[not]{dS@$d^{S}$}%
$d^{S}$ denotes differentiation on $S$.
\begin{proposition}\label{prop:cov}
	We have the following identity at $s_{0}$:
	\begin{equation}\label{eq:cova1}
		\left[\n^{\mathcal{D}_{Z,<a}\left(Y,F\right)}_{s_{0}},i_{Z}\vert_{\mathcal{D}_{Z,<a}\left(Y,F\right)}\right]\vert_{s=s_{0}}=P_{Z\vert_{Y_{s_{0}}},<a} i_{d^{S}Z\vert _{s=s_{0}}}P_{Z_{s_{0},<a}}.
\end{equation}
\end{proposition}
\begin{proof}
	We take $\varphi\in \mathcal{D}_{Z\vert_{Y_{s_{0}}},<a}\left(Y,F\right)$. Since $P_{Z,<a}$ 
	and $i_{Z}$ commute, near $s_{0}$, we get
	\begin{equation}\label{eq:cova2}
i_{Z}P_{Z,<a}\varphi=P_{Z,<a}i_{Z}\varphi.
\end{equation}
At $s$ close to $s_{0}$, equation (\ref{eq:cova2}) can be written in 
the form,
\begin{equation}\label{eq:cova3}
\left(i_{Z}P_{Z,<a}\varphi\right)_{s}=P_{Z\vert_{Y_{s}},<a,}i_{Z\vert_{Y_{s}}}\varphi.
\end{equation}

Moreover, since $P_{Z\vert_{Y_{s}},<a}$ identifies $\mathcal{D}_{Z\vert_{Y_{s_{0}}},<a}\left(Y,F\right)$ with 
$\mathcal{D}_{Z\vert_{Y_{s}},<a}\left(Y,F\right)$ near $s_{0}$, for $s$ close enough to $s_{0}$, 
there is a  unique $\rho_{s}\in\mathcal{D}_{Z\vert_{Y_{s_{0}}},<a}\left(Y,F\right)$ such that
\begin{equation}\label{eq:cova4}
\left( i_{Z}P_{Z,<a}\varphi \right) _{s}=P_{Z\vert_{Y_{s}},<a}\rho_{s}.
\end{equation}

Since $P_{Z,<a}\varphi$ is flat with respect to 
$\n^{\mathcal{D}_{Z,<a}\left(Y,F\right)}_{s_{0}}$, by (\ref{eq:cova4}), 
 we get
\begin{equation}\label{eq:cova7}
\left[\n^{\mathcal{D}_{Z,<a}\left(Y,F\right)}_{s_{0}},i_{Z}\right]P_{Z,<a}\varphi=P_{Z\vert_{Y_{s_{0}}},<a}d^{S}\rho\vert_{s=s_{0}}.
\end{equation}

By (\ref{eq:cova3}), we can rewrite (\ref{eq:cova4}) in the form,
\begin{equation}\label{eq:cova5}
 P_{Z\vert_{Y_{s}},<a}\left( i_{Z\vert_{Y_{s}}}\varphi -\rho_{s}\right)=0.
\end{equation}
By (\ref{eq:cova5}), since $\mathcal{D}_{Z\vert_{Y_{s_{0}}},<a}\left(Y,F\right)$ is preserved 
by $i_{Z\vert_{Y_{s_{0}}}}$, we get 
\begin{equation}\label{eq:cova6}
\left(i_{Z}\varphi-\rho\right)_{s_{0}}=0.
\end{equation}

By (\ref{eq:cova5}), if $\psi\in \mathcal{D}_{-Z\vert_{Y_{s_{0}}},<a}\left(Y,\Ft\right)$, then
\begin{equation}\label{eq:cova6a1}
\left\langle i_{Z\vert_{Y_{s}}}\varphi-\rho_{s},P_{-Z\vert_{Y_{s}},<a}\psi  
\right\rangle=0.
\end{equation}
Clearly, as $s\to s_{0}$, 
\begin{equation}\label{eq:cova6a2}
\frac{1}{s-s_{0}} \left( i_{Z\vert_{Y_{s}}}\varphi-\rho_{s}\right) \to 
 \left( i_{d^{S}Z\vert_{s=s_{0}}}\varphi-d^{S}\rho\vert_{s=s_{0}} 
 \right),
\end{equation}
the convergence taking place in 
$\mathcal{D}_{\mathcal{Y}^{*}_{s}}\left(Y,F\right)$. This is indeed 
obvious for the first term in the left-hand side, and also easy for 
the second term which varies in the fixed vector space 
$\mathcal{D}_{Z\vert_{Y_{s_{0}}},<a}\left(Y,F\right)$.

By (\ref{eq:cova6a1}), (\ref{eq:cova6a2}), we get
\begin{equation}\label{eq:cova6a3}
\left\langle   
i_{d^{S}Z\vert_{s=s_{0}}}\varphi-d^{S}\rho\vert_{s=s_{0}} 
,\psi\right\rangle=0.
\end{equation}
By (\ref{eq:cova7}), (\ref{eq:cova6a3}), we get (\ref{eq:cova1}). The proof of our proposition is complete. 
\end{proof}
\subsection{The case where $F$ is flat on $M$}%
\label{subsec:flam}
We assume now that  $\n^{F}$ is a flat connection on $F$ over the full $M$. Then  
$H\left(Y,F\vert_{Y}\right)$ is a flat vector bundle on $S$, the flat 
connection being the
\index[terms]{Gauss-Manin connection}%
Gauss-Manin connection. In this case, 
(\ref{eq:co6}) induces an isomorphism of smooth line bundles, and 
$\nu$ inherits a corresponding flat connection $\n^{\nu}$.

The Gauss-Manin connection is defined as follows. Let $F\Omega\left(M,F\right)$ be the decreasing filtration on $ \Omega\left(M,F\right)$ by the degree on $S$, so that if  $q\in \N$, elements of $F^{q}\Omega\left(M,F\right)$ are forms in which $\Omega^{q}\left(S,\C\right)$ appear as a factor. Let $\alpha\in C^{\infty}\left(M, \Lambda^{p}\left(T^{*}Y\right) \otimes _{\R}F\right)$ be such that $d^{Y} \alpha=0$, and let $\left[\alpha\right]$ be the corresponding section of $H\left(Y,F\right)$ on $S$. Recall that 
\index[not]{ga@$\gamma$}%
$\gamma$ is the embedding of the generic fiber  $Y\to M$. Let $\beta\in \Omega^{p}\left(M,F\right)$ be a lift of $\alpha$, i.e.,   $\gamma^{*}\beta=\alpha$. Then $\gamma^{*}d^{M}\beta=0$, so that $d^{M} \beta\in F^{1}\Omega^{p+1}\left(M,F\right)$. Let $\left[d^{M}\beta\right]$ be the image of $d^{M}\beta$ in 
$$\left(F^{1}/F^{2}\right)\Omega^{p+1}\left(M,F\right)=C^{\infty }\left(M,\pi^{*}T^{*}S \otimes \Lambda^{p}\left(\TsX\right) \otimes _{\R}F\right).$$
 Using the fact that $d^{M,2}=0$, we find that $d^{Y}\left[d^{M}\beta\right]=0$. One can prove  that $\left[d^{M}\beta\right]$ ultimately defines the Gauss-Manin differential  $\n^{H\left(Y,F\right)}\left[\alpha\right]$ on $H^{p}\left(Y,F\right)$. 

  In the trivialization of Subsection \ref{subsec:smotr}, $F$ and 
$d^{Y}$ do not depend any more on $s$, and only $Z\vert_{Y_{s}}$ varies with $s$.  Then $\n^{H\left(Y,F\right)}$ is simply defined via the action of $d^{S}$. In this case, the nontrivial fact on the Gauss-Manin connection is that it does not depend on the choice of the product structure.

We use the same notation as in the previous subsections. In particular the $C^{1}$ section 
\index[not]{tniZ@$\tau_{\nu}\left(i_{Z}\right)$}%
$\tau_{\nu}\left(i_{Z}\right)$ of $\nu$ was defined in Theorem \ref{thm:smo2}.
\begin{theorem}\label{thm:flat}
	For $\mathcal{O}$ small enough,   
		 $P_{Z\vert_{Y_{s_{0}}},<a}P_{Z\vert_{Y_{s}},<a}P_{Z\vert_{Y_{s_{0}}},<a}$ 
		 is a $C^{1}$ family of automorphisms of 
		 $\left(D_{Z\vert_{Y_{s_{0}}},<a}\left(Y,F\right),d^{Y}\vert_{\mathcal{D}_{Z\vert_{Y_{s_{0}}},<a}\left(Y,F\right)}\right)$ which is 
		 $C^{1}$ $d^{Y}$-homotopic to the identity. Moreover, 
		$\tau_{\nu}\left(i_{Z}\right)$ is a flat nonzero section of 
		$\nu$, and as such it is a smooth section of $\nu$.
\end{theorem}
\begin{proof}
We use the notation in the proof of Theorem \ref{thm:smo2}. The key point is that in our 
trivialization, $d^{Y}$ does not depend on the base point 
$s\in\mathcal{O}$.  Also recall that $Y=Y_{s_{0}}$. 

If $s\in \mathcal{O}$ is close to 
$s_{0}$, 
$P_{Z\vert_{Y_{s}},<a}$ induces an isomorphism of complexes 
$\left( 
\mathcal{D}_{Z\vert_{Y_{s_{0}}},<a}\left(Y,F\right),d^{Y}\vert_{\mathcal{D}_{Z\vert_{Y_{s_{0}}},<a}\left(Y,F\right)} \right) \to 
\left( \mathcal{D}_{Z\vert_{Y_{s}},<a}\left(Y,F\right),d^{Y}\vert_{\mathcal{D}_{Z\vert_{Y_{s}},<a}\left(Y,F\right)} \right) $. By our construction of the flat connection $\n^{\mathcal{D}_{Z,<a}\left(Y,F\right)}_{s_{0}}$, $d^{Y}\vert_{\mathcal{D}_{Z,<a}\left(Y,F\right)}$ is parallel with respect to this connection. Since the above morphism is $C^{1}$, it 
induces a $C^{1}$ identification of the corresponding cohomology 
groups,
\begin{equation}\label{eq:trav1}
\psi_{s_{0},s}: H\mathcal{D}_{Z\vert_{Y_{s_{0}}},<a}\left(Y,F\right) \simeq 
H\mathcal{D}_{Z\vert_{Y_{s}},<a}\left(Y,F\right).
\end{equation}

By (\ref{eq:doma22}), we get
\begin{equation}\label{eq:trav3}
P_{Z,>a}=\left[d^{Y},i_{Z}L_{Z}^{-1}P_{Z,>a}\right].
\end{equation}
By (\ref{eq:trav3}), we deduce that
\begin{multline}\label{eq:trav4}
P_{Z\vert_{Y_{s'_{0}}},<a}P_{Z\vert_{Y_{s_{0}},<a}}-
P_{Z\vert_{Y_{s'_{0}}},<a}P_{Z,<a}P_{Z\vert_{Y_{s_{0}}},<a}\\
=\left[d^{Y},P_{Z\vert_{Y_{s'_{0}}},<a}i_{Z}L_{Z}^{-1}P_{Z,>a}
P_{Z\vert_{Y_{s_{0}}},<a}\right].
\end{multline}

We will now prove that the term in the commutator with $d^{Y}$ in the right-hand side of (\ref{eq:trav4}) is properly $C^{1}$ in $s\in\mathcal{O}$. Observe that
\begin{equation}\label{eq:ny-1}
L_{Z}^{-1}P_{Z,>a}=\frac{1}{2i\pi}\int_{z\in S_{a}}^{}\frac{1}{z}\left(L_{Z}-z\right)^{-1}P_{Z,>a}dz.
\end{equation}
We will rewrite (\ref{eq:ny-1}) in the form
\begin{equation}\label{eq:ny-2}
	L_{Z}^{-1}P_{Z,>a}=\frac{1}{2i\pi}\int_{z\in S_{a}}^{}\frac{1}{z}\left(\left(L_{Z}-z\right)^{-1}- \left(L_{Z}-z\right)^{-1}P_{Z,<a}\right) dz.
\end{equation}
Using Theorems \ref{thm:smo2-a} and \ref{thm:smo2}, we find that the term
 $$P_{Z\vert_{Y_{s'_{0}}},<a}i_{Z}L_{Z}^{-1}P_{Z,>a}
P_{Z\vert_{Y_{s_{0}}},<a}$$
 is indeed properly $C^{1}$ in $s\in\mathcal{O}$.

By (\ref{eq:trav4}), we have proved the  first part of our theorem.
We also conclude that for $s'_{0},s\in \mathcal{O}$ close to $s_{0}$,
\begin{equation}\label{eq:trav2}
\psi_{s_{0},s'_{0}}=\psi_{s,s'_{0}}\psi_{s_{0},s}.
\end{equation}

By (\ref{eq:trav2}),  on $\mathcal{O}$, $H\mathcal{D}_{Z,<a}\left(Y,F\right)$ has 
been canonically trivialized. In particular, $H\mathcal{D}_{Z,<a}\left(Y,F\right)$ is equipped 
with a canonical flat connection. We will now prove that this 
connection is the Gauss-Manin connection.

By Theorem \ref{thm:smo2},  $P_{Z,<a}:\left( 
\Omega\left(Y,F\right),d^{Y}\right) \to 
\left( \mathcal{D}_{Z,<a}\left(Y,F\right),d^{Y}\vert_{\mathcal{D}_{Z,<a}\left(Y,F\right)} \right) $ is a $C^{1}$ quasi-isomorphism. Also  $$P_{Z\vert_{Y_{s_{0}}},<a}P_{Z,<a}:\left( 
\Omega\left(Y,F \right) ,d^{Y}\right)\to \left( \mathcal{D}_{Z\vert_{Y_{s_{0}}},<a}\left(Y,F\right),d^{Y}\vert_{\mathcal{D}_{Z\vert_{Y_{s_{0}}},<a}\left(Y,F\right)} \right)$$ 
is a $C^{1}$ quasi-isomorphism. 

Let $f\in 
\Omega\left(Y,F\right) $ be a 
closed form. Its cohomology class $\left[f\right]$ can be viewed as a 
flat section of $H\left(Y,F\right)$ near $s_{0}$ with respect to the 
Gauss-Manin connection, so that  $P_{Z,<a}\left[f\right]$ 
is a flat section of $H\mathcal{D}_{Z,<a}\left(Y,F\right)$. Therefore the natural flat 
connection on $H\mathcal{D}_{Z,<a}\left(Y,F\right)$ is the Gauss-Manin connection.

By the results of Subsection \ref{subsec:derfl}, near $s_{0}$, 
$\mathcal{D}_{Z,<a}\left(Y,F\right)$ is equipped with a
flat connection $\n^{\mathcal{D}_{Z,<a}\left(Y,F\right)}_{s_{0}}$. 
On $V_{a}$, $\vartheta_{<a}=\det 
\mathcal{D}_{Z,<a}\left(Y,F\right)$ is a $C^{1}$  line bundle, that we equip with 
the induced flat connection $\n^{\vartheta_{<a}}_{s_{0}}$. By 
(\ref{eq:zomb18}) and what we just saw, we have the identification of 
flat line bundles on $V_{a}$, 
\begin{equation}\label{eq:trav5}
\vartheta_{<a} \simeq \nu.
\end{equation}
We will now show that $\tau_{\nu}\left(i_{Z}\right)$ is a flat 
section of $\nu$.

 On $V_{a}$,  $\tau_{<a}\left(i_{Z}\right)$  is a $C^{1}$ nonzero 
section of $\vartheta_{<a}$. 
On $\mathcal{O}$, the line bundle 
$\nu$ is identified with $\nu_{s_{0}}$, and this identification is 
just induced by the Gauss-Manin connection  $\n^{\nu}$.
  We only need to show that  
$\tau_{\nu}\left(i_{Z}\right)$ is flat at $s_{0}$. 

 Let $\alpha$ be a smooth one form on $Y$ such that $i_{Z}\alpha=1$. 
 Note that $\alpha$  depends on $s\in \mathcal{O}$.
Over $V_{a}$, set
\begin{equation}\label{eq:trav6}
	\alpha_{<a}=P_{Z,<a}\alpha P_{Z,<a}.
\end{equation}
Then $\alpha_{<a}$ is $C^{1}$, and moreover it is a homotopy with 
respect to $i_{Z}\vert_{\mathcal{D}_{Z,<a}\left(Y,F\right)}$, i.e.,
\begin{equation}\label{eq:trav6x1}
\left[i_{Z}\vert_{\mathcal{D}_{Z,<a}\left(Y,F\right)},\alpha_{<a}\right]=1.
\end{equation}
By (\ref{eq:clev15}), we get
\begin{equation}\label{eq:trav7}
\frac{\n^{\vartheta_{<a}}_{s_{0}}\tau_{<a}\left(i_{Z}\right)}{\tau_{<a}\left(i_{Z}\right)}=\Trs\left[\alpha_{<a}\left[\n^{\mathcal{D}_{Z,<a}\left(Y,F\right)}_{s_{0}},i_{Z}\right]\right].
\end{equation}

In   (\ref{eq:romb9z1}),  
$R_{Z,F}\left(\sigma\right),R_{Z,F,<a}\left(\sigma\right),R_{Z,F,>a}\left(\sigma\right)$ depend also on $s\in V_{a}$.
By the properties we gave of the vector bundle $\mathcal{D}_{Z,<a}\left(Y,F\right)$, 
for $\Re \sigma$ large enough, $R_{Z,F,<a}\left(\sigma\right)$ is a  
holomorphic function of $\sigma$ that is also a
$C^{1}$ function of $s\in V_{a}$.   In the sequel operators like $d^{Y}, 
i_{Z}$ are always restricted to $\mathcal{D}_{Z,<a}\left(Y,F\right)$. We get
\begin{multline}\label{eq:trav9}
d^{S}\log 
R_{Z,F,<a}\left(\sigma\right)\\
=\Trs^{\mathcal{D}_{Z,<a}\left(Y,F\right)}\left[N\left[\n^{\mathcal{D}_{Z,<a}\left(Y,F\right)}_{s_{0}},L_{Z}\vert_{\mathcal{D}_{Z,<a}\left(Y,F\right)}\right]\left(L_{Z}+\sigma\right)^{-1}\right].
\end{multline}
Since $d^{Y}$ is parallel with respect to 
$\n^{\mathcal{D}_{Z,<a}\left(Y,F\right)}_{s_{0}}$,  we get
\begin{equation}\label{eq:malo1}
\left[\n^{\mathcal{D}_{Z,<a}\left(Y,F\right)}_{s_{0}},L_{Z}\vert_{\mathcal{D}_{Z,<a}\left(Y,F\right)}\right]=-\left[d^{Y},\left[\n^{\mathcal{D}_{Z,<a}\left(Y,F\right)}_{s_{0}},i_{Z}\right]\right].
\end{equation}
Also
\begin{equation}\label{eq:malo2}
\left[d^{Y},N\right]=-d^{Y}.
\end{equation}
Since supertraces vanish on supercommutators, by (\ref{eq:malo1}), 
(\ref{eq:malo2}), we obtain,
\begin{equation}\label{eq:malo3}
d^{S}\log 
R_{Z,F,<a}\left(\sigma\right)=
-\Trs^{\mathcal{D}_{Z,<a}\left(Y,F\right)}\left[d^{Y}\left[\n^{\mathcal{D}_{Z,<a}\left(Y,F\right)}_{s_{0}},i_{Z}\right]\left(L_{Z}+\sigma\right)^{-1}\right].
\end{equation}

Since $i_{Z}^{2}=0$, we get
\begin{equation}\label{eq:malo4}
\left[i_{Z},\left[\n^{\mathcal{D}_{Z,<a}\left(Y,F\right)}_{s_{0}},i_{Z}\right]\right]=0.
\end{equation}
By  (\ref{eq:trav6x1}), (\ref{eq:malo4}), we obtain
\begin{equation}\label{eq:malo5}
\left[\n^{\mathcal{D}_{Z,<a}\left(Y,F\right)}_{s_{0}},i_{Z}\right]=\left[i_{Z},\alpha_{<a}\left[\n^{\mathcal{D}_{Z,<a}\left(Y,F\right)}_{s_{0}},i_{Z}\right]\right].
\end{equation}
Since supertraces vanish on supercommutators, by (\ref{eq:malo3}), 
(\ref{eq:malo5}), we get
\begin{equation}\label{eq:malo6}
d^{S}\log 
R_{Z,F,<a}\left(\sigma\right)=
-\Trs^{\mathcal{D}_{Z,<a}\left(Y,F\right)}\left[\alpha_{<a}\left[\n^{\mathcal{D}_{Z,<a}\left(Y,F\right)}_{s_{0}},i_{Z}\right]L_{Z}\left(L_{Z}+\sigma\right)^{-1}\right].
\end{equation}
We can rewrite (\ref{eq:malo6}) in the form
\begin{multline}\label{eq:malo7}
d^{S}\log 
R_{Z,F,<a}\left(\sigma\right)=
-\Trs^{\mathcal{D}_{Z,<a}\left(Y,F\right)}\left[\alpha_{<a}\left[\n^{\mathcal{D}_{Z,<a}\left(Y,F\right)}_{s_{0}},i_{Z}\right]\right]\\
+\sigma 
\Trs^{\mathcal{D}_{Z,<a}\left(Y,F\right)}\left[\alpha_{<a}\left[\n^{\mathcal{D}_{Z,<a}\left(Y,F\right)}_{s_{0}},i_{Z}\right]\left(L_{Z}+\sigma\right)^{-1}\right].
\end{multline}

By (\ref{eq:trav7}), (\ref{eq:malo7}), we obtain,
\begin{multline}\label{eq:malo8}
\frac{\n^{\vartheta_{<a}}_{s_{0}}\left[\tau_{<a}\left(i_{Z}\right)R_{Z,F,<a}\left(\sigma\right)\right]}{\tau_{<a}\left(i_{Z}\right)R_{Z,F,<a}\left(\sigma\right)}\\
=
\sigma 
\Trs^{\mathcal{D}_{Z,<a}\left(Y,F\right)}\left[\alpha_{<a}\left[\n^{\mathcal{D}_{Z,<a}\left(Y,F\right)}_{s_{0}},i_{Z}\right]\left(L_{Z}+\sigma\right)^{-1}\right].
\end{multline}

By \cite[Corollary 4.5]{DangShen20}, for $\Re\,\sigma\gg 1$, 
$R_{Z,F}\left(\sigma\right)$ is a $C^{1}$ function of $s$, and moreover, for
$r>0$ small enough, we get
\begin{equation}\label{eq:malo9}
d^{S}\log 
R_{Z,F}\left(\sigma\right)=\sigma\Tr_{\mathrm{s,int}}\left[\alpha 
i_{d^{S}Z}e^{-r\left(L_{Z}+\sigma\right)}\left(L_{Z}+\sigma\right)^{-1}\right].
\end{equation}
Formally, equation (\ref{eq:malo9}) is the exact analogue 
of the finite dimensional (\ref{eq:malo7}),  except that the trivial equation
\begin{equation}\label{eq:malo10}
	\mathrm{Tr^{reg}_{s,int}}\left[\alpha i_{d^{S}Z}\right]=0
\end{equation}
has been used.

 By (\ref{eq:romb9z1}), (\ref{eq:malo8}), and  (\ref{eq:malo9}), we 
 conclude that for $r>0$ small enough and 
$\sigma\gg 1$, 
\begin{multline}\label{eq:trav10}
\frac{\n^{\vartheta_{<a}}_{s_{0}}\left[\tau_{<a}\left(i_{Z}\right)
R^{-1}_{Z,F,>a}\left(\sigma\right)\right]}{\tau_{<a}\left(i_{Z}\right)
R_{Z,F,>a}^{-1}\left(\sigma\right)} \\
=\sigma 
\Trs^{\mathcal{D}_{Z,<a}\left(Y,F\right)}\left[\alpha_{<a}\left[\n^{\mathcal{D}_{Z,<a}\left(Y,F\right)}_{s_{0}},i_{Z}\right]\left(L_{Z}+\sigma\right)^{-1}\right]\\
-\sigma\Tr_{\mathrm{s,int}}\left[\alpha 
i_{d^{S}Z}e^{-r\left(L_{Z}+\sigma\right)}\left(L_{Z}+\sigma\right)^{-1}\right].
\end{multline}

The left-hand side of (\ref{eq:trav10}) extends to a holomorphic 
function near $\sigma=0$.   More 
precisely, we have the identity,
\begin{equation}\label{eq:trav11}
\frac{\n^{\vartheta_{<a}}_{s_{0}}\left[\tau_{<a}\left(i_{Z}\right)
R^{-1}_{Z,F,>a}\left(\sigma\right)\right]}{\tau_{<a}\left(i_{Z}\right)
R^{-1}_{Z,F,>a}\left(\sigma\right)}\vert_{\sigma=0}=\frac{\n^{\vartheta_{<a}}_{s_{0}}\tau_{\nu}\left(i_{Z}\right)}{\tau_{\nu}\left(i_{Z}\right)}.
\end{equation}

The two terms in the right-hand side of (\ref{eq:trav10}) are meromorphic 
near
$\sigma=0$. To compute the value of (\ref{eq:trav11}) at $\sigma=0$, 
it is enough to compute the constant term of the corresponding 
expansion of (\ref{eq:trav10}) at $\sigma=0$. If $A\left(\sigma\right) $ is meromorphic at 
$\sigma=0$, we denote by $A^{0}$ the constant term.

Recall that $P_{Z,0}$ is the spectral projector associated with the eigenvalue $0$. Clearly,
\begin{multline}\label{eq:trav12}
\left[\sigma 
\Trs^{\mathcal{D}_{Z,<a}\left(Y,F\right)}\left[\alpha_{<a}\left[\n^{\mathcal{D}_{Z,<a}\left(Y,F\right)}_{s_{0}},i_{Z}\right]\left(L_{Z}+\sigma\right)^{-1}\right]\right]^{0}\\
=
\left[\sigma\Trs^{\mathcal{D}_{Z,0}\left(Y,F\right)}\left[P_{Z,0}\alpha_{<a}\left[\n^{\mathcal{D}_{Z,<a}\left(Y,F\right)}_{s_{0}},i_{Z}\right]P_{Z,0}\left(L_{Z}+\sigma\right)^{-1}P_{Z,0}\right]\right]^{0}.
\end{multline}
Put
\begin{equation}\label{eq:trav12a1}
\alpha_{0}=P_{Z,0}\alpha P_{Z,0}.
\end{equation}
Since $\left[i_{Z},\alpha\right]=1$ and 
$\left[i_{Z},P_{Z,0}\right]=0$,  we get the identity of 
operators acting on $\mathcal{D}_{Z,0}\left(Y,F\right)$,
\begin{equation}\label{eq:trav12a2}
\left[i_{Z},\alpha_{0}\right]=1.
\end{equation}
Using (\ref{eq:rond3}), (\ref{eq:malo4}), and (\ref{eq:trav12a2}), we have the identity of operators acting on 
$\mathcal{D}_{Z,0}\left(Y,F\right)$,
\begin{align}\label{eq:trav12a3}
&\left[i_{Z},P_{Z,0}\alpha_{<a}\left[\n^{\mathcal{D}_{Z,<a}\left(Y,F\right)}_{s_{0}},i_{Z}\right]P_{Z,0}\right]=P_{Z,0}\left[\n^{\mathcal{D}_{Z,<a}\left(Y,F\right)}_{s_{0}},i_{Z}\right]P_{Z,0},\\
&\left[i_{Z},P_{Z,0}\alpha_{0}\left[\n^{\mathcal{D}_{Z,<a}\left(Y,F\right)}_{s_{0}},i_{Z}\right]P_{Z,0}\right]=P_{Z,0}\left[\n^{\mathcal{D}_{Z,<a}\left(Y,F\right)}_{s_{0}},i_{Z}\right]P_{Z,0}. \nonumber 
\end{align}
By (\ref{eq:trav12a3}), we deduce that
\begin{multline}\label{eq:trav12a4}
\left[i_{Z},\alpha_{0}P_{Z,0}\left( \alpha_{<a}-\alpha_{0} \right) 
\left[\n^{\mathcal{D}_{Z,<a}\left(Y,F\right)}_{s_{0}},i_{Z}\right]P_{Z,0}\right]\\
=P_{Z,0}\left( \alpha_{<a}-\alpha_{0} \right) 
\left[\n^{\mathcal{D}_{Z,<a}\left(Y,F\right)}_{s_{0}},i_{Z}\right]P_{Z,0}.
\end{multline}
Since supertraces vanish on supercommutators and 
$\left[i_{Z},\left(L_{Z}+\sigma\right)^{-1}\right]=0$, by 
 (\ref{eq:trav12a4}), we obtain,
\begin{multline}\label{eq:trav12a5}
\Trs^{\mathcal{D}_{Z,0}\left(Y,F\right)}\left[P_{Z,0}\alpha_{<a}\left[\n^{\mathcal{D}_{Z,<a}\left(Y,F\right)}_{s_{0}},i_{Z}\right]P_{Z,0}\left(L_{Z}+\sigma\right)^{-1}P_{Z,0}\right]\\
=\Trs^{\mathcal{D}_{Z,0}\left(Y,F\right)}\left[\alpha_{0}P_{Z,0}\left[\n^{\mathcal{D}_{Z,<a}\left(Y,F\right)}_{s_{0}},i_{Z}\right]P_{Z,0}\left(L_{Z}+\sigma\right)^{-1}P_{Z,0}\right].
\end{multline}

By proceeding as in (\ref{eq:trav12}), for $r>0$ small enough, 
\begin{multline}\label{eq:trav13}
\left[\sigma 
\mathrm{Tr_{s,int}}\left[\alpha 
i_{d^{S}i_{Z}}e^{-r\left(L_{Z}+\sigma\right)}\left(L_{Z}+\sigma\right)^{-1}\right]\right]^{0}\\
=
\left[\sigma\Trs^{\mathcal{D}_{Z,0}\left(Y,F\right)}\left[P_{Z,0}\alpha\left[d^{S},i_{Z}\right]P_{Z,0}e^{-r\left(L_{Z}+\sigma\right)}\left(L_{Z}+\sigma\right)^{-1}P_{Z,0}\right]\right]^{0}.
\end{multline}
Using the  same arguments as in 
(\ref{eq:trav12a1})--(\ref{eq:trav12a5}), we 
get
\begin{multline}\label{eq:trav13a1}
\Trs^{\mathcal{D}_{Z,0}\left(Y,F\right)}\left[P_{Z,0}\alpha\left[d^{S},i_{Z}\right]P_{Z,0}e^{-r\left(L_{Z}+\sigma\right)}\left(L_{Z}+\sigma\right)^{-1}P_{Z,0}\right]\\
=\Trs^{\mathcal{D}_{Z,0}\left(Y,F\right)}\left[\alpha_{0}P_{Z,0}\left[d^{S},i_{Z}\right]P_{Z,0}e^{-r\left(L_{Z}+\sigma\right)}\left(L_{Z}+\sigma\right)^{-1}P_{Z,0}\right].
\end{multline}

By (\ref{eq:trav10}), (\ref{eq:trav11}), and (\ref{eq:trav12a5})--(\ref{eq:trav13a1}), we obtain
\begin{multline}\label{eq:trav14}
\frac{\n^{\vartheta_{<a}}_{s_{0}}\tau_{\nu}\left(i_{Z}\right)}{\tau_{\nu}\left(i_{Z}\right)}
=\Biggl[\sigma\Trs^{\mathcal{D}_{Z,0}\left(Y,F\right)}\Biggl[\alpha_{0}P_{Z,0}\Biggl(\left[\n^{\mathcal{D}_{Z,<a}\left(Y,F\right)}_{s_{0}},i_{Z}\right] \\
-\left[d^{S},i_{Z}\right]e^{-r\left(L_{Z}+\sigma\right)}\Biggr)P_{Z,0}
\left(L_{Z}+\sigma\right)^{-1}P_{Z,0}\Biggr]\Biggr]^{0}.
\end{multline}

We will calculate (\ref{eq:trav14}) at $s=s_{0}$. Using Proposition 
\ref{prop:cov} and (\ref{eq:trav14}), at $s=s_{0}$, we obtain
\begin{multline}\label{eq:trav15}
\Trs^{\mathcal{D}_{Z,0}\left(Y,F\right)}\Biggl[\alpha_{0}P_{Z,0}\left(\left[\n^{\mathcal{D}_{Z,<a}\left(Y,F\right)}_{s_{0}},i_{Z}\right]-\left[d^{S},i_{Z}\right]e^{-r\left(L_{Z}+\sigma\right)}\right)\\
P_{Z,0}\left(L_{Z}+\sigma\right)^{-1}P_{Z,0}\Biggr]
=\Trs^{\mathcal{D}_{Z,0}\left(Y,F\right)}\left[\alpha_{0}P_{Z,0}i_{\left[d^{S},Z\right]}P_{Z,0}\frac{1-e^{-r\left(L_{Z}+\sigma\right)}}{L_{Z}+\sigma}P_{Z,0}\right].
\end{multline}
Since the function $\frac{1-e^{-r x}}{x}$ extends holomorphic at 
$x=0$,  at $s=s_{0}$, (\ref{eq:trav15}) is holomorphic near $\sigma=0$, and its value at $0$ is given by
\begin{multline}\label{eq:trav16}
\Trs^{\mathcal{D}_{Z,0}\left(Y,F\right)}\Biggl[\alpha_{0}P_{Z,0}\left(\left[\n^{\mathcal{D}_{Z,<a}\left(Y,F\right)}_{s_{0}},i_{Z}\right]-\left[d^{S},i_{Z}\right]e^{-r\left(L_{Z}+\sigma\right)}\right)\\
P_{Z,0}
\left(L_{Z}+\sigma\right)^{-1}P_{Z,0}\Biggr]_{\sigma=0}
=\Trs^{\mathcal{D}_{Z,0}\left(Y,F\right)}\left[\alpha_{0}P_{Z,0}i_{\left[d^{S},Z\right]}P_{Z,0}\frac{1-e^{-r L_{Z}}}{L_{Z}}P_{Z,0}\right].
\end{multline}

By (\ref{eq:trav14}), (\ref{eq:trav16}), at $s_{0}$, we get
\begin{equation}\label{eq:trav17}
\frac{\n^{\vartheta_{<a}}\tau_{\nu}\left(i_{Z}\right)}{\tau_{\nu}\left(i_{Z}\right)}=0,
\end{equation}
which says that $\tau_{\nu}\left(i_{Z}\right)$ is a flat 
section of $\nu$. The proof of our theorem is complete. 
\end{proof}
\section{Reidemeister metrics and the Fried conjecture}%
\label{sec:rsmfr}
The purpose of this Section is to provide a refurbished version of the Fried conjecture in full generality. We assume that our flat bundle is unimodular, which means that $\det F$ carries a flat metric.  Our conjecture basically says that the norm of the section $\tau_{\nu}\left(i_{Z}\right)$ with respect to the Reidemeister metric on  $\nu$ is equal to $1$. 

This short section is organized as follows. In Subsection \ref{subsec:reme}, we recall the definition of the Reidemeister metric $\left\Vert \,\right\Vert_{\nu}$ on $\nu$. In our context, this metric does not depend on the choice of a metric on $\det F$.

In Subsection \ref{subsec:refr}, we give a general formulation of the Fried conjecture. 
\subsection{The Reidemeister metric}%
\label{subsec:reme}
Let $Y$ be a compact connected manifold, and let $\left(F,\n^{F}\right)$ be a complex vector bundle on $Y$ equipped with a flat connection. Let $\n^{\det F}$ be the induced connection on $\det F$. In the sequel, we assume that $F$ is unimodular, i.e., $\det F$ carries a flat metric, this metric being unique up to a constant. Equivalently, if $m=\mathrm{rank}F$, the holonomy representation in $\mathrm{GL}\left(m,\C\right)$ factors through matrices whose module of their determinant is equal to $1$. A metric $g^{F}$ on $F$ will be said to be unimodular if the induced metric $g^{\det F}$ on $\det F$ is flat.

Set
\begin{equation}\label{eq:bode1}
\nu=\det H\left(Y,F\right).
\end{equation}

Let $g^{F}$ be a Hermitian metric on $F$, that induces a  metric $g^{\det F}$ on $\det F$.

When $g^{F}$ is flat, given a triangulation $K$ of $Y$, Reidemeister \cite{Reidemeister35} defined a metric $\left\Vert \,\right\Vert_{\nu}$ on $\nu$, and he showed that this metric does not depend on the triangulation $K$. This construction was extended by Bismut-Zhang \cite{BismutZhang92} and Müller  \cite{Muller93b} when $g^{\det F}$ is flat, in which case these authors proved a corresponding version of the Cheeger-Müller theorem on the equality of the Reidemeister and the Ray-Singer metrics.
\begin{proposition}\label{prop:inde}
The metric $\left\Vert \,\right\Vert_{\nu}$ depends only on $g^{\det F}$. More precisely, if $g^{\det F}$ is a flat metric on $\det F$, if $\lambda>0$ and if $g^{\det F\prime}=\lambda g^{\det F}$, if $\left\Vert \,\right\Vert'_{\nu}$ is the corresponding metric on $\nu$, then
\begin{equation}\label{eq:bode2}
\frac{\left\Vert \,\right\Vert'_{\nu}}{\left\Vert \,\right\Vert_{\nu}}=\lambda^{\frac{1}{2}\chi\left(Y\right)}.
\end{equation}
In particular, if $\chi\left(Y\right)=0$, the metric $\left\Vert \,\right\Vert_{\nu}$ is independent of any parameter.
\end{proposition}

Let $\mathrm{Diff}\left(Y\right)$ be the group of diffeomorphisms of $Y$.  Let $\mathrm{Diff}^{F}\left(Y\right)$ be the group of diffeomorphisms of $Y$ that lift to $\left(F,\n^{F}\right)$. Then $\mathrm{Diff}^{F}\left(Y\right)$ acts on $H\left(Y,F\right)$.

Observe that $\psi_{*}g^{\det F}$ is also a flat metric on $\det F$, so that
\begin{equation}\label{eq:bode2a1}
\psi_{*}g^{\det F}=\lambda\left(\psi\right)g^{\det F},
\end{equation}
and $\lambda\left(\psi\right)$ does not depend on $g^{\det F}$. 
\begin{definition}\label{def:actc}
If $\psi\in \mathrm{Diff}^{F}\left(Y\right)$, set
\begin{equation}\label{eq:bode2a1}
\rho\left(\psi \right) =\psi\vert_{\nu}.
\end{equation}
\end{definition}
Then $\rho: \mathrm{Diff}^{F}\left(Y\right)\to \C^{*}$ is a character.
\begin{proposition}\label{prop:val}
The following identity holds:
\begin{equation}\label{eq:bode2a2}
\left\vert \rho\left(\psi\right)\right\vert=\lambda\left(\psi\right)^{-\frac{1}{2}\chi\left(Y\right)}.
\end{equation}
In particular, if $\chi\left(Y\right)=0$, $\rho$ takes its values in $\mathrm{U}\left(1\right)$.
\end{proposition}
\begin{proof}
	Let $\left\Vert \,\right\Vert_{\nu},\left\Vert \,\right\Vert'_{\nu}$ be the metrics on $\nu$ associated with $g^{\det F}, \psi_{*}g^{\det F}$. Tautologically,  $\psi:\left(\nu,\left\Vert \,\right\Vert_{\nu} \right) \to \left(\nu,\left\Vert \,\right\Vert'_{\nu}\right)$ is an isometry. By combining Proposition \ref{prop:inde} and the above, our proposition follows.
\end{proof}

\subsection{A modified Fried conjecture}%
\label{subsec:refr}
We make the same assumptions as in Section \ref{sec:detli} and we use the corresponding notation. In particular, $Y$ is a compact manifold such that $\chi\left(Y\right)=0$. Also, we assume that $F$ is unimodular.

As we saw in Subsection \ref{subsec:reme}, the metric $\left\Vert \,\right\Vert_{\nu}$ on $\nu$ is unambiguously defined.

Recall that $\tau_{\nu}\left(i_{Z}\right)$ is a nonzero section of $\nu$.
\begin{conjecture} 
The following identity holds:
\begin{equation}\label{eq:bode3}
\left\Vert \tau_{\nu}\left(i_{Z}\right)\right\Vert_{\nu}=1.
\end{equation}
\end{conjecture}

Let us consider the family situation of Subsection \ref{subsec:flam}. By Theorem \ref{thm:flat}, $\tau_{\nu}\left(i_{Z}\right)$ is a flat section of $\nu$. We assume that $F$ is fiberwise unimodular. This just means that $\det F \otimes \det \overline{F}$ is a line bundle which is the pull back of a line bundle $\eta$ on $S$. Given a positive section $m$ of $\eta^{-1}$, this defines a fiberwise flat metric on $\det F$.
\begin{proposition}\label{prop:flatn}
The metric $\left\Vert \,\right\Vert_{\nu}$ is flat on $S$.
\end{proposition}
\begin{proof}
We may trivialize the fibration $\pi:M\to S$ on a small open $\mathcal{O} \subset S$ and we use the invariance of the metric $\left\Vert \,\right\Vert_{\nu}$ under scaling. 
\end{proof}
 
By Theorem \ref{thm:flat}, we know that $\tau_{\nu}\left(i_{Z}\right)$ is flat on $S$.

Then (\ref{eq:bode3}) is compatible with the above results.
\section{Declarations}%
\label{sec:decl}
\subsection{Funding}%
\label{subsec:fun}
Shu Shen was partially supported by ANR Grant ANR-20-CE40-0017.
\subsection{Conflicts of interest/Competing interests}%
\label{subsec:coin}
The authors have no relevant financial or non-financial interests to disclose.

\printindex[terms]
\printindex[not]

\bibliography{Bismut,Others}

\end{document}